\documentclass{amsart}

\usepackage{amssymb,amsthm}
\usepackage{stmaryrd} 
\usepackage{mathrsfs} 
\usepackage{colonequals}

\usepackage{longtable}
\usepackage{booktabs}

\usepackage{comment}

\usepackage[utf8]{inputenc}
\usepackage[OT2,T1]{fontenc}
\usepackage[english]{babel}
\DeclareFontFamily{U}{wncy}{}
\DeclareFontShape{U}{wncy}{m}{n}{<->wncyr10}{}
\DeclareSymbolFont{mcy}{U}{wncy}{m}{n}
\DeclareMathSymbol{\Sha}{\mathord}{mcy}{"58}

\usepackage{csquotes}

\usepackage[protrusion=false]{microtype}

\usepackage{xcolor}

\usepackage{enumerate}

\usepackage[
	colorlinks, citecolor=darkgreen,
	linkcolor=darkblue,
	pagebackref,
	unicode=true,
	pdfauthor={Timo Keller and Michael Stoll},
	pdftitle={},
]{hyperref}

\renewcommand{\epsilon}{\varepsilon}
\renewcommand{\phi}{\varphi}

\def\ol#1{\overline{#1}}
\def\ul#1{\underline{#1}}
\def\Alphabet{A,B,C,D,E,F,G,H,I,J,K,L,M,N,O,P,Q,R,S,T,U,V,W,X,Y,Z}
\def\alphabet{a,b,c,d,e,f,g,h,i,j,k,l,m,n,o,p,q,r,s,t,u,v,w,x,y,z}
\def\endpiece{xxx}
\def\makeAlphabet[#1]{\expandafter\makeA#1,xxx,}
\def\makealphabet[#1]{\expandafter\makea#1,xxx,}
\def\makeA#1,{\def\temp{#1}\ifx\temp\endpiece\else%
	\mkbb{#1}\mkfrak{#1}\mkbf{#1}\mkcal{#1}\mkscr{#1}\mkbs{#1}\expandafter\makeA\fi}%
\def\makea#1,{\def\temp{#1}\ifx\temp\endpiece\else\mkfrak{#1}\mkbf{#1}\mkbs{#1}\expandafter\makea\fi}%
\def\mkbb#1{\expandafter\def\csname bb#1\endcsname{\mathbb{#1}}}
\def\mkfrak#1{\expandafter\def\csname fr#1\endcsname{\mathfrak{#1}}}
\def\mkbf#1{\expandafter\def\csname b#1\endcsname{\mathbf{#1}}}
\def\mkcal#1{\expandafter\def\csname c#1\endcsname{\mathcal{#1}}}
\def\mkscr#1{\expandafter\def\csname s#1\endcsname{\mathscr{#1}}}
\def\mkbs#1{\expandafter\def\csname bs#1\endcsname{{\boldsymbol{#1}}}}
\def\makeop[#1]{\xmakeop#1,xxx,}
\def\mkop#1{\expandafter\def\csname #1\endcsname{{\operatorname{#1}}}} %
\def\xmakeop#1,{\def\temp{#1}\ifx\temp\endpiece\else\mkop{#1}\expandafter\xmakeop\fi}%
\def\makeup[#1]{\xmakeup#1,xxx,}
\def\mkup#1{\expandafter\def\csname #1\endcsname{{\mathrm{#1}\,}}} %
\def\xmakeup#1,{\def\temp{#1}\ifx\temp\endpiece\else\mkup{#1}\expandafter\xmakeup\fi}%
\makeAlphabet[\Alphabet]
\makealphabet[\alphabet]
\makeop[Hom,Aut,End,Mor,SL,GL,H,ord,Irr,Ell,Gal,Cl,Pic,NS,Gal,d,Re,Im,res,Symb,Ev,Char,Ram,SU,BSD,Res,cork,Tam,id] 
\makeup[Spec,Proj,dR,new,old,AJ,tr,ker,im,coker]

\newcommand{\Q}{{\mathbf{Q}}}
\newcommand{\Qbar}{{\overline{\Q}}}
\newcommand{\Z}{{\mathbf{Z}}}
\newcommand{\Zhat}{{\hat{\Z}}}
\newcommand{\F}{{\mathbf{F}}}

\newcommand{\T}{{\mathbf{T}}}

\newcommand{\R}{{\mathbf{R}}}
\renewcommand{\C}{{\mathbf{C}}}
\newcommand{\mfr}{{\mathfrak{m}}}
\renewcommand{\P}{{\mathbf{P}}}

\newcommand{\Jac}{\operatorname{Jac}}
\newcommand{\PGL}{\operatorname{PGL}}
\newcommand{\PSL}{\operatorname{PSL}}
\newcommand{\GSp}{\operatorname{GSp}}

\newcommand{\Ann}{\operatorname{Ann}}
\newcommand{\GalQ}{{\Gal(\Qbar|\Q)}}
\newcommand{\GalQp}{{\Gal(\Qbar_p|\Q_p)}}
\newcommand{\GalQl}{{\Gal(\Qbar_\ell|\Q_\ell)}}
\newcommand{\et}{{\mathrm{\acute{e}t}}}
\renewcommand{\H}{{\operatorname{H}}}
\newcommand{\Het}{\H_\et}

\providecommand{\Derived}[1]{{\mathscr{D}#1}}

\newcommand{\pr}{\operatorname{pr}}

\newcommand{\Mat}{\operatorname{Mat}}
\newcommand{\Lie}{\operatorname{Lie}}
\newcommand{\Frob}{\operatorname{Frob}}
\newcommand{\Sel}{\operatorname{Sel}}

\newcommand{\Nm}{\mathbf{N}}
\newcommand{\Nabs}{\mathscr{N}}
\newcommand{\disc}{\operatorname{disc}}
\newcommand{\Tr}{\operatorname{Tr}}
\newcommand{\rk}{\operatorname{rk}}
\newcommand{\Lrk}{L\operatorname{-rk}}
\newcommand{\lcm}{\operatorname{lcm}}
\renewcommand{\Re}{\operatorname{Re}}

\newcommand{\resultant}{\operatorname{res}}

\newcommand{\charpol}{\operatorname{charpol}}
\newcommand{\cond}{\operatorname{cond}}

\renewcommand{\max}{{\operatorname{max}}}

\renewcommand{\ss}{{\mathrm{ss}}}
\newcommand{\hhat}{{\hat{h}}}
\newcommand{\an}{{\mathrm{an}}}
\newcommand{\nr}{{\mathrm{nr}}}
\newcommand{\tame}{{\mathrm{t}}}
\newcommand{\wild}{{\mathrm{w}}}

\newcommand{\tors}{{\mathrm{tors}}}

\newcommand{\Sym}{\operatorname{Sym}}

\newcommand{\Reg}{\operatorname{Reg}}

\newcommand{\Div}{\operatorname{Div}}

\newcommand{\diag}{\operatorname{diag}}

\newcommand{\Frac}{\operatorname{Frac}}
\newcommand{\isoto}{\xrightarrow{\sim}}

\newcommand{\iso}{\simeq}
\newcommand{\isom}{\cong}

\newcommand{\dual}{{\vee}}
\newcommand{\inj}{\hookrightarrow}
\newcommand{\surj}{\twoheadrightarrow}
\newcommand{\defeq}{\colonequals}

\providecommand{\gen}[1]{\langle{#1}\rangle}
\providecommand{\ideal}[1]{\langle{#1}\rangle}
\newcommand{\bpi}{\boldsymbol{\pi}}

\newcommand{\To}{\longrightarrow}

\newcommand{\LMFDBLabel}[1]{\textnormal{\href{https://www.lmfdb.org/EllipticCurve/Q/#1/}{\texttt{#1}}}}

\theoremstyle{plain}
\newtheorem{theorem}{Theorem}[subsection]
\newtheorem{lemma}[theorem]{Lemma}
\newtheorem{corollary}[theorem]{Corollary}
\newtheorem{proposition}[theorem]{Proposition}
\newtheorem{conjecture}[theorem]{Conjecture}

\theoremstyle{definition}
\newtheorem{definition}[theorem]{Definition}
\newtheorem{definition-proposition}[theorem]{Definition-Proposition}
\newtheorem{algorithm}[theorem]{Algorithm}
\newtheorem{question}[theorem]{Question}
\newtheorem{assumption}[theorem]{Assumption}

\theoremstyle{remark}
\newtheorem{example}[theorem]{Example}
\newtheorem{examples}[theorem]{Examples}
\newtheorem{remark}[theorem]{Remark}

\usepackage{cleveref}
\crefname{theorem}{Theorem}{Theorems}
\crefname{lemma}{Lemma}{Lemmata}
\crefname{corollary}{Corollary}{Corollaries}
\crefname{proposition}{Proposition}{Propositions}
\crefname{definition}{Definition}{Definitions}
\crefname{definition-proposition}{Definition-Proposition}{Definition-Propositions}
\crefname{conjecture}{Conjecture}{Conjectures}
\crefname{question}{Question}{Questions}
\crefname{example}{Example}{Examples}
\crefname{examples}{Examples}{Examples}
\crefname{algorithm}{Algorithm}{Algorithms}
\crefname{remark}{Remark}{Remarks}
\crefname{assumption}{Assumption}{Assumptions}

\usepackage{tikz-cd}
\usepackage[all]{xy}

\definecolor{darkgreen}{rgb}{0,0.5,0}
\definecolor{darkblue}{rgb}{0,0,0.8}
\makeatletter
\let\@wraptoccontribs\wraptoccontribs
\makeatother

\begin{document}
\title[Verification of strong BSD for modular abelian surfaces]{Complete verification of strong BSD\\ for many modular abelian surfaces over $\Q$}

\author{Timo Keller}
\address{Rijksuniveriteit Groningen, Bernoulli Institute, Bernoulliborg, Nijenborgh 9, 9747 AG Groningen, The Netherlands}
\address{Leibniz Universität Hannover, Institut für Algebra, Zahlentheorie und Diskrete Mathematik, Welfengarten 1, 30167 Hannover, Germany}
\address{Department of Mathematics, Chair of Computer Algebra, Universität Bay\-reuth, Germany)}
\email{math@kellertimo.de}
\urladdr{\url{https://www.timo-keller.de}}

\author{Michael Stoll}
\address{Department of Mathematics, Chair of Computer Algebra, Universität Bay\-reuth, Germany}
\email{Michael.Stoll@uni-bayreuth.de}
\urladdr{\url{https://www.mathe2.uni-bayreuth.de/stoll/}}

\contrib[with an appendix joint with]{Sam Frengley}


\date{\today}

\begin{abstract}
  We develop the theory and algorithms necessary to be able to verify the strong
  Birch--Swinnerton-Dyer Conjecture for absolutely simple modular abelian varieties
  over~$\Q$. We apply our methods to all 28 Atkin--Lehner quotients of~$X_0(N)$ of genus~$2$,
  all 97 genus~$2$ curves from the LMFDB whose Jacobian is of this type and six
  further curves originally found by Wang.
  We are able to verify the strong BSD~Conjecture unconditionally and exactly
  in all these cases; this is the first time that strong~BSD has been confirmed
  for absolutely simple abelian varieties of dimension at least~$2$.
  We also give an example where we verify that the order of the Tate--Shafarevich group
  is~$7^2$ and agrees with the order predicted by the BSD~Conjecture.
\end{abstract}

\maketitle

\numberwithin{subsection}{section}
\numberwithin{equation}{section}

\tableofcontents

	
	\section{Introduction}
	
	\subsection{Background}
	
	The \emph{Conjecture of Birch and Swinnerton-Dyer} (\enquote{BSD} for short),
	originally formulated based on extensive computations by Birch
	and Swinnerton-Dyer~\cite{BSD1965} in the 1960s for elliptic curves over~$\Q$,
	is one of the most important open conjectures in number theory. For example,
	it is one of the seven \enquote{Millennium Problems}, for whose solution the Clay
	Foundation is offering a million dollars each. It relates in a surprising
	way analytic invariants of an elliptic curve~$E$, which are obtained
	via its $L$-series from its local properties (essentially the number
	of points modulo~$p$ on~$E$, for all prime numbers~$p$), to global
	arithmetic invariants like the rank of the Mordell--Weil group~$E(\Q)$,
	its regulator, and the rather mysterious Tate--Shafarevich group~$\Sha(E/\Q)$.
	The conjecture has been generalized to cover all abelian varieties over
	all algebraic number fields. It consists of two parts, which we will
	explain for the case of an abelian variety~$A$ of dimension~$g$ over~$\Q$.
	
	One attaches to~$A$ its $L$-function $L(A/\Q,s)$, which is defined by an Euler
	product over all prime numbers~$p$. If $A$ is the Jacobian variety of
	a curve~$X$ of genus~$g$, the Euler factor at~$p$ for a prime~$p$ of
	good reduction is determined by the number of $\F_{p^n}$-points on the mod~$p$
	reduction of~$X$ for $n \leq g$.
	It follows from the Weil conjectures for
	varieties over finite fields that the Euler product converges for
	$\Re(s) > \frac{3}{2}$ to a holomorphic function. A standard conjecture predicts
	that $L(A/\Q,s)$ extends to an entire function; this is known when $A$ is
	\emph{modular}, i.e., occurs as an isogeny factor of the Jacobian~$J_0(N)$
	of one of the modular curves~$X_0(N)$. By the Modularity Theorem
	of Wiles and others~\cite{Wiles1995,TaylorWiles1995,BCDT},
	this is always the case when $A$ is an elliptic curve over~$\Q$
	(this is now a special case of Serre's Modularity Conjecture~\cite{KhareWintenberger2010}).
	
	We now introduce the relevant global invariants of~$A$.
	By the Mordell--Weil Theorem, the abelian group $A(\Q)$ of rational points
	on~$A$ is finitely generated, so it splits as $A(\Q) \isom A(\Q)_\tors \oplus \Z^r$,
	where $A(\Q)_\tors$ is the finite \emph{torsion subgroup} and $r$ is
	a nonnegative integer, the \emph{rank of $A(\Q)$}. There is a natural
	positive definite quadratic form~$\hat{h}$ on $A(\Q) \otimes_\Z \R \cong \R^r$,
	the \emph{canonical height}, turning $A(\Q)/A(\Q)_\tors$ into a lattice
	in a euclidean vector space. The squared covolume of this lattice
	(equivalently, the determinant of the Gram matrix of~$\hat{h}$ with respect
	to a lattice basis) is the \emph{regulator}~$\Reg_{A/\Q}$. The final global
	arithmetic invariant of~$A$ that we need is the \emph{Tate--Shafarevich group}~$\Sha(A/\Q)$.
	It can be defined as the localization kernel
	\[ \Sha(A/\Q) = \ker\Big(\H^1(\Q,A) \to \bigoplus_{v} \H^1(\Q_v,A)\Big) \]
	in Galois cohomology; here $\Q_v$ denotes the completion of~$\Q$ with respect
	to a place~$v$ and the direct sum is over all places of~$\Q$.
	Geometrically, $\Sha(A/\Q)$ is the group of equivalence classes of
	everywhere locally trivial $A/\Q$-torsors. This group is conjectured to be finite,
	but this is not known in general; for example, it is not known for a single elliptic
	curve with (algebraic or analytic) rank at least~$2$.
	
	We also need some local invariants. To each prime~$p$,
	one associates the \emph{Tamagawa number}~$c_p(A)$; this is the number of
	connected components of the special fiber at~$p$ of the Néron model~$\sA/\Z$
	of~$A$ that are fixed by Frobenius and equals~$1$ for all primes of good reduction.
	Let $(\omega_1, \ldots, \omega_g)$ be the pull-back to~$\H^0(A, \Omega^1)$ of a basis of
	the free $\Z$-module $\H^0(\sA, \Omega^1)$ of rank~$g$. Then the
	\emph{real period} of~$A$ is the volume of~$A(\R)$ measured
	using $|\omega_1 \wedge \dots \wedge \omega_g|$:
	$\Omega_A = \int_{A(\R)} |\omega_1 \wedge \dots \wedge \omega_g|$.
	
	The \emph{weak BSD} or \emph{BSD rank conjecture} says that $L(A/\Q,s)$
	has an analytic continuation to a neighborhood of $s = 1$ and
	\[ r_\an \defeq \ord_{s=1} L(A/\Q, s) = r \,. \]
	The order of vanishing of~$L(A/\Q, s)$ at $s = 1$ is also called the
	\emph{analytic rank} of~$A/\Q$.
	
	We will from now on assume that $A$ is principally polarized,
	for example the Jacobian variety of a curve. In particular, $A \isom A^\dual$,
	where $A^\dual$ is the dual abelian variety. Then the
	\emph{strong BSD conjecture} says that in addition $\Sha(A/\Q)$ is finite and
	\[ L^*(A/\Q,1) \defeq \lim_{s \to 1} (s-1)^{-r} L(A/\Q, s)
	= \frac{\Omega_A \prod_p c_p(A) \cdot \Reg_{A/\Q} \#\Sha(A/\Q)}{(\#A(\Q)_\tors)^2} \,.
	\]
	
	Since all the other invariants of~$A$ can (usually) be computed at least
	numerically, we define the \emph{analytic order of Sha} to be
	\[ \#\Sha(A/\Q)_\an \defeq \frac{L^*(A/\Q,1)}{\Omega_A \Reg_{A/\Q}} \cdot
	\frac{(\#A(\Q)_\tors)^2}{\prod_p c_p(A)} \,.
	\]
	Assuming the BSD rank conjecture, strong BSD can then be phrased as
	\enquote{$\Sha(A/\Q)$ is finite and $\#\Sha(A/\Q) = \#\Sha(A/\Q)_\an$.}
	
	Even the weak BSD conjecture for elliptic curves over~$\Q$ is wide open
	in general (this is the Clay Millennium Problem mentioned above).
	However, the strong BSD conjecture has been verified for many \enquote{small}
	elliptic curves; see below. In this article, we verify the
	strong BSD conjecture for the first time in dimension greater than~$1$,
	namely for a number of abelian \emph{surfaces}~$A/\Q$,
	in a situation where it cannot be reduced to BSD for some elliptic curves.
	Concretely, this means that $A$ is absolutely simple.
	
	Recall that an abelian variety~$A$ of positive dimension over~$\Q$ is
	\emph{absolutely simple} if $A_{\Qbar}$ is not isogenous to a product of
	at least two abelian varieties of positive dimension.
	An abelian variety of dimension~$g$ whose endomorphism ring~$\End_\Q(A)$
	is isomorphic to an order~$\cO$ in a totally real number field~$F$ of degree
	$[F:\Q] = g$ is said to have \emph{real multiplication} (RM).
	
	Absolutely simple abelian varieties with real multiplication over~$\Q$ are \emph{modular}
	(see~\cref{RM abelian varieties are modular} for references).
	This means that $A$ can be obtained as an isogeny factor of some~$J_0(N)$, where
	$J_0(N)$ denotes the Jacobian variety of the modular curve~$X_0(N)$. These
	isogeny factors correspond to (Galois orbits of) newforms in~$S_2(\Gamma_0(N))$;
	see~\cref{RM abelian varieties are modular} below. (Note that we reserve the
	term \enquote{modular} for $\GL_2$-type abelian varieties here as opposed to the
	more general property of being \enquote{automorphic}; see~\cite[§\,9.1]{BCGP2021}.)
	
	\subsection{General results}
	
	While the BSD~conjecture is wide open in general, there are some cases where
	parts of it are known to be true. Assume that $A/\Q$ is an absolutely simple
	abelian variety of dimension~$g$ with real multiplication by an order~$\cO$
	in a totally real number field of degree~$g$. Then $A$ is of $\GL_2$-type;
	in particular, for each prime ideal $\frp$ of~$\cO$, the common kernel~$A[\frp]$
	of all elements of~$\frp$ acting on~$A$ is a $2$-dimensional vector space
	over~$\cO/\frp$, hence induces a Galois representation into $\GL_2(\cO/\frp)$.
	In this situation there is a newform~$f$ of weight~$2$ and some level~$N$
	with $q$-expansion coefficients that generate an order commensurable with~$\cO$ and such that
	$A$ is an isogeny factor of~$J_0(N)$; furthermore,
	\[ L(A/\Q, s) = \prod_{\sigma \colon \cO \inj \R} L(f^\sigma, s) \,, \]
	where $\sigma$ acts on the $q$-expansion coefficients. Since it is known
	that $L(f^\sigma, s)$ is an entire function, the same is true for~$L(A/\Q, s)$.
	So for $A/\Q$ with RM we can at least speak of the analytic rank~$r_\an$ and
	the leading coefficient $L^*(A/\Q,1)$ of the $L$-function at $s = 1$.
	The parity of the order of vanishing of~$L(f^\sigma, s)$ at $s = 1$ does not depend
	on~$\sigma$ (it is determined by the eigenvalue~$\epsilon_N$ of~$f^\sigma$ under
	the Fricke involution, which is the same for all~$f^\sigma$),
	and the order of vanishing itself does not depend on~$\sigma$
	when $\ord_{s=1} L(f^\sigma,s) \le 1$ for some~$\sigma$,
	so in this case we have that $r_\an = g \cdot \ord_{s=1} L(f,s)$,
	where $\ord_{s=1} L(f,s) \defeq \ord_{s=1} L(f^\sigma,s)$ for any~$\sigma$;
	see~\cite[Cor.~V.1.3]{GrossZagier1986}. We call
	$\ord_{s=1} L(f,s)$ the \emph{$L$-rank of~$A$} in this case and abbreviate it as $\Lrk A$.
	If the BSD rank conjecture holds for~$A$, then the $L$-rank of~$A$ is the same as the rank
	of~$A(\Q)$ as an $\cO$-module.
	
	Based on work of Gross--Zagier~\cite{GrossZagier1986}
	relating the canonical height of Heegner points to~$L^{(g)}(A/K,1)$
	for suitable imaginary quadratic fields~$K$, Kolyvagin~\cite{Kolyvagin1988}
	(for modular elliptic curves) and Kolyvagin--Logachëv~\cite{KolyvaginLogachev}
	(for modular abelian varieties in general) were able to show that the BSD rank
	conjecture holds under the assumption that the $L$-rank is $0$ or~$1$,
	that in this case, $\Sha(A/\Q)$ is finite (this is the only case where we know
	finiteness), and that $\#\Sha(A/\Q)_\an$ is a rational number.
	
	The rational number~$\#\Sha(A/\Q)_\an$ can be computed when $A$ is
	an elliptic curve, and we show in this paper how to do that when $A$ is
	a modular abelian surface. To complete the verification of the strong BSD~Conjecture
	for~$A$, it remains to determine $\#\Sha(A/\Q)$ and to check that the two numbers
	agree.
	
	This involves
	showing that $\Sha(A/\Q)[p]$ is trivial for all primes $p \notin S$, where
	$S$ is an explicit finite set of primes, and then determining $\#\Sha(A/\Q)[p^\infty]$ for
	the finitely many $p \in S$. When $A$ is an elliptic curve,
	a suitable set~$S$ (or even an explicit annihilator of~$\Sha(A/\Q)$) can be extracted from Kolyvagin's
	work and subsequent refinements; see below. We show in this paper
	how to obtain a suitable set~$S$ when $A$ is a modular abelian surface.
	
	The remaining task is to determine $\Sha(A/\Q)[p^\infty]$ for a given
	prime~$p$. This is always possible in theory (assuming that $\Sha(A/\Q)[p^\infty]$
	is finite), since one can compute the $p^n$-Selmer group of~$A$
	for $n = 1,2,\ldots$, which is defined as
	\[ \Sel_{p^n}(A/\Q) = \ker\Bigl(\H^1(\Q, A[p^n]) \to \bigoplus_v \H^1(\Q_v, A)\Bigr) \]
	and sits in an exact sequence
	\[ 0 \To A(\Q)/p^n A(\Q) \To \Sel_{p^n}(A/\Q) \To \Sha(A/\Q)[p^n] \To 0 \,. \]
	Since we know $A(\Q)$, this gives us~$\Sha(A/\Q)[p^n]$,
	and as soon as $\Sha(A/\Q)[p^n] = \Sha(A/\Q)[p^{n+1}]$, we have determined
	$\Sha(A/\Q)[p^\infty] = \Sha(A/\Q)[p^n]$ (and if $\Sha(A/\Q)[p] = 0$, then
	$\Sha(A/\Q)[p^\infty] = 0$ as well). For the computability of the
	Selmer group in theory see~\cite{Stoll2012} for elliptic curves
	and~\cite{BPS} in general. In practice, there are fairly tight
	limits on~$p^n$, since the computation requires the knowledge of the
	class and unit groups of number fields of degree growing quickly with~$p^n$,
	for which no really efficient algorithms are available so far.
	
	If one has a conjecturally tight bound on~$\#\Sha(A/\Q)[p^\infty]$,
	then another approach is to try and get a lower bound that agrees
	with the upper bound. If the upper bound is nontrivial, this involves
	showing the existence of nontrivial elements of~$\Sha(A/\Q)$ in some way.
	One possibility for this is
	\enquote{visibility},
	which uses another related abelian variety~$B$, for which one can construct
	a nontrivial map $B(\Q) \to \H^1(\Q, A)$, whose image one can show
	to contain nontrivial elements of~$\Sha(A/\Q)$ under suitable conditions.
	This is used for the example in Appendix~\ref{sec:7-torsion in Sha}.
	
	We now give a short overview of what has been done so far regarding
	the verification of the strong BSD~Conjecture in concrete cases.
	
	\subsection{Exact verification of strong BSD for elliptic curves}
	
	In the case of elliptic curves, the various ingredients mentioned above
	have been worked out, made explicit and been improved to an extent that
	it was possible to verify the strong BSD conjecture for all elliptic
	curves~$E$ over~$\Q$ of rank $\le 1$ and conductor $N < 5000$;
	see~\cite{GJPST,Miller2011,MillerStoll2013,CreutzMiller2012,LawsonWuthrich2016}.
	
	An explicit finite set~$S$ of primes such that $\Sha(E/\Q)[p] = 0$ for $p \notin S$
	can be obtained using Kolyvagin's work and refinements building on it~\cite{Jetchev2008,Cha2005}.
	The size of~$\Sha(E/\Q)[p^\infty]$ for $p \in S$ can be obtained by several methods,
	e.g., using Iwasawa theory and $p$-adic $L$-functions~\cite{SteinWuthrich2013} or by
	performing descents~\cite{Cremona1997,SchaeferStoll,CFOSS2008,CFOSS2009,CreutzMiller2012,MillerStoll2013,CFOSS2015,Stoll2012,Creutz2014,BPS}.
	
	\subsection{Numerical verification of strong BSD for higher-dimensional abelian varieties}
	
	Compared to the case of elliptic curves, considerably less has been done
	regarding the verification of the BSD conjectures for higher-dimensional
	abelian varieties~$A$ over~$\Q$. If $A$ is not absolutely simple,
	then $A$ splits up to isogeny (and possibly after base-change to
	an algebraic number field) as a product of abelian varieties of
	lower dimension. Since the validity of strong BSD is invariant
	under isogenies~\cite{Tate1968} and Weil restriction~\cite{Milne1972},
	this reduces the verification of
	strong BSD for~$A$ to cases of lower dimension. We will therefore
	assume that $A$ is absolutely simple in the following.
	
	In~\cite{FLSSSW}, all factors in the BSD formula except for the order of the
	Tate--Shafarevich group (only its $2$-torsion is computed)
	and the analytic order of Sha are determined exactly in the $L$-rank~$0$
	cases and numerically to high precision in the $L$-rank~$1$ cases
	for the Jacobian varieties of $29$ genus~$2$ curves over~$\Q$ such that the
	Jacobians are absolutely simple, of $\GL_2$-type and have level $N \le 200$.
	Their work also includes results on three curves whose Jacobians
	are Weil restrictions of elliptic curves over~$\Q(\sqrt{-3})$.
	
	More recently, van Bommel~\cite{vanBommel2019} has done computations
	similar to those in~\cite{FLSSSW}
	for various (in general non-modular) Jacobians of hyperelliptic curves of
	genus $\leq 5$. He did not provably compute the
	regulator or the torsion subgroup, which means that the approximate value
	of~$\#\Sha(A)_\an$ that he computes may be off by a square rational factor.
	Van Bommel also provides an algorithm for the computation of the real period,
	which corrects the version described in~\cite{FLSSSW} in the case
	when some of the special fibers of a minimal regular model of the curve
	have multiple components.
	
	However, it was still an open problem to provably compute
	$\#\Sha(A/\Q)_\an$ as an exact rational number when $A$ has positive rank
	and to determine $\#\Sha(A/\Q)$.
	See for example William Stein's blog post~\cite{SteinBlog} for the former.
	
	We now verify the strong BSD~Conjecture \emph{unconditionally} and exactly for
	all the curves in~\cite{FLSSSW} with absolutely simple Jacobian
	and all genus~$2$ curves in the LMFDB~\cite{lmfdb} with absolutely simple modular Jacobian.
	
	\subsection{New general results in this paper}
	
	We note that, compared to elliptic curves, a number of additional
	difficulties show up when trying to verify strong BSD for higher-dimensional
	modular abelian varieties. By the Modularity Theorem, every elliptic curve~$E$
	over~$\Q$ of conductor~$N$ is the target of a nontrivial morphism $X_0(N) \to E$.
	This makes it fairly easy to compute Heegner points on~$E$. Also, elliptic
	curves are given explicitly by a Weierstrass equation, and a variety of
	algorithms are available for them. There is no comparable explicit representation
	of a general (modular) abelian variety of higher dimension. We can, however,
	work with curves~$X$ and their Jacobians. In particular for hyperelliptic curves,
	a variety of algorithms exist. However, in general there is only a dominant
	homomorphism from~$J_0(N)$ to the Jacobian in question and no nontrivial
	morphism from $X_0(N)$ to the curve~$X$. When there is such a morphism,
	the relevant computations are much simpler; we have dealt with this case
	for surfaces first, and the results are described in~\cite{KellerStoll2022}.
	In the other cases, the required arguments are much more subtle; for example,
	it is quite nontrivial to obtain a formula for the canonical height of a
	Heegner point on the Jacobian of~$X$ from the Gross--Zagier formula.
	
	In this paper, we overcome these difficulties and devise general methods
	to verify strong BSD exactly for absolutely simple modular
	Jacobians $J/\Q$ of $L$-rank $0$ and $1$ and apply them to several examples.
	We denote the curve $J$~is the Jacobian of by~$X$.
	Many of our results and algorithms apply to any dimension or at least to hyperelliptic
	Jacobians. We note that an abelian variety (assumed to be absolutely simple) is automatically
	a Jacobian when it is principally polarized and its dimension is $2$ or~$3$.
	
	More specifically, given such a Jacobian~$J$ and/or an attached newform~$f$, we
	do the following. (The numbers link to the corresponding sections.)
	\begin{enumerate}[(2)]
		\item[\eqref{sec:computation-of-the-residual-galois-representations}]
		We determine the (projective) images of the associated mod-$\frp$
		Galois representations for all maximal ideals~$\frp$ of the
		endomorphism ring; in particular, we determine which of them are reducible.
		\item[\eqref{sec:computation-of-the-heegner-points-and-index}]
		We develop an efficient algorithm for the computation of Heegner points, their
		canonical heights, and Heegner indices. This involves the
		computation of Petersson norms of newforms of weight~$2$. We also
		provide the refined information of the Heegner index as a characteristic
		ideal of the endomorphism ring.
		\item[\eqref{sec:Sha_an}]
		We derive explicit formulas for $\#\Sha(J/\Q)_\an$ and $\#\Sha(J/K)_\an$
		(where $K$ is a Heegner field). When the $L$-rank is~$1$, this involves
		some fairly nontrivial arguments.
		\item[\eqref{sec:finite support}]
		We give an explicit upper bound for the set of primes dividing~$\#\Sha(J/\Q)$
		and for the primes dividing~$\#\Sha(J/K)$, where $K$ is a Heegner field.
		\item[\eqref{sec:descent}]
		We perform $\frp$-isogeny descents in some cases where the mod-$\frp$ Galois
		representation is reducible to get an upper bound on the $\frp$-Selmer group
		and thus show that $\Sha(J/\Q)[\frp] = 0$.
		\item[\eqref{sec:Iwasawa theory}]
		We provide a feasible algorithm for the computation of $p$-adic $L$-functions
		in our setting.
		\item[\eqref{sec:examples}]
		We provide an algorithm that, combining the above algorithms,
		verifies strong BSD for the Jacobian $J$ (absolutely simple and modular)
		of a given genus $2$ curve $C$ of level $N$ or returns
		at least a small finite set of primes $\frp$ for which $\Sha(J/\Q)[\frp^\infty]$
		needs to be computed to finish the verification.
	\end{enumerate}
	We also improve van Bommel's algorithm for the determination of the real period
	so that it does not rely on a gcd~computation with real numbers;
	see~\cref{computation of Omega_J}.
	
	Our methods and algorithms generalize to RM abelian varieties over~$\Q$
	of arbitrary dimension, provided one can compute Mordell--Weil groups and canonical heights and
	one has an analogue of~\cref{exclude sub-line}. In joint work in progress
	with Pip Goodman and John Voight, we are planning to treat modular abelian varieties
	over totally real number fields.
	
	We then use the algorithms we developed to verify strong~BSD exactly
	for the first time for a number of absolutely simple abelian surfaces.
	The specific examples are described below.

	\subsection{Examples} \label{S:data}
	
	All computations were carried out with Magma~\cite{Magma}.
	The code to reproduce our computations can be found at
	\begin{center}
		\url{https://github.com/TimoKellerMath/strongBSDgenus2}.
	\end{center}
	The \texttt{README.md} file contains short descriptions of the Magma files and
	which sections in this paper they belong to.
	
	Using the methods and algorithms developed in this paper,
	we verified strong BSD completely for the Jacobians of the following genus~$2$
	curves over $\Q$ (whose Jacobians are modular and absolutely simple).
	\begin{enumerate}[(a)]
		\item The \href{https://www.lmfdb.org/}{LMFDB}~\cite{lmfdb} currently (as of September~2024) lists exactly~$97$ genus~$2$ curves
		with absolutely simple Jacobian of $\GL_2$-type; they all have level $\leq 1000$.
		By their completeness statement, this comprises all such examples with absolute value
		of their discriminant at most~$10^6$ and \enquote{small} coefficients.
		We will refer to these as the \emph{LMFDB examples}.
		Note that there are more newforms of weight $2$
		with real quadratic coefficients of level $\leq 1000$ contained in the LMFDB;
		our algorithms would at least give an upper bound on the size of the
		Tate--Shafarevich group of their associated modular abelian variety
		given a Jacobian in their isogeny class.
		Some of the examples mentioned below provide such a Jacobian
		for additional newforms.
		\item The $28$ \enquote{Hasegawa curves} from~\cite{Hasegawa1995} that
		have absolutely simple Jacobian.
		These are all quotients of~$X_0(N)$ by a subgroup of Atkin--Lehner involutions.
		Because of this, these examples are easier to deal with
		(compare~\cref{same heights}, which shows that the computation of
		Heegner points is simpler in this case), which is
		why we treated them first, before extending the theory and algorithms
		to the general case.
		See~\cite{KellerStoll2022} for an overview of the results.
		Note that $X_0(161)/\langle{w_7, w_{23}}\rangle$ is the only curve on this list
		whose Jacobian is not isogenous to the Jacobian of one of the LMFDB examples.
		(We check this by comparing the associated newforms.)
		Hence strong BSD for the other~$27$ Hasegawa curves
		follows from isogeny invariance and the validity of strong~BSD for the LMFDB examples.
		\item The four \enquote{Wang curves} from~\cite{FLSSSW} that are neither
		Hasegawa curves nor have Jacobian isogenous to that of a curve
		in the LMFDB. They are the curves labeled
		65A, 117B, 125B and~175 in~\cite{FLSSSW}.
		\item Sam Frengley's example of a curve with $N = 3200$ and $\#\Sha(J/\Q) = 7^2$.
	\end{enumerate}
	Note that there is some overlap between the first two sets:
	$21$ of the Hasegawa curves are in the LMFDB.
	In total, the LMFDB, Hasegawa, and Wang examples comprise the
	Jacobians of $108$ isomorphism classes of curves,
	whose Jacobians fall into $95$ distinct isogeny classes. Including the last example,
	we therefore have verified the strong BSD conjecture completely for $96$
	isogeny classes of absolutely simple modular abelian surfaces.
	The distribution of the $L$-ranks in our examples is as follows:
	\begin{enumerate}[(a)]
		\item There are $36$~isomorphism classes of $L$-rank~$0$ and $61$ of $L$-rank $1$.
		These belong to $31$ and~$59$ isogeny classes, respectively.
		\item There are $6$~isomorphism classes of $L$-rank~$0$ and $22$ of $L$-rank~$1$.
		Out of the $7$~Hasegawa non-LMFDB examples, there are
		$4$~isomorphism classes of $L$-rank~$0$ and $3$ of $L$-rank~$1$.
		The latter include the single isogeny class not represented by LMFDB curves.
		All these examples belong to distinct isogeny classes. The abundance of $L$-rank~$1$ examples
		is explained by the fact that one often quotients out by the Fricke involution~$w_N$,
		hence the sign in the functional equation shows that the $L$-rank is~odd.
		\item All $4$~isogeny classes have $L$-rank~$0$.
	\end{enumerate}
	
	In total, there are $44$~isomorphism classes and $35$~isogeny classes
	of $L$-rank~$0$ and $64$~isomorphism classes and $60$~isogeny classes of $L$-rank~$1$.
	
	\subsection*{Completeness of our data}
	
	We consider all the Hasegawa examples (already contained in~\cite{KellerStoll2022}) and
	all Wang-only examples from~\cite{FLSSSW}. The LMFDB examples comprise all absolutely simple
	modular genus~$2$ Jacobians of curves with \enquote{small} coefficients and of level $N \le 1000$.
	The smallest level for which there exists a pair of conjugate newforms that is not related to
	one of these examples is \href{https://www.lmfdb.org/ModularForm/GL2/Q/holomorphic/?weight=2&char_order=1&dim=2}{$43$}.
	There is an \href{https://github.com/edgarcosta/ModularAbelianSurfaces}{ongoing project} that
	attempts to produce, for a given weight~$2$ newform~$f$ with real quadratic coefficients,
	a genus-$2$ curve over~$\Q$ with Jacobian isogenous to $A_f$. Our code can prove strong~BSD
	for many of these examples automatically, especially by using~\cite{KellerYin2024}
	to skip many descent computations.
	
	\subsection{Structure of the paper}
	
	We give an overview of the paper; more details are given at the beginning of each section.
	In Section~\ref{sec:computation-of-the-residual-galois-representations}, we give algorithms to
	determine whether the residual Galois representations attached to $f$ are irreducible or not.
	In Section~\ref{sec:computation-of-the-heegner-points-and-index}, we compute (a multiple of)
	the Heegner index, which is used in the following two sections:
	In Section~\ref{sec:Sha_an}, we compute $\#\Sha(J/\Q)_\an \in \Q_{> 0}$ exactly.
	In Section~\ref{sec:finite support}, we give a description of a finite set~$S$ of prime ideals~$\frp$
	such that $\Sha(J/\Q)[\frp] = 0$ for $\frp \notin S$; this strongly depends on the
	determination of the residual Galois representations and the Heegner index.
	In Section~\ref{sec:descent}, we perform isogeny descents to prove $\Sha(J/\Q)[\frp] = 0$
	for several $\frp \in S$.
	In Section~\ref{sec:Iwasawa theory}, we show how results from Iwasawa theory and the
	computation of $p$-adic $L$-functions can be used to prove an upper bound on~$\#\Sha(J/\Q)[\frp^\infty]$.
	In Appendix~\ref{sec:7-torsion in Sha}, we prove strong~BSD for an example of Sam Frengley,
	where $\#\Sha(J/\Q) = 7^2$. In all our other examples, $\#\Sha(J/\Q) \in \{1,2,4\}$.
	We also exhibit examples $J/\Q$ for which $p^2 \mid \Sha(J/\Q)_\an$ with
	$p \in \{3,5,7\}$ and prove the $\ell$-part strong BSD for them except for $\ell \in \{2,p\}$,
	where we only get an upper bound.
	These examples are obtained as quadratic twists~$J^K$ of some of our
	main examples, where $K$ is a suitable Heegner field.

	\subsection{Terms and notation}
	
	We denote canonical isomorphisms by~$\iso$ and arbitrary, not necessarily canonical
	isomorphisms by~$\isom$. We fix an embedding $\Qbar \inj \C$ once and for all.
	
	We use boldface $\bpi$ to denote the area of the unit disk to avoid confusion with
	our use of the letter~$\pi$ to denote an isogeny in most of the paper.
	
	
	\section{Computation of the residual Galois representations}
	\label{sec:computation-of-the-residual-galois-representations}
	
	The purpose of this section is to generalize several results on the image of mod-$p$
	Galois representations of elliptic curves over~$\Q$ (mainly from~\cite{Serre1972}
	and~\cite{Cojocaru2005}) to modular abelian varieties over~$\Q$ of higher dimension.
	
	Let $A$ be an abelian variety of dimension $g \ge 1$. Let $\cO$ be an order in a totally
	real number field~$F$ of degree~$g$ over~$\Q$.
	Recall that $A$ has \emph{real multiplication
		by~$\cO$ over $\Q$} if  $\End_\Q(A) \isom \cO$, where $\End_\Q(A)$ denotes the ring
	of $\Q$-defined endomorphisms of~$A$. Then $A$ is of $\GL_2$-type in the
	following sense.
	Let $\frp$ be a non-zero prime ideal of~$\cO$ lying above a rational prime~$p$.
	We denote its finite residue field~$\cO/\frp$ by~$\F_\frp$ and call $[\F_\frp : \F_p]$
	the \emph{degree~$\deg\frp$} of~$\frp$; $\F_\frp$ is isomorphic to~$\F_{p^{\deg\frp}}$.
	If $\frp$ is \emph{regular}, i.e., its local ring is a discrete valuation ring, or
	equivalently, $\frp$ does not divide the conductor ideal
	$\frf(\cO_F/\cO) = \{a \in \cO : a \cO_F \subseteq \cO\}$ of~$\cO$ in~$\cO_F$, then
	$A[\frp^n](\Qbar)$ is free of rank~$2$ over~$\cO/\frp^n$ for all $n \ge 1$.
	
	We then obtain $2$-dimensional Galois representations
	\[ \rho_{\frp^n,A} \colon \GalQ \to \Aut_{\cO/\frp^n}(A[\frp^n](\Qbar)) \isom \GL_2(\cO/\frp^n) \,. \]
	In a similar way, we have $2g$-dimensional Galois representations
	\[ \rho_{p^n,A} \colon \GalQ \to \Aut_{\Z/p^n\Z}(A[p^n](\Qbar)) \isom \GL_{2g}(\Z/p^n\Z) \,. \]
	Since the Galois action preserves the Weil pairing, the image of~$\rho_{p,A}$ lies
	in the general symplectic group $\GSp_{2g}(\F_p)$.
	
	We define the \emph{$\frp$-adic Tate module $T_\frp A \defeq \varprojlim_n A[\frp^n](\Qbar)$};
	it is free of rank~$2$ over the completion~$\cO_\frp$.
	We also define $V_\frp A = F_\frp \otimes_{\cO_\frp} T_\frp A$; this is a $2$-dimensional
	vector space over~$F_\frp$. There is also the standard $p$-adic Tate-module~$T_p A$,
	which is a free module of rank~$2g$ over~$\Z_p$
	and the associated vector space~$V_p A = \Q_p \otimes_{\Z_p} T_p A$.
	We obtain the \emph{$\frp$-adic Galois representation}
	\[ \rho_{\frp^\infty,A} \colon \GalQ \to \Aut_{F_\frp}(V_\frp A) \isom \GL_2(F_\frp) \]
	and the \emph{$p$-adic Galois representation}
	\[ \rho_{p^\infty,A} \colon \GalQ \to \Aut_{\Q_p}(V_p A) \isom \GL_{2g}(\Q_p) \,. \]
	As before, the image of~$\rho_{p^\infty,A}$ is contained in~$\GSp_{2g}(\Q_p)$.
	
	If $A$ is understood, we omit it from the notation and write $\rho_\frp$ etc.
	
	We heavily exploit that we can work with $2$-dimensional Galois representations instead
	of $2g$-dimensional ones in the following. For example, there is an easy classification
	of (maximal) subgroups of~$\GL_2(\F_\frp)$, whereas the
	subgroups of~$\GSp_{2g}(\F_p)$ are more complicated.
	
	The goal of this section is to determine the image $G_\frp$ of the mod-$\frp$
	Galois representation
	\[ \rho_\frp \colon \GalQ \to \GL_{\cO/\frp}(A[\frp](\Qbar)) \isom \GL_2(\F_\frp) \]
	in the case when $g = 2$, so $\cO$ is an order in a real quadratic number field.
	In particular, we want to decide whether $\rho_\frp$ is irreducible
	as an $\F_\frp[\GalQ]$- or $\F_p[\GalQ]$-representation and whether its image
	in~$\GL_2(\F_\frp)$ is as large as allowed by the extra endomorphisms coming from~$\cO$.
	
	We will state our results for general~$g$ if this is easily possible,
	but in some cases we assume $g = 2$ to simplify the statements and algorithms.
	
	Let $f \in S_2(\Gamma_0(N))$ be a newform, i.e., a normalized eigenform for
	the action of the Hecke algebra~$\T_\Z$ on the new subspace of~$S_2(\Gamma_0(N))$.
	The Fourier coefficients of~$f$ generate an order~$\Z[f]$ in a totally real
	number field~$\Q(f)$. Let $I_f \defeq \Ann_{\T_\Z}(f)$ be the annihilator of~$f$;
	then $\T_\Z/I_f \iso \Z[f]$, where the Hecke operator~$T_n$ is mapped to
	the Fourier coefficient~$a_n(f)$.
	The Hecke algebra also acts via $\Q$-defined endomorphisms on~$J_0(N)$,
	and so we can define an abelian variety $A_f$ over~$\Q$ as $A_f \defeq J_0(N)/I_f J_0(N)$.
	Then $\dim A_f = [\Q(f) : \Q]$ and $\End_\Q(A_f) \iso \T_\Z/I_f \iso \Z[f]$.
	Acting by~$\GalQ$ on the Fourier coefficients of~$f$, we obtain a Galois orbit
	of conjugate newforms $f^\sigma$, which has size $[\Q(f) : \Q]$.
	More generally, if $\alpha \colon \Z[f] \inj \R$ is an
	embedding of~$\Z[f]$ into~$\R$, then $f^\alpha$
	denotes the newform with (real) coefficients $\alpha(a_n(f))$.
	The abelian variety~$A_f$ only depends on the Galois orbit of~$f$.
	
	We give a short summary of the contents of this section. After recalling basic results about modular
	abelian varieties and their Galois representations in Section~\ref{ssec:preliminaries},
	we determine the maximal possible image of~$\rho_\frp$ in
	Section~\ref{sec:maximal image} and state the classification of its maximal subgroups
	in Section~\ref{sec:classification-of-the-maximal-subgroups-of-psl2pfr}.
	This is eventually used to show that $\rho_\frp$ has maximal image
	for all~$\frp$ outside an explicit finite set
	by excluding the possibility that the image is contained
	in one of the maximal subgroups. For fixed~$\frp$, we give an algorithm
	that returns a set of types of maximal subgroups that could contain the
	image of~$\rho_\frp$ in Section~\ref{sec:criteria for fixed frp}.
	In Section~\ref{ssec:rho_p explicit} we show how $\rho_\frp$
	can be determined explicitly for a given prime ideal~$\frp$.
	We then give some criteria for when the image of the decomposition group at~$p$
	is contained in a Cartan subgroup in Section~\ref{sec:determination-of-the-image-of-ip}. Together with some
	results on the image of inertia at primes $\ell \neq p$, which we recall
	in Section~\ref{sec:image of inertia not p}, this provides the input for
	an algorithm that determines a (small and explicit) finite set~$S$ of prime
	ideals~$\frp$ such that $\rho_\frp$ is irreducible for all $\frp \notin S$
	in Section~\ref{sec:irreducibility-for-almost-all-pfr}.
	To approach the goal of determining an analogous set with respect to
	$\rho_\frp$ with maximal image, we first describe a method that allows
	us to eliminate two further types of maximal subgroups (other than
	Borel subgroups, which correspond to reducible representations) in
	Section~\ref{sec:excluding-L-Ex}. To deal with maximal images, we need
	to exclude that the given newform~$f$ has complex multiplication, so we
	provide an algorithm that checks that in Section~\ref{sec:proving-non-cm}.
	We then derive an algorithm that computes a small explicit finite set
	of prime ideals~$\frp$ such that $\rho_\frp$ has maximal image for all~$\frp$
	not in this set in Section~\ref{sec:maximal-image-for-almost-all-pfr}.
	Finally, we provide a table giving the types of all representations~$\rho_\frp$
	attached to our LMFDB examples.

	\subsection{Preliminaries} \label{ssec:preliminaries}
	
	We begin by stating the correspondence between (absolutely simple) abelian
	varieties with real multiplication over~$\Q$ and weight-$2$ newforms for~$\Gamma_0(N)$.
	
	Recall that $L(A/\Q, s)$ denotes the $L$-series of~$A$ and $L(f, s)$ denotes the
	$L$-series of~$f$
	and that $L(A/\Q,s)$ is defined as
	\begin{align*}
		L(A/\Q,s) &= \prod_p\frac{1}{\det(1 - \Frob_p^{-1} p^{-s} \mid \Het^1(A \otimes \Qbar, \Q_\ell)^{I_p})},
	\end{align*}
	where for each Euler factor at $p$ one chooses a prime $\ell \neq p$ for the
	$\ell$-adic cohomology group; this is well-defined because the Euler factors
	are independent of~$\ell$.
	The product converges for $\Re(s) > \frac{3}{2}$ to a holomorphic function.
	The $L$-function associated to $f = \sum_n a_n q^n$ with coefficients $a_n \in \C$ is
	\begin{align*}
		L(f,s) &= \sum_{n \ge 1}\frac{a_n}{n^s} = \prod_p\frac{1}{1 - a_pp^{-s} + \epsilon(p)p^{1-2s}},
	\end{align*}
	where $\epsilon(p)$ is $1$ if $p \nmid N$ and $0$ otherwise. Since $L(f,s)$ is the Mellin
	transform of~$f$ and $f$ is a cusp form, $L(f, s)$ is holomorphic on the whole complex plane.
	
	\begin{theorem}[Characterization of modular abelian varieties over $\Q$]
		\label[theorem]{RM abelian varieties are modular}
		Let $A/\Q$ be an absolutely simple abelian variety. The following are equivalent.
		\begin{enumerate}[\upshape(i)]
			\item \label{item:RM}
			$A$ has real multiplication over~$\Q$.
			\item \label{item:isogeny factor}
			There is some~$N$ such that $A$ is an isogeny factor of~$J_0(N)$.
			\item \label{item:L-series}
			There is some~$N$ and a newform $f \in S_2(\Gamma_0(N))$ such that
			\[ L(A/\Q,s) = \prod_{\alpha \colon \Z[f] \inj \R} L(f^\alpha,s) \,. \]
		\end{enumerate}
		The number~$N$ in~\eqref{item:L-series} is uniquely determined; we call it the
		\emph{level~$N_A$} of~$A/\Q$. The statement in~\eqref{item:isogeny factor}
		holds for the same~$N$ and its multiples.
		If these equivalences hold, then
		$\End_\Q^0(A) \defeq \End_\Q(A) \otimes_\Z \Q$ is isomorphic to~$\Q(f)$
		and the conductor of~$A/\Q$ equals $N_A^{\dim{A}}$.
		Furthermore, $A$ is of \emph{$\GL_2$-type}, i.e.,
		the $p$-adic Tate modules $V_p A$ are free modules of rank~$2$ over the
		completion of~$\Q(f)$ at~$p$; if $\frp$ is a prime ideal of~$\Q(f)$ above
		the rational prime~$p$, then $V_\frp A$ is a vector space of dimension~$2$
		over the completion~$\Q(f)_\frp$, a local field.
	\end{theorem}
	
	\begin{proof}
		The equivalence of~\eqref{item:isogeny factor} and~\eqref{item:L-series} is a well-known
		characterization of modular abelian varieties following from the Eichler--Shimura relation
		and Faltings' Isogeny Theorem.
		The equivalence of~\eqref{item:RM} and~\eqref{item:isogeny factor} can be found
		as~\cite[Thm.~5]{Serre1987} as a consequence of Serre's Modularity Conjecture
		for absolutely simple $2$-dimensional residual odd Galois representations (which is
		formulated in the same paper); this conjecture is now a theorem~\cite{KhareWintenberger2010}.
		See also~\cite{Ribet2004}. The remaining statements are well-known.
	\end{proof}
	
	We will use these equivalences tacitly. Note that in the literature,
	sometimes more general modular abelian varieties are considered,
	which are quotients of~$J_1(N)$ and which can have complex multiplication.
	
	However, when $A$ is an \emph{absolutely simple} abelian surface with CM,
	then $A$ cannot be of $\GL_2$-type over~$\Q$, as the following result shows.
	We thank Pip Goodman for pointing it out to us.
	
	\begin{proposition} \label[proposition]{not CM}
		Let $A/\Q$ be an absolutely simple abelian surface with CM.
		Then $\End_\Q^0(A) =\Q$; in particular $A/\Q$ is not of $\GL_2$-type.
	\end{proposition}
	
	\begin{proof}
		Let $E = \End_{\Qbar}^0(A) = \End_{\Qbar}(A) \otimes_\Z \Q$
		be the geometric endomorphism algebra of $A$.
		By~\cite[Proposition 30]{ShimuraAbVarCM}, the minimal field over which the endomorphisms of $A$ are defined is $E^*$, the reflex field of $E$ (note that the base field is just $\Q$ here).
		In particular, $E^*|\Q$ is Galois,
		and the absolute Galois group~$\GalQ$ acts on~$\End_{\Qbar}^0(A)$ through~$\Gal(E^*|\Q)$,
		so we obtain an embedding
		\[ i \colon \Gal(E^*|\Q) \hookrightarrow \Aut(\End_{\Qbar}^0(A)) \,. \]
		We now consult Examples~8.4\,(2) in \emph{loc.~cit.}
		In Example~(C) there the reflex field is not Galois, and in Example~(A)
		the CM-type is not primitive, which means that $A$ is not absolutely simple.
		So both these cases cannot occur, and by Example~(B) it follows that $E^* = E$;
		in particular, the map~$i$ above is an isomorphism. This finally implies that
		\[ \End^0_\Q(A) = \End_{\Qbar}^0(A)^{\Gal(E^*/\Q)} = E^{\Gal(E/\Q)} = \Q \,. \qedhere \]
	\end{proof}
	
	\begin{remark} \label[remark]{RM ab var}
		In the situation of \cref{RM abelian varieties are modular}, $A$ is isogenous to~$A_f$
		(by Faltings' Isogeny Theorem) and therefore
		the Galois representations on~$V_\frp A$ and on~$V_\frp A_f$ are isomorphic
		(similarly for $V_p A$ and~$V_p A_f$). When $A$ and/or~$f$ are clear from the context,
		we write $\rho_{\frp^\infty}$ and $\rho_{p^\infty}$ for these representations,
		which depend only on the Galois orbit of~$f$.
		The fact that $A$ and~$A_f$ are isogenous also implies that the semi-simplifications
		of~$\rho_{\frp,A}$ and of~$\rho_{\frp,f} \defeq \rho_{\frp,A_f}$ are isomorphic
		when $\frp$ is a regular prime of both $\Z[f] \iso \End_{\Q}(A_f)$ and~$\End_{\Q}(A)$
		(and similarly for $\rho_{p,A}$ and~$\rho_{p,f}$).
		
		Note that the canonical isomorphism $\Z[f] \iso \T_\Z/I_f \iso \End_{\Q}(A_f)$
		induces a canonical identification of~$\Q(f)$ with~$\End^0_{\Q}(A_f)$, which in
		turn is isomorphic to~$\End^0_{\Q}(A)$ via the isogeny between $A$ and~$A_f$.
		Fixing the isogeny, this identifies $\cO = \End_{\Q}(A)$ with an order in the
		totally real number field~$\Q(f)$. If $\Z[f]$ is contained in~$\cO$ under this
		identification (e.g., when $\cO$ is the maximal order), then the Fourier
		coefficient~$a_n$ of~$f$, which is the image
		of the Hecke operator $T_n \in \T_{\Z}$ in~$\Z[f]$, can be interpreted as an
		element of~$\cO$, i.e., an endomorphism of~$A$.
		We make use of this to get the correct identifications of $\sigma$-isotypic
		components when dealing with the Gross--Zagier formula for the height of a Heegner point in
		Section~\ref{sec:computing-the-height-of-a-heegner-point-using-gross-zagier}.
		
		Write $F = \Q(f)$ and $\cO_F$ for the maximal order of~$F$, and let $\cO \subseteq \cO_F$
		be any order in~$F$. Recall
		the conductor ideal of~$\cO$ in~$\cO_F$,
		\[ \frf(\cO_F/\cO) = \{a \in \cO : a \cO_F \subseteq \cO\} \,; \]
		it is the largest ideal of~$\cO$ that is also an ideal of~$\cO_F$.
		If $A$ is an abelian variety such that $\End_\Q(A) \isom \cO$, then one can
		check that the isogenous abelian variety $A' \defeq A/A[\frf(\cO_F/\cO)]$
		has $\End_\Q(A') \isom \cO_F$. So by working with~$A'$ instead of with~$A$
		(or~$A_f$), we can assume that the endomorphism ring is the maximal order.
	\end{remark}
	
	\begin{definition}
		Let $p$ be a prime. We write
		\[ \chi_{p^n} \colon \GalQ \to \Aut(\mu_{p^n}(\Qbar)) \iso (\Z/p^n\Z)^\times \]
		for the mod-$p^n$ cyclotomic character and
		\[ \chi_{p^\infty} \colon \GalQ \to \Aut(\mu_{p^\infty}(\Qbar)) \iso \Z_p^\times \]
		for the $p$-adic cyclotomic character.
	\end{definition}
	
	\begin{definition}
		If $\frp$ is a maximal ideal in an order~$\cO$ of a number field, we write $p(\frp)$
		for the characteristic of the finite field~$\F_\frp = \cO/\frp$.
	\end{definition}
	
	\begin{theorem}[Characteristic polynomials of Frobenii of a modular Galois representation]
		\label[theorem]{characteristic polynomials of rho_pfr}
		Let $f \in S_2(\Gamma_0(N))$ be a newform with Fourier coefficients~$a_\ell$
		and coefficient field~$\Q(f)$ a totally real field of degree~$g$.
		
		Associated to $f$,
		there is a strictly compatible system of $\frp$-adic Galois representations~$\rho_{\frp^\infty}$, unramified outside~$N p(\frp)$. For all $\ell \nmid N p(\frp)$
		the characteristic polynomial of~$\rho_{\frp^\infty}(\Frob_\ell)$ equals
		\[ \charpol(f, \ell; T) \defeq \det\bigl(T - \rho_{\frp^\infty}(\Frob_\ell)\bigr)
		= T^2 - a_\ell T + \ell \in \Z[f][T] \,.
		\]
		One has
		\[ \det \circ \rho_{\frp^\infty} = \chi_{p^\infty} \]
		for $\frp \nmid N$.
		In particular, $\rho_{\frp^\infty}$ is odd.
		The determinant of the $p$-adic Galois representation
		\[ \rho_{p^\infty} \colon \GalQ \to \GL_{2g}(\Z_p) \]
		is $\chi_{p^\infty}^g$.
	\end{theorem}
	
	\begin{proof}
		This is well-known and shown for more generally for weight $k \geq 2$ in~\cite{Deligne1971} (and for weight $2$ earlier by Shimura).
	\end{proof}
	
	If $\frp$ is a regular prime of~$\Z[f]$ (or~$\cO$) not dividing $N\ell$,
	then we write $\charpol(f, \ell, \frp; T)$ for the characteristic polynomial
	of the image of~$\Frob_\ell$ under~$\rho_{\frp,f}$ (or $\rho_{\frp,A}$);
	it is the image of~$\charpol(f, \ell; T)$ in~$\F_\frp[T]$.
	
	Magma can compute the Fourier coefficients~$a_\ell$ of a newform~$f$ and its coefficient
	ring~$\Z[f]$ efficiently. This will be crucial for computing the image
	of~$\rho_\frp$, because the only access to elements of the absolute Galois group of~$\Q$
	we have is via Frobenius elements, and we can reconstruct the characteristic
	polynomials of the Frobenii acting on~$A[\frp](\Q)$, which uniquely determine
	their semi-simple part. The fact that we know only the characteristic polynomials
	also means that we do not have direct access to
	the unipotent part of~$\rho_p(\Frob_\ell)$ via~$a_\ell$ alone.
	
	Let $p$ be a prime. We fix an embedding of~$\Qbar$ into~$\Qbar_p$; this
	determines a decomposition group $D_p = \GalQp \inj \GalQ$ and its inertia
	subgroup $I_p = \Gal(\Qbar_p|\Q_p^\nr)$. The inertia subgroup has a descending
	filtration by its (normal) higher ramification subgroups, the first of which
	is the wild ramification subgroup~$I_p^\wild$, the unique (hence normal) Sylow pro-$p$
	subgroup of~$I_p$. The quotient~$I_p/I_p^\wild$ is the tame inertia group~$I_p^\tame$,
	which is canonically isomorphic to the pro-cyclic $p'$-group
	\[ \varprojlim_{\Nm_{\F|}} \F^\times \isom \Zhat^{(p')}(1) \iso \prod_{\ell \neq p} \Z_\ell(1) \,, \]
	where the limit is taken over all finite fields~$\F$ of characteristic~$p$,
	the transition maps in the projective limit are the field norms and
	$\Z_\ell(1)$ is the Galois module $\varprojlim_n \mu_{\ell^n}(\Qbar_p)$.
	(See~\cite[\S\,{1.3}]{Serre1972}. Note that Serre uses~$I_p$ to denote~$I_p^\wild$,
	$I_t$ to denote~$I_p^\tame$ and~$I$ to denote~$I_p$.)
	
	\begin{lemma} \label[lemma]{Frob on Iptame}
		The absolute Galois group of the residue field~$\F_p$, which is canonically
		isomorphic to~$\Zhat\Frob_p$, acts on~$I_p^\tame$ via conjugation. One has
		\[ \Frob_p x \Frob_p^{-1} = x^p \quad \text{for $x \in I_p^\tame$} \,. \]
		Note that we have written $I_p^\tame$ multiplicatively here.
	\end{lemma}
	
	\begin{proof}
		See~\cite[Theorem~7.5.3]{NSW2.3}.
	\end{proof}
	
	See~\cite[\S\,{1.7}]{Serre1972} for the following definition.
	
	\begin{definition} \label[definition]{fundamental chars}
		Let $k \ge 1$. We define the character $\psi_k$ of~$I_p^\tame$ via the canonical
		projection from the projective limit as
		\[ \psi_k \colon I_p^\tame \isoto \varprojlim_{\Nm_{\F|}} \F^\times \surj \F_{p^k}^\times \,. \]
		One has $\Nm_{\F_{p^k}|\F_p} \circ \psi_k = \psi_1 = \chi_p$.
		
		The $k$ \emph{fundamental characters of level $k$} are the powers~$\psi_k^{p^n}$
		for $0 \le n < k$ (equivalently, $\psi_k$ followed by the $k$ automorphisms of~$\F_{p^k}$).
	\end{definition}

	\subsection{General set-up and notation}
	\label{sec:set-up}
	
	In the following, $f$ will always denote a newform of weight~$2$, level~$N$ and
	trivial nebentypus. We let $\cN(N)$ denote the set of such newforms;
	$\cN(N, g)$ denotes the subset consisting of forms whose Galois orbit
	has size~$g$.
	The Fourier coefficients of~$f$ will be denoted~$a_n$
	(or~$a_n(f)$ if we want to make the dependence on~$f$ explicit); they generate
	the coefficient ring~$\Z[f]$, which is an order in a totally real number field
	$F = \Q(f)$ (of degree~$g$ when $f \in \cN(N, g)$).
	
	Further, $A$ will denote an abelian variety over~$\Q$ that is
	$\Q$-isogenous to~$A_f$ (e.g., $A = A_f$ or $A = A' = A_f/A_f[\frf(\cO_F/\Z[f])]$
	as in~\cref{RM ab var}) and has endomorphism ring~$\cO$. Let $\frp$ be a regular
	prime ideal of~$\cO$; then $\frp$ is a maximal ideal of~$\cO$; we write
	$\F_\frp = \cO/\frp$ for its residue class field and $p(\frp)$ for the
	residue characteristic (i.e., the characteristic of~$\F_\frp$).
	Then $\rho_\frp = \rho_{\frp,A}$ is the Galois representation on~$A[\frp]$;
	its semi-simplification~$\rho_\frp^\ss$ is independent of the choice of~$A$
	(as long as $\frp$ is regular). Since we are mostly interested in determining
	when $\rho_\frp$ is irreducible (or has maximal image), knowing~$\rho_\frp^\ss$
	is usually enough, and so we suppress the dependency on~$A$ in the notation.
	(Note that when $\frp$ is a prime ideal of~$\cO_F$ that does not correspond
	to a regular prime ideal of~$\cO$, then $\rho_{\frp,A'}$ will be reducible,
	since there is an isogeny $A' \to A$ whose kernel has nontrivial intersection
	with~$A'[\frp]$.)

	\subsection{Determination of the maximal image}
	\label{sec:maximal image}
	
	One of our goals is to show that the image of~$\rho_\frp$ is as large as
	possible for all but finitely
	many prime ideals~$\frp$ (with a small explicit set of possible exceptions).
	The first step is to determine what this maximal image is. Then we will
	show that it suffices to consider the image in~$\PGL_2(\F_\frp)$.
	
	To show that the projective image is maximal, we have to exclude the
	possibility that it is contained in one of the maximal subgroups of the
	maximal projective image, so we need a classification of these maximal
	subgroups. This will be done in
	Section~\ref{sec:classification-of-the-maximal-subgroups-of-psl2pfr} below.
	
	\begin{definition} \label[definition]{def: Gfrpmax}
		We write $G_\frp \defeq \rho_\frp(\GalQ) \subseteq \GL_2(\F_\frp)$ for the image of~$\rho_\frp$.
		\cref{characteristic polynomials of rho_pfr} implies that $\det(G_\frp) = \F_p^\times$,
		since $\det \circ \rho_\frp$ is the mod-$p$ cyclotomic character. We set
		\[ G_\frp^{\max} \defeq \{M \in \GL_2(\F_\frp) : \det(M) \in \F_p^\times\} \,; \]
		then $G_\frp \subseteq G_\frp^{\max}$.
	\end{definition}
	
	\begin{definition}
		We write $\P \colon \GL_2(\F_\frp) \to \PGL_2(\F_\frp)$ for the canonical surjection.
		For a subgroup~$G$ of~$\GL_2(\F_\frp)$, we write $\P G$ for its image in~$\PGL_2(\F_\frp)$.
		We call $\P G$ the \emph{projective image of~$G$}. We also say that $\P{G_\frp}$
		is the \emph{projective image of~$\rho_\frp$}.
		We write $\overline{\det} \colon \PGL_2(\F_\frp) \to \F_\frp^\times/\F_\frp^{\times 2}$
		for the homomorphism induced by the determinant.
	\end{definition}
	
	We write $\PSL_2(\F_\frp) \defeq \P \SL_2(\F_\frp)$ for the quotient of~$\SL_2(\F_\frp)$ by its
	center~$\{\pm I_2\}$. Note that this is not the same as the group of $\F_\frp$-points
	of the algebraic group~$\PSL_2$. When $p$ is odd, $\PSL_2(\F_\frp)$ has index~$2$
	in~$\PGL_2(\F_\frp)$ and is the kernel of~$\overline{\det}$.
	
	\begin{lemma} \label[lemma]{proj im of Gmax}
		We have that $\P G_\frp^{\max} = \PGL_2(\F_\frp)$ when $\deg \frp$ is odd
		and $\P G_\frp^{\max} = \PSL_2(\F_\frp)$ when $\deg \frp$ is even.
	\end{lemma}
	
	\begin{proof}
		The image of~$G_\frp^{\max}$
		in~$\PGL_2(\F_\frp)$ consists of the elements $\gamma \in \PGL_2(\F_\frp)$ such
		that $\overline{\det}(\gamma)$ is in the image of~$\F_p^\times$. The latter is
		trivial if and only if $\deg \frp$ is even (or $p = 2$, in which case
		$\PSL_2(\F_\frp) = \PGL_2(\F_\frp)$). (See also~\cite[\S\,5.2]{Ribet1976}.)
	\end{proof}
	
	\begin{proposition} \label[proposition]{red_to_proj}
		Let $G \le G_\frp^{\max}$ be a subgroup such that
		$\det(G) = \F_p^\times$ and $\P G = \P G_\frp^{\max}$.
		Then $G = G_\frp^{\max}$.
	\end{proposition}
	
	\begin{proof}
		First assume that $\#\F_\frp > 3$ and
		let $\Derived H$ denote the derived subgroup of a group~$H$. Since the center
		of~$\GL_2(\F_\frp)$ is abelian, the assumption $\P G = \P G_\frp^{\max}$
		implies that $\Derived G = \Derived G_\frp^{\max} = \SL_2(\F_\frp)$,
		so $\SL_2(\F_\frp) \le G$. The second equality follows from
		\[ \SL_2(\F_\frp) = \Derived\SL_2(\F_\frp) \le \Derived G_\frp^{\max}
		\le \Derived \GL_2(\F_\frp) \le \SL_2(\F_\frp) \,,
		\]
		where the first equality follows from the fact that $\PSL_2(\F_\frp)$
		is non-solvable simple (since non-abelian when $\#\F_\frp > 3$)
		by~\cite[Theorem~8.4]{LangAlgebra}.
		Since both groups map onto~$\F_p^\times$ under the determinant, we then have exact sequences
		\[ 1 \to \SL_2(\F_\frp) \to G \stackrel{\det}{\to} \F_p^\times \to 1 \quad\text{and}\quad
		1 \to \SL_2(\F_\frp) \to G_\frp^{\max} \stackrel{\det}{\to} \F_p^\times \to 1 \,,
		\]
		so $\#G = \#\SL_2(\F_\frp) \#\F_p^\times = \#G_\frp^{\max}$, whence the claim.
		
		The two cases $\F_\frp = \F_2$ or~$\F_3$ can be checked by an easy computation.
	\end{proof}

	\subsection{Classification of the maximal subgroups of $\P{G_\frp^\max}$}
	\label{sec:classification-of-the-maximal-subgroups-of-psl2pfr}
	
	By~\cref{proj im of Gmax}, $\P{G_\frp^\max} = \PGL_2(\F_\frp)$ when $\deg \frp$
	is odd, and $\P{G_\frp^\max} = \PSL_2(\F_\frp)$ when $\deg \frp$ is even.
	By~\cref{red_to_proj}, we know that $G_\frp = G_\frp^\max$ if and only if
	$\P{G_\frp} = \P{G_\frp^\max}$, which is equivalent to $\P{G_\frp} \not\subseteq \Gamma$
	for every maximal subgroup~$\Gamma$ of~$\P{G_\frp^\max}$. In this section,
	we recall the classification of these maximal subgroups.
	
	We begin with the case $\deg \frp$ even, where $\P{G_\frp^\max} = \PSL_2(\F_\frp)$.
	
	\begin{theorem}[Maximal subgroups of~$\P{G_\frp^\max}$, $\deg \frp$ even]
		\label[theorem]{subgroups of PSL_2(q) even}
		Let $p$ be a prime and let $q = p^{2e}$ be an even power of~$p$.
		The maximal subgroups of~$\PSL_2(\F_q)$ are as follows.
		\begin{enumerate}[\upshape(i)]
			\item \label[item]{class.Be}
			(Borel) The stabilizer of a point of~$\P^1(\F_q)$. It has order $q (q-1)/2$ when
			$q$ is odd and $q(q-1)$ when $q$ is even.
			\item \label[item]{class.Se}
			(Sub-line) The stabilizer $\PGL_2(\F_{q'}) \cap \PSL_2(\F_q)$ of a sub-line $\P^1(\F_{q'})$,
			where $q = {q'}^\ell$ with a prime~$\ell$ (in particular, $\ell \mid 2e$).
			\item \label[item]{class.De}
			(Dihedral) Stabilizers of a pair of points in~$\P^1(\F_q)$ (normalizer
			of a split Cartan subgroup, order $q-1$ for $q \neq 9$ odd and $2(q-1)$ for $q$ even)
			or of a pair of $\F_q$-conjugate points in~$\P^1(\F_{q^2})$ (normalizer of a nonsplit
			Cartan subgroup, order $q+1$ for $q \neq 9$ odd and $2(q+1)$ for $q$ even).
			\item \label[item]{class.Ee}
			(Exceptional) Subgroups isomorphic to $S_4$ (when $e = 1$ and $3 < p \equiv \pm 3 \bmod 8$),
			or~$A_5$ (when $e = 1$ and $p \equiv \pm 3 \bmod 10$).
		\end{enumerate}
	\end{theorem}
	
	\begin{proof}
		See~\cite[Corollary~2.2]{King2005}, taking into account that $q$ is an even power of~$p$.
	\end{proof}
	
	When $q = 4$, the sub-line and normalizer of a split Cartan case are
	in the same conjugacy class. 
	When $q = 9$, the normalizers of split Cartan subgroups
	are contained in exceptional subgroups of type~$A_5$, and the normalizers of nonsplit Cartan
	subgroups are contained in sub-line stabilizers.
	
	When $\deg \frp$ is odd, we have $\P{G_\frp^\max} = \PGL_2(\F_\frp)$.
	Since $\det(G_\frp) = \F_p^\times$ contains elements of~$\F_\frp^\times$
	that are non-squares in this case, we also know that $\P{G_\frp}$ is not
	contained in~$\PSL_2(\F_\frp)$.
	
	\begin{theorem}[Maximal subgroups of~$\P{G_\frp^\max}$, $\deg \frp$ odd]
		\label[theorem]{subgroups of PSL_2(q) odd}
		Let $p \neq 2$ be a prime and let $q = p^{2e+1}$ be an odd power of~$p$.
		The maximal subgroups of~$\PGL_2(\F_q)$ different from~$\PSL_2(\F_q)$ are as follows.
		\begin{enumerate}[\upshape(i)]
			\item \label[item]{class.Bo}
			(Borel) The stabilizer of a point of~$\P^1(\F_q)$. It has order $q (q-1)$.
			\item \label[item]{class.So}
			(Sub-line) The stabilizer $\PGL_2(\F_{q'})$ of a sub-line $\P^1(\F_{q'})$,
			where $q = {q'}^\ell$ with a prime~$\ell$ (in particular, $\ell \mid 2e+1$).
			\item \label[item]{class.Do}
			(Dihedral) Stabilizers of a pair of points in~$\P^1(\F_q)$ (normalizer
			of a split Cartan subgroup, order $2(q-1)$, when $q > 5$) or of a pair of
			$\F_q$-conjugate points in~$\P^1(\F_{q^2})$ (normalizer of a nonsplit
			Cartan subgroup, order $2(q+1)$).
			\item \label[item]{class.Eo}
			(Exceptional) Subgroups isomorphic to $S_4$
			(when $e = 0$ and $3 < p \equiv \pm 3 \bmod 8$),
			and if $q = 3$, $A_4$.
		\end{enumerate}
	\end{theorem}
	
	\begin{proof}
		See~\cite[Corollary~2.3]{King2005}, which excludes $q = 3$.  For $q = 3$, Magma computes that $\PGL_2(\F_3) \isom S_4$ has $3$ maximal subgroups, $S_3, D_4, A_4$, which correspond to the Borel, normalizer of nonsplit Cartan, and exceptional maximal subgroup case, respectively.
	\end{proof}
	
	\begin{definition}
		We say that $\rho_\frp$ or $G_\frp$ is \emph{Borel}, \emph{sub-line}, \emph{dihedral},
		or \emph{exceptional} when $\P{G_\frp}$ is contained in a maximal subgroup of~$\P{G_\frp^\max}$
		of the corresponding type. In the dihedral case, we distinguish between \emph{split}
		and \emph{nonsplit}, according to the Cartan subgroup involved. We say that $\rho_\frp$
		or~$G_\frp$ is \emph{reducible} or \emph{irreducible}, if the action of~$G_\frp$ on $\F_\frp^2$ is,
		and we say that it is \emph{maximal}, if $G_\frp = G_\frp^\max$.
	\end{definition}
	
	The action of~$G_\frp$ is reducible if and only if $\rho_\frp$ is Borel.
	In the sub-line case, the invariant sub-line can be the image of a nontrivial invariant
	subspace of~$A[\frp]$ considered as an $\F_p$-vector space. In this case (if $\rho_\frp$
	is not also Borel), $\rho_\frp$ is irreducible as a $2$-dimensional $\F_\frp$-representation,
	but reducible as a $2 (\deg \frp)$-dimensional $\F_p$-representation.
	See Section~\ref{sec:excluding-L-Ex} below for a more detailed discussion.

	\subsection{Irreducibility and maximality criteria for fixed $\frp$}
	\label{sec:criteria for fixed frp}
	
	In this section we collect some criteria that allow us to verify that $\rho_\frp$
	is irreducible or maximal for a given prime ideal~$\frp$, using information from
	the characteristic polynomials of~$\rho_\frp(\Frob_\ell)$ for $\ell \nmid Np$.
	Recall from~\cref{characteristic polynomials of rho_pfr}
	that the characteristic polynomial of~$\rho_\frp(\Frob_\ell)$ has the form
	\[ T^2 - \bar{a}_\ell T + \bar{\ell} \,, \]
	where $x \mapsto \bar{x}$ denotes the reduction homomorphism $\Z[f] \to \F_\frp$.
	
	We define some invariants associated to elements of~$\PGL_2(\F_\frp)$;
	see~\cite[\S\,{2}]{Serre1972}.
	For $\F$ a finite field of odd characteristic, we define the Legendre symbol
	for $a \in \F$ as usual:
	\[ \Bigl(\frac{a}{\F}\Bigr) =
	\begin{cases} 0 & \text{if $a = 0$,} \\
		1 & \text{if $a = b^2$ for some $b \in \F^\times$,} \\
		-1 & \text{otherwise.}
	\end{cases}
	\]
	
	\begin{lemma}
		Let $\F$ be a finite field.
		\begin{enumerate}[\upshape(1)]
			\item The function
			\[ \GL_2(\F) \to \F, \quad M \mapsto \frac{\Tr(M)^2}{\det(M)} \]
			descends to a function $u \colon \PGL_2(\F) \to \F$.
			\item Assume that $\F$ has characteristic $p \neq 2$. The function
			\[ \GL_2(\F) \to \{0,1,-1\}, \quad M \mapsto \Bigl(\frac{\Tr(M)^2 - 4 \det(M)}{\F}\Bigr) \]
			descends to $\Delta \colon \PGL_2(\F) \to \{0,1,-1\}$.
		\end{enumerate}
	\end{lemma}
	
	\begin{proof}
		Since $\Tr^2$ and $\det$ are both homogeneous of degree $2$ and $\det(M) \neq 0$,
		the existence of~$u$ follows. Similarly, $\Tr(M)^2 - 4 \det(M)$ is well-defined
		up to multiplication with a non-zero square, so the Legendre symbol is well-defined.
	\end{proof}
	
	We now assume that the characteristic of~$\F$ is odd.
	
	If $u(g) \neq 0$ (equivalently, $\Tr(M) \neq 0$,
	where $M$ is a lift of~$g$ to~$\GL_2(\F)$), then $\Delta(g) = \left(\frac{u(g)(u(g)-4)}{\F}\right)$.
	If $u(g) = 0$, then $\Delta(g) = \left(\frac{-\det(M)}{\F}\right) \neq 0$,
	so $\Delta(g) = 0$ is equivalent to $u(g) = 4$.
	
	We note that $\Delta(g)$ gives the square class of the discriminant of the characteristic
	polynomial of any lift~$M$ of~$g$ to~$\GL_2(\F)$. This implies that $\Delta(g) \neq 0$
	if and only if $M$ has distinct eigenvalues (and hence is semi-simple).
	The eigenvalues are in~$\F$ when
	$\Delta(g) = 1$ and in the quadratic extension of~$\F$ and conjugate when $\Delta(g) = -1$.
	It follows that the elements of any Borel subgroup of $\PGL_2(\F)$ have $\Delta \neq -1$.
	We therefore obtain the following.
	
	\begin{corollary}
		Let $\frp$ be a prime ideal of odd residue characteristic. If $\Delta(g) = -1$
		for some $g \in \P{G_\frp}$, then $\rho_\frp$ is irreducible.
	\end{corollary}
	
	\begin{proof}
		If $\rho_\frp$ were reducible, then $\P{G_\frp}$ would be contained in a Borel subgroup,
		and so $\Delta(\P{G_\frp}) \subseteq \{0, 1\}$, contradicting the assumption.
	\end{proof}
	
	This gives a method to prove the irreducibility of~$\rho_\frp$ by computing
	\[ \Delta(\P \rho_\frp(\Frob_\ell)) = \Bigl(\frac{a_\ell^2 - 4 \ell}{\F_\frp}\Bigr) \]
	for a number of primes $\ell \nmid N p(\frp)$. If we obtain the value~$-1$ for one such~$\ell$,
	this shows that $\rho_\frp$ is irreducible.
	
	\medskip
	
	We can use the invariant~$u(g)$ to obtain information on the order of~$g$.
	(See~\cite[§\,2.6\,iii]{Serre1972}.)
	
	\begin{proposition}\label[proposition]{detection-of-elements-of-small-order}
		Let $\F$ be a finite field of characteristic $p$ and let $g \in \PGL_2(\F)$.
		\begin{enumerate}[\upshape(1)]
			\item $g$ is unipotent $\iff u(g) = 4 \iff \Delta(g) = 0$.
			\item If $p \neq 2$: $\ord(g) = 2 \iff u(g) = 0$.
			\item If $p \neq 3$: $\ord(g) = 3 \iff u(g) = 1$.
			\item If $p \neq 2$: $\ord(g) = 4 \iff u(g) = 2$.
			\item If $p \neq 5$: $\ord(g) = 5 \iff u(g)^2 - 3 u(g) + 1 = 0$.
		\end{enumerate}
	\end{proposition}
	
	\begin{proof}
		Let $g = \P M$ for some $M \in \GL_2(\F)$. If $\Delta(g) \neq 0$ (equivalently,
		$u(g) \neq 4$), then $M$ is semi-simple by the discussion above, and so, up to scaling,
		we can diagonalize $M$ over~$\bar{\F}$ as $M \sim \diag(1, \zeta)$, where $\zeta$ is some
		root of unity of order~$\ord(g)$. Then $u(g) = (1 + \zeta)^2/\zeta = \zeta + 2 + \zeta^{-1}$.
		Two values of~$u$ agree if and only if the corresponding values of~$\zeta$ are either
		equal or inverses of each other. Claim~(1) follows from the discussion above,
		and the others follow by observing that the condition on~$p$ ensures that
		the corresponding $u(g)$ is not equal to~$4$ and by matching roots of unity~$\zeta$
		with~$u$: $\zeta = -1 \iff u = 0$, $\ord(\zeta) = 3 \iff \zeta + \zeta^{-1} = -1 \iff u = 1$,
		$\ord(\zeta) = 4 \iff \zeta + \zeta^{-1} = 0 \iff u = 2$, and the two values
		of $\zeta + 2 + \zeta^{-1}$ for a fifth root of unity~$\zeta$ are the roots of
		$u^2 - 3 u + 1$.
	\end{proof}
	
	We can use this to show that $\rho_\frp$ is not exceptional, since the elements
	of~$S_4$ have order at most~$4$ and the elements of~$A_5$ have order $5$ or at most~$3$.
	So if we can find an element $g \in \P{G_\frp}$ of order at least~$5$, then
	$\P{G_\frp} \not\subseteq S_4$, and if we can find an element of order~$4$ or at least~$6$,
	then $\P{G_\frp} \not\subseteq A_5$.
	\medskip
	
	We now want to rule out the other possible maximal subgroups.
	
	If $\rho_\frp$ is dihedral, then $\P{G_\frp}$ is contained in the normalizer~$N(C)$
	either of a split or of a nonsplit Cartan subgroup~$C$. The elements of $N(C) \setminus C$
	have order~$2$, hence $u = 0$. If $C$ is split, then the nontrivial elements
	of~$C$ have $\Delta = 1$; if $C$ is nonsplit, its nontrivial elements have $\Delta = -1$.
	So if we find elements in~$\P{G_\frp}$ with
	$\Delta = 1$ and $u \neq 0$ and also elements with $\Delta = -1$ and $u \neq 0$,
	then $\P{G_\frp}$ cannot be dihedral.
	
	For the sub-line case, we restrict to $\deg \frp \le 2$. If $\rho_\frp$ is sub-line,
	we must then have $\deg \frp = 2$, and $\P{G_\frp} \subseteq \PGL_2(\F_p)$
	(up to conjugation in~$\PGL_2(\F_\frp)$). Since clearly $u(g) \in \F_p$ for each element
	$g \in \PGL_2(\F_p) \subset \PSL_2(\F_\frp)$, we can exclude the sub-line case
	when we find an element~$g \in \P{G_\frp}$ such that $u(g) \in \F_\frp \setminus \F_p$.
	Without the restriction on $\deg \frp$, we can similarly exclude the sub-line
	case when we find $g \in \P{G_\frp}$ such that $\F_p(u(g)) = \F_\frp$.
	It is also the case that the discriminant of the characteristic polynomial of any
	element is in~$\F_p$ (up to squares in~$\F_\frp^\times$) and therefore a square
	in~$\F_\frp$, so that $\Delta \in \{0, 1\}$.
	So, similar to the Borel case, this case can also be ruled out as soon as we find an element
	with $\Delta = -1$. (This last argument is specific to $\deg \frp = 2^n$ for some~$n$.)
	
	Assuming that we already know that the image is not exceptional, we can
	therefore prove that it is maximal by considering primes $\ell \nmid Np$,
	computing
	\[ \Delta(\ell) := \Delta(\P\rho_\frp(\Frob_\ell)) = \Bigl(\frac{a_\ell^2 - 4 \ell}{\F_\frp}\Bigr) \]
	until we have found one~$\ell$ such that $\Delta(\ell) = -1$ and $\frp \nmid a_\ell$
	and another~$\ell$ such that $\Delta(\ell) = 1$ and $\frp \nmid a_\ell$.
	(Recall that $u(\P\rho_\frp(\Frob_\ell)) \neq 0 \iff \frp \nmid a_\ell$).
	
	\medskip
	
	We obtain the following algorithm that returns a set of possible types of subgroups
	of~$\P{G_\frp^\max}$ that can contain~$\P{G_\frp}$. If this set is empty, then
	$\rho_\frp$ has maximal image. We assume $g = 2$ here since the discussion
	of the sub-line case above was assuming $\deg \frp \le 2$, and $g = 2$ is our main
	case of interest. The algorithm can be modified to work for general~$g$ if desired.
	
	We use the symbols $S_4$ and~$A_5$ to denote subgroups isomorphic to the respective
	groups, $R$ (\enquote{reducible}) for a Borel subgroup, $L$ for a sub-line stabilizer, and $N_s$
	and~$N_{ns}$ for the normalizers of a split or nonsplit Cartan subgroup.
	
	\begin{algorithm} \label[algorithm]{image for fixed frp} \strut
		
		\noindent\textsc{Input:} A newform~$f \in \cN(N, 2)$.
		A prime ideal~$\frp$ of the maximal order~$\cO$ of~$\Z[f]$. A bound~$B$.
		
		\noindent\textsc{Output:} A subset of $\{R, L, N_s, N_{ns}, S_4, A_5\}$ such that
		if a type is not in the set, then $\P{G_\frp}$ is not contained in a maximal subgroup
		of~$\P{G_\frp^\max}$ of this type.
		
		\begin{enumerate}[1.]
			\item {[Initialize]}
			Set $S \defeq \{R, L, N_s, N_{ns}, S_4, A_5\}$. Set $p \defeq p(\frp)$.
			\item {[Degree~1]}
			If $\deg \frp = 1$:
			\begin{enumerate}[a.]
				\item Remove $L$ and~$A_5$ from~$S$.
				\item If $p \in \{2, 3\}$, then remove $N_s$ and~$S_4$ from~$S$.
			\end{enumerate}
			\item {[Degree~2]} \label{xxx}
			If $\deg \frp = 2$:
			\begin{enumerate}[a.]
				\item If $p = 2$, then remove $N_s$, $S_4$ and~$A_5$ from~$S$.
				\item If $p = 3$, then remove $N_s$, $N_{ns}$ and~$S_4$ from~$S$.
				\item \label{xxxc} If $p \ge 5$ and $p^2 \nmid N$, then remove $N_{ns}$ from~$S$.
				\item If $p \not\equiv \pm 3 \bmod 10$, then remove $A_5$ from~$S$.
			\end{enumerate}
			\item {[$S_4$ possible?]}
			If $p \not\equiv \pm 3 \bmod 8$, then remove $S_4$ from~$S$.
			\item {[Loop over primes]}
			For each prime~$\ell \le B$ such that $\ell \nmid Np$:
			\begin{enumerate}[a.]
				\item Compute the image~$u(\ell)$ of $a_\ell^2/\ell$ in~$\F_\frp$.
				\item If $p \neq 2$, then compute $\Delta(\ell) \defeq \left(\frac{a_\ell^2 - 4 \ell}{\F_\frp}\right)$.
				\item If $u(\ell) \notin \{0,1,2,4\}$, then remove $S_4$ from~$S$.
				\item If $u(\ell) \notin \{0,1,4\}$ and $u(\ell)^2 - 3 u(\ell) + 1 \neq 0$,
				then remove $A_5$ from~$S$.
				\item If $\deg \frp = 2$ and $u(\ell) \notin \F_p$, then remove~$L$ from~$S$.
				\item If $p = 2$, $\deg \frp = 1$ and $u(\ell) = 1$, then remove~$R$ from~$S$.
				\item If $p = 2$ and $\deg \frp = 2$, then remove~$R$ from~$S$. \\
				If in addition $u(\ell) = 1$, then remove~$N_{ns}$ from~$S$.
				\item If $p \neq 2$ and $\Delta(\ell) = -1$, then remove~$R$ and~$L$ from~$S$. \\
				If in addition $u(\ell) \neq 0$, then remove~$N_s$ from~$S$.
				\item If $p \neq 2$, $\Delta(\ell) = 1$ and $u(\ell) \neq 0$, then remove~$N_{ns}$ from~$S$.
				\item If $S = \emptyset$, then return~$\emptyset$.
			\end{enumerate}
			\item Return $S$.
		\end{enumerate}
	\end{algorithm}
	
	The correctness of the algorithm follows from the
	classification results in Section~\ref{sec:classification-of-the-maximal-subgroups-of-psl2pfr}
	and the discussion in this section.
	The fact that $N_{ns}$ can be excluded when $\deg \frp = 2$ in Step~\ref{xxxc}
	follows from~\cref{no Nns} below.
	
	If the image is indeed maximal, then $\P{G_\frp} = \P{G_\frp^\max}$ contains elements with
	$u \neq 0$ and $\Delta = 1$, with $u \neq 0$ and $\Delta = -1$, of
	order $\ge 6$ (when $p \ge 5$) and of order~$4$.
	Chebotarëv's density theorem then guarantees that suitable primes~$\ell$ exist
	to rule out all the possible types. Using an effective version of the density
	theorem would give an explicit bound~$B$ for the primes~$\ell$ that have to be
	considered in the algorithm to be able to decide whether $\rho_\frp$ has maximal image.
	This bound will be too large to be useful in practice, however.

	\subsection{Explicit computation} \label{ssec:rho_p explicit}
	
	When~\cref{image for fixed frp} returns a non-empty set of types, we can try
	to determine the image explicitly as follows. We assume that we have given a
	curve~$X$ of genus~$2$ over~$\Q$ whose Jacobian~$J$ is isogenous to~$A_f$;
	we will determine the image of~$\rho_{J,\frp}$ (which has the same semi-simplification
	as~$\rho_{\frp,f}$), assuming that $\frp$ is a regular prime of~$\End_{\Q}(J)$.
	
	We compute (using Magma, say) the big period matrix associated to~$J$, which
	allows us to write $J(\C) = \C^2/\Lambda$ for some (numerically) explicit
	lattice~$\Lambda$. We can also determine the action of~$\End(J)$ on~$\C^2$
	and so we can approximate numerically the points in~$J(\C)[\frp]$. We represent
	these points by (numerical) divisors on~$X$, which we then recognize as divisors
	supported in algebraic points (this will work when the precision is sufficiently
	large). We then verify that the algebraic points on~$J$ we obtain are indeed
	in~$J[\frp]$. Knowing the points explicitly as algebraic points allows us to
	determine the Galois action.
	
	Since Magma can easily determine the torsion subgroup of~$J(\Q)$
	using the algorithm described in~\cite[\S\,11]{Stoll1999}, we can at least
	deduce that the representation associated to some prime ideal~$\frp$ with $p(\frp) = p$
	is reducible if $J(\Q)[p]$ is nontrivial.
	
	\begin{examples} \label[example]{ex reducible 1}
		We give two examples for such an explicit computation.
		\begin{enumerate}[(1)]
			\item \label{ex_red_125}
			$A = J_0(125)^+$ with endomorphism ring the maximal order of $\Q(\sqrt{5})$.
			The representation $\rho_{\ideal{\sqrt{5}}}$ is reducible.
			$A[\sqrt{5}]$ has constituents $\mu_5^{\otimes 2}$ and~$\mu_5^{\otimes 3}$.
			In this case, $5 \nmid [\Q(A[\sqrt{5}]) : \Q]$, so
			\[ A[\sqrt{5}] \isom \mu_5^{\otimes 2} \oplus \mu_5^{\otimes 3} \,. \]
			\item \label{ex_red_147}
			$A = J_0(147)^{\langle w_3,w_{49} \rangle}$
			with endomorphism ring the maximal order of $\Q(\sqrt{2})$.
			There is a prime $\frp \mid 7$ in $\mathrm{End}(A)$
			such that $\rho_\frp$ is reducible. We find that
			its irreducible constituents are $\mu_7^{\otimes 3}$ and~$\mu_7^{\otimes 4}$.
			As $[\Q(A[\frp]) : \Q] = 7 \cdot (7-1)$ and
			$A[7](\Q(\sqrt{-7})) = A[7](\Q(\mu_7^{\otimes 3})) = 0$, one has a nonsplit
			short exact sequence of Galois modules
			\[ 0 \to \mu_7^{\otimes 4} \to A[\frp] \to \mu_7^{\otimes 3} \to 0 \,.\]
		\end{enumerate}
	\end{examples}

	\subsection{The image of inertia at~$p$}
	\label{sec:determination-of-the-image-of-ip}
	
	Our next goal will be to prove that $\rho_\frp$ is irreducible (or even
	maximal) for all but finitely many~$\frp$, with a small explicit set of possible
	exceptions. To this end, we need to study the representations~$\rho_\frp$
	more carefully, so that we can extract some uniform statements. We begin
	by considering the action of the inertia group at the prime~$p$.
	Recall the definitions and the notations $I_p$, $I_p^\wild$, $I_p^\tame$
	from Section~\ref{ssec:preliminaries}. Also recall the fundamental
	characters~$\psi_k$ from~\cref{fundamental chars}.
	In the following, $g$ is arbitrary again.
	
	We now consider~$\rho_\frp|_{I_p}$, where $p = p(\frp)$ is the residue characteristic
	of~$\frp$. We have the following result.
	
	\begin{theorem} \label[theorem]{action of Ip}
		Assume $p^2 \nmid N$. Exactly one of the following two statements is true.
		\begin{enumerate}[\upshape(i)]
			\item \label{I_p:i} $\rho_\frp|_{I_p}$ has, up to conjugation, the form
			\[ \begin{pmatrix} \chi_p^n & * \\ 0 & \chi_p^{1-n} \end{pmatrix} \]
			for some $n \in \{0, 1\}$.
			\item \label{I_p:ii} After extending to the quadratic extension of~$\F_\frp$
			when $\deg \frp$ is odd, $\rho_\frp|_{I_p}$ has, up to conjugation, the form
			\[ \begin{pmatrix} \psi_2 & 0 \\ 0 & \psi_2^{p} \end{pmatrix} \,. \]
		\end{enumerate}
	\end{theorem}
	
	\begin{proof}
		By \cite[Prop.~1]{Serre1987}, the claim is true up to the exponents of
		the characters. (Note that according to \emph{loc.\ cit.}, fundamental characters
		of level $> 2$ cannot occur.) By~\cite[Cor.~3.4.4]{Raynaud1974}, as extended
		via~\cite[Lemma~4.9]{LarsonVaintrob2014} to the semistable case, the characters
		must be among $\chi_p^0$ and~$\chi_p^1$ in the first case,
		and among $\psi_2^0$, $\psi_2^1$, $\psi_2^p$ and~$\psi_2^{p+1}$ in the
		second case. The condition that $\det \circ \rho_\frp = \chi_p = \psi_2^{p+1}$,
		together with the fact that the characters are conjugate in the second case,
		then fixes the exponents. See also~\cite[Thm.~3.6]{Lombardo2016}.
	\end{proof}
	
	\begin{definition}
		In case~\eqref{I_p:i}, we say that $\rho_\frp|_{I_p}$ has \emph{level~1},
		and in case~\eqref{I_p:ii}, $\rho_\frp|_{I_p}$ has \emph{level~2}.
	\end{definition}
	
	\begin{corollary} \label[corollary]{image of Ip lower bound}
		Assume $p^2 \nmid N$.
		When $\rho_\frp|_{I_p}$ has level~1, then
		$\P{\rho_\frp(I_p)}$ contains a cyclic subgroup of order~$p-1$
		of a split Cartan subgroup.
		In the case of level~2, $\P{\rho_\frp(I_p)}$ is cyclic of order~$p+1$.
	\end{corollary}
	
	\begin{proof}
		This follows immediately from~\cref{action of Ip}.
	\end{proof}
	
	\begin{corollary} \label[corollary]{no Nns}
		If $\deg \frp$ is even and $p > 3$ with $p^2 \nmid N$, then $\P{G_\frp}$
		cannot be contained in the normalizer of a nonsplit Cartan subgroup.
	\end{corollary}
	
	\begin{proof}
		By~\cref{image of Ip lower bound}, $\P{G_\frp}$ contains elements
		of order $\ge p - 1 > 2$. Since $\deg \frp$ is even, no quadratic extension
		is necessary in~\cref{action of Ip} in the level~$2$ case, so
		in all cases, we find elements of order $> 2$ in a \emph{split}
		Cartan subgroup.
		Since such elements lie in a unique Cartan subgroup (which is the
		centralizer of the element) and the normalizer of a nonsplit Cartan
		subgroup contains only one Cartan subgroup, the claim follows.
	\end{proof}
	
	\begin{lemma} \label[lemma]{tame}
		Let $\chi \colon \GalQp \to \F^\times$ be a one-dimensional character
		of order prime to~$p$. Then $\chi|_{I_p}$ is a power of~$\chi_p$.
	\end{lemma}
	
	\begin{proof}
		Since $\chi$ has order prime to~$p$ and $I_p^\wild$
		is a pro-$p$ group, its image under~$\chi$ must be trivial, so $\chi$ is at
		most tamely ramified. Since $\Q_p^\nr(\mu_p)$ is the maximal abelian
		tamely ramified extension of~$\Q_p^\nr$, $\chi|_{I_p}$ must factor through
		$I_p \to I_p^\tame \to \F_p^\times$, which implies the claim.
	\end{proof}
	
	\begin{corollary} \label[corollary]{shape of reducible rep}
		Let $\frp$ be a regular prime of~$\Z[f]$ of residue characteristic~$p$
		and assume that $\rho_\frp$ is reducible and $p^2 \nmid N$. Then there is a
		character~$\epsilon$ of~$\GalQ$ with values in~$\F_\frp^\times$ and conductor~$d$
		such that $d^2 \mid N$ and (with respect to a suitable basis)
		\[ \rho_\frp = \begin{pmatrix} \epsilon \chi_p^n & * \\ 0 & \epsilon^{-1} \chi_p^{1-n} \end{pmatrix} \]
		with $n = 0$ or $n = 1$.
	\end{corollary}
	
	\begin{proof}
		The semi-simplification of~$\rho_\frp$ splits as a direct sum $\pi_1 \oplus \pi_2$
		of one-dimensional characters $\pi_1, \pi_2 \colon \GalQ \to \F_\frp^\times$.
		By~\cref{tame}, $\pi_1|_{I_p} = \chi_p^n$ for some~$n$, and so $\pi_2|_{I_p} = \chi_p^{1-n}$.
		We can therefore write
		$\pi_1 = \epsilon \chi_p^n$ with some character~$\epsilon$ that is unramified at~$p$.
		Since $\chi_p = \det \circ \rho_\frp = \pi_1 \cdot \pi_2$, it follows that
		$\pi_2 = \epsilon^{-1} \chi_p^{1-n}$. We then have
		\[ d^2 = \cond(\epsilon) \cond(\epsilon^{-1}) \mid \cond(\rho_\frp) \mid N \,. \]
		Since $p^2 \nmid N$, we have $n \in \{0, 1\}$ by~\cref{action of Ip}.
	\end{proof}
	
	\begin{remark} \label[remark]{refinement of irred}
		One can use this to refine~\cref{image for fixed frp} by potentially eliminating
		type~$R$ in more cases. For each prime $\ell \nmid Np$, compare the reduction
		of~$a_\ell$ mod~$\frp$ with all elements of the form
		$\epsilon(\ell) \ell^n + \epsilon(\ell)^{-1} \ell^{1-n}$
		for the finitely many possible characters~$\epsilon$ and
		$n \in \{0, 1\}$,
		and let $S_\ell$ be the set of compatible pairs~$(\epsilon, n)$.
		Then one takes the intersection of the sets~$S_\ell$ for several~$\ell$.
		If the intersection is empty, then $\rho_\frp$ must be irreducible.
	\end{remark}
	
	\begin{examples} \label[example]{ex reducible}
		We give two examples that illustrate~\cref{shape of reducible rep}.
		In order to determine whether $G_\frp$ has a nontrivial unipotent part,
		we determine explicit generators of~$A[\frp](\Qbar)$, which allows us to
		find $[\Q(A[\frp]) : \Q]$ (and, in fact, to determine the Galois action);
		see Section~\ref{ssec:rho_p explicit}.
		\begin{enumerate}[(1)]
			\item \label{ex_red_39}
			$A = J_0(39)^{w_{13}}$ with endomorphism ring the maximal order of $\Q(\sqrt{2})$.
			There is exactly one prime $\frp \mid 7$ such that $\rho_\frp$
			is reducible. Since $39$ is squarefree, $\rho_\frp^\ss \isom 1 \oplus \chi_7$.
			We find that $A[\frp](\Q) \isom \Z/7\Z$, so we have a short exact sequence
			\[ 0 \to \Z/7\Z \to A[\frp] \to \mu_7 \to 0 \,, \]
			which turns out to be nonsplit, since $[\Q(A[\frp]) : \Q]$ is divisible by~$7$.
			\item \label{ex_red_87}
			$A = J_0(87)^{w_{29}}$ with endomorphism ring the maximal order of $\Q(\sqrt{5})$.
			The representation $\rho_{\ideal{\sqrt{5}}}$ is reducible.
			Similarly as in Example~\eqref{ex_red_39}, the constituents are $\Z/5\Z$ and~$\mu_5$.
			Since $A[\sqrt{5}](\Q) \isom \Z/5\Z$, we have the exact sequence
			\[ 0 \to \Z/5\Z \to A[\frp] \to \mu_5 \to 0 \,, \]
			which is again nonsplit.
		\end{enumerate}
	\end{examples}
	
	\medskip
	
	We now consider the case that $\rho_\frp(\GalQp)$ is contained in the
	normalizer of a Cartan subgroup of~$\GL_2(\F_\frp)$.
	
	\begin{lemma} \label[lemma]{Frobp acts nontrivial iff level greater than 1}
		Assume $\frp$ is a regular prime ideal lying above a rational prime $p > 3$
		such that $p^2 \nmid N$ and
		that $\rho_\frp(\GalQp)$ is contained in the normalizer~$N(C)$ of a
		Cartan subgroup~$C$ of~$\GL_2(\F_\frp)$. Then the following are equivalent:
		\begin{enumerate}[\upshape(i)]
			\item $\rho_\frp|_{I_p}$ has level~1.
			\item $\rho_\frp(\GalQp) \subseteq C$.
			\item $\P\rho_\frp(\GalQp) \subseteq \P C$ has order~$p-1$ and $C$ is split.
		\end{enumerate}
	\end{lemma}
	
	\begin{proof}
		Clearly, (iii) implies~(ii) (since $\P^{-1} (\P C) = C$).
		If (ii) holds, then $\rho_\frp(\GalQp)$ is abelian of order prime
		to~$p$, so \cref{tame} implies that $\rho_\frp|_{I_p}$ has level~1.
		
		To show that (i) implies~(iii), first note that the image of the wild inertia
		group~$I_p^\wild$ must be trivial since $I_p^\wild$ is a pro-$p$ group
		and the order of~$N(C)$ is prime to~$p$ (since $p > 2$).
		So $\rho_\frp|_{I_p}$ factors through~$I_p^\tame$, which is pro-cyclic,
		and hence $\rho_\frp(I_p)$ is a cyclic group, which has order~$p-1$,
		since $\rho_\frp|_{I_p}$ has level~1 and $p^2 \nmid N$. Let $\Frob_p$ be any lift of the
		$p$-Frobenius on~$\ol{\F}_p$ to~$\GalQp$. \cref{Frob on Iptame} implies
		that conjugating by~$\rho_\frp(\Frob_p)$ has the effect of taking
		$p$th powers on~$\rho_\frp(I_p)$. Since $p \equiv 1 \bmod \#\rho_\frp(I_p)$,
		this action is trivial, so the image of~$\Frob_p$ commutes with the
		image of~$I_p$. This shows that $\rho_\frp(\GalQp)$ is abelian.
		Since $\#\rho_\frp(\GalQp)$ contains elements of order $p - 1 > 3$, this implies that
		$\rho_\frp(\GalQp) \le C$. (This is where we use that $\rho_\frp(\GalQp) \le N(C)$:
		all the abelian subgroups of~$N(C)$ containing elements of order $\ge 3$
		are contained in~$C$.) On the other hand, it follows from the
		discussion above that $\P\rho_\frp(I_p)$
		has order~$p - 1$; both together imply statement~(iii).
	\end{proof}
	
	\begin{corollary} \label[corollary]{not in nonsplit Cartan}
		Assume $\frp$ is a regular prime ideal lying above a rational prime $p > 3$
		such that $p^2 \nmid N$ and that $G_\frp$ is contained in a Cartan subgroup~$C$.
		Then $C$ is split; in particular, $\rho_\frp$ is reducible.
	\end{corollary}
	
	\begin{proof}
		The assumption implies that $\rho_\frp(\GalQp) \subseteq G_\frp \subseteq C$.
		The claim then follows from the implication \enquote{(ii) $\Rightarrow$ (iii)}
		in~\cref{Frobp acts nontrivial iff level greater than 1}.
	\end{proof}
	
	\begin{corollary}
		\label[corollary]{projective image of inertia contained in Cartan for dihedral}
		Assume $p > 3$ and $p^2 \nmid N$. If $\P{G_\frp}$ is dihedral, then $\P{\rho_\frp(I_p)}$
		is cyclic and contained in the corresponding Cartan subgroup~$C$.
	\end{corollary}
	
	\begin{proof}
		The first statement follows as in the proof
		of~\cref{Frobp acts nontrivial iff level greater than 1}.
		The second statement follows from the classification of~$\rho_\frp(I_p)$.
	\end{proof}

	\subsection{The image of inertia at a prime $\ell \neq p$}
	\label{sec:image of inertia not p}
	
	We now consider the image $\rho_\frp(I_\ell)$ of the inertia subgroup
	at a prime $\ell \neq p$.
	
	\begin{lemma} \label[lemma]{inertia unipotent if semistable}
		Let $\frp$ be a regular prime ideal of~$\Z[f]$ of residue characteristic~$p$
		and let $\ell \neq p$ be a prime. If $\ell^2 \nmid N$, then the image
		$\rho_\frp(I_\ell)$ of the inertia subgroup at~$\ell$ consists of unipotent
		elements.
	\end{lemma}
	
	\begin{proof}
		If $\ell \nmid N$, then $\rho_\frp$ is unramified at~$\ell$ and so the image
		of inertia is trivial. Otherwise, since the prime-to-$p$ part of the Artin conductor of~$\rho_\frp$
		divides the prime-to-$p$ part of~$N$, it follows from $v_\ell(N) = 1$ that the image
		of wild inertia is trivial and that the image of inertia has a one-dimensional
		fixed subspace; see~\cite[p.~181]{Serre1987}. So
		\[ \rho_\frp|_{I_\ell} \sim \begin{pmatrix} 1 & * \\ 0 & \chi_p \end{pmatrix} \,. \]
		Since $\chi_p$ is unramified at~$\ell$, this implies that
		$\rho_\frp(I_\ell)$ is unipotent.
		(See also~\cite[Section~2]{Ribet1997}; note that the definition of the conductor
		is purely local.)
	\end{proof}
	
	\begin{corollary} \label[corollary]{ramification when in Cartan normalizer}
		Let $\frp$ be a regular prime ideal of~$\Z[f]$ of residue characteristic~$p$
		and assume that $G_\frp$ is contained in the normalizer~$N(C)$ of a Cartan
		subgroup~$C$. Then $\rho_\frp$ is unramified at all primes $\ell \neq p$
		such that $\ell^2 \nmid N$.
	\end{corollary}
	
	\begin{proof}
		This follows from~\cref{inertia unipotent if semistable}, since $N(C)$ contains
		no nontrivial unipotent elements.
	\end{proof}
	
	\begin{corollary} \label[corollary]{conductor of character on Cartan normalizer}
		Let $\frp$ be a regular prime ideal of~$\Z[f]$ of residue characteristic~$p > 3$
		and assume that $G_\frp$ is contained in the normalizer~$N(C)$ of a Cartan
		subgroup~$C$. Since $C$ has index~$2$ in~$N(C)$, we obtain a quadratic character
		(which can be trivial)
		\[ \epsilon_\frp \colon \GalQ \stackrel{\rho_\frp}{\to} N(C) \to N(C)/C \simeq \{\pm 1\} \,; \]
		let $d$ be its conductor. Then the odd part of~$d^2$
		divides~$N$.  Moreover, if $4 \nmid N$, then $d^2 \mid N$.
	\end{corollary}
	
	\begin{proof}
		Note that the odd part of~$d$ is squarefree.
		Let $\ell \neq p$ be an odd prime. If $\ell^2 \nmid N$, then $\rho_\frp$ is
		unramified at~$\ell$ by~\cref{ramification when in Cartan normalizer}, and
		so $\ell \nmid d$. This shows the claim except for powers of $2$ or~$p$.
		If $p^2 \nmid N$, then $\rho_\frp(I_p) \subseteq C$
		by~\cref{projective image of inertia contained in Cartan for dihedral} (here we use $p > 3$)
		and so $\epsilon_\frp$ is unramified at~$p$. This takes care of the power of~$p$.
		If $4 \nmid N$, \cref{ramification when in Cartan normalizer} applies as well
		to show that $2 \nmid d$. This completes the proof.
	\end{proof}
	
	\begin{question}
		Can the stronger statement be extended to the case $4 \mid N$?
	\end{question}
	
	\begin{corollary} \label[corollary]{rho restricted to K still irreducible}
		Let $p > 2$ be prime and let $K|\Q$ be an imaginary quadratic extension
		such that $N$ and~$D_K$ are coprime and $K$ is
		not equal to $\Q(\sqrt{-1})$, $\Q(\sqrt{-2})$, or $\Q(\sqrt{-3})$.
		If $\rho_\frp$ is irreducible, then the restriction $\rho_\frp|_{G_K}$ is still irreducible.
	\end{corollary}
	
	\begin{proof}
		Suppose that $\rho_\frp$ is irreducible, but $\rho_\frp|_{G_K}$ is reducible.
		Then the quadratic character on~$G_\Q$ associated to~$K|\Q$ induces a nontrivial
		quadratic character~$\epsilon$ on the image of~$\rho_\frp$: $G_K$ fixes a one-dimensional
		subspace~$V$, and $\epsilon$ is given by the action on $\{V, \rho_\frp(\sigma) V\}$
		for $\sigma \in G_\Q \setminus G_K$. Since $\rho_\frp$ is ramified only at
		primes dividing~$Np$ and $K|\Q$ is ramified exactly at the primes dividing~$D_K$, it
		follows from the condition that $N$ and~$D_K$ are coprime that the conductor~$|D_K|$
		of~$\epsilon$ is a power of~$p$. Since $D_K \neq -3, -4, -8$, we must have $p > 3$
		and $D_K = -p$.
		
		Since $\rho_\frp$ fixes an unordered pair of complementary one-dimensional
		subspaces, but does not fix the two subspaces individually, it must be dihedral
		(compare~\cref{subgroups of PSL_2(q) even,subgroups of PSL_2(q) odd}). Then
		$\epsilon$ is the character as in~\cref{conductor of character on Cartan normalizer},
		so~\cref{conductor of character on Cartan normalizer} implies that $|D_K|^2 = p^2 \mid N$
		(here we use that $p > 3$), contradicting again the coprimality of $N$ and~$D_K$.
	\end{proof}

	\subsection{Explicit irreducibility for almost all $\frp$}
	\label{sec:irreducibility-for-almost-all-pfr}
	
	It is known that $\rho_\frp$ is irreducible for all but finitely many~$\frp$.
	Lombardo~\cite[Thm.~1.4]{Lombardo2016} gives an explicit (but very large) bound
	for the reducible primes (actually, for the primes such that the image is not maximal).
	What we need, however, is to determine as exactly as possible the finite set
	of primes~$\frp$ such that $\rho_\frp$ is reducible in each concrete case.
	
	In view of~\cref{shape of reducible rep}, we make the following definition.
	
	\begin{definition}
		We write~$d_\max$ for the largest positive integer~$d$ such that $d^2 \mid N$.
	\end{definition}
	
	We then obtain the following criterion. (Compare~\cite[§\,3.1]{Dieulefait2002},
	where similar criteria are used to show that the mod-$p$ Galois representations
	associated to an abelian surface with minimal endomorphism ring are maximal
	for almost all primes~$p$.)
	
	\begin{proposition} \label[proposition]{criterion for irred}
		Let $\ell \nmid N$ be a prime and let $m(\ell)$ be the order of~$\ell$
		in~$(\Z/d_\max\Z)^\times$. Let $\frp$ be a regular prime of~$\Z[f]$ of
		residue characteristic~$p$. If $p^2 \nmid N$ and
		\[ \frp \nmid \ell \cdot \resultant_T(T^2 - a_\ell T + \ell, T^{m(\ell)} - 1) \in \Z[f] \,, \]
		then $\rho_\frp$ is irreducible.
	\end{proposition}
	
	We note that when $\varphi(d_\max) = 1$ (which is the case, e.g., when $N$ is squarefree)
	or $\ell \equiv 1 \bmod d_\max$, the resultant simplifies to~$\ell + 1 - a_\ell$.
	
	\begin{proof}
		Let us prove the contrapositive. Thus suppose that $\rho_\frp$ is reducible.
		Let $\epsilon$ be the character in~\cref{shape of reducible rep}.
		It can be considered as a Dirichlet character of conductor
		$d \mid d_\max$ with values in~$\F_\frp^\times$.
		Since $p^2 \nmid N$, we have $n = 0$ or $n = 1$ in~\cref{shape of reducible rep},
		and by symmetry (and since what we do depends only on the semi-simplification of~$\rho_\frp$), we
		can assume $n = 0$. If $\ell = p$, then $\frp \mid \ell$, and we are done.
		Hence suppose $\ell \nmid Np$.
		Then $\rho_\frp(\Frob_\ell)$ is well-defined and has $\epsilon(\ell)$ as an eigenvalue,
		which is also a root of unity of order dividing~$m(\ell)$. So the
		characteristic polynomial $T^2 - \bar{a}_\ell T + \bar{\ell} \in \F_\frp[T]$
		of~$\rho_\frp(\Frob_\ell)$
		has at least one root in common with $T^{m(\ell)} - 1$. This is equivalent to
		\[ \resultant_T(T^2 - \bar{a}_\ell T + \bar{\ell}, T^{m(\ell)} - 1) = 0 \,. \]
		Since the resultant is compatible with ring homomorphisms, this implies
		that $\frp$ divides the resultant in the statement.
	\end{proof}
	
	\begin{remark}
		One can alternatively consider the condition
		\[ p \mid \ell \cdot
		\resultant_T\bigl(\charpol(\rho_{p^\infty}(\Frob_\ell))(T), T^{m(\ell)} - 1\bigr) \in \Z \]
		with the same $m(\ell)$ as above.
		Using the fact that the resultant is the product of all differences of roots
		of the first and of the second polynomial, together with the Weil conjectures
		for the characteristic polynomial of~$\Frob_\ell$ on~$T_p A$, we see that
		the resultant above is an integer~$R$ satisfying
		\[ 0 < (\sqrt{\ell} - 1)^{2gm(\ell)} < R < (\sqrt{\ell}+1)^{2gm(\ell)} \,. \]
		In particular, it is non-zero, so taking just one $\ell \nmid N$ gives a
		relatively small bound for the set of primes~$p$ such that $\rho_\frp$ is reducible
		for some $\frp \mid p$. In practice, one takes several~$\ell$ and uses the gcd
		of the resultants to obtain reasonably sharp bounds.
	\end{remark}
	
	This leads to the following algorithm. Recall the notation $p(\frp)$ for the residue
	characteristic of the prime ideal~$\frp$.
	
	\begin{algorithm} \label[algorithm]{Dieulefait irreducibility algorithm} \strut
		
		\noindent\textsc{Input:} A newform~$f \in \cN(N, g)$. A bound~$B$.
		
		\noindent\textsc{Output:} A finite set~$S$ of primes of the maximal order~$\cO$ of~$\Z[f]$
		such that $\rho_\frp$ is irreducible for all $\frp \not\in S$, or \enquote{failure}.
		
		\begin{enumerate}[1.]
			\item {[Maximal conductor]}
			Compute $d_\max \defeq \prod_{p \mid N} p^{\lfloor v_p(N)/2 \rfloor}$.
			\item {[Initialization]} Let $I \defeq \langle 0 \rangle$ as an ideal of~$\cO$.
			\item {[Loop over primes]} For all primes $\ell \le B$ such that $\ell \nmid N$:
			\begin{enumerate}[a.]
				\item Compute the order $m(\ell)$ of~$\ell$ in~$(\Z/d_\max\Z)^\times$.
				\item Set $I \defeq I
				+ \langle \ell \cdot \resultant_T(T^2 - a_\ell T + \ell, T^{m(\ell)} - 1) \rangle$.
			\end{enumerate}
			\item {[Result]} If $I = \langle 0 \rangle$, then output \enquote{failure}, else output
			$\{\frp : \text{$\frp \mid I$ or $p(\frp)^2 \mid N$}\}$.
		\end{enumerate}
	\end{algorithm}
	
	We can then use~\cref{image for fixed frp} on the regular odd primes in~$S$ to try
	to show that $\rho_\frp$ is irreducible even though \cref{Dieulefait irreducibility algorithm}
	was unable to prove that. One can also use the idea from~\cref{refinement of irred}.
	
	\begin{remark}
		If we take $B$ sufficiently large in~\cref{Dieulefait irreducibility algorithm},
		then by the Chebotarëv density theorem, the set~$S$ that the algorithm returns
		will contain only the prime ideals~$\frp$ such that $p(\frp)^2 \mid N$ or
		that the image of~$\rho_\frp$ consists of elements with one eigenvalue in the
		image of a character of conductor dividing~$d_\max$. This is compatible with
		the image being contained in the normalizer of a split Cartan subgroup (but not
		in the Cartan subgroup itself).
	\end{remark}
	
	\begin{example}
		Let $A = J_0(35)^{w_7}$ be the Jacobian of the modular curve
		quotient~$X_0(35)/\langle w_7 \rangle$, where $w_7$ denotes the Atkin--Lehner involution.
		$A$ corresponds to the
		\href{http://www.lmfdb.org/ModularForm/GL2/Q/holomorphic/35/2/a/b/}%
		{Galois orbit of size~$2$ of newforms of level~$35$},
		weight~$2$ and trivial nebentypus. Let $f$ be one of these two newforms. Then $\Z[f]$
		is the maximal order of~$\Q(\sqrt{17})$. Since the level $N = 35$ is squarefree,
		$d_\max = 1$.
		
		Let $\beta = (1 + \sqrt{17})/2$ be a generator of~$\Z[f]$. We then have $a_2 = -\beta$
		(for one of the two conjugate newforms), so the resultant is
		\[ \resultant_T(T^2 + \beta T + 2, T - 1) = 3 + \beta \,, \]
		which is an element of norm~$2^3$. This shows that $\rho_\frp$ is irreducible
		for all prime ideals~$\frp$ with odd residue characteristic. We also have $a_3 = \beta - 1$
		and the corresponding resultant is
		\[ \resultant_T(T^2 - (\beta - 1) T + 3, T - 1) = 5 - \beta \,, \]
		whose norm is~$2^4$ and which is not divisible by~$2$, so it generates a power
		of one of the two prime ideals above~$2$ in~$\Z[f]$; explicitly,
		\[ \langle 5 - \beta \rangle = \langle \beta + 1 \rangle^4 = \frp^4 \,. \]
		Write $\langle 2 \rangle = \frp \frp'$. Then we can deduce that
		$\rho_{\frp'}$ is irreducible.
		
		Magma computes that $A(\Q)_\tors \isom \Z/16\Z$. This shows that $\rho_{\frp}$
		must be reducible (and that it is nontrivial, since $A[\frp](\Q) \isom \Z/2\Z$
		and not $(\Z/2\Z)^2$).
	\end{example}
	
	\begin{example}
		We now consider $A = J_0(125)^+$, which corresponds to a
		\href{http://www.lmfdb.org/ModularForm/GL2/Q/holomorphic/125/2/a/a/}{Galois orbit of newforms}
		with coefficient ring the maximal order of~$\Q(\sqrt{5})$. Let $\alpha = (1 + \sqrt{5})/2$.
		Then we can pick one of the newforms~$f$ such that $a_2 = -\alpha$ and $a_3 = \alpha - 2$.
		Here $d_\max = 5$ and so $m(2) = m(3) = 4$. We find that
		\[ \begin{array}{r@{{}\defeq{}}l@{{}={}}l}
			r_2 & \resultant_T(T^2 + \alpha T + 2, T^4 - 1)     & 15 + 5 \alpha \quad\text{and} \\
			r_3 & \resultant_T(T^2 - (\alpha-2) T + 3, T^4 - 1) & 90 - 15 \alpha \,.
		\end{array}
		\]
		The ideal of~$\Z[f]$ generated by $2 r_2$ and~$3 r_3$ has norm~$5^2$. So the
		algorithm shows that $\rho_\frp$ is irreducible for all primes
		$\frp \neq \langle \sqrt{5} \rangle$.
		An explicit computation shows that $\rho_{\langle \sqrt{5} \rangle}$ is reducible;
		see~\cref{ex reducible 1}, \eqref{ex_red_125}.
	\end{example}

	\subsection{Excluding the sub-line and exceptional cases for almost all $\frp$}
	\label{sec:excluding-L-Ex}
	
	In the following, we require that the coefficient field of the newform is of degree $2$,
	and so $\dim A = 2$ as well. This simplifies the discussion of the sub-line case.
	
	Our next goal will be to show that $\rho_\frp$ has maximal image for all~$\frp$
	outside a small explicit finite set. Starting from the result of the previous
	subsection, it remains to show that for all but a few explicit~$\frp$, $\P{G_\frp}$
	is not contained in the stabilizer of a sub-line (when $\deg \frp = 2$), in a maximal
	subgroup of type $S_4$ or~$A_5$, or in the normalizer of a Cartan subgroup.
	
	We begin with the sub-line case and assume that $\rho_\frp$ is irreducible
	and $p(\frp)$ is odd.
	In this case, $G_\frp \subseteq \P^{-1}(\PGL_2(\F_p)) \cap G_\frp^\max$,
	which is a group containing~$\GL_2(\F_p)$ with index~$2$. More precisely,
	let $\alpha \in \F_\frp^\times \setminus \F_p^\times$ be such that $\alpha^2 \in \F_p^\times$;
	then the group is the union of~$\GL_2(\F_p)$ and $\alpha \GL_2(\F_p)$.
	So $G_\frp$ stabilizes the union of a $2$-dimensional $\F_p$-subspace~$U$
	(that is not a $1$-dimensional $\F_\frp$-subspace) of~$A[p]$ and~$\alpha U$.
	If $G_\frp \subseteq \GL_2(\F_p)$ (for some choice of $\F_\frp$-basis of~$A[p]$),
	then these two $\F_p$-subspaces are fixed individually, and so $\rho_p$ is
	the direct sum of two copies of a $2$-dimensional $\GalQ$-representation over~$\F_p$;
	in particular, $\rho_p$ is reducible. Conversely, if $\rho_p$ is reducible
	(but $\rho_\frp$ is not), then $\P{G_\frp}$ must be contained in the stabilizer
	of the sub-line that is the image of a $2$-dimensional $\F_p$-subspace fixed by~$G_\frp$.
	So, by showing that $\rho_\frp$ is irreducible and $\P{G_\frp}$ is not contained
	in the stabilizer of a sub-line, we also show that $A$ has no nontrivial isogenies
	of degree a power of~$p$.
	
	In any case, we know from the discussion in Section~\ref{sec:criteria for fixed frp}
	that the element $u(\P\rho_\frp(\Frob_\ell)) = a_\ell^2/\ell \in \F_\frp$ must be in~$\F_p$
	for all $\ell \nmid Np$ if $\P{G_\frp}$ is contained in the stabilizer of a sub-line.
	This is equivalent to $p \mid a_\ell^2 - {a_\ell^\sigma}^2$, where $\sigma$ denotes
	the nontrivial automorphism of~$\Z[f]$. Since in our setting, there always are $\frp$ of degree~$2$
	such that $\P{G_\frp}$ does not consist entirely of elements~$g$ such that $u(g) \in \F_p$,
	there will be a set of primes~$\ell \nmid N$ of positive density such that
	$a_\ell^2 \neq {a_\ell^\sigma}^2$. (If $f$ has an inner twist by a quadratic
	character, then $a_\ell = \pm a_\ell^\sigma$ for almost all~$\ell$. But in this
	case $A_f$ is isogenous to the Weil restriction of an elliptic curve over the
	quadratic field associated to the character, and so $A_f$ is not absolutely simple.)
	We can replace $a_\ell^2 - {a_\ell^\sigma}^2$
	by $(a_\ell^2 - {a_\ell^\sigma}^2)/\sqrt{\disc(\Z[f])} \in \Z$ (note that the
	prime divisors of~$\disc(\Z[f])$ are always ramified, so the corresponding primes
	of~$\Z[f]$ have degree~$1$). This leads to the following algorithm.
	
	\begin{algorithm} \label[algorithm]{exclude sub-line} \strut
		
		\noindent\textsc{Input:} A newform~$f \in \cN(N, 2)$. A bound~$B$.
		
		\noindent\textsc{Output:} A finite set~$S$ of prime ideals of the
		maximal order~$\cO$ of~$\Z[f]$ such that $\P{G_\frp}$ is not contained
		in the stabilizer of a sub-line for all $\frp \notin S$, or \enquote{failure}.
		
		\begin{enumerate}[1.]
			\item {[Initialize]} Set $R \defeq 0 \in \Z$.
			\item {[Loop over primes]}
			For all primes $\ell \le B$ such that $\ell \nmid N$:
			\begin{enumerate}[a.]
				\item Set $R \defeq \gcd(R, \ell \cdot (a_\ell^2 - {a_\ell^\sigma}^2)/\sqrt{\disc(\Z[f])})$.
				\item If $R = 1$, then exit the loop.
			\end{enumerate}
			\item {[Result]}
			If $R = 0$, then return \enquote{failure}, \\
			else return $\{\frp : \text{$p(\frp) \mid R$ and $\deg \frp = 2$}\}$.
		\end{enumerate}
	\end{algorithm}
	
	By the discussion above, the algorithm will not return \enquote{failure}
	when $B$ is sufficiently large.
	
	We now consider the case of exceptional image. Our analysis of the projective image
	of the inertia group~$I_p$ allows us to find elements of order at least~$6$ in
	$\P{\rho_\frp(I_p)} \leq \P{G_\frp}$ in most cases, which implies
	that $\P{G_\frp}$ is not contained in $S_4$ or~$A_5$.
	
	\begin{proposition} \label[proposition]{image not exceptional}
		If $p \ge 7$ is a prime such that $p^2 \nmid N$, then $\P{G_\frp}$ is not exceptional
		for all $\frp \mid p$.
	\end{proposition}
	
	\begin{proof}
		Since we assume $p^2 \nmid N$, it follows from~\cref{image of Ip lower bound}
		that $\P{G_\frp}$ contains elements of order $p-1$ or~$p+1$.
		So if $p \ge 7$, there are always elements of order at least~$6$ in~$\P{G_\frp}$,
		which implies that $\P{G_\frp}$ cannot be contained in a group isomorphic
		to~$S_4$ or to~$A_5$.
	\end{proof}
	
	So we only have to consider prime ideals~$\frp$ of residue characteristic~$p$
	such that $p \le 5$ or $p^2 \mid N$. By the classification in
	Section~\ref{sec:classification-of-the-maximal-subgroups-of-psl2pfr},
	we can also exclude $p = 3$ when $\deg \frp = 1$ and $p = 2$.
	We can then run~\cref{image for fixed frp} on these finitely many~$\frp$ to reduce
	the set of possibly exceptional primes further.

	\subsection{Proving non-CM} \label{sec:proving-non-cm}
	
	Recall that our goal is to show that $G_\frp = G_\frp^\max$ for all $\frp$
	outside an explicit small finite set. Now if $f$ has~CM, then $G_\frp$ will
	\emph{always} be contained in the normalizer of a Cartan subgroup, so in this
	case our task is impossible. Note that by~\cref{not CM}, in our
	case of interest when $g = 2$, if the associated abelian surface is absolutely
	simple, $f$ cannot have~CM. However, since it may be of interest in other
	situations, we describe a suitable algorithm.
	
	We recall the relevant definition (see~\cite[§\,3]{Ribet1977}, specialized
	to the case of interest).
	
	\begin{definition} \label[definition]{def: CM form}
		Let $f$ be a newform of weight~$2$, level~$N$ and trivial nebentypus, and
		let $\epsilon$ be a Dirichlet character.
		We say that $f$ has \emph{CM by~$\epsilon$}, if $f \otimes \epsilon = f$.
		We say that \emph{$f$ has~CM}, if $f$ has~CM by some nontrivial Dirichlet
		character~$\epsilon$.
	\end{definition}
	
	\begin{remark}
		The Dirichlet character~$\epsilon$ is then necessarily quadratic:
		Since $f$ has totally real coefficients, $\epsilon$ has to take real values.
	\end{remark}
	
	If $f$ has CM by~$\epsilon$, then $\epsilon(\ell) a_\ell = a_\ell$ for
	all $\ell \nmid N \cond(\epsilon)$, so $a_\ell = 0$ whenever $\epsilon(\ell) = -1$.
	This can be used to show that $f$ does \emph{not} have CM by~$\epsilon$,
	by exhibiting a prime $\ell \nmid N$ such that $\epsilon(\ell) = -1$ and $a_\ell \neq 0$.
	
	The idea is then to first determine a finite set of possibilities for the
	conductor~$D$ of~$\epsilon$ and then to check for each of the finitely many
	possible characters~$\epsilon$ of conductor~$D$ that $f$ does not have CM
	by~$\epsilon$ using this approach.
	
	\begin{theorem} \label[theorem]{twist-of-newform-by-Dirichlet-character}
		Let $f$ be a newform of weight~$2$, level~$N$ and trivial nebentypus, and let $\epsilon$
		be a quadratic Dirichlet character of conductor~$D$. Then the twist $f \otimes \epsilon$ of~$f$
		by~$\epsilon$ is a normalized eigenform of level dividing~$\lcm(N,D^2)$
		and trivial nebentypus.
	\end{theorem}
	
	\begin{proof}
		See~\cite[Proposition~3.64]{ShimuraIATAF}, using that $\epsilon^2$ is trivial.
	\end{proof}
	
	\begin{proposition} \label[proposition]{newforms-from-Hecke-characters}
		If $f$ is a CM form of level~$N$, it is the newform associated to a Hecke character~$\psi$
		of some conductor~$\mfr$ of an imaginary quadratic number field~$K$ of discriminant~$-\Delta_K$.
		One has $N = \Delta_K M$ with $M$ the absolute norm $\Nabs(\mfr) \defeq \#\cO_K/\mfr$
		and $\langle \Delta_K \rangle \mid \mfr$. In particular, $\Delta_K^2 \mid N$.
	\end{proposition}
	
	\begin{proof}
		\cite[Theorem~1.4 and Corollary~1.5]{Schuett2009}
	\end{proof}
	
	In the situation of~\cref{newforms-from-Hecke-characters}, the CM~character~$\epsilon$
	is the quadratic character associated to~$K$.
	
	This leads to the following algorithm.
	
	\begin{algorithm} \label[algorithm]{non-CM} \strut
		
		\noindent\textsc{Input:} A newform $f$ of weight~$2$, level~$N$ and trivial nebentypus. A bound~$B$.
		
		\noindent\textsc{Output:} \enquote{non-CM} or \enquote{no result}.
		
		\begin{enumerate}[1.]
			\item Let $S$ be the set of all negative fundamental discriminants~$-\Delta$
			such that $\Delta^2 \mid N$.
			\item For each $\Delta \in S$, do the following.
			\begin{enumerate}[a.]
				\item For all primes~$\ell \le B$ such that $\ell$ is inert in~$\Q(\sqrt{-\Delta})$, do:
				\begin{enumerate}[(i)]
					\item if $a_\ell \neq 0$, then continue with the next~$\Delta$.
				\end{enumerate}
				\item Return \enquote{no result}.
			\end{enumerate}
			\item Return \enquote{non-CM}.
		\end{enumerate}
		
	\end{algorithm}
	
	\begin{remark}
		If $f$ does have CM, then this algorithm will eventually return \enquote{no result}.
		Then one has a candidate character~$\epsilon$ (the character associated to~$-\Delta$),
		and one can try to verify that $f \otimes \epsilon = f$ (this is a finite computation,
		as the space of cusp forms of weight~$2$, level $\lcm(N, \Delta^2)$ and trivial nebentypus
		is of finite dimension and can be computed).
	\end{remark}
	
	The following result
	guarantees that there will be enough primes~$\ell$ as above when $f$ has no~CM.

	\begin{theorem}\label[theorem] {density-of-ap-eq-0}
		Let $f$ be a non-CM form. The set of primes~$\ell$ such that $a_\ell = 0$ has density~$0$.
	\end{theorem}
	
	\begin{proof}
		See \cite[p.~174]{Serre1981}.
	\end{proof}
	
	\begin{example}
		If $N$ is squarefree, then \cref{newforms-from-Hecke-characters} implies that
		$f$ has no~CM, since the only possible value of~$\Delta$ would lead to $-\Delta_K = -1$,
		which is not a discriminant of an imaginary quadratic number field.
	\end{example}
	
	\begin{remark}
		The newforms of weight~$2$, level~$800$, trivial nebentypus and with coefficients
		in~$\Q(\sqrt{5})$ in the Galois orbit with LMFDB
		label~\href{https://www.lmfdb.org/ModularForm/GL2/Q/holomorphic/800/2/a/j/}{800.2.a.j}
		have CM by~$\Q(\sqrt{-5})$.
		Computing its endomorphism ring using Magma, we see that it has nontrivial
		idempotents, so $A_f$ is not absolutely simple as predicted by~\cref{not CM}.
	\end{remark}

	\subsection{Maximal image for almost all $\frp$}
	\label{sec:maximal-image-for-almost-all-pfr}
	
	Let $f$ as usual be a newform of weight~$2$, level~$N$ and trivial nebentypus.
	We now assume that $f$ does not have~CM.
	In this section, we will describe an algorithm inspired by~\cite{Cojocaru2005}
	that finds a small finite set of primes~$\frp$ such that for all $\frp \notin S$,
	the representation~$\rho_\frp$ has maximal image.
	
	Using the algorithms we have described so far (and assuming $g = 2$,
	as that is required for some of these algorithms), we can
	determine a finite set~$S$ of prime ideals such that for all $\frp \notin S$,
	the representation~$\rho_\frp$ is either maximal or irreducible and dihedral.
	So if $\P{G_\frp}$ is not maximal for $\frp \notin S$, it is contained in
	the normalizer~$N(C)$ of a Cartan subgroup~$C$. It therefore remains to find a finite
	set of prime ideals such that $\P{G_\frp}$ is not dihedral for $\frp$ outside
	this set.
	
	So assume that $\rho_\frp$ is irreducible and dihedral,
	with $\P{G_\frp} \subseteq N(C)$ as above. We can also assume that $p^2 \nmid N$,
	as this excludes only finitely many~$p$, which we can consider separately.
	By~\cref{not in nonsplit Cartan}, $\P{G_\frp}$ is not contained in~$C$.
	So the character~$\epsilon_\frp$ defined in~\cref{conductor of character on Cartan normalizer}
	is nontrivial.
	
	\begin{proposition} \label[proposition]{not dihedral when semistable}
		If $N$ is squarefree, the residue characteristic of~$\frp$ is not~$2$ and
		$\rho_\frp$ is irreducible, then $\rho_\frp$ is not dihedral.
	\end{proposition}
	
	\begin{proof}
		Let us prove the contrapositive. Thus suppose that $\rho_\frp$ is dihedral.
		Then by the discussion
		above, $\epsilon_\frp$ is nontrivial. On the other hand, the conductor~$d > 1$
		of~$\epsilon_\frp$
		satisfies $d^2 \mid N$ by~\cref{conductor of character on Cartan normalizer}.
		Thus $N$ is not squarefree.
	\end{proof}
	
	So in the semi-stable case, we already know that $\rho_\frp$ is maximal for all
	$\frp \notin S$. If $N$ is not squarefree, then \cref{conductor of character on Cartan normalizer}
	provides us with a finite set of possibilities for~$\epsilon_\frp$, and we
	can try to rule each of them out for all prime ideals outside a finite set.
	
	\begin{lemma} \label[lemma]{epsilon-p-inert-then-aell}
		Assume that $\frp$ has odd residue characteristic~$p$ and that $\rho_\frp$
		is dihedral, with associated character~$\epsilon_\frp$. If $\ell \nmid Np$
		is a prime such that $\epsilon_\frp(\ell) = -1$, then $\frp \mid a_\ell$.
	\end{lemma}
	
	\begin{proof}
		Since $\epsilon_\frp(\ell) = -1$ by assumption,
		we have $\rho_\frp(\Frob_\ell) \in N(C) \setminus C$.
		Then $\P\rho_\frp(\Frob_\ell)$ has order~$2$, which implies that
		$a_\ell \equiv \Tr(\rho_\frp(\Frob_\ell)) = 0 \bmod \frp$.
	\end{proof}
	
	We make use of this as follows. For each quadratic character~$\epsilon$ of
	conductor~$d$ such that $d^2 \mid N$ (or $d \mid 8 d_0$ with $d_0$ odd such that
	$d_0^2 \mid N$ if $4 \mid N$), find some prime $\ell = \ell(\epsilon) \nmid N$ such
	that $a_\ell \neq 0$ and $\epsilon(\ell) = -1$
	(there are many such primes by~\cref{density-of-ap-eq-0}). Then for all~$\frp$ such that
	$\frp \nmid \ell a_\ell$ (which are all but finitely many), it follows that
	$\epsilon_\frp \neq \epsilon$. (In practice, it makes sense to use several such primes~$\ell$
	to cut the set of possible exceptions down further.) So replacing $S$ by the union of~$S$
	with the finitely many finite sets $S_\epsilon = \{\frp : \frp \mid \ell a_{\ell(\epsilon)}\}$,
	we obtain the desired finite set~$S$ of prime ideals such that for $\frp \notin S$,
	$\rho_\frp$ has maximal image.
	
	\begin{algorithm} \label[algorithm]{superset of non-max primes} \strut
		
		\noindent\textsc{Input:} A non-CM newform $f \in \cN(N, 2)$. A bound~$B$.
		
		\noindent\textsc{Output:} A finite set~$S$ of prime ideals of the maximal order~$\cO$ of~$\Z[f]$
		such that for all $\frp \notin S$, $\rho_\frp$ has maximal image, or \enquote{failure}.
		
		\begin{enumerate}[1.]
			\item {[Initialize]}
			Let $S$ be the union of
			\[ \{\frp : \text{$p(\frp) \in \{2,3,5\}$ or $p(\frp)^2 \mid N$}\} \]
			and the finite sets returned by
			\cref{Dieulefait irreducibility algorithm,exclude sub-line}
			(run on $f$ with the bound~$B$). \\
			Return \enquote{failure} when one of these algorithms failed.
			\item {[Possible conductors]}
			Set $\cD \defeq \{d \in \Z_{>0} : \text{$d^2 \mid N$ and $2 \nmid d$}\}$. \\
			If $4 \mid N$, set $\cD \defeq \cD \cup \{4d : d \in \cD\} \cup \{8d : d \in \cD\}$.
			\item {[Loop over characters]}
			For each $d \in \cD$ and each quadratic Dirichlet character~$\epsilon$ of
			conductor~$d$, do:
			\begin{enumerate}[a.]
				\item {[Initialize]}
				Set $I \defeq \langle 0 \rangle$ as an ideal in~$\cO$.
				\item {[Loop over primes]}
				For each prime $\ell \le B$ such that $\ell \nmid N$ and $\epsilon(\ell) = -1$,
				set $I \defeq I + \langle \ell a_\ell \rangle$.
				\item {[Failure]} If $I = \langle 0 \rangle$, then return \enquote{failure}.
				\item {[Record  prime ideals]}
				Set $S \defeq S \cup \{\frp : I \subseteq \frp\}$.
			\end{enumerate}
			\item {[Refine]}
			Run \cref{image for fixed frp} on each $\frp \in S$ and remove $\frp$
			from~$S$ when the result is the empty set.
			\item Return $S$.
		\end{enumerate}
	\end{algorithm}
	
	If $B$ is sufficiently large, then the algorithm will not return \enquote{failure},
	and by the discussion above, $S$ will satisfy the specification.
	
	We can use the information obtained from the various algorithms
	to provide a list of possible types of maximal subgroups that could
	contain~$G_\frp$ for those primes~$\frp$ that are in the set returned
	by~\cref{superset of non-max primes}.

	\subsection{The image of the $\frp$-adic Galois representation}\label{sec:the-image-of-the-pfr-adic-galois-representation}
	
	For~\cref{Kato divisibility}, we also need information about the image of $\rho_{\frp^\infty}|_{\Gal(\Qbar|\Q(\mu_{p^\infty}))}$, and for~\cref{Howard Euler system}, about the image of $\rho_{\frp^\infty}$.
	
	\begin{proposition} \label[proposition]{p-adic image}
		Let $\cO$ be the ring of integers of an unramified extension of~$\Z_p$.
		Let $G \subseteq \SL_2(\cO)$ be a closed subgroup.
		\begin{enumerate}[\upshape(i)]
			\item If $p > 3$ and $G$ surjects onto $\SL_2(\cO/p)$,
			then $G = \SL_2(\cO)$.
			
			\item If $p = 3$ and $G$ surjects onto $\SL_2(\cO/3^2)$, then $G = \SL_2(\cO)$.
		\end{enumerate}
	\end{proposition}
	
	\begin{proof}
		See~\cite[Lemma~IV.23.3]{Serre1998}, noting that the proof works for~$\cO$
		instead of~$\Z_p$; for $p = 3$, our claim follows from the proof given there.
	\end{proof}
	
	\begin{proposition} \label[proposition]{mod-9 image}
		Assume that $\cO/3 \in \{\F_3,\F_{3^2}\}$.
		Let $\rho\colon \GalQ \to  \GL_2(\cO)$ be a continuous homomorphism with mod-$3^n$ reduction $\rho_{3^n}$.
		If $\rho_{3}$ is surjective and the number of characteristic polynomials
		of elements of $\rho_{3^2}(\GalQ)$ with constant term $1$ is larger than
		the number of characteristic polynomials of elements of $\rho_{3}(\GalQ)$ with constant term $1$,
		then $\rho_{3^2}$ is surjective.
	\end{proposition}
	
	\begin{proof}
		This is a Magma computation, looping over all subgroups of~$\SL_2(\cO/3^2)$
		that surject onto $\SL_2(\cO/3)$ and computing characteristic polynomials.
	\end{proof}
	
	\begin{proposition} \label[proposition]{p-adic image GL}
		Let $\cO$ be the ring of integers of an unramified extension of~$\Z_p$. Let
		\[ G \subseteq G_{\frp^\infty}^\max \defeq \{M \in \GL_2(\cO) : \det(M) \in \Z_p^\times\} \]
		be a closed subgroup with $\det(G) = \Z_p^\times$.
		
		If $p > 3$ and the image of~$G$ in $\GL_2(\F_\frp)$ contains~$\SL_2(\F_p)$, then $G = G_{\frp^\infty}^\max$.
	\end{proposition}
	
	\begin{proof}
		This follows from the proof of~\cite[Theorem~4.22]{Lombardo2016}.
	\end{proof}
	
	We do not need \cref{p-adic image GL} for the examples in this article, but it is useful for further examples.

	\subsection{Examples}
	\label{sec:examples-rho-frp}
	
	The following table contains the result of running our algorithms on all
	absolutely simple Jacobians with real multiplication of genus~$2$ curves
	over~$\Q$ that are contained in the~LMFDB. (The genus~$2$ curves in
	the~LMFDB have discriminant bounded by~$10^6$; since the conductor of
	the Jacobian is the square of the level~$N$ and divides the discriminant,
	this implies that $N \le 1000$.) We add information on the
	isogeny classes coming from Hasegawa or Wang curves that are not also
	represented by an LMFDB curve.
	
	The entry \enquote{$N$} gives the level and the letter~$x$ of the isogeny
	class of the curve in the~LMFDB (the LMFDB~label of the isogeny class
	is then $N^2.x\dots$). For the Hasegawa and Wang curves not representing
	an isogeny class of an LMFDB curve, we use the label from~\cite{FLSSSW}.
	The entry~\enquote{$p^2 \mid N$} lists the
	primes at which the Jacobian does not have semi-stable reduction.
	The third entry \enquote{$\disc(\cO)$} gives the discriminant of the endomorphism ring
	of the Jacobian.
	The next entry lists the prime ideals~$\frp$ such that $\rho_\frp$
	is reducible and gives the splitting of~$\rho_\frp^{\ss}$ into characters.
	We use~$\epsilon_d$ to denote the quadratic character associated to
	the quadratic extension of discriminant~$d$.
	Since $\rho_\frp^{\ss}$ is the same for isogenous Jacobians,
	we list each isogeny class only once.
	The primes are given as~\enquote{$p$} when $\frp = \langle p \rangle$
	is of degree~$2$, as~\enquote{$\frp_p'$} or~\enquote{$\frp_p''$}
	when $p = p(\frp)$ is split and as~\enquote{$\frp_p$} when $p$ is ramified.
	The last entry lists the prime ideals~$\frp$ such that $\rho_\frp$
	is irreducible, but~\cref{superset of non-max primes}
	does not prove that~$\rho_\frp$ is maximal. In these cases, we
	have determined the isomorphism type of~$\P{G_\frp}$ by a direct
	computation; we give it in the table.
	
	{\small
		\begin{longtable}{r@{}lccll}
			\toprule
			$N$ & & $p^2 \mid N$ & $\disc(\cO)$ & reducible $\rho_\frp$ & irreducible non-maximal $\rho_\frp$ \\
			\midrule
			\endfirsthead
			\toprule
			$N$ & & $p^2 \mid N$ & $\disc(\cO)$ & reducible $\rho_\frp$ & irreducible non-maximal $\rho_\frp$ \\
			\midrule
			\endhead
			\bottomrule
			\endfoot
			\bottomrule
			\endlastfoot
			
			23&a & & $5$ & $\frp_{11}'$: $\mathbf{1} \oplus \chi_{11}$ & $3$: $A_5$  \\
			29&a & & $8$ & $\frp_{7}'$: $\mathbf{1} \oplus \chi_{7}$ & $\frp_{2}$: $N(C_{ns})$  \\
			31&a & & $20$ & $\frp_{5}$: $\mathbf{1} \oplus \chi_{5}$ &  \\
			35&a & & $17$ & $\frp_{2}'$: $\mathbf{1} \oplus \chi_{2}$  & $\frp_{2}''$: $N(C_{ns})$ \\
			39&a & & $8$ & $\frp_{2}$: $\mathbf{1} \oplus \chi_{2}$;\quad $\frp_{7}'$: $\mathbf{1} \oplus \chi_{7}$ &  \\
			51&a & & $17$ & $\frp_{2}'$: $\mathbf{1} \oplus \chi_{2}$  & $\frp_{2}''$: $N(C_{ns})$;\quad $3$: $A_5$ \\
			65&a & & $12$ & $\frp_{2}$: $\mathbf{1} \oplus \chi_{2}$;\quad $\frp_{3}$: $\mathbf{1} \oplus \chi_{3}$ &  \\
			67&a & & $5$ &  & $3$: $A_4$ \\
			67&c & & $5$ & $\frp_{11}'$: $\mathbf{1} \oplus \chi_{11}$ &  \\
			73&a & & $13$ & $\frp_{3}'$: $\mathbf{1} \oplus \chi_{3}$ &  \\
			73&b & & $5$ &  & $3$: $A_4$ \\
			77&b & & $5$ & $2$: $\mathbf{1} \oplus \chi_{2}$ &  \\
			85&a & & $8$ & $\frp_{2}$: $\mathbf{1} \oplus \chi_{2}$ &  \\
			85&b & & $12$ & $\frp_{2}$: $\mathbf{1} \oplus \chi_{2}$;\quad $\frp_{3}$: $\mathbf{1} \oplus \chi_{3}$ &  \\
			87&a & & $5$ & $\frp_{5}$: $\mathbf{1} \oplus \chi_{5}$ &  \\
			88&b & [ 2 ] & $17$ & $\frp_{2}'$: $\mathbf{1} \oplus \chi_{2}$  & $\frp_{2}''$: $N(C_{ns})$ \\
			93&a & & $5$ &  & $3$: $A_5$ \\
			103&a & & $5$ &  & $3$: $A_4$ \\
			107&a & & $5$ &  & $3$: $A_5$ \\
			115&b & & $5$ &  & $3$: $A_5$ \\
			123&b & & $8$ & $\frp_{7}'$: $\mathbf{1} \oplus \chi_{7}$ & $\frp_{2}$: $N(C_{ns})$;\quad $3$: $A_5$  \\
			125&a & [ 5 ] & $5$ & $\frp_{5}$: $\chi_{5}^{2} \oplus \chi_{5}^{3}$ & $3$: $A_5$  \\
			129&a & & $8$ &  & $\frp_{2}$: $N(C_{ns})$;\quad $3$: $A_5$ \\
			133&c & & $5$ &  & $3$: $A_4$ \\
			133&d & & $13$ & $\frp_{3}'$: $\mathbf{1} \oplus \chi_{3}$ &  \\
			133&e & & $5$ & $\frp_{5}$: $\mathbf{1} \oplus \chi_{5}$ &  \\
			135&c & [ 3 ] & $52$ & $\frp_{3}'$: $\mathbf{1} \oplus \chi_{3}$ &  \\
			147&a & [ 7 ] & $8$ & $\frp_{2}$: $\mathbf{1} \oplus \chi_{2}$;\quad $\frp_{7}'$: $\chi_{7}^{3} \oplus \chi_{7}^{4}$ &  \\
			165&a & & $8$ & $\frp_{2}$: $\mathbf{1} \oplus \chi_{2}$ &  \\
			167&a & & $5$ &  &  \\
			176&a & [ 2 ] & $17$ & $\frp_{2}'$: $\mathbf{1} \oplus \chi_{2}$ & $\frp_{2}''$: $N(C_{ns})$  \\
			177&a & & $5$ &  & $3$: $A_5$ \\
			188&a & [ 2 ] & $5$ &  & $2$: $N(C_{ns})$;\quad $3$: $A_5$ \\
			191&a & & $5$ &  & $3$: $A_5$ \\
			193&a & & $5$ &  & $3$: $A_4$ \\
			205&a & & $5$ &  & $3$: $A_5$ \\
			207&b & [ 3 ] & $8$ & $\frp_{2}$: $\mathbf{1} \oplus \chi_{2}$  & $3$: $A_5$ \\
			209&a & & $8$ &  & $\frp_{2}$: $N(C_{ns})$ \\
			211&a & & $5$ & $\frp_{5}$: $\mathbf{1} \oplus \chi_{5}$ & $3$: $A_5$  \\
			213&a & & $5$ &  & $3$: $A_5$ \\
			221&a & & $5$ &  & $3$: $A_5$ \\
			223&a & & $8$ &  & $\frp_{2}$: $N(C_{ns})$ \\
			227&a & & $5$ &  & $2$: $N(C_{ns})$ \\
			245&a & [ 7 ] & $8$ & $\frp_{7}'$: $\chi_{7}^{3} \oplus \chi_{7}^{4}$ & $\frp_{2}$: $N(C_{ns})$;\quad $3$: sub-line  \\
			250&a & [ 5 ] & $5$ & $\frp_{5}$: $\chi_{5}^{2} \oplus \chi_{5}^{3}$ & $2$: $N(C_{ns})$  \\
			261&c & [ 3 ] & $20$ &  &  \\
			275&a & [ 5 ] & $5$ & $\frp_{5}$: $\mathbf{1} \oplus \chi_{5}$ & $3$: $A_5$  \\
			275&b & [ 5 ] & $13$ & $\frp_{3}'$: $\epsilon_{5} \oplus \epsilon_{-3 \cdot 5}$ &  \\
			287&a & & $5$ &  &  \\
			289&a & [ 17 ] & $13$ & $\frp_{3}'$: $\epsilon_{17} \oplus \epsilon_{-3 \cdot 17}$;\quad $\frp_{17}'$: $\chi_{17}^{3} \oplus \chi_{17}^{14}$ &  \\
			299&a & & $5$ &  &  \\
			303&a & & $8$ &  & $\frp_{2}$: $N(C_{ns})$ \\
			313&a & & $5$ &  &  \\
			321&a & & $5$ &  &  \\
			334&a & & $5$ &  &  \\
			357&a & & $8$ &  & $\frp_{2}$: $N(C_{ns})$ \\
			358&a & & $5$ & $\frp_{5}$: $\mathbf{1} \oplus \chi_{5}$ &  \\
			375&a & [ 5 ] & $5$ & $\frp_{5}$: $\chi_{5}^{2} \oplus \chi_{5}^{3}$ & $3$: $A_5$  \\
			376&b & [ 2 ] & $5$ &  & $2$: $N(C_{ns})$ \\
			376&e & [ 2 ] & $5$ &  & $2$: $N(C_{ns})$ \\
			383&a & & $5$ &  & $3$: $A_5$ \\
			389&a & & $8$ &  & $\frp_{2}$: $N(C_{ns})$ \\
			457&a & & $5$ &  &  \\
			461&a & & $5$ &  & $3$: $A_5$ \\
			491&a & & $5$ &  &  \\
			499&a & & $5$ &  & $3$: $A_5$ \\
			523&a & & $5$ &  & $3$: $A_4$ \\
			533&a & & $8$ &  & $\frp_{2}$: $N(C_{ns})$ \\
			599&a & & $5$ &  &  \\
			621&a & [ 3 ] & $8$ &  & $\frp_{2}$: $N(C_{ns})$;\quad $3$: $A_5$ \\
			621&c & [ 3 ] & $5$ &  & $3$: $A_5$ \\
			637&a & [ 7 ] & $5$ &  & $3$: $A_4$ \\
			640&a & [ 2 ] & $5$ & $2$: $\mathbf{1} \oplus \chi_{2}$ &  \\
			640&b & [ 2 ] & $5$ & $2$: $\mathbf{1} \oplus \chi_{2}$ &  \\
			647&a & & $5$ &  & $3$: $A_5$ \\
			677&a & & $5$ &  & $3$: $A_5$ \\
			683&a & & $5$ &  & $2$: $N(C_{ns})$ \\
			689&a & & $5$ &  & $2$: $N(C_{ns})$ \\
			752&a & [ 2 ] & $5$ &  & $2$: $N(C_{ns})$;\quad $3$: $A_5$ \\
			752&f & [ 2 ] & $5$ &  & $2$: $N(C_{ns})$ \\
			752&j & [ 2 ] & $5$ &  & $2$: $N(C_{ns})$ \\
			783&a & [ 3 ] & $5$ &  &  \\
			799&a & & $5$ &  & $3$: $A_5$ \\
			809&a & & $5$ &  &  \\
			837&b & [ 3 ] & $8$ &  & $\frp_{2}$: $N(C_{ns})$ \\
			841&a & [ 29 ] & $5$ & $\frp_{29}'$: $\chi_{29}^{5} \oplus \chi_{29}^{24}$ &  \\
			845&a & [ 13 ] & $5$ &  & $3$: $A_5$ \\
			877&a & & $5$ &  & $2$: $N(C_{ns})$ \\
			887&a & & $5$ &  & $3$: $A_5$ \\
			929&a & & $5$ &  & $3$: $A_5$ \\
			\midrule
			\multicolumn{6}{c}{Hasegawa curve isogeny class not in the LMFDB} \\
			\midrule
			161& & & $5$ &  & $3$: $A_5$ \\
			\midrule
			\multicolumn{6}{c}{\enquote{Wang only} curve isogeny classes not in the LMFDB} \\
			\midrule
			65&A & & $8$ & $\frp_{2}$: $\mathbf{1} \oplus \mathbf{1}$;\quad $\frp_{7}'$: $\mathbf{1} \oplus \chi_{7}$ & $3$: $A_5$  \\
			117&B & [ 3 ] & $8$ & $\frp_{2}$: $\mathbf{1} \oplus \chi_{2}$ &  \\
			125&B & [ 5 ] & $5$ & $\frp_{5}$: $\mathbf{1} \oplus \chi_{5}$ & $3$: $A_5$  \\
			175& & [ 5 ] & $5$ & $\frp_{5}$: $\mathbf{1} \oplus \chi_{5}$ & \\
		\end{longtable}
	}
	
	
	\section{Computation of Heegner points and the Heegner index}
	\label{sec:computation-of-the-heegner-points-and-index}
	
	For this section, we fix the following set-up. Let $f \in \cN(N, g)$
	be a newform of level~$N$ with Galois orbit of size~$g$, so that its
	coefficient ring~$\Z[f]$ is an order in the totally real number
	field~$\Q(f)$ of degree~$g$. Its Fourier coefficients are $a_n = a_n(f) \in \Z[f]$.
	We denote the set of embeddings $\Q(f) \inj \R$
	by~$\Sigma$, and we write $f^\sigma$ for the modular form with real
	coefficients obtained from~$f$ by applying $\sigma \in \Sigma$ to
	its coefficients. Recall that $I_f$ denotes the annihilator of~$f$
	in the integral Hecke algebra~$\bT$, and that we have morphisms
	of abelian varieties
	\[ A_f^\dual = J_0(N)[I_f] \stackrel{\iota_f}{\inj} J_0(N)
	\stackrel{\pi_f}{\surj} J_0(N)/I_f J_0(N) = A_f \,; \]
	the composition $\lambda_f = \pi_f \circ \iota_f \colon A_f^\dual \to A_f$
	is a polarization of~$A_f^\dual$; it is the polarization induced by
	the canonical principal polarization of~$J_0(N)$ as the Jacobian
	of~$X_0(N)$. (If $\lambda \colon A \to A^\dual$ is the polarization
	coming from $L \in \NS(A)$
	and $\phi \colon B \subseteq A$ is an abelian subvariety, then
	$\phi^\dual \circ \lambda \circ \phi \colon B \to B^\dual$ is the polarization
	on~$B$ coming from $\phi^* L$. See~\cite[Cor.~2.4.6\,(d)]{BirkenhakeLange}.)
	Note that $\iota_f$ is the composition of~$\pi_f^\dual$ with
	the inverse of the canonical polarization of~$J_0(N)$.
	We write
	\begin{equation} \label{Eq:def d_f}
		d_f \defeq d_1 \cdots d_g ,
	\end{equation}
	where $(d_1, \ldots, d_g)$ is the type of~$\lambda_f$; then $\deg \lambda_f = d_f^2$
	(see \cite[Thm.~3.6.1 and Cor.~3.6.2]{BirkenhakeLange}).
	The number~$d_f$ is sometimes called the \emph{modular degree} of~$A_f$;
	see for example~\cite[\S\,3.3]{AgasheStein2005}.
	
	We further assume that we are given a (nice) curve~$X$ whose Jacobian~$J$
	is isogenous to~$A_f$ via an isogeny $\pi \colon A_f \to J$.
	We denote $\End_\Q(J)$ by~$\cO$.
	The isogeny~$\pi$ induces an isomorphism of endomorphism algebras (where
	$\End^0_\Q(A) \defeq \End_\Q(A) \otimes_\Z \Q$)
	\begin{equation} \label{Eq:iso End0}
		\pi^0_* \colon \Q(f) = \End^0_\Q(A_f) \stackrel{\iso}{\To} \End^0_\Q(J) = \Frac(\cO)\,, \qquad
		\phi \longmapsto \pi \phi \pi^{-1} ,
	\end{equation}
	which we use to identify $\End^0_\Q(J) = \Q(f)$.
	In particular, $\cO$ is identified with an order of~$\Q(f)$.
	Then for any $\gamma \in \Z[f] \cap \cO$,
	it follows that
	\begin{equation} \label{Eq:pi equivariant}
		\pi \circ \gamma = \gamma \circ \pi .
	\end{equation}
	
	We write $\pi_J = \pi \circ \pi_f \colon J_0(N) \to J$. We then
	have a commutative diagram
	\begin{equation} \label{Eqn:AJ diagram}
		\SelectTips{cm}{}
		\xymatrix@!C{J_0(N) \ar@{=}[d] \ar@{->>}[r]_-{\pi_f} \ar@/^10pt/[rr]^{\pi_J}
			& A_f \ar@{->}[r]_{\pi}
			& J \ar@{<-}[d]^{\lambda} \\
			J_0(N)^\dual
			& A_f^\dual \ar@{_(->}[l]_-{\pi_f^\dual} \ar[u]_{\lambda_f}
			& J^\dual \ar@{->}[l]_-{\pi^\dual} \ar@/^10pt/[ll]^{\pi_J^\dual} }
	\end{equation}
	with a polarization $\lambda = \pi \circ \lambda_f \circ \pi^\dual = \pi_J \circ \pi_J^\dual$
	of~$J^\dual$.
	Then by pre-composing $\lambda$ with the canonical principal polarization~$\lambda_J$
	of~$J$, we obtain an element
	\begin{equation} \label{Eq:def alpha}
		\alpha \defeq \lambda \circ \lambda_J
		= \pi \circ \lambda_f \circ \pi^\dual \circ \lambda_J
		\in \cO \subseteq \Q(f) .
	\end{equation}
	By~\cite[Prop.~12.12]{MilneAV}, $(\deg \pi)^2 (\deg \lambda_f) = \deg \alpha = \Nm(\alpha)^2$,
	which implies that
	\begin{equation} \label{Eq:Nm alpha}
		\Nm(\alpha) = d_f \cdot \deg \pi \,.
	\end{equation}
	
	In practice, we start with the curve~$X$ of genus~$g$ and we know that its Jacobian~$J$
	has real multiplication (and is absolutely simple). We then need to find
	the corresponding newform~$f \in \cN(N, g)$. We first determine~$N$.
	When $X$ is a quotient of~$X_0(N)$, then we know~$N$ by construction.
	In general, we find $N$ as the square root of the conductor of~$J$, which can be
	computed up to finitely many choices of power of~$2$ at worst; in our LMFDB
	examples, the conductor has been determined exactly and is available in the LMFDB.
	Given~$N$, we then compare the traces of the
	Fourier coefficients at primes $\ell \nmid N$ of the various candidate~$f$
	with the corresponding coefficient of the $L$-function of~$J$.
	This quickly leaves only one candidate, which must then be the correct~$f$
	(up to the Galois action). We now assume that $f$ is fixed.
	
	One of the ingredients we need in order to prove that $\Sha(J/\Q)[\frp] = 0$
	for all except an explicit finite set of prime ideals~$\frp$ of~$\cO$ is the
	\emph{Heegner index}, whose definition we now recall.
	(See Sections~\ref{sec:Sha_an} and~\ref{sec:finite support} below for why the
	Heegner index is important.) Let $K$ be a
	\emph{Heegner field} for~$f$; this is an imaginary quadratic field such that
	all prime divisors of~$N$ split in~$K$ and such that the $L$-series
	$L(f/K, s) = L(f, s) L(f \otimes \epsilon_K, s)$ (with $\epsilon_K$ the
	quadratic character corresponding to~$K$) vanishes to first order at~$s = 1$.
	The first condition
	implies that $\cO_K$ contains ideals~$\frn$ of norm~$N$ such that $\cO_K/\frn$
	is a cyclic group of order~$N$. Then the natural map $\C/\cO_K \to \C/\frn^{-1}$
	corresponds to a cyclic isogeny of degree~$N$ between two elliptic curves
	with CM by~$\cO_K$ and so defines a point in~$X_0(N)$, which is known
	to be defined over the Hilbert class field~$H$ of~$K$. More generally,
	let $\fra$ be some ideal of~$\cO$; then we can consider $\C/\fra \to \C/\fra \frn^{-1}$.
	We obtain $h_K$ points $x_{[\fra]} \in X_0(N)(H)$ in this way, where
	$h_K$ denotes the class number of~$K$ and the point depends only on the
	ideal class of~$\fra$. These points form an orbit under~$\Gal(H|K)$;
	their formal sum~$\bx_K$ is the \emph{Heegner cycle} on~$X_0(N)$
	associated to~$K$ and~$\frn$; it is defined over~$K$.
	Let $\infty \in X_0(N)(\Q)$ denote the cusp at infinity.
	Then $y_K = [\bx_K - h_K \cdot (\infty)] \in J_0(N)(K)$ is a \emph{Heegner point}
	associated to~$K$. By varying~$\frn$ in the construction, we may get
	different Heegner points, but they all agree up to sign and adding
	a torsion point. (See also~\cite{Gross1984}.) So in the following,
	we will consider Heegner points up to sign and modulo torsion.
	
	We then obtain a point $y_{K,\pi} = \pi_J(y_K) \in J(K)$.
	By~\cite{GrossZagier1986}, the $\cO$-span of $y_{K,\pi} \in J(K)$
	is a rank $g = \dim{J}$ subgroup of finite index of~$J(K)$
	(which does not depend on the choice of the Heegner point).
	This index is the \emph{Heegner index}; we denote it by
	\begin{equation} \label{Eq:Def I_Kpi}
		I_{K,\pi} \defeq (J(K) : \cO \cdot y_{K,\pi}) .
	\end{equation}
	Considering the characteristic ideal
	$\cI_{K,\pi} \defeq \Char_\cO(J(K)/\cO y_{K,\pi})$ gives refined information;
	this refinement is helpful for our intended application,
	because we can study the summands of
	$\Sha(J/\Q)[p^\infty] = \bigoplus_{\frp \mid p} \Sha(J/\Q)[\frp^\infty]$ individually.
	In the same way, we have $y_{K,A_f} = \pi_f(y_K) \in A_f(K)$, and we
	set $I_K \defeq (A_f(K) : \Z[f] y_{K,A_f})$ and
	$\cI_K \defeq \Char_{\Z[f]}(A_f(K)/\Z[f] y_{K,A_f})$.
	
	For an abelian variety $A/\Q$ and a quadratic number field~$K$,
	we denote the quadratic twist of~$A$ by the quadratic character associated
	to~$K$ by~$A^K$. Then $A^K$ is isomorphic to~$A$ over~$K$.
	The natural map $A(\Q) \times A^K(\Q) \to A(K)$ has finite
	kernel and cokernel killed by~$2$: The kernel is the diagonally
	embedded $A(\Q)[2]$, and the image contains~$2 A(K)$.
	
	When $\Lrk J = 0$, then $J(K)$ is essentially~$J^K(\Q)$; more precisely,
	the image of~$J^K(\Q)$ in~$J(K)$ contains~$2 J(K)$ up to torsion, and so we can
	identify~$2 y_{K,\pi}$ up to torsion with a rational point on~$J^K$.
	When $\Lrk J = 1$, then $2 J(K)$ is contained in~$J(\Q)$ up to torsion,
	and we can identify~$2 y_{K,\pi}$ up to torsion with a
	rational point on~$J$. (See also~\cite[Lemma~2.1]{Miller2010}.)
	This simplifies the computations, since certain
	algorithms (for example, computing canonical heights on~$J$) are so far
	only implemented when the base field is~$\Q$.
	
	The aim of this section is to explain how we can compute the
	Heegner index~$I_{K,\pi}$ (or the corresponding ideal~$\cI_{K,\pi})$.
	
	The first step is to determine a Heegner field~$K$. This is
	explained in Section~\ref{sec:Heegner fields}.
	In order to determine the $\cO$-span of~$y_{K,\pi}$, we need to
	determine~$\cO = \End_\Q(J)$ and its action on~$J(K)$ (we can determine
	generators of~$J(K)$ from generators of~$J(\Q)$ and of~$J^K(\Q)$,
	which Magma can usually compute). Section~\ref{sec:Computing-the-endomorphism-ring}
	explains how to do that. Then, of course, we need to find
	the Heegner point~$y_{K,\pi}$ on~$J$.
	
	One approach is to compute the $j$-invariant morphism $X_0(N) \to \bP^1_\Q$
	given by sending the point representing an isogeny $E \to E'$ to~$j(E)$
	as an algebraic map. Then, given the $h_K$ different $j$-invariants
	of elliptic curves with CM by~$\cO_K$,
	we can lift them to the corresponding points in~$X_0(N)(H)$
	and thus get an algebraic description of the Heegner cycle.
	However, this turns out to be too slow even for
	moderately large~$N$. Therefore, we do not give more details here.
	
	Instead, we use an analytic approach. We start with the
	$h_K$ reduced integral binary quadratic forms whose roots with positive imaginary
	part map to the points in the support of the Heegner cycle on~$X_0(N)$
	under the uniformization map $\bH \to X_0(N)(\C)$, where $\bH$ denotes
	the upper half plane. These quadratic forms can easily be determined
	using the built-in Magma function~\texttt{HeegnerForms}.
	Via the uniformization map $\bH \to X_0(N)(\C)$,
	we obtain the set of $h_K$ points in the Heegner cycle. If the curve~$X$
	is a quotient of~$X_0(N)$, we can then map the Heegner cycle directly to~$X$
	and try to recognize it as a divisor defined over~$K$. We can then obtain the point
	on~$J$ given by the Heegner cycle.
	If $X$ is not a quotient of~$X_0(N)$, we do the following.
	We first use the Abel--Jacobi map $X_0(N)(\C) \to J_0(N)(\C)$
	and the map $\pi_{f,\C} \colon J_0(N)(\C) \to A_f(\C) \isom \C/\Lambda_f$
	to map the Heegner cycle to~$A_f(\C)$. We then compute the point $y_K \in A_f(\C)$
	(by taking a sum in~$\C/\Lambda_f$). Then we use an explicit numerical
	representation of the
	isogeny $\pi_\C \colon A_f(\C) \to J(\C)$ to map the Heegner point from $A_f(\C)$ to~$J(\C)$.
	Finally, we recognize the image as a point defined over the Heegner field~$K$.
	This is explained in Section~\ref{sec:Computing-Heegner-points-analytically}.
	However, $J(K) \subset J(\C)$ is dense. Hence we must prove that we have found the correct point.
	We do this by determining its canonical height (which is well-defined
	since it does not change when adding a torsion point or changing the sign)
	via the Gross--Zagier formula.
	Since there are only finitely many points with bounded height,
	knowing the height is sufficient to cut the possibilities
	down to finitely many candidates up to sign and torsion; in practice,
	there is only one candidate. So we check that the point we have
	computed has the correct height and that no other point (up to sign
	and torsion) has the same height up to the numerical precision used
	in the computation. In principle, we could use this approach to
	bypass the analytic computation of the Heegner point altogether
	and just recognize it from its height, but using both approaches
	provides an additional level of confirmation that our results are correct.
	A further benefit of computing the Heegner point (and its image
	under a generator of the endomorphism ring of~$J$) is that this provides
	us with generators of a finite-index subgroup of~$J(K)$. So in order to
	determine~$J(K)$ (which is necessary for the computation of the Heegner
	index~$I_{K,\pi}$), it then suffices to saturate the known subgroup,
	which means that we do not have to search for points first. This can
	save a considerable amount of time.
	
	Since the Gross--Zagier formula involves the Petersson norm of~$f^\sigma$
	for the various embeddings $\sigma \colon \Z[f] \inj \R$,
	we need a way to compute these Petersson norms; see
	Section~\ref{sec:computing-the-Petersson-norm}.
	To apply the Gross--Zagier formula, we project $y_K$, viewed as an
	element of $J_0(N)(K) \otimes_\Z \R$, to its various $\sigma$-components~$y_{K,\sigma}$
	(which have the property that $\Q(f)$ acts on them via the embedding~$\sigma$).
	The formula then gives an expression for~$\hhat(y_{K,\sigma})$,
	where $\hhat \colon J_0(N)(K) \otimes_\Z \R \to \R$ is the normalized canonical
	height on~$J_0(N)$ associated to twice the theta divisor.
	This is discussed in
	Section~\ref{sec:computing-the-height-of-a-heegner-point-using-gross-zagier}.
	Finally, we have to relate the height~$\hhat_J(y_{K,\pi})$ with respect
	to twice the theta divisor on~$J$ to the heights~$\hhat(y_{K,\sigma})$;
	this is done in Section~\ref{sec:Comparing-canonical-heights}.

	\subsection{Computational representation of Diagram~\eqref{Eqn:AJ diagram}}
	\label{sec:period matrices}
	
	For our computations, we need to represent $A_f^\dual$, $A_f$, $J$
	as complex tori and the isogenies between them. This is done as follows.
	
	\begin{definition}
		Let $A$ be an arbitrary abelian variety over~$\C$, of dimension~$g$. Associated
		to a $\C$-basis $\ul{\omega} = (\omega_1, \ldots, \omega_g)$ of~$\H^0(A, \Omega^1)$
		and a $\Z$-basis $\ul{\gamma} = (\gamma_1\, \ldots, \gamma_{2g})$ of the integral homology
		$\H_1(A(\C), \Z)$, there is the \emph{period matrix}
		\[ \Pi_A \defeq \Pi_{A,\ul{\omega},\ul{\gamma}}
		\defeq \Bigl(\int_{\gamma_j} \omega_i\Bigr)_{i,j} \in \C^{g \times 2g} \,. \]
		Its $2g$ columns generate a lattice~$\Lambda$, and $A(\C) \isom \C/\Lambda$ via
		$x \mapsto \bigl(\int_0^x \omega_i\bigr)_i + \Lambda$.
		We also write $\Lambda_A$ for~$\Lambda$ to indicate the associated abelian variety.
	\end{definition}
	
	\begin{definition}
		Let $A$ and~$B$ be two abelian varieties over~$\C$ of the same dimension~$g$,
		and let $\Pi_A$ and~$\Pi_B$ be associated period matrices. If $\phi \colon A \to B$
		is an isogeny, then there are uniquely determined matrices $\alpha_\phi \in \GL_g(\C)$
		and $M_\phi \in \Z^{2g \times 2g}$ such that
		\[ \alpha_\phi \cdot \Pi_A = \Pi_B \cdot M_\phi \,. \]
		We call $(\alpha_\phi, M_\phi)$ the \emph{pair of matrices associated to~$\phi$}.
	\end{definition}
	
	We observe that, given $\Pi_A$ and~$\Pi_B$, each of $M_\phi$ and~$\alpha_\phi$ can be
	determined from the other.
	
	Note that $\alpha_\phi$ is the matrix of the $\C$-linear map $\omega \mapsto \phi^* \omega$
	with respect to the bases of the spaces of holomorphic differentials used for
	$\Pi_A$ and~$\Pi_B$, and $M_\phi$ is the matrix of the $\Z$-linear map $\gamma \mapsto \phi_* \gamma$
	on the homology bases.
	
	If $A$, $B$ and~$\phi$ are defined over~$\Q$ and we use $\Q$-bases of $\H^0(A, \Omega^1)$
	and~$\H^0(B, \Omega^1)$, then $\alpha_\phi \in \GL_g(\Q)$ (since $\phi^*$ is a $\Q$-linear
	map). Similarly, if we use $\Z$-bases of $\H^0(\sA, \Omega^1)$ and $\H^0(\sB, \Omega^1)$,
	where $\sA$ and~$\sB$ are the Néron models of $A$ and~$B$ over~$\Z$, then
	$\alpha_\phi \in \Z^{g \times g} \cap \GL_g(\Q)$.
	
	\begin{definition} \label[definition]{c of an isogeny}
		Let $A$ and~$B$ be abelian varieties defined over~$\Q$, with Néron models
		$\sA$ and~$\sB$ over~$\Z$, respectively. Let $\phi \colon A \to B$ be
		an isogeny defined over~$\Q$. Then we set
		\[ c_\phi \defeq \bigl(\H^0(\sA, \Omega^1) : \phi^* \H^0(\sB, \Omega^1)\bigr) \in \Z_{\ge 1} \,. \]
	\end{definition}
	
	If $\Pi_A$ and~$\Pi_B$ are computed using $\Z$-bases of $\H^0(\sA, \Omega^1)$
	and~$\H^0(\sB, \Omega^1)$, then $c_\phi = |\det \alpha_\phi|$.
	
	\begin{definition} \label[definition]{Manin constant}
		Let $f \in \cN(N, g)$. We define $S_2(f, \Z)$ to be the $\Z$-sublattice of
		the $\C$-span of $f$ and its Galois conjugates in~$S_2(\Gamma_0(N))$ that
		consists of forms whose $q$-expansions have integral coefficients.
		Under the natural identification of $S_2(\Gamma_0(N))$ with
		$\H^0(X_0(N), \Omega^1) \iso \H^0(J_0(N), \Omega^1)$, the image of~$S_2(f, \Z)$
		contains $\pi_f^* \H^0(\sA_f,\Omega^1)$ (where $\sA_f$ is the Néron model
		of~$A_f$ over~$\Z$). The index
		\[ c_f \defeq \bigl(S_2(f,\Z) : \pi_f^* \H^0(\sA_f,\Omega^1)\bigr) \in \Z_{\ge 1} \]
		is the \emph{Manin constant} of~$\pi_f$.
	\end{definition}
	
	See~\cite[Def.~3.3 and Thm.~3.4]{ARS}.
	
	\begin{proposition}
		The Manin constant~$c_f$ is divisible only by primes~$p$ such that $p^2 \mid N$
		or $p = 2$ and the conductor of~$\Z[f]$ is even.
		
		In particular, $c_f = 1$ if $N$ is squarefree and the conductor of~$\Z[f]$ is odd.
	\end{proposition}
	
	\begin{proof}
		This is~\cite[Cor.~3.7]{ARS} for odd primes and~\cite[Thm.~5.19]{Cesnavicius2018}
		for $p = 2$.
	\end{proof}
	
	It has been conjectured (see~\cite[Conj.~5.2]{Cesnavicius2018} and the text
	preceding~it) that $c_f$ is always~$1$, but~\cite[Thm.~5.10]{Cesnavicius2018} gives
	a counterexample in dimension~$24$ (with $N = 431$ odd and $2 \mid c_f$).
	
	Magma can compute a period matrix~$\Pi_{A_f^\dual}$ of~$A_f^\dual$ with respect to
	a $\Z$-basis of~$S_2(f, \Z)$ and some homology basis. Magma also computes the matrix~$I$
	of the intersection pairing (inside the homology of~$J_0(N)(\C)$) on the first homology
	of~$A_f^\dual(\C)$. Then $\Pi_{A_f} \defeq \Pi_{A_f^\dual} \cdot I^{-1}$ is a period matrix
	for~$A_f$ (with respect to the same basis of~$S_2(f, \Z)$),
	and $(I_g, I)$ is the pair of matrices associated to the polarization~$\lambda_f$.
	Let now $J$ be the Jacobian of a curve of genus~$2$ over~$\Q$ such that there is an isogeny
	$\pi \colon A_f \to J$ as in~\eqref{Eqn:AJ diagram}. Magma can compute a period
	matrix~$\Pi_J$ for~$J$ with respect to a certain $\Q$-basis~$B$ of~$\H^0(J, \Omega^1)$
	(if $J$ is the Jacobian of a genus~$2$ curve $y^2 = f(x)$, then $B$ corresponds
	to the differentials $dx/y$ and~$x\,dx/y$ on the curve)
	and a symplectic homology basis. We can then find the associated pair of
	matrices~$(\alpha_\pi, M_\pi)$. The algorithm~\cite[Algorithm~13]{vanBommel2019}
	determines the \enquote{compensation factor} (called $W$ in \emph{loc.~cit.})
	\begin{equation} \label{Eq:compensation factor}
		C = \bigl(\H^0(\sJ, \Omega^1) : \langle B \rangle_{\Z}\bigr)
	\end{equation}
	(where $\sJ$ is the Néron model of~$J$ over~$\Z$ and the index of two commensurable
	$\Z$-lattices in~$\H^0(J, \Omega^1)$ is in general a positive rational number).
	Combining these computations gives the following.
	
	\begin{lemma} \label[lemma]{computation of c_f c_pi}
		\[ c_f \cdot c_\pi = C \cdot |\det \alpha_\pi| \,. \]
	\end{lemma}
	
	In our LMFDB examples, the compensation factor~$C$ (with respect
	to a minimal Weierstrass model) is always~$1$, and $c_f c_\pi$
	divides the degree of the isogeny~$\pi$.
	
	For later applications, we want to compute the sizes of the kernel and cokernel
	of the map~$\pi_\R \colon A_f(\R) \to J(\R)$ induced by the isogeny~$\pi$
	on the groups of real points. We note that we can obtain the action of
	complex conjugation~$\tau$ on $A(\C) \isom \H_1(A(\C), \Z) \otimes \R/\Z$ by solving
	$\overline{\Pi_A} = \Pi_A \cdot M_{A,\tau}$ for $M_{A,\tau} \in \Z^{2g \times 2g}$.
	We obtain $\ker \pi \isom M_\pi^{-1} \Z^{2g}/\Z^{2g} \iso \Z^{2g}/M_\pi \Z^{2g}$,
	and we can find its $\tau$-invariant part~$\ker \pi_\R$ using~$M_{A_f,\tau}$.
	The group~$\pi_0(J(\R))$ of connected components of~$J(\R)$ is isomorphic to
	\[ \ker(1+\tau \mid \Lambda_J)/(1 - \tau) \Lambda_J
	\isom \ker(I_{2g} + M_{J,\tau} \mid \Z^{2g})/(I_{2g} - M_{J,\tau}) \Z^{2g} \]
	and similarly for~$A_f$, so
	\[ \coker \pi_\R \isom
	\frac{\ker(I_{2g} + M_{J,\tau} \mid \Z^{2g})}%
	{(I_{2g} - M_{J,\tau}) \Z^{2g} + \ker(I_{2g} + M_{J,\tau} \mid M_\pi \Z^{2g})} \,. \]
	When considering the quadratic twist~$\pi^K$ for an imaginary quadratic field~$K$,
	then we have to replace $M_{J,\tau}$ by~$-M_{J,\tau}$ to obtain the twisted action
	of~$\tau$ on~$J^K(\C) \isom J(\C)$.
	
	We can use a similar idea to compute~$\Omega_J$ from the period
	matrix~$\Pi_J$ and the compensation factor~$C$ from~\eqref{Eq:compensation factor},
	as follows.
	
	\begin{lemma} \label[lemma]{computation of Omega_J}
		Let $T \in \GL_{2g}(\Z)$ be such that $(M_{J,\tau} + I_{2g}) \cdot T = (\tilde{M} \mid 0)$
		with $\tilde{M} \in \Z^{2g \times g}$. Then
		\[ \Omega_J = C \cdot |\det (\Pi_J \cdot \tilde{M})| \,. \]
	\end{lemma}
	
	\begin{proof}
		It is well-known that $\Omega_J$ is the covolume of the lattice given by
		integrating a Néron basis over the $\C|\R$-trace of~$\H_1(X(\C), \Z)$
		(see \cite[\S\,3.5]{FLSSSW} or~\cite[Def.~11]{vanBommel2019}).
		The $\Z$-lattice generated by the columns of the matrix $M_{J,\tau} + I_{2g}$
		corresponds to the $\C|\R$-trace of~$\H_1(X(\C), \Z)$
		(w.r.t.\ the homology basis used to compute~$\Pi_J$).
		This lattice is known to have rank~$g$, so there exists a unimodular matrix~$T$
		as in the statement, and multiplying on the right by~$T$ preserves the lattice.
		So the columns of~$\tilde{M}$ give a basis of the $\C|\R$-trace of~$\H_1(X(\C), \Z)$,
		and the result follows (taking into account the factor~$C$ arising from
		changing the basis of differentials used in the computation of~$\Pi_J$ to
		a Néron basis).
	\end{proof}
	
	This improves over the method currently implemented in Magma (which is based
	on~\cite[Algorithm~13]{vanBommel2019}) in that it uses an exact computation
	to find the correct integral linear combination of $g \times g$ minors of~$\Pi_J$
	instead of relying on a \enquote{real gcd} computation with numerical approximations.
	(The approximation step is in the computation of~$M_{J,\tau}$, but here we
	know that the entries are integers, so we can simply round.)
	
	We now want to determine the endomorphism $\alpha \in \cO = \End J$ that
	was defined in~\eqref{Eq:def alpha}. Since $\Pi_J$ is
	computed with respect to a symplectic homology basis, we obtain the $M$-matrix
	of the canonical polarization of~$J$ as the matrix~$I'$ of the standard
	symplectic pairing. Then $M \defeq M_\pi \cdot I \cdot M_\pi^\top \cdot I'$ gives
	the action of $\alpha \in \cO$ on the lattice associated to~$J$, and its action
	on the tangent space can be recovered from that. (Recall that $I$
	denotes the matrix of the intersection pairing for~$A_f^\dual(\C)$.)
	We can (and do) \enquote{optimize}~$\alpha$ by post-composing $\pi$ with an
	automorphism $\epsilon \in \End_\Q(J)^\times$ (this has the effect of
	multiplying $\alpha$ by~$\epsilon^2$) in the sense that we minimize
	the images $\alpha^\sigma \in \R_{>0}$ under the real embeddings of~$\cO$
	(in practice, we minimize the trace of~$\alpha$).
	This leads to potentially smaller Heegner points on~$J$, which simplifies
	some of the computations.

	\subsection{Determining Heegner fields}
	\label{sec:Heegner fields}
	
	To be able to use the results of~\cite{GrossZagier1986} and some other
	results that require the discriminant of the Heegner field to be odd, we
	restrict to odd discriminants in the following.
	
	We find a Heegner field~$K$ by enumerating the odd discriminants~$-D$ of imaginary
	quadratic number fields with the property that all prime divisors of~$N$
	split completely in~$\cO_K$ (this can be checked easily by computing Legendre
	symbols). The condition $\ord_{s=1} L(f/K, s) = 1$ is equivalent to
	$L(f \otimes \epsilon_K, 1) \neq 0$ when $\Lrk J = 1$ and to
	$L'(f \otimes \epsilon_K, 1) \neq 0$ when $\Lrk J = 0$.
	Using modular symbols as described in~\cite[\S\,{2.8}]{Cremona1997}, we can decide
	whether $L(f \otimes \epsilon_K, 1) = 0$ or not. The non-vanishing
	of~$L'(f \otimes \epsilon_K, 1)$ can be proved by computing it to a high enough
	precision using Dokchitser's Magma implementation~\cite{Dokchitser2004}.
	Alternatively and in practice (because the evaluation of the twisted $L$-value
	can take fairly long when $N$ is large),
	we can compute the Heegner point for a
	given~$K$; if it is non-torsion, then $K$ is a suitable Heegner field.

	\subsection{Computing the endomorphism ring and its action on Mordell--Weil groups}
	\label{sec:Computing-the-endomorphism-ring}
	
	We need to determine the endomorphism ring~$\cO$ of the Jacobian~$J$
	and how it acts on the Mordell--Weil group~$J(\Q)$ or~$J^K(\Q)$, or,
	more generally, on~$J(L)$ for some number field~$L$.
	For this, we compute a numerical approximation to the big period matrix
	as in Section~\ref{sec:period matrices}; potential endomorphisms can be guessed
	from this information. To verify that the presumed endomorphism ring
	is the correct one, we can use data from the LMFDB~\cite{lmfdb}.
	(Alternatively, one could use~\cite{Lombardo2019}.)
	This shows that the numerical endomorphisms are close to
	actual endomorphisms and thus gives us
	a representation of~$\cO$ as a subring of a matrix algebra over~$\Z$,
	together with its action on the complex torus $\C^2/\Lambda \isom J(\C)$.
	To compute the action of~$\cO$ on~$J(L)$, we use an improved version
	of Magma's \texttt{(To/From)AnalyticJacobian} to convert between
	points in~$J(\C)$ in Mumford representation and representatives in~$\C^2$.
	(The improvement also handles points at infinity and Weierstrass points.)
	For a generator~$\gamma$ of~$\cO$ and each generator~$x$ of~$J(L)$, we map
	$x$ to~$\C^2/\Lambda$, apply $\gamma$ to the image, and map back to~$J(\C)$.
	We then recognize the coefficients of the Mumford representation as elements
	of~$L$ using Magma's \texttt{MinimalPolynomial} and check that the coefficients
	we recognize really define a point in~$J(L)$. We then write the resulting
	point as a linear combination of the generators of~$J(L)$. In this way,
	we obtain a matrix giving the action of~$\gamma$ on~$J(L)$ with respect
	to the chosen generators. We can bound the height of $\gamma \cdot x$
	by $\max_\sigma |\gamma^\sigma|^2$ times the height of~$x$, so there are
	only finitely many candidates for $\gamma \cdot x$, which allows us to
	determine $\gamma \cdot x$ exactly by computing with sufficient precision.
	As an additional check, we verify that the matrix we obtain has the same minimal
	polynomial as~$\gamma$.

	\subsection{Computing Heegner points analytically}
	\label{sec:Computing-Heegner-points-analytically}
	
	Recall that $\pi$ denotes the isogeny $A_f \to J$ and that we want
	to compute the Heegner point
	\[ y_{K,\pi} = \pi_J(y_K) = \pi(y_{K,A_f}) \in J(K) \,. \]
	In this section we explain how to find~$y_{K,\pi}$ explicitly.
	Also recall that $\bpi$ denotes the real number giving the area
	of the unit disk, to avoid confusion with the isogeny~$\pi$.
	
	We obtain a computational representation of the isogeny
	$\pi \colon A_f \to J$ as described in Section~\ref{sec:period matrices}.
	Since we know that an isogeny has to exist, we can be sure that
	what we obtain indeed describes an actual isogeny.
	
	To find the Heegner point on~$J$,
	we first determine the integral binary quadratic \enquote{Heegner forms}
	associated to~$K$ and~$N$. These are representatives of the $h_K$ classes
	of positive definite binary quadratic forms of discriminant~$D_K$ such
	that their roots $\tau \in \bH$ map to the points in the Heegner cycle~$\bx_K$.
	Let $(f_1, \ldots, f_g)$ be the $\Z$-basis of~$S_2(f, \Z)$ that
	is used for the computation of the big period matrix of~$A_f$ (it can be
	obtained via the Magma function \texttt{qIntegralBasis}).
	We then compute the period integrals
	\[ P(\tau, j) \defeq 2 \bpi i \int_{i \infty}^\tau f_{j}(z)\,dz
	= \int_0^{e^{2 \bpi i \tau}} \sum_{n \ge 1} a_n(f_j) q^n \, \frac{dq}{q}
	= \sum_{n \geq 1} \frac{a_n(f_j)}{n} e^{2 \bpi i \tau n} \in \C \]
	for each of these roots~$\tau$ and each $1 \le j \le g$ to the desired precision.
	(We pick our Heegner forms in such a way that $\tau$ has imaginary part as large as possible.
	We can use the bound $|a_n^\sigma| \le \sqrt{3} n$
	(see~\cite[Lemma~2.9]{GJPST}, where a bound $|a_n^\sigma| \le n$ is claimed,
	but their argument bounding $|a_{p^m}^\sigma|$ is not correct for powers of $2$ or~$3$)
	and the representation of~$f_j$ as a linear combination of the~$f^\sigma$
	to determine the number of terms we need.
	The points $y_\tau = \bigl(P(\tau, j)\bigr)_{j} \in \C^g/\Lambda_f \iso A_f(\C)$
	then are of the form $[x_\tau - (\infty)]$ projected to~$A_f$, where $\infty$
	is the cusp at infinity and $x_\tau$ runs through the points in the support
	of the Heegner cycle~$\bx_K$. In particular, we have $\sum_\tau y_\tau = y_{K,A_f}$.
	
	We then use the matrix~$\alpha_\pi$ associated to the isogeny~$\pi$ to map~$y_K$
	or all the points~$y_\tau$ to~$\C^g/\Lambda \isom J(\C)$. We apply the
	numerical inverse of the Abel--Jacobi map to find the Mumford representation
	of this or these points as points on~$J$. We then try to recognize the
	coefficients in the Mumford representation as elements of~$K$
	(for~$\pi(y_{K, A_f}) = y_{K,\pi}$) or of~$H$ (for~$\pi(y_\tau)$) and check that
	this really gives rise to a point in~$J(K)$ or~$J(H)$.
	
	Using the action of~$\cO$ on~$J(K)$ that we have determined in
	Section~\ref{sec:Computing-the-endomorphism-ring}, we can then determine the
	$\cO$-span of~$y_{K,\pi}$ and from it the ideal
	\[ \cI_{K,\pi} = \Char_\cO(J(K)/\cO y_{K,\pi}) \]
	and the index $I_{K,\pi} = (J(K) : \cO y_{K,\pi})$.
	The corresponding index for~$A_f$ is $I_K = (A_f(K) : \Z[f] \cdot y_{K,A_f})$.
	We can express it in terms of~$I_{K,\pi}$ via
	\[ I_K = \frac{(\cO_{\Q(f)} : \Z[f])}{(\cO_{\Q(f)} : \cO)}
	\cdot \frac{\#A_f[\pi](K)}{\#\bigl(J(K)/\pi(A_f(K))\bigr)} \cdot I_{K,\pi} \,. \]
	The first factor on the right takes care of the fact that $\Z[f]$ and
	$\cO = \End_{\Q}(J)$ can be different orders in~$\Q(f)$;
	it can be computed easily as $\sqrt{\disc \Z[f]/\disc \cO}$. The second factor
	captures the effect of the isogeny~$\pi$. Note that the second factor
	can be multiplicatively bounded from above by~$\#A_f[\pi](K) \mid \deg \pi$,
	which gives a multiplicative upper bound for~$I_K$ as well.
	(We may get a better bound than $\deg \pi$ from bounding
	$\#A_f(K)_\tors$ using the coefficients of the $L$-series of~$A_f/K$.)
	
	We will need $I_{K,\pi}$ in Section~\ref{SS: Sha_an/K}
	for the computation of~$\Sha(J/\Q)_\an$ when $\Lrk J = 1$ and in
	Section~\ref{ssec:finite support} for the determination of an explicit
	finite support of~$\Sha(J/\Q)$. We will need (a multiplicative
	upper bound for)~$I_K$ in Section~\ref{sec:Iwasawa theory}.

	\subsection{Computing the Petersson norm of a newform}
	\label{sec:computing-the-Petersson-norm}
	
	Let $\sigma \colon \Q(f) \inj \R$ be an embedding. For the Gross--Zagier formula,
	we need the Petersson norms of the conjugates $f^\sigma$ for the various
	possible~$\sigma$. We identify $X_0(N)(\C)$ with $\Gamma_0(N) \backslash \bH^*$
	(where $\bH^* = \bH \cup \bP^1(\Q)$ is the upper half-plane together with the cusps)
	and use the normalization
	\[ \|f^\sigma\|^2 = \int_{X_0(N)(\C)} |f^\sigma(x + yi)|^2\, dx \wedge dy \]
	for the Petersson norm of $f^\sigma \in S_2(\Gamma_0(N), \C)$ as
	in~\cite[(5.1)]{GrossZagier1986}.
	(Sometimes this is normalized differently by dividing by the volume~$\mu(X_0(N)(\C))$
	to make it independent of the choice of~$N$.)
	We compute the Petersson norm by relating it to the
	\emph{symmetric square $L$-function} $L(\Sym^2 f^\sigma, s)$.
	
	If an $L$-function $L(\cX, s)$ has an Euler product expansion, we write it as
	\[ L(\cX, s) = \prod_\ell L_\ell(\cX, \ell^{-s})^{-1} \,, \]
	where $L_\ell(\cX, T) \in R[T]$ (with $R$ the coefficient ring of the $L$-function)
	is the Euler polynomial at~$\ell$.
	
	We define the symmetric square $L$-function $L(\Sym^2 f, s)$ as the $L$-function
	associated to the strictly compatible system~$(\Sym^2 \rho_{\frp^\infty, f})$
	of $\frp$-adic Galois representations. For a prime $\ell \nmid N$, write
	\[ L_\ell(f, T) = 1 - a_\ell T + \ell T^2 = (1 - \alpha_\ell T) (1 - \beta_\ell T) \,; \]
	then
	\begin{equation} \label{E:LSym2}
		\begin{array}{r@{{}={}}l}
			L_\ell(\Sym^2 f, T)
			& (1 - \alpha_\ell^2 T) (1 - \alpha_\ell \beta_\ell T) (1 - \beta_\ell^2 T) \\[6pt]
			& (1 - \ell T) \bigl((1 + \ell T)^2 - a_\ell^2  T\bigr) .
		\end{array}
	\end{equation}
	We define the \emph{imprimitive symmetric square $L$-function}
	$\tilde{L}(\Sym^2 f, s)$ by this formula
	for the Euler polynomial at \emph{all} primes (then we take $\alpha_\ell = a_\ell$
	and $\beta_\ell = 0$ when $\ell \mid N$); compare~\cite[p.~110]{CoatesSchmidt1987}.
	(The difference is whether we take $I_\ell$-coinvariants
	before ($\tilde{L}(\Sym^2 f, s)$) or after ($L(\Sym^2 f, s)$)
	applying $\Sym^2$ when defining the Euler polynomials.)
	This imprimitive version is what Shimura denotes~$D(s)$ in~\cite{Shimura1975}.
	
	We thank user334725 on MathOverflow~\cite{410483} for pointers to the
	relevant literature.
	
	\begin{proposition} \label[proposition]{Petersson L-tilde}
		Let $f \in S_2(\Gamma_0(N), \C)$ be a normalized eigenform.
		Then the Petersson norm of~$f$ is given by
		\[ \|f\|^2 = \frac{N}{8 \bpi^{3}} \cdot \tilde{L} (\Sym^2 f, 2) \,. \]
	\end{proposition}
	
	\begin{proof}
		Denote the Fourier coefficients of $f$ by~$a_n$. We set
		\[ D(f, s) \defeq \sum_{n \ge 1} \frac{a_n^2}{n^s} \,. \]
		By~\cite[Satz~6]{Petersson1949} (and taking into account the different normalization;
		compare also~\cite[Eq.~(2.5)]{Shimura1976}),
		\begin{align}
			\|f\|^2
			&= [\Gamma(1) : \Gamma_0(N)] \, \frac{\bpi}{3} \, \frac{1}{(4\bpi)^2} \, \res_{s=2} D(f, s) \nonumber \\
			&= N \prod_{\ell \mid N} \Bigl(1 + \frac{1}{\ell}\Bigr) \frac{1}{48 \bpi} \, \res_{s=2} D(f, s) .
			\label{Pn_formula1}
		\end{align}
		By~\cite[Eq.~(0.4)]{Shimura1975} (see \cite[Lemma~1]{Shimura1976} for the relation between
		the Euler factors), we have the following equality, where the superscript~$N$
		means that we leave out the Euler factors coming from prime divisors of~$N$.
		\[ D(f, s) = \frac{\zeta^N(s-1)}{\zeta^N(2s - 2)} \, \tilde{L}(\Sym^2 f, s) \,. \]
		Taking the residue at $s = 2$ on both sides, we obtain
		\[ \res_{s=2} D(f, s)
		= \frac{6}{\bpi^2} \prod_{\ell \mid N} \Bigl(1 + \frac{1}{\ell}\Bigr)^{-1} \, \tilde{L}(\Sym^2 f, 2) \,, \]
		which gives the desired result when used in~\eqref{Pn_formula1}.
	\end{proof}
	
	So we need to compute $\tilde{L}(\Sym^2 f, 2)$. However, we cannot directly
	do that since $\tilde{L}(\Sym^2 f, 2)$ does not in general satisfy a suitable
	functional equation (which is needed to obtain a reasonably fast converging
	series for the value via a Mellin transform).
	We can, however, compute~$L(\Sym^2 f, 2)$, if we know its Euler
	factors at primes dividing~$N$. So we need
	to determine these Euler factors; combining this with~\eqref{E:LSym2} will also
	tell us what the correction factor $\tilde{L}(\Sym^2 f, 2)/L(\Sym^2 f, 2)$ is.
	
	For a prime~$\ell$, we set
	\[ C_\ell \defeq \frac{L_\ell(\Sym^2 f, \ell^{-2})}{\tilde{L}_\ell(\Sym^2 f, \ell^{-2})} \,. \]
	
	\begin{corollary} \label[corollary]{Petersson formula}
		Let $f$ and~$C_\ell$ be as above. Then
		\[ \|f\|^2 = \frac{N}{8 \bpi^3} \prod_{\ell^2 \mid N} C_\ell \cdot L(\Sym^2 f, 2) \,. \]
		In particular,
		\[ \|f\|^2 = \frac{N}{8 \bpi^3} \cdot L(\Sym^2 f, 2) \]
		when the level~$N$ is squarefree.
	\end{corollary}
	
	\begin{proof}
		This follows from~\cref{Petersson L-tilde} and the definition of~$C_\ell$,
		together with the fact that $C_\ell = 1$ unless $\ell^2 \mid N$, which will be
		shown in~\cref{L vs L-tilde 1} below.
	\end{proof}
	
	Alternative algorithms for computing the Petersson inner product are described
	in~\cite{Cohen2013} and have been implemented in Pari.
	
	Using the formula in~\cref{Petersson formula}, we can compute $\|f^\sigma\|^2$
	using~\cite{Dokchitser2004}
	for all $\sigma \in \Sigma$ if we can determine the Euler factors $L_\ell(\Sym^2 f, T)$
	for the primes $\ell \mid N$. We will do that in the following subsection.

	\subsection{Euler factors of the symmetric square $L$-function}
	\label{sec:LSym2}
	
	In this section, we explain how to find the Euler factors of~$L(\Sym^2 f, s)$
	at primes~$\ell$ dividing the level~$N$ of~$f$. In~\cite[\S\,1]{CoatesSchmidt1987},
	analogous statements are shown for the $L$-function of an elliptic curve, and
	\cite{Schmidt1988} has similar results stated in a somewhat different language.
	
	Recall that when $\ell \nmid N$,
	$L_\ell(\Sym^2 f, T) = \tilde{L}_\ell(\Sym^2 f, T)$ (and hence \hbox{$C_\ell = 1$}),
	and we can write down $\tilde{L}_\ell(\Sym^2 f, T)$ easily in terms of~$L_\ell(f, T)$.
	We now consider the case $v_\ell(N) = 1$.
	
	\begin{lemma} \label[lemma]{L vs L-tilde 1}
		Assume that $v_\ell(N) = 1$. Then
		\[ L_\ell(\Sym^2 f, T) = \tilde{L}_\ell(\Sym^2 f, T) = 1 - T \,. \]
		In particular, $C_\ell = 1$ whenever $\ell^2 \nmid N$.
	\end{lemma}
	
	\begin{proof}
		Fix some $\frp \nmid N \ell$ and set $p = p(\frp)$ and $V_\frp = V_\frp(A_f)$.
		Since $v_\ell(N) = 1$, $V_\frp$ has a one-dimensional quotient of $I_\ell$-coinvariants.
		Since $\chi_{p^\infty} = \det \circ \rho_{\frp^\infty}$ is trivial on~$I_\ell$, it follows that
		$\rho_{\frp^\infty}(I_\ell)$ lands in a unipotent subgroup
		(compare~\cref{inertia unipotent if semistable}). The image is nontrivial,
		since otherwise $\ell \nmid N$. Since $\ell \neq p$ and the unipotent subgroups
		of~$\GL_2(\Q(f)_\frp)$ are pro-$p$ groups, whereas the wild inertia at~$\ell$
		is a pro-$\ell$ group, it follows that $\rho_{\frp^\infty}|_{I_\ell}$ factors
		through the tame inertia group $I_\ell^\tame$.
		By~\cref{Frob on Iptame}, conjugating by any lift~$\Frob_\ell$ of the Frobenius
		automorphism to~$I_\ell^\tame$ has the effect of raising to the $\ell$th power.
		Using that $\det \rho_{\frp^\infty}(\Frob_\ell) = \chi_{p^\infty}(\Frob_\ell) = \ell$,
		it follows that (with respect to a suitable $\Q(f)_\frp$-basis),
		$\rho_{\frp^\infty}(\Frob_\ell) = \pm \left(\begin{smallmatrix} \ell & 0 \\ 0 & 1 \end{smallmatrix}\right)$
		and $\rho_{\frp^\infty}|_{I_\ell} \subseteq \left(\begin{smallmatrix} 1 & * \\ 0 & 1 \end{smallmatrix}\right)$.
		This implies that $L_\ell(f, T) = 1 \mp T$, and so $\tilde{L}_\ell(\Sym^2 f, T) = 1 - T$.
		We also see that
		\[ \Sym^2 \rho_{\frp^\infty}|_{\GalQl}
		= \begin{pmatrix} \chi_{p^\infty}^2 & * & * \\ 0 & \chi_{p^\infty} & * \\ 0 & 0 & 1 \end{pmatrix}
		\]
		with a one-dimensional $I_\ell$-coinvariant quotient, on which $\Frob_\ell$
		acts trivially. This shows that $L_\ell(\Sym^2 f, T) = 1 - T$ as well.
		(This is analogous to~\cite[Lemma~1.2]{CoatesSchmidt1987}.)
		
		In particular, $C_\ell = 1$, which, together with the discussion preceding
		this lemma, gives the last claim.
	\end{proof}
	
	We now consider the case $\ell^2 \mid N$. We first note that $L(\Sym^2 f, s)$ does not
	change under quadratic twists.
	
	\begin{lemma} \label[lemma]{LSym2 qu twist}
		Let $\tilde{f}$ be a quadratic twist of~$f$. Then
		\[ L(\Sym^2 \tilde{f}, s) = L(\Sym^2 f, s) \,. \]
	\end{lemma}
	
	\begin{proof}
		We consider the Euler factor at~$\ell$. Fix some $\frp \nmid N \ell$ and set $p = p(\frp)$
		and $V_\frp = V_\frp(A_f)$. Let $\epsilon$ be the quadratic character such that
		$\tilde{f} = f \otimes \epsilon$. Since the canonical group homomorphism
		$\GL(V_\frp) \to \GL(\Sym^2 V_\frp)$ is trivial on $\pm \id$,
		it follows that $\Sym^2 (\rho_{\frp^\infty} \otimes \epsilon) = \Sym^2 \rho_{\frp^\infty}$,
		which, upon restricting to~$\GalQl$, directly translates into
		$L_\ell(\Sym^2 \tilde{f}, T) = L_\ell(\Sym^2 f, T)$. The claim follows.
	\end{proof}
	
	The argument in the proof together with the fact that
	$\rho_{\frp^\infty}(I_\ell) \subseteq \SL(V_\frp)$
	shows that the action of~$I_\ell$ on~$\Sym^2 V_\frp$
	depends only on the projective image $\P \rho_{\frp^\infty}(I_\ell) \subseteq \PSL(V_\frp)$.
	We will see that the dimension of the $I_\ell$-coinvariants of~$\Sym^2 V_\frp$
	depends on whether this projective image is abelian or not.
	
	\begin{lemma} \label[lemma]{PIell abelian}
		Let $k$ be a field of characteristic zero and let $V$ be a two-dimensional
		$k$-vector space. Let $G \subseteq \SL(V)$ be such that $\P G$ is not unipotent.
		\begin{enumerate}[\upshape(1)]
			\item If $\Sym^2 V$ has a nontrivial $G$-invariant quotient, then $G$ is abelian.
			\item If $\P G$ is abelian, then $\bar{V} \defeq V \otimes_k \bar{k}$ has a basis $e_1$, $e_2$ consisting of
			simultaneous eigenvectors for the elements of~$G$. The $G$-coinvariant
			space of~$\Sym^2 \bar{V} = \langle e_1^2, e_1 e_2, e_2^2 \rangle$
			is one-dimensional and is isomorphic to the direct summand $\bar{k} \cdot e_1 e_2$.
		\end{enumerate}
	\end{lemma}
	
	\begin{proof}
		We can assume without loss of generality that $k$ is algebraically closed.
		\begin{enumerate}[(1)]
			\item If $\Sym^2 V$ has a nontrivial $G$-invariant quotient, then
			$\Sym^2 V^*$ has a non-zero $G$-invariant element, which is a quadratic
			form~$q$ on~$V$. Then $G$ must fix the zero set of~$q$ in~$\P^1(k)$.
			This zero set can have either one or two elements.
			
			In the first case,
			$G$ fixes a point in~$\P^1$, hence is contained in a Borel subgroup,
			so, after fixing a suitable basis, the associated representation~$\rho$
			has the form
			$\left(\begin{smallmatrix} \chi & \alpha \\ 0 & \chi^{-1} \end{smallmatrix}\right)$
			with a character~$\chi$ such that $\chi^2 \neq \mathbf{1}$ (recall that $\P G$
			is not unipotent). Then (with the columns giving the action on $X^2$, $XY$, $Y^2$,
			when $X, Y$ is the given basis with
			$\rho(g) X = \chi(g) X$, $\rho(g) Y = \alpha(g) X + \chi^{-1}(g) Y$)
			\[ \Sym^2 \rho = \begin{pmatrix}
				\chi^2 & \alpha \chi & \alpha^2 \\
				0 &  \mathbf{1} & 2 \alpha \chi^{-1} \\
				0 &           0 & \chi^{-2}
			\end{pmatrix} ,
			\]
			and this has nontrivial $G$-coinvariants only when $\alpha = 0$,
			which implies that $G$ is abelian.
			
			In the second case, $G$ is contained in the normalizer of a Cartan
			subgroup, so its elements are (with respect to a suitable basis $(e_1, e_2)$)
			either of the form
			$\left(\begin{smallmatrix} a & 0 \\ 0 & a^{-1} \end{smallmatrix}\right)$
			(which fix $e_1 \cdot e_2$) or of the form
			$\left(\begin{smallmatrix} 0 & -a \\ a^{-1} & 0 \end{smallmatrix}\right)$.
			However, the elements of the second form send $e_1 \cdot e_2$ to its
			negative, so (noting that there must be elements of the first form
			with $a^2 \neq 1$, again since $\P G$ is not unipotent, so that neither
			$e_1^2$ nor $e_2^2$ can be fixed by~$G$) such elements
			cannot be present in~$G$, which again implies that $G$ is abelian.
			\item If $\P G$ is abelian, then so is~$G$. The representation on~$V$ then
			splits as a sum of two characters $\chi$ and~$\chi^{-1}$ such that
			$\chi^2 \neq 1$. Let $e_1$ and~$e_2$ be corresponding eigenvectors.
			Then $\Sym^2 V$ splits as $\chi^2 \oplus \mathbf{1} \oplus \chi^{-2}$,
			with the $G$-action on $e_1 e_2$ being trivial. This shows the claim.
			\qedhere
		\end{enumerate}
	\end{proof}
	
	This leads to the following classification.
	
	\begin{lemma} \label[lemma]{L vs L-tilde}
		Let $\ell$ be a prime such that $\ell^2 \mid N$. Then $\tilde{L}_\ell(\Sym^2 f, T) = 1$.
		Let $\tilde{f}$ be a quadratic twist of~$f$ whose level~$\tilde{N}$ is
		(multiplicatively) minimal. Fix a regular prime ideal $\frp \nmid N \ell$ of~$\Z[f]$.
		\begin{enumerate}[\upshape(1)]
			\item $\P \rho_{\frp^\infty}(I_\ell)$ is trivial if and only if $\ell \nmid \tilde{N}$.
			In particular, using~\eqref{E:LSym2},
			\[ C_\ell = \frac{(\ell - 1)((\ell + 1)^2 - a_\ell(\tilde{f})^2)}{\ell^3} \,. \]
			\item $\P \rho_{\frp^\infty}(I_\ell)$ is nontrivial and unipotent if and only if
			$v_\ell(\tilde{N}) = 1$.
			In this case, $L_\ell(\Sym^2 f, T) = 1 - T$ and the conductor exponent of~$\Sym^2 f$
			at~$\ell$ is~$2$.
			In particular, $C_\ell = \frac{\ell^2 - 1}{\ell^2}$.
			\item \label{LLcase3}
			$\P \rho_{\frp^\infty}(I_\ell)$ is abelian and not unipotent if and only if
			$v_\ell(\tilde{N}) = 2$. In this case,
			$L_\ell(\Sym^2 f, T) = 1 \mp \ell T$, with the negative sign if and only if
			$\rho_{\frp^\infty}(\GalQl)$ is abelian, and the conductor exponent of~$\Sym^2 f$
			at~$\ell$ is~$2$.
			In particular, $C_\ell = \frac{\ell \mp 1}{\ell}$ (with the same sign).
			\item \label{LLcase4}
			$\P \rho_{\frp^\infty}(I_\ell)$ is non-abelian if and only if $v_\ell(\tilde{N}) \ge 3$.
			In this case, we have $L_\ell(\Sym^2 f, T) = 1$, and the
			conductor exponent of~$\Sym^2 f$ at~$\ell$ is at least~$4$
			and at most $2v_\ell(\tilde{N}) - 1$.
			In particular, $C_\ell = 1$.
		\end{enumerate}
	\end{lemma}
	
	\begin{proof}
		When $\ell^2 \mid N$, then the space of $I_\ell$-coinvariants is trivial,
		hence so is its symmetric square. This means that $\tilde{L}_\ell(\Sym^2 f, T) = 1$.
		By~\cref{LSym2 qu twist}, we have that $L(\Sym^2 f, s) = L(\Sym^2 \tilde{f}, s)$.
		It suffices to show the \enquote{only if} direction of the equivalences at the beginning
		of each statement, since the consequences exhaust all possibilities disjointly.
		\begin{enumerate}[(1)]
			\item If $\P \rho_{\frp^\infty}(I_\ell)$ is trivial, then
			$\rho_{\frp^\infty}|_{I_\ell}$ is of the form
			$(\mathbf{1} \oplus \mathbf{1}) \otimes \epsilon$ with a
			quadratic character~$\epsilon$. Twisting by~$\epsilon$ makes
			the representation unramified at~$\ell$, so $\ell \nmid \tilde{N}$
			(using that $\tilde{N}$ is minimal). The statement on~$C_\ell$
			then follows.
			\item If the projective image is nontrivial and unipotent,
			then there is a quadratic character~$\epsilon$ such that
			$\rho_{\frp^\infty} \otimes \epsilon|_{I_\ell}$ is unipotent and nontrivial,
			which implies that $v_\ell(\tilde{N}) = 1$. The statement on the Euler
			factor then follows from~\cref{L vs L-tilde 1}.
			\item We assume that $\P \rho_{\frp^\infty}(I_\ell)$ is abelian, but not unipotent.
			By~\cref{PIell abelian}, $\Sym^2 V_\frp$ has a one-dimensional
			$I_\ell$-coinvariant space. Also, because $\rho_{\frp^\infty}|_{I_\ell}$
			is a sum of two characters of order coprime to~$\ell$, the representation
			factors through the tame inertia group, hence there is no wild part
			in the conductor. This shows that $v_\ell(\tilde{N}) = 2$ and that
			the conductor exponent of~$\Sym^2 f = \Sym^2 \tilde{f}$ is $3 - 1 = 2$.
			Also by~\cref{PIell abelian}, the coinvariant space corresponds
			to the tensor product of the two one-dimensional representations
			in the splitting of~$V_\frp$. Frobenius either fixes each of these
			two one-dimensional spaces, in which case its action on the tensor
			product is by $\det(\rho_{\frp^\infty}(\Frob_\ell)) = \ell$; then
			$L_\ell(\Sym^2 f, T) = 1 - \ell T$, and $\rho_{\frp^\infty}(\GalQl)$
			is abelian. Or else $\Frob_\ell$ swaps the two spaces; then it
			acts by the negative of the determinant (compare the proof
			of~\cref{PIell abelian}), so $L_\ell(\Sym^2 f, T) = 1 + \ell T$,
			and $\rho_{\frp^\infty}(\GalQl)$ is non-abelian.
			\item We assume that $\P \rho_{\frp^\infty}(I_\ell)$ is non-abelian.
			This implies that it is not unipotent.
			By~\cref{PIell abelian}, the $I_\ell$-coinvariant space
			of~$\Sym^2 V_\frp$ is trivial. This shows that $L_\ell(\Sym^2 f, T) = 1$.
			Also, $\P \rho_{\frp^\infty}|_{I_\ell}$
			cannot factor through the tame inertia group since the latter
			is abelian. So there must be wild ramification at~$\ell$ both
			in~$\P \rho_{\frp^\infty}$ and in~$\Sym^2 \rho_{\frp_\infty}$.
			As the tame parts of the conductor exponents of these two are given
			by $\dim V_\frp = 2$ and $\dim \Sym^2 V_\frp = 3$, respectively,
			it follows that $\ell^3 \mid \tilde{N}$ and the conductor exponent~$c$
			of~$\Sym^2 f$ is at least~$4$.
			To obtain the claimed upper bound,
			we observe that when $\rho \colon G \to \GL(V)$ is a $2$-dimensional
			representation, then the codimension of the invariant subspace
			of~$\Sym^2 \rho$ is at most twice the codimension of the invariant
			subspace of~$\rho$. Then~\cite[Eq.~(1.2.1)]{Serre1987} implies
			that the wild part $c - 3$ of~$c$
			is at most twice the wild part $v_\ell(\tilde{N}) - 2$ of the
			conductor exponent of~$\tilde{f}$. This gives the desired bound.
			\qedhere
		\end{enumerate}
	\end{proof}
	
	Given~$f$, a choice of~$\tilde{f}$ can be obtained from the LMFDB.
	Alternatively, the conductor~$d$ of the twisting character~$\epsilon$ must
	satisfy $d^2 \mid N$, so we can check the finitely many possibilities for~$\epsilon$
	and compare the resulting levels to find~$\tilde{f}$.
	
	Which of the two possibilities for the Euler factor in case~\eqref{LLcase3}
	is correct and what the correct choice of conductor exponent is in
	case~\eqref{LLcase4} can be checked by trying all possibilities
	and determining which one is compatible with the functional equation.
	Using the function \texttt{SymmetricPower} that Magma provides for
	constructing symmetric power $L$-functions seems to result in fairly slow
	code. Instead, we compute the relevant number of coefficients ourselves
	and use this coefficient sequence when constructing the $L$-series, which
	is then used for testing the functional equation and evaluating at $s = 2$.

	\subsection{Computing the height of a Heegner point using the Gross--Zagier formula}
	\label{sec:computing-the-height-of-a-heegner-point-using-gross-zagier}
	
	To state the Gross--Zagier formula, we need to introduce some more notation.
	Let $K$ be a Heegner field for~$f$. Recall that $H$ denotes the Hilbert class
	field of~$K$. The Heegner cycle $\bx_K$ on~$X_0(N)$ and the Heegner point
	$y_K = [\bx_K - h_K \cdot (\infty)] \in J_0(N)(K)$ have been defined in the introduction
	to this section. Recall that $I_f = \Ann_{\bT}(f)$.
	
	The action of~$\bT$ (or its quotient $\End_\Q(J_0(N))$) on~$J_0(N)(K)$ extends
	to a linear action on the real vector space $J_0(N)(K) \otimes_\Z \R$. Since
	the center~$Z$ of~$\End_\Q(J_0(N))$ is an order in a totally real étale $\Q$-algebra,
	we obtain a canonical decomposition
	\[ J_0(N)(K) \otimes_\Z \R = \bigoplus_{\sigma \colon Z \inj \R} J_0(N)(K)_\sigma \]
	into isotypical linear subspaces. If $\sigma$ factors through~$\Z[f]$, then
	\[ J_0(N)(K)_\sigma \subseteq A_f^\dual(K) \otimes_\Z \R = J_0(N)(K)[I_f] \otimes_\Z \R \,, \]
	and by the Heegner hypothesis (which implies that $A_f^\dual(K)$ has rank~$1$
	as a $\Z[f]$-module) it follows that \hbox{$\dim J_0(N)(K)_\sigma = 1$}.
	We will abuse notation slightly and also write $J_0(N)(K)_\sigma$
	when $\sigma \in \Sigma$, implicitly pre-composing with the projection
	$Z \to \Z[f]$. We write $y_{K,\sigma} \in J_0(N)(K)_\sigma$ for the components
	of~$y_K$ with respect to this composition and set
	\[ y_K^f \defeq \sum_{\sigma \in \Sigma} y_{K,\sigma} \in A_f^\dual(K) \otimes_\Z \R \,; \]
	Then $\lambda_f(y_K^f) = y_{K,A_f}$; compare the diagram~\eqref{Eqn:AJ diagram}.
	Note that $\omega \in \Q(f)$ acts on~$y_{K,\sigma}$ as
	$\omega \cdot y_{K,\sigma} = \omega^\sigma y_{K,\sigma}$.
	Explicitly, when $(b_j)_{1 \le j \le g}$ is a $\Z$-basis of~$\Z[f]$
	and $(b_j^*)_{1 \le j \le g}$ is its dual basis in $\Z[f] \otimes_{\Z} \R$
	with respect to the trace form, we have
	\begin{equation} \label{Eq:y_sigma vs y}
		y_{K,\sigma} = \sum_{j=1}^g (b_j \cdot y_K) \otimes {b_j^*}^\sigma .
	\end{equation}
	
	The normalized canonical height on~$J_0(N)(K)$ (with respect to twice the
	theta divisor) induces a positive definite quadratic
	form on~$J_0(N)(K) \otimes_\Z \R$, which (by abuse of notation) we also
	denote~$\hhat$. Since the endomorphisms are self-adjoint with respect
	to the height pairing (this is because they are fixed under the Rosati
	involution; see~\cite[Section~5.5]{BirkenhakeLange} and recall that
	the endomorphism ring is totally real), it follows that the $\sigma$-components
	are pairwise orthogonal under the height pairing.
	
	Recall that we write $L(f/K, s)$ for the $L$-function of~$f$ base-changed to~$K$.
	This is the same as $L(f, \mathbf{1}, s)$ for the trivial character
	$\mathbf{1} \colon \Gal(H|K) \to \C^\times$ in the notation of~\cite{GrossZagier1986}.
	
	\begin{theorem}[Gross--Zagier formula] \label[theorem]{height of Heegner point}
		With the notation introduced above and assuming that $D_K$ is odd, we have
		\[ \hhat(y_{K,\sigma}) = L'(f^\sigma/K, 1) \frac{u_K^2 \sqrt{-D_K}}{16 \bpi^2 \|f^\sigma\|^2} \,. \]
		Here $u_K \defeq \#\cO_K^\times/\Z^\times$,
		which equals $1$ for $D_K < -4$, $2$ for $D_K = -4$ and $3$ for $D_K = -3$.
	\end{theorem}
	\begin{proof}
		This is a reformulation of~\cite[Theorem~I.6.3]{GrossZagier1986},
		taking into account that our~$\hhat$ is $1/2h_K$~times the height used there
		(see~\cite[Eq.~(I.6.4)]{GrossZagier1986}), where $h_K$ is the class number of $K$.
		Note that Gross and Zagier
		assume that the Heegner discriminant is odd; see~\cite[\S\,I.3]{GrossZagier1986}.
	\end{proof}
	
	To evaluate this formula, we need the Petersson norm
	from~\ref{sec:computing-the-Petersson-norm}, and we need to
	evaluate~$L'(f^\sigma/K, 1)$. By the Artin formalism of $L$-functions,
	\[ L(f/K, s) = L(f,s) L(f \otimes \epsilon_K, s) \]
	with $f \otimes \epsilon_K$ the twist of $f$ by (the Kronecker character
	associated to)~$K$. Its first derivative at $s = 1$ is
	\[ L'(f/K, 1) = L(f,1) L'(f \otimes \epsilon_K, 1) + L'(f, 1) L(f \otimes \epsilon_K, 1) \,. \]
	Since $K$ is a Heegner field by assumption, $L(f/K, s)$ vanishes
	to first order at~$s = 1$. This implies that exactly one of the two terms in the sum
	is non-zero; which one it is can be decided by considering the
	action of the Fricke involution~$w_N$ on~$f$: if $w_N \cdot f = f$, then
	$L(f, 1) = 0$; otherwise, $w_N \cdot f = -f$, and $L'(f, 1) = 0$.
	The special values of the $L$-functions of newforms and their derivatives can
	be computed to arbitrary precision using Tim Dokchitser's Magma
	implementation~\cite{Dokchitser2004}.
	It provides a \texttt{TensorProduct} function for $L$-functions, which however
	tends to be slow in our use case. So we construct the tensor product $L$-function
	$L(f \otimes \epsilon_K, s)$ \enquote{by hand} for performance reasons, explicitly giving the Euler factors.
	
	We finally obtain a formula for~$\hhat(y_K^f)$.
	
	\begin{corollary} \label[corollary]{height via GZ}
		Let $K$ be a Heegner field for~$f$ and let $y_K^f \in A_f(K)^\dual$ be an associated
		Heegner point. Then
		\[ \hhat(y_K^f) = \frac{u_K^2 \sqrt{-D_K} \bpi}{2 N \prod_{\ell^2 \mid N} C_\ell}
		\sum_{\sigma \in \Sigma} \frac{L'(f^\sigma/K, 1)}{L(\Sym^2 f^\sigma, 2)} \,.
		\]
	\end{corollary}
	
	\begin{proof}
		Since the $y_{K,\sigma}$ are orthogonal with respect to the height pairing,
		we have $\hhat(y_K^f) = \sum_\sigma \hhat(y_{K,\sigma})$. Now combine
		\cref{height of Heegner point} and~\cref{Petersson formula}.
	\end{proof}

	\subsection{Comparing canonical heights}
	\label{sec:Comparing-canonical-heights}
	
	Our goal in this section is to determine~$\hhat_J(y_{K,\pi})$ (so that
	we can either use that to identify~$y_{K,\pi}$ up to a sign and adding
	torsion assuming there is an essentially unique point of that height,
	or to verify that our computation of~$y_{K,\pi}$ is correct).
	Recall the diagram~\eqref{Eqn:AJ diagram}, in particular the
	endomorphism~$\alpha$ of~$J$ defined in~\eqref{Eq:def alpha}.
	Note that $\lambda = \pi_J \circ \pi_J^\dual$ equals $\alpha$
	composed with the inverse of the canonical
	polarization~$\lambda_J$ of~$J$ induced by the theta divisor.
	
	We freely use standard facts about height pairings on abelian varieties;
	see for example~\cite[§\,9]{BombieriGubler}.
	We denote by $\langle {-}, {-} \rangle_J$ the height pairing on~$J$
	(such that $\hhat_J(x) = \langle x, x \rangle_J$).
	By~\cite[Prop.~9.3.6]{BombieriGubler} (noting that our $\hhat_J$
	is twice their~$\hhat_\theta$), it satisfies
	\[ \langle x, x' \rangle_J = \hhat_{\sP}(\lambda_J(x), x') \,, \]
	where $\hhat_{\sP}({-})$ is the canonical height on $J^\dual \times J$
	associated to the Poincaré bundle~$\sP$.
	Similarly, we obtain the canonical height associated to a polarization
	$\lambda \colon J^\dual \to J$ as $\hhat_\lambda(x) = \hhat_{\sP}(x, \lambda(x))$.
	If $\phi \colon A \to J$ is a homomorphism and $\lambda_{A^\dual} \colon A^\dual \to A$
	is a polarization such that
	\[ \lambda = {\phi^\dual}^* \lambda_{A^\dual}
	= \phi \circ \lambda_{A^\dual} \circ \phi^\dual
	= \alpha \circ \lambda_J^{-1} \]
	with $\alpha \in \End_\Q(J)$, then by functoriality of heights, we have for $x \in J(\Qbar)$
	\begin{align*}
		\hhat_{\lambda_{A^\dual}}(\phi^\dual(\lambda_J(x)))
		&= \hhat_{{\phi^\dual}^* \lambda_{A^\dual}}(\lambda_J(x))
		= \hhat_{\lambda}(\lambda_J(x))
		= \hhat_{\sP}(\lambda_J(x), \lambda(\lambda_J(x))) \\
		&= \hhat_{\sP}(\lambda_J(x), \alpha(x))
		= \langle x, \alpha(x) \rangle_J .
	\end{align*}
	
	\begin{proposition} \label[proposition]{height relation}
		For each $\sigma \in \Sigma$, write
		$y_{K,\pi,\sigma} = \pi_J(y_{K,\sigma}) \in J(K) \otimes_\Z \R$
		for the $\sigma$-component of~$y_{K,\pi}$. Then
		$\hhat_J(y_{K,\pi,\sigma}) = \alpha^\sigma \hhat(y_{K,\sigma})$, and so
		\[ \hhat_J(y_{K,\pi}) = \sum_{\sigma \in \Sigma} \alpha^\sigma \hhat(y_{K,\sigma}) \,. \]
	\end{proposition}
	
	\begin{proof}
		Chasing $y_{K,\sigma}$ through the diagram~\eqref{Eqn:AJ diagram} and
		taking into account the definition of~$\alpha \in \cO \subseteq \Q(f)$, we see that
		(identifying $A_f^\dual(K)$ with its image under the inclusion~$\iota_f$)
		\begin{align*}
			\pi^\dual\bigl(\lambda_J(y_{K,\pi,\sigma})\bigr)
			&= \pi^\dual\bigl((\lambda_J \circ \pi_J)(y_{K,\sigma})\bigr)
			= \pi_J^\dual\bigl((\lambda_J \circ \pi_J)(y_{K,\sigma})\bigr) \\
			&= \alpha \cdot y_{K,\sigma}
			= \alpha^\sigma y_{K,\sigma} .
		\end{align*}
		By the discussion preceding the proposition, we have
		\begin{align*}
			\hhat\bigl(\pi^\dual\bigl(\lambda_J(y_{K,\pi,\sigma})\bigr)\bigr)
			&= \langle y_{K,\pi,\sigma}, \alpha \cdot y_{K,\pi,\sigma} \rangle_J
			= \langle y_{K,\pi,\sigma}, \alpha^\sigma y_{K,\pi,\sigma} \rangle_J \\
			&= \alpha^\sigma \langle y_{K,\pi,\sigma}, y_{K,\pi,\sigma} \rangle_J
			= \alpha^\sigma \hhat_J(y_{K,\pi,\sigma}) .
		\end{align*}
		Therefore,
		\[ \hhat_J(y_{K,\pi,\sigma})
		= \frac{\hhat\bigl(\pi^\dual\bigl(\lambda_J(y_{K,\pi,\sigma})\bigr)\bigr)}{\alpha^\sigma}
		= \frac{\hhat(\alpha^\sigma y_{K,\sigma})}{\alpha^\sigma}
		= \alpha^\sigma \hhat(y_{K,\sigma}) \,. \qedhere \]
	\end{proof}
	
	When $X$ is a quotient of~$X_0(N)$, this gives a particularly simple
	formula.
	
	\begin{corollary} \label[corollary]{same heights}
		Assume that $\pi_X \colon X_0(N) \to X$ is a finite covering of curves of degree~$n$
		and that $\pi_J \colon J_0(N) \to J$ is induced by~$\pi_X$ via Albanese
		functoriality. Then
		\[ \hhat_J(y_{K,\pi}) = n \hhat(y_K^f) \,. \]
	\end{corollary}
	
	\begin{proof}
		In this case, $\alpha = \lambda_J \circ \pi_J \circ \pi_J^\dual$
		is multiplication by $\deg \pi_X = n$, so $\alpha^\sigma = n$
		for all $\sigma \in \Sigma$. Now use~\cref{height relation}.
	\end{proof}
	
	In the general case, we can determine~$\alpha$ as described in
	Section~\ref{sec:period matrices}. We record the final general formula
	for the height of~$y_{K,\pi}$.
	
	\begin{corollary} \label[corollary]{height on J via GZ}
		With the notation introduced so far, we have
		\[ \hhat_J(y_{K,\pi}) =
		\frac{u_K^2 \sqrt{-D_K} \bpi}{2 N \prod_{\ell^2 \mid N} C_\ell}
		\sum_{\sigma \in \Sigma} \alpha^\sigma \frac{L'(f^\sigma/K, 1)}{L(\Sym^2 f^\sigma, 2)} \,.
		\]
	\end{corollary}
	
	\begin{proof}
		Combine~\cref{height of Heegner point,height relation}.
	\end{proof}
	
	\begin{remark}
		In a similar way as in the proof of~\cref{height relation}, we obtain
		the formula
		\[ \langle \beta \cdot y_{K,\pi}, \gamma \cdot y_{K,\pi} \rangle_J
		= \sum_{\sigma \in \Sigma} \alpha^\sigma \beta^\sigma \gamma^\sigma \hhat(y_{K,\sigma})
		\]
		for arbitrary $\beta, \gamma \in \cO$.
		This allows us to compute the height pairing matrix~$M$ for a $\Z$-basis
		of~$\cO y_{K,\pi}$ and from this the regulator $\Reg_{\cO y_{K,\pi}} = \det M$. Then the Heegner
		index is given by
		\[ I_{K,\pi} = \#J(K)_\tors \sqrt{\frac{\Reg_{\cO y_{K,\pi}}}{\Reg_{J(K)}}} \,. \]
	\end{remark}
	
	
	\section{Computing the analytic order of $\Sha$}
	\label{sec:Sha_an}
	
	Recall that $J$ is an absolutely simple and principally polarized
	abelian variety over $\Q$ of dimension $g$ of $\GL_2$-type with associated
	newform $f \in S_2(\Gamma_0(N))$, and $A_f$ is the modular abelian variety
	associated to~$f$. In particular, $A_f$ and~$J$ are isogenous.
	
	For an abelian variety~$J$ over a number field~$F$, we define
	the \emph{Tamagawa product} to be
	\[ \Tam(J/F) \defeq \prod_v c_v(J/F) \,, \]
	where $v$ runs through the finite places of~$F$.
	When $J$ is the Jacobian variety of an explicitly given curve,
	the Tamagawa numbers~$c_v(C/F)$ (which are~$1$ for all places of good reduction)
	and hence the Tamagawa product~$\Tam(J/F)$ can be computed. For Jacobians
	of genus~$2$ curves in the LMFDB~\cite{lmfdb}, this information is also available
	in the~LMFDB. For the Tamagawa number at~$2$ in the example
	in Appendix~\ref{sec:7-torsion in Sha}, we compute a regular model by hand.
	
	We now describe how to compute the \emph{analytic order of the Tate--Shafarevich group}
	\begin{equation} \label{Eq:Sha_an}
		\#\Sha(J/\Q)_\an \defeq \frac{L^{(r)}(J/\Q,1)}{r! \, \Omega_J \Reg_{J/\Q}}
		\cdot \frac{(\#J(\Q)_\tors)^2}{\Tam(J/\Q)}
	\end{equation}
	as an exact positive rational number, assuming that $\Lrk J \in \{0,1\}$.
	
	Note that we can provably verify that $\Lrk J \in \{0, 1\}$ and determine
	$\Lrk J$ in this case.
	The Fricke involution~$w_N$ sends $f$ to $f$ or~$-f$. In the first
	case, the analytic order of~$L(f,s)$ is odd, and in the second case, it is even.
	In the even case, we can show that $L(f,1) \neq 0$, and in the odd case that $L'(f,1) \neq 0$
	by computing the respective value numerically to a high enough precision.

	\subsection{Comparing the real periods of $A_f$ and $J$} \label{ssec:comparing real periods}
	
	Let $A$ be an abelian variety over~$\Q$ of dimension~$g$ with
	Néron model~$\sA$ over~$\Z$. We say that a $\Q$-basis of~$\H^0(A, \Omega^1)$
	is a \emph{Néron basis} for~$A$ if it is a $\Z$ basis of the image of~$\H^0(\sA, \Omega^1_{\sA/\Z})$.
	Let $(\omega_1, \ldots, \omega_g)$ be a Néron basis for~$A$. Then
	$\omega_A \defeq \omega_1 \wedge \dots \wedge \omega_g$ is a generator of the free $\Z$-module
	of rank~$1$ $\H^0(\sA, \Omega^g_{\sA/\Z})$. Recall that the \emph{real period}
	of~$A$ is
	\[ \Omega_A \defeq \int_{A(\R)} |\omega_A| = \Bigl| \int_{A(\R)} \omega_A \Bigr| \,. \]
	
	Let $B$ be another abelian variety over~$\Q$ of dimension~$g$ with Néron
	model~$\sB$ over~$\Z$, and let $\pi \colon A \to B$ be an isogeny.
	Since by the Néron mapping property, $\pi$ uniquely extends to the Néron models,
	one has $\pi^*\omega_B = n_\pi \cdot \omega_A$ with an integer~$n_\pi$.
	By the above, $|n_\pi| = c_\pi$, where $c_\pi$ is defined in~\cref{c of an isogeny}.
	We now compare $\Omega_A$ and~$\Omega_B$.
	
	\begin{lemma} \label[lemma]{behavior of real periods under isogeny}
		Let $\pi\colon A \to B$ be an isogeny of abelian varieties of dimension~$g$ over~$\Q$.
		Denote by $\pi_\R$ the induced morphism $A(\R) \to B(\R)$ on the real Lie groups. Then
		\[
		\frac{\Omega_B}{\Omega_A} = \frac{\#\coker\pi_\R \cdot c_\pi}{\#\ker\pi_\R} \in \Q_{>0} \,.
		\]
		Here, $c_\pi$ divides~$e(\pi)^g$, where $e(\pi)$ is the exponent of~$\ker \pi$,
		$\#\ker\pi_\R$ divides $\deg\pi$, and $\#\coker\pi_\R$ divides the
		number $\#\pi_0(B(\R))$ of connected components of~$B(\R)$,
		which divides $2^g$.
	\end{lemma}
	
	\begin{proof}
		The isogeny $\pi$ induces a short exact sequence of real Lie groups
		\[
		0 \To (\ker \pi)(\R) \To A(\R) \stackrel{\pi_\R}{\To} \pi(A(\R)) \To 0 \,.
		\]
		This gives, for $\omega \in \H^0(B, \Omega^g)$,
		\[ \int_{A(\R)} \pi^* \omega = \#\ker \pi_\R \cdot \int_{\im(\pi_\R)} \omega
		= \frac{\#\ker \pi_\R}{\#\coker \pi_\R} \int_{B(\R)} \omega \,,
		\]
		where the second equality uses that $\omega$ is translation-invariant
		(compare~\cite[Lemma~5.13]{Jorza2005}). Hence,
		\begin{align*}
			\frac{\Omega_B}{\Omega_A} &= \frac{\int_{B(\R)}|\omega_B|}{\int_{A(\R)}|\omega_A|}
			= \frac{c_\pi \cdot \int_{B(\R)}|\omega_B|}{\int_{A(\R)}|\pi^*\omega_B|} \\
			&= \frac{c_\pi \cdot \#\coker \pi_\R \cdot \int_{B(\R)}|\omega_B|}{\#\ker \pi_\R \cdot\int_{B(\R)}|\omega_B|}
			= \frac{c_\pi \cdot \#\coker \pi_\R}{\#\ker \pi_\R}.
		\end{align*}
		
		One has $\#\ker\pi_\R \mid \deg\pi$ because $\ker \pi_\R \subseteq \ker \pi$.
		Let $\pi' \colon B \isoto A/\ker \pi \to A/A[e(\pi)] \isoto A$ be the isogeny such that $\pi' \circ \pi$ is
		multiplication by~$e(\pi)$.  Then
		\[ n_{\pi'} n_\pi \cdot \omega_A
		= n_{\pi'} \cdot \pi^* \omega_B
		= \pi^* (n_{\pi'} \cdot \omega_B)
		= \pi^* {\pi'}^* \omega_A
		= (\pi' \circ \pi)^* \omega_A
		= [e(\pi)]^* \omega_A
		= e(\pi)^g \cdot \omega_A \,,
		\]
		so $c_\pi = |n_\pi|$ divides $e(\pi)^g$.
		$\pi_\R$ is a topological
		covering map, so its image is open and closed, i.e., a union of connected components.
		This implies that $\#\coker \pi_\R$ divides $\#\pi_0(B(\R))$. Since the trace
		map $B(\C) \to B(\R)$ has image the connected component~$B(\R)^0$ of the origin,
		it follows that $\pi_0(B(\R))$ is killed by~$2$. This implies that $\pi_0(B(\R))$
		is isomorphic to $B(\R)[2]/B(\R)^0[2]$ of order
		dividing $4^g/2^g = 2^g$~\cite[proof of Lemma~3.10]{Schaefer1996}.
	\end{proof}
	
	Note that we can determine $\#\ker \pi_\R$ and $\#\coker \pi_\R$ explicitly
	if we have a suitable computational representation of the isogeny~$\pi$; see
	Section~\ref{sec:period matrices}.
	
	\begin{remark}
		See \cite[Lemma~5.13]{Jorza2005} for a similar statement over arbitrary completions
		of global fields.
	\end{remark}

	\subsection{Computing $L(J/\Q, 1)/\Omega_J$} \label{ssec: LJ/Om}
	
	We now consider the isogeny $\pi \colon A_f \to J$.
	The formula for~$\#\Sha(J/\Q)_\an$ in the case of $L$-rank~$0$
	contains the factor~$L(J/\Q, 1)/\Omega_J$.
	In this section we explain how this quotient can be computed as a rational number.
	We will then also use this later applied to a rank~$0$ quadratic twist of~$J$
	when dealing with the $L$-rank~$1$ case.
	By~\cref{behavior of real periods under isogeny}, we have
	\[ \frac{L(J/\Q, 1)}{\Omega_J}
	= \frac{L(A_f/\Q, 1)}{\Omega_{A_f}} \cdot \frac{\Omega_{A_f}}{\Omega_J}
	= \frac{L(A_f/\Q, 1)}{\Omega_{A_f}} \cdot \frac{1}{c_\pi} \cdot \frac{\#\ker \pi_\R}{\#\coker \pi_\R} \,,
	\]
	and $\Omega_{A_f} = c_f \cdot \Omega'_{A_f}$, where $\Omega'_{A_f}$ is the
	volume computed with respect to a $\Z$-basis of~$S_2(f, \Z)$ instead of a
	Néron basis. This gives
	\begin{equation} \label{Eq:LJ/Om}
		\frac{L(J/\Q, 1)}{\Omega_J}
		= \frac{L(A_f/\Q, 1)}{\Omega'_{A_f}} \cdot \frac{1}{c_f c_\pi} \cdot \frac{\#\ker \pi_\R}{\#\coker \pi_\R} .
	\end{equation}
	The quotient $LR(A_f) \defeq L(A_f/\Q, 1)/\Omega'_{A_f}$ is what Magma calls
	the \texttt{LRatio} of~$A_f$. This Magma function computes $LR(A_f) \in \Q_{\ge 0}$
	directly using modular symbols, but this computation is very slow and needs
	lots of memory when the level~$N$ of~$f$ is not very small (this seems to be caused
	by a computation of an integral homology basis of the ambient modular symbols
	space). The computation runs in reasonable time for $N \le 1000$, which is enough
	for the $L$-rank zero case, but becomes infeasible for example when $N = 67 \cdot 7^2$,
	which is the first relevant level of a suitable quadratic twist in the $L$-rank one case.
	
	So we use a numerical method instead. We compute numerical approximations
	to~$L(A_f/\Q, 1)$ and to~$\Omega'_{A_f}$ or~$\Omega_J$ and recognize the quotient
	as a rational number of small height. (See~\cite{Wuthrich2018} for a similar
	approach in the context of elliptic curves.) To do that reliably, we need a bound for
	the denominator of this quotient.
	
	\begin{proposition} \label[proposition]{denom bound A_f}
		\[ LR(A_f) = \frac{m}{\#\pi_0(A_f(\R)) \cdot \#A_f(\Q)_\tors} \qquad
		\text{for some $m \in \Z_{\ge 0}$.}
		\]
	\end{proposition}
	
	\begin{proof}
		By~\cite[Prop.~4.6]{AgasheStein2005}, the denominator of $\#\pi_0(A_f(\R)) \cdot LR(A_f)$
		divides the order~$n$ of the image in~$A_f$ of the difference of the cusps represented
		by $0$ and~$\infty$. This image is a rational torsion point, so $n \mid \#A_f(\Q)_\tors$.
	\end{proof}
	
	\begin{corollary} \label[corollary]{denom bound J}
		Let $g$ denote the dimension of $A_f$ and of~$J$. Then
		\[ \frac{L(J/\Q, 1)}{\Omega_J} = \frac{m}{4^g \cdot c_f c_\pi \cdot \#J(\Q)_\tors} \qquad
		\text{for some $m \in \Z_{\ge 0}$.}
		\]
	\end{corollary}
	
	\begin{proof}
		Note that $\#A_f(\Q)_\tors$ divides $\#\ker \pi_\R \cdot \#J(\Q)_\tors$ and
		that $\#\pi_0(A_f(\R))$ and $\#\coker \pi_\R$ both divide~$2^g$.
		The claim then follows from~\eqref{Eq:LJ/Om} and~\cref{denom bound A_f}.
	\end{proof}
	
	Since we can determine $\#J(\Q)_\tors$ (an upper bound obtained from the $L$-series
	coefficients as in~\cite[\S\,3.5]{AgasheStein2005} would be enough) and we can
	compute~$c_f c_\pi$ by~\cref{computation of c_f c_pi}, it suffices to compute
	$L(J/\Q, 1)$ and~$\Omega_J$ to sufficient precision so that the resulting
	approximation to $4^g \cdot c_f c_\pi \cdot \#J(\Q)_\tors \cdot L(J/\Q, 1)/\Omega_J$
	has error $< 1/2$. We then round to the nearest integer to obtain the numerator~$m$
	in~\cref{denom bound J}. In practice, we use higher precision and check that
	the error is as small as can be expected.
	
	\subsection{The case of $L$-rank $0$}  \label{ssec:Shaan rk 0}
	
	We obtain the following formula.
	
	\begin{proposition} \label[proposition]{formula for Sha_an rk 0}
		Assume that $L(f, 1) \neq 0$. Then
		\[ \#\Sha(J/\Q)_\an =\frac{L(J/\Q, 1)}{\Omega_J}
		\cdot \frac{(\#J(\Q)_\tors)^2}{\Tam(J/\Q)} \in \Q_{>0} \,.
		\]
	\end{proposition}
	
	\begin{proof}
		This is \eqref{Eq:Sha_an} for $r = 0$.
	\end{proof}
	
	Note that all quantities in the formula in~\cref{formula for Sha_an rk 0} can be computed
	explicitly: for the first factor see Section~\ref{ssec: LJ/Om},
	for the torsion subgroup see~\cite[\S\,11]{Stoll1999},
	and for the Tamagawa product see the beginning of this section.

	\subsection{The case of $L$-rank $1$: Computing $\#\Sha(J/K)_\an$}
	\label{SS: Sha_an/K}
	
	In the following we keep assuming that $J$ is a Jacobian.
	In particular, $J$ is principally polarized.
	
	When the $L$-rank is~$1$, we first find a Heegner field~$K$ and compute the analytic
	order of~$\Sha$ for~$J/K$ exactly from the BSD formula
	\[ L^*(J/K,1) = \#\Sha(J/K) \cdot \frac{\Omega_{J/K} \Reg'_{J/K}}{\sqrt{|D_K|}^g}
	\cdot \frac{\Tam(J/K)}{(\#J(K)_\tors)^2} \,. \]
	Here the period $\Omega_{J/K}$ is defined as
	\begin{equation} \label{eq:def Omega_JK}
		\Omega_{J/K} = \int_{J(\C)} |\omega \wedge \ol{\omega}| ,
	\end{equation}
	where $\omega$ is a generator of the free rank~$1$ $\Z$-module of top Néron
	differentials on~$J$ (this works since $J/K$ is base-changed from an abelian
	variety over~$\Q$). Note that this is $2^g$~times the covolume of the period
	lattice (which is generated by the columns of the big period matrix~$\Pi_J$,
	if it is computed with respect to a Néron basis of the invariant $1$-forms).
	
	The regulator $\Reg'_{J/K}$ is computed with respect to heights over~$K$.
	We will write~$\Reg_{J/K}$ to denote the regulator with respect to the normalized
	height; we then have that $\Reg'_{J/K} = [K : \Q]^{\rk J(K)} \cdot \Reg_{J/K}$.
	See~\cite{Tate1968}. (In the literature, the formula is often stated without
	making precise what \enquote{the regulator} and \enquote{the period} are, which can lead to
	confusion. See the answers to the Math Overflow question at~\cite{MO139575}
	and~\cite{CremonaBSD} for a discussion.) We deduce that
	\begin{align} \label{eq:ShaAKan}
		\#\Sha(J/K)_\an
		&= \frac{(\#J(K)_\tors)^2}{\Tam(J/K)}
		\cdot \frac{L^*(J/K,1) \sqrt{|D_K|}^g}{\Omega_{J/K} [K : \Q]^{\rk J(K)}\Reg_{J/K}} \\
		&= \frac{(\#J(K)_\tors)^2}{\Tam(J/K) \cdot u_K^{2g}}
		\cdot \frac{\prod_\sigma 8 \bpi^2 \|f^\sigma\|^2}{\Omega_{J/K}}
		\cdot \frac{\prod_\sigma \hhat(y_{K,\sigma})}{\Reg_{J/K}},
		\label{eq:ShaAKan2}
	\end{align}
	where $\sigma$ runs through the $g$ embeddings $\sigma\colon \Q(f) \inj \R$.
	The second equality
	follows from the Gross--Zagier formula~\cref{height of Heegner point} with
	\[ L^*(J/K, 1) = \frac{L^{(g)}(J/K,1)}{g!} = \prod_\sigma L'(f/K^\sigma,1) \,, \]
	where $L(f/K^\sigma, s) = L(f^\sigma, s) L(f^\sigma \otimes \epsilon_K, s)$
	with the quadratic character~$\epsilon_K$ associated to~$K|\Q$.
	
	Note that all primes~$p$ of bad reduction for~$J/\Q$ split as~$\frp\bar{\frp}$
	in~$K$ by the Heegner condition. This implies that $K_{\frp} \iso \Q_p \iso K_{\bar\frp}$
	and so in particular that $c_\frp(J/K) = c_p(J/\Q) = c_{\bar{\frp}}(J/K)$. Therefore
	the Tamagawa product over~$K$ is the square of the Tamagawa product over~$\Q$,
	\begin{equation} \label{Eq:Tam JK}
		\Tam(J/K) = \Tam(J/\Q)^2 .
	\end{equation}
	We now describe in a series of lemmas how to determine the last two factors
	in~\eqref{eq:ShaAKan2}. Combining the results
	gives the explicit formula in~\cref{formula for Sha J/K an}.
	
	\begin{lemma} \label[lemma]{Petersson to integral}
		Let $f \in \cN(N,g)$. Then the Petersson norm~$\|f\|^2$ satisfies
		\[ 8\bpi^2\|f\|^2 = \|\omega_f\|^2 \defeq \int_{X_0(N)(\C)}\omega_f \wedge \ol{i \omega_f} \]
		with $\omega_f = 2\bpi i f(z) \,dz$.
	\end{lemma}
	\begin{proof}
		See~\cite[\S\,1.6]{GrossZagier1986}.
	\end{proof}
	
	We now want to relate the product of the Petersson norms to the complex
	period of~$A_f$. Extending diagram~\eqref{Eqn:AJ diagram} to the left, we obtain
	the following diagram, where $B_f = I_f J_0(N)$ is the kernel of~$\pi_f$.
	\begin{equation} \label{Eqn:AJ diagram 2}
		\SelectTips{cm}{}
		\xymatrix@!C{B_f \ar@{^(->}[r]^-{\iota_f} \ar[d]_{\lambda'_f}
			& J_0(N) \ar@{=}[d] \ar@{->>}[r]^-{\pi_f}
			& A_f \\
			B_f^\dual
			& J_0(N)^\dual \ar@{->>}[l]_-{\iota_f^\dual}
			& A_f^\dual \ar@{_(->}[l]_-{\pi_f^\dual} \ar[u]_{\lambda_f} }
	\end{equation}
	Also recall the definition of~$d_f$ from~\eqref{Eq:def d_f}.
	
	\begin{lemma} \label[lemma]{intersection of W with B}
		Let $W_g$ be the image of~$\Sym^g X_0(N)$ in~$J_0(N)$ (with respect to some
		base divisor of degree~$g$). Let $B_f = I_f J_0(N) = \ker \pi_f$. Then the
		intersection number $W_g \cdot B_f$ equals~$d_f$.
	\end{lemma}
	
	We thank Jakob Stix and Yusuf Mustopa for help with the proof.
	
	\begin{proof}
		Let $m \defeq \dim J_0(N) = g(X_0(N))$.
		By~\cite[Thm.~V.1.3]{ACGH} (with $(r, d, n) \leftarrow (0, g, m)$)
		the class of~$W_g$ is $1/(m-g)! \cdot \theta^{m-g}$, where $\theta$ is the
		class of the theta divisor on~$J_0(N)$. We write $d_f'$ for the product
		$d'_1 \cdots d'_g$, where $(d'_1, \ldots, d'_g)$ is the type of~$\lambda_f'$
		in diagram~\eqref{Eqn:AJ diagram 2} above.
		Then
		\begin{align*}
			{d_f'}^2 &= \deg \lambda_f'
			= \#(B_f \cap \ker \iota_f^\dual)
			= \#(\ker \pi_f \cap A_f^\dual)
			= \deg \lambda_f
			= d_f^2 ,
		\end{align*}
		so $d_f' = d_f$. This implies that the intersection number is
		(compare \cite[Thm.~3.6.3]{BirkenhakeLange})
		\begin{align*}
			W_g \cdot B_f &= \frac{\theta|_{B_f}^{m-g}}{(m-g)!}
			= \frac{(m-g)! \cdot d_f'}{(m-g)!}
			= d_f'
			= d_f .
			\qedhere
		\end{align*}
	\end{proof}
	
	We denote the Abel--Jacobi morphism $X_0(N) \inj J_0(N)$ with respect to the
	cusp~$\infty$ by~$\iota$. Since
	$\iota^* \colon \H^0(J_0(N), \Omega^1) \isoto \H^0(X_0(N), \Omega^1)$
	is an isomorphism, we can identify the differentials~$\omega_{f^\sigma}$
	with holomorphic (hence invariant) \hbox{$1$-forms} on~$J_0(N)$, which we also
	denote by~$\omega_{f^\sigma}$. The map $\pi_f \colon J_0(N) \to A_f$
	induces an injective homomorphism
	$\pi_f^* \colon \H^0(A_f(\C), \Omega^1) \to \H^0(J_0(N)(\C), \Omega^1)$
	whose image is the subspace spanned by the~$\omega_{f^\sigma}$.
	We write $\omega_{A_f,\sigma}$ for the uniquely determined preimage
	of~$\omega_{f^\sigma}$ under this map.
	
	\begin{lemma} \label[lemma]{prod omega to int A_f}
		With the notation introduced so far,
		\[
		\prod_\sigma \|\omega_{f^\sigma}\|^2
		= d_f \cdot \int_{A_f(\C)}\bigwedge_\sigma(\omega_{A_f,\sigma} \wedge \ol{i \omega_{A_f,\sigma}}) \,.
		\]
	\end{lemma}
	
	\begin{proof}
		To simplify notation, fix a numbering $\sigma_1, \ldots, \sigma_g$
		of the embeddings $\sigma \colon \Q(f) \inj \R$ and write $\omega_j$ for~$\omega_{f^{\sigma_j}}$.
		
		We first show that
		\[
		\prod_\sigma \|\omega_{f^\sigma}\|^2
		= \int_{W_g(\C)} \omega_1 \wedge \ol{i \omega_1} \wedge \dots \wedge \omega_g \wedge \ol{i \omega_g} \,,
		\]
		where $W_g$ is as in~\cref{intersection of W with B} with base divisor $g \cdot \infty$.
		Consider the composition $X_0(N)^g \stackrel{\iota^g}{\to} J_0(N)^g \stackrel{s}{\to} J_0(N)$
		with the first morphism $\iota \times \ldots \times \iota$ and $s$ the summation morphism.
		This morphism has degree~$g!$ above its image~$W_g$ since it factors through the $g$-fold
		symmetric power of~$X_0(N)$, which is birational to~$W_g$ via~$s$. This gives
		\begin{align*}
			\int_{W_g(\C)} \omega_1 \wedge \ol{i \omega_1} \wedge \dots &\wedge \omega_g \wedge \ol{i \omega_g} \\
			&= \frac{1}{g!} \int_{X_0(N)(\C)^g}
			(\iota^g)^* s^* (\omega_1 \wedge \ol{i \omega_1} \wedge \dots \wedge \omega_g \wedge \ol{i \omega_g}).
		\end{align*}
		
		Now for any invariant $1$-form~$\omega$ on an abelian variety, we have
		that $s^* \omega = \sum_{k=1}^g \pr_k^* \omega$ with $\pr_k$ the $k$th projection
		$J_0(N)^g \to J_0(N)$; see~\cite[§\,1.5\,(9)]{BirkenhakeLange}. This implies
		\begin{align*}
			(\iota^g)^*s^*\Bigg(\bigwedge_{j=1}^g \omega_j \wedge \ol{i \omega_j}\Bigg)
			&= (\iota^g)^*\Bigg(\bigwedge_{j=1}^g \Big(\sum_{k=1}^g \pr_k^*\omega_j\Big)
			\wedge \Big(\sum_{k=1}^g \pr_k^* \ol{i \omega_j}\Big)\Bigg) .
		\end{align*}
		We expand the right hand side. Terms containing two factors $\pr_k^* \omega_j$
		or two factors $\pr_k^* \ol{i \omega_j}$ with the same~$k$ vanish since the wedge product of
		two holomorphic differentials on a curve vanishes.
		So we are left with a sum of terms of the form
		\begin{align*}
			\pm \frac{1}{g!} \int_{X_0(N)(\C)^g} \pr_1^*(\omega_{j_1} \wedge \ol{i \omega_{j'_1}}) &\wedge \dots
			\wedge \pr_g^*(\omega_{j_g} \wedge \ol{i \omega_{j'_g}}) \\
			&= \pm \frac{1}{g!} \prod_{k=1}^g \int_{X_0(N)(\C)} \omega_{j_k} \wedge \ol{i \omega_{j'_k}}
		\end{align*}
		with $\{j_1, \ldots, j_g\} = \{j'_1, \ldots, j'_g\} = \{1, \ldots, g\}$.
		Now when $j_k \neq j'_k$ for some~$k$, then the corresponding integral vanishes since
		the $f^\sigma$ are pairwise orthogonal with respect to the Petersson inner product.
		All the remaining terms have a positive sign (all relevant permutations are even)
		and differ only in the ordering of the factors; in particular, there are exactly~$g!$
		such terms. So we obtain
		\begin{align*}
			\int_{W_g(\C)} \bigwedge_\sigma (\omega_{f^\sigma} \wedge \ol{i \omega_{f^\sigma}})
			&= \int_{W_g(\C)} \omega_1 \wedge \ol{i \omega_1} \wedge \dots \wedge \omega_g \wedge \ol{i \omega_g} \\
			&= \prod_{j=1}^g \int_{X_0(N)(\C)} \omega_{j} \wedge \ol{i \omega_{j}}
			= \prod_\sigma \|\omega_{f^\sigma}\|^2
		\end{align*}
		as desired.
		
		We now consider $\pi_f|_{W_g} \colon W_g \to A_f$. Since $\dim W_g = \dim A_f = g$
		and by~\cref{intersection of W with B}, $W_g$ meets generic cosets of~$B_f = \ker \pi_f$
		transversally in $W_g \cdot B_f = d_f$ points, we finally obtain that
		\begin{align*}
			\prod_\sigma \|\omega_{f^\sigma}\|^2
			&= \int_{W_g(\C)} \omega_1 \wedge \ol{i \omega_1} \wedge \dots \wedge \omega_g \wedge \ol{i \omega_g} \\
			&= \int_{W_g(\C)} \bigwedge_\sigma (\pi_f^* \omega_{A_f,\sigma} \wedge \ol{i \pi_f^* \omega_{A_f,\sigma}}) \\
			&= d_f \cdot \int_{A_f(\C)} \bigwedge_\sigma (\omega_{A_f,\sigma} \wedge \ol{i \omega_{A_f,\sigma}}) .
			\qedhere
		\end{align*}
	\end{proof}
	
	We now relate the integral on the right hand side in~\cref{prod omega to int A_f}
	to the period~$\Omega_{A_f/K}$.
	
	\begin{lemma} \label[lemma]{relation period over omegafsigma and Neron differential}
		One has
		\begin{align*}
			\int_{A_f(\C)} \bigwedge_\sigma (\omega_{A_f,\sigma} \wedge \ol{i \omega_{A_f,\sigma}})
			&= \frac{\disc \Z[f]}{c_f^2} \cdot \int_{A_f(\C)}|\omega_{A_f} \wedge \ol{\omega_{A_f}}| \\
			&= \frac{\disc \Z[f]}{c_f^2} \cdot \Omega_{A_f/K} ,
		\end{align*}
		where $\omega_{A_f}$ is a top Néron differential on $A_f$, i.e., a generator
		of~$\H^0(\sA_f,\Omega^g)$ with $\sA_f/\Z$ the Néron model of~$A_f$,
		$\Omega_{A_f/K} \defeq \int_{A_f(\C)}|\omega_{A_f} \wedge \ol{\omega_{A_f}}|$,
		and $c_f$ is the Manin constant from~\cref{Manin constant}.
		
		Over $\R$, one has
		\[
		\int_{A_f(\R)} \bigwedge_\sigma \omega_{A_f,\sigma}
		= \pm \frac{\sqrt{\disc \Z[f]}}{c_f} \cdot \Omega_{A_f/\Q} \,.
		\]
	\end{lemma}
	
	\begin{proof}
		Let $(f_j)_{j=1}^g$ be a $\Z$-basis of~$S_2(f, \Z)$.
		Then $f = \sum_j b_j f_j$ for some $b_j \in \Z[f]$, which form a
		$\Z$-basis of~$\Z[f] = \Z[a_n(f) : n \ge 1]$. The matrix $A = (b_j^\sigma)_{\sigma,j}$
		then is such that $A \cdot (f_j)_j^\top = (f^\sigma)_\sigma^\top$, and it satisfies
		\[ \det(A)^2 = \det(b_i^\sigma)_{i,\sigma}^2 = \disc \Z[f] \,. \]
		By the definition of~$c_f$, we have that
		\[ |\pi_f^* \omega_{A_f}| = c_f \cdot |\omega_{f_1} \wedge \dots \wedge \omega_{f_g}| \,, \]
		so
		\[ |\omega_{A_f} \wedge \ol{\omega_{A,f}}|
		= c_f^2 \cdot \bigl| \bigwedge_j (\omega_{A,j} \wedge \ol{i \omega_{A,j}}) \bigr| \,,
		\]
		where $\pi_f^* \omega_{A,j} = \omega_{f_j}$. We also have that
		\[ \bigwedge_\sigma (\omega_{A_f,\sigma} \wedge \ol{i \omega_{A_f,\sigma}})
		= (\det A)^2 \bigwedge_j (\omega_{A,j} \wedge \ol{i \omega_{A,j}}) \,.
		\]
		Combining these gives the result.
		
		The formula for $A_f(\R)$ follows in the same way, the only difference being
		that we do not take wedge products with conjugate differentials, hence
		we get the square root of the factor.
	\end{proof}
	
	We need to compare the periods $\Omega_{J/K}$ and~$\Omega_{A_f/K}$.
	
	\begin{lemma} \label[lemma]{J A period comparison}
		One has
		\[ \frac{\Omega_{J/K}}{\Omega_{A_f/K}} = \frac{c_\pi^2}{\deg \pi} \,, \]
		where~$c_\pi$ is as in~\cref{c of an isogeny}.
	\end{lemma}
	
	\begin{proof}
		Note that the top Néron differentials $\omega_{A_f}$ and~$\omega_J$ on $A_f$
		and~$J$ are related by $\pi^*\omega_J = \pm c_\pi \cdot \omega_{A_f}$. Hence
		\begin{align*}
			\Omega_{J/K}
			&= \int_{J(\C)}|\omega_J \wedge \ol{\omega_J}|
			= \frac{1}{\deg \pi}\int_{A_f(\C)}|\pi^*(\omega_J  \wedge \ol{\omega_J})| \\
			&= \frac{1}{\deg \pi}\int_{A_f(\C)}|(c_\pi \cdot \omega_{A_f}) \wedge \ol{(c_\pi \cdot \omega_{A_f})}|
			= \frac{c_\pi^2}{\deg \pi} \Omega_{A_f/K}. \qedhere
		\end{align*}
	\end{proof}
	
	Combining \cref{Petersson to integral,prod omega to int A_f,%
		relation period over omegafsigma and Neron differential,J A period comparison}
	yields the following explicit expression for the second factor in~\eqref{eq:ShaAKan2}.
	
	\begin{corollary} \label[corollary]{third factor}
		One has
		\[
		\frac{\prod_\sigma 8\bpi^2\|f^\sigma\|^2}{\Omega_{J/K}}
		= \frac{\deg \pi \cdot d_f \cdot \disc \Z[f]}{(c_f c_\pi)^2} \in \Q_{>0} \,.
		\]
	\end{corollary}
	
	We now consider the third (and last) factor in~\eqref{eq:ShaAKan2}.
	
	\begin{lemma} \label[lemma]{last factor}
		One has
		\[ \frac{\prod_\sigma\hhat(y_{K,\sigma})}{\Reg_{J/K}}
		= \frac{I_{K,\pi}^{2}}{(\#J(K)_\tors)^2 \cdot \Nm(\alpha) \cdot \disc\End_\Q(J)}
		\in \Q_{>0} \,,
		\]
		where $\alpha \in \End_\Q J$ is defined in~\eqref{Eq:def alpha} and $I_{K,\pi}$
		is the Heegner index of~$J$ with respect to the chosen isogeny
		$\pi \colon A_f \to J$; see~\eqref{Eq:Def I_Kpi}.
	\end{lemma}

	\begin{proof}
		Since $\alpha^\sigma\hhat(y_{K,\sigma}) = \hhat_J((\pi \circ \lambda_f)(y_{K,\sigma}))$
		(see~\cref{height relation}),
		\[
		\Nm(\alpha)\prod_\sigma\hhat(y_{K,\sigma})
		= \prod_\sigma \alpha^\sigma\hhat(y_{K,\sigma})
		= \prod_\sigma\hhat_J((\pi \circ \lambda_f)(y_{K,\sigma})) \,.
		\]
		Here $\sigma$ runs through the embeddings $\End_{\Q}(J) \inj \R$.
		Now
		\[
		\Reg_J(\End_{\Q}(J)\cdot y_{K,\pi}) = \det(\langle b_i\cdot y_{K,\pi}, b_j\cdot y_{K,\pi}\rangle_J)
		\]
		with $(b_i)_{i=1}^g$ a $\Z$-basis of $\End_{\Q}(J)$.
		
		But $(\pi \circ \lambda_f)(y_{K,\sigma}) = \sum_{j=1}^g b_j \cdot y_{K,\pi} \otimes b_j^{*,\sigma}$,
		where $(b_j^*)_{j = 1}^g$ is the dual basis of~$\End_{\Q}(J) \otimes_\Z \R$
		with respect to the trace pairing $(a,b) \mapsto \Tr_{\End_{\Q}(J)/\Z}(ab)$ of~$\End_{\Q}(J)$;
		see~\eqref{Eq:y_sigma vs y}.
		Using this and the fact that the $(\pi \circ \lambda_f)(y_{K,\sigma})$ are
		orthogonal in pairs with respect to the height pairing, we find that
		\begin{align*}
			\prod_\sigma \hhat_J((\pi \circ \lambda_f)(y_{K,\sigma}))
			&= \det(\langle (\pi \circ \lambda_f)(y_{K,\sigma_1}),
			(\pi \circ \lambda_f)(y_{K,\sigma_2})\rangle_{J}) \\[-9pt]
			&= \Reg_J(\End_{\Q}(J)\cdot y_{K,\pi}) \cdot \det(b_j^{*,\sigma})^2 \\
			&= \Reg_J(\End_{\Q}(J)\cdot y_{K,\pi}) \cdot \det(b_j^{\sigma})^{-2}\\
			&= \Reg_J(\End_{\Q}(J)\cdot y_{K,\pi}) \cdot (\disc\End_{\Q}(J))^{-1}.
		\end{align*}
		Using that $\Reg_J(\End_{\Q}(J)\cdot y_{K,\pi}) = I_{K,\pi}^2 \Reg_{J/K}/(\#J(K)_\tors)^2$,
		we finally obtain
		\begin{align*}
			\frac{\prod_\sigma\hhat(y_{K,\sigma})}{\Reg_{J/K}}
			&= \frac{\prod_\sigma \hat{h}_J((\pi \circ \lambda_f)(y_{K,\sigma}))}{\Nm(\alpha) \cdot \Reg_{J/K}}
			= \frac{\Reg_J(\End_{\Q}(J)\cdot y_{K,\pi})}{\Nm(\alpha) \cdot \Reg_{J/K} \cdot \disc\End_{\Q}(J)} \\
			&= \frac{I_{K,\pi}^2}{(\#J(K)_\tors)^2 \cdot \Nm(\alpha) \cdot \disc\End_{\Q}(J)} .
			\qedhere
		\end{align*}
	\end{proof}
	
	We can now compute $\#\Sha(J/K)_\an \in \Q_{>0}$ exactly, as follows.
	
	\begin{corollary} \label[corollary]{formula for Sha J/K an}
		\[ \#\Sha(J/K)_\an
		= \frac{1}{(c_f c_\pi)^2} \cdot \frac{\disc \Z[f]}{\disc \End_{\Q} (J)}
		\cdot \Bigl(\frac{I_{K,\pi}}{\Tam(J/\Q) \cdot u_K^g}\Bigr)^2 \,.
		\]
	\end{corollary}
	
	\begin{proof}
		This follows from using \cref{third factor,last factor} in~\eqref{eq:ShaAKan2},
		noting that $\Nm(\alpha) = d_f \cdot \deg \pi$
		by~\eqref{Eq:Nm alpha}.
	\end{proof}
	
	Since $\Z[f]$ and $\End_\Q(J)$ both are sub-orders of the ring of integers of~$\Q(f)$,
	the quotient of their discriminants is a square. So $\#\Sha(J/K)_\an$ is a square;
	this is consistent with the fact that $J$ over~$K$ is even in the sense
	of~\cite{PoonenStoll1999} since the only bad places are primes that split in~$K$,
	and the curve $J$ is the Jacobian of is simultaneously deficient or not at both
	places above a bad prime~$p$ of~$J$ over~$\Q$.
	
	\begin{remark}
		Assuming $\Z[f] = \End_{\Q}(J)$ and $u_K = 1$, all invariants on the
		right hand side of~\cref{formula for Sha J/K an}
		are orders of finite $\cO$-modules in a natural way.
		It is natural to ask for a refined BSD formula over $\cO$, namely
		whether the element in the Grothendieck group
		of finite $\cO$-modules corresponding to the right hand side
		of~\cref{formula for Sha J/K an} equals that defined by~$\Sha(J/K)$.
	\end{remark}
	
	\subsection{Periods of quadratic twists} \label{ssect: periods}
	
	In order to compute $\#\Sha(J/\Q)_\an$ in the $L$-rank~$1$ case,
	we also need to compute~$\#\Sha(J^K/\Q)_\an$. We can do this as described
	in Section~\ref{ssec:Shaan rk 0}. This requires the computation of the
	quotient $L(J^K/\Q, 1)/\Omega_{J^K}$. We do that as explained
	in Section~\ref{ssec: LJ/Om}. The computation of~$\Omega_{J^K}$ requires
	the period matrix of~$J^K$. We explain in this section how we can obtain
	this period matrix easily from that of~$J$. See~\cref{twisted period matrix}
	below for a slightly more general version. We also need to determine
	$c_{f^K} c_{\pi^K}$ to obtain the bound for the denominator. We show
	in~\cref{c_f c_pi under twist} that it is the same as $c_f c_\pi$.
	
	The following result is a generalization of \cite[Lemma~3.1 and Cor.~2.6]{Pal2012}
	from elliptic curves to more general abelian varieties. We first state a local
	version. We will frequently use the embedding
	\[ \H^0(\sA^K, \Omega^1) \inj \H^0(A^K, \Omega^1) \inj \H^0(A_K, \Omega^1) \]
	induced by the isomorphism $A^K_K \iso A_K$, where $\sA^K$ is a Néron model
	of the quadratic twist~$A^K$.
	
	\begin{lemma} \label[lemma]{twisted differentials}
		Let $A$ be an abelian variety over~$\Q_p$, where $p$ is an odd prime;
		assume that $A$ has good reduction. Let $K|\Q_p$ be a ramified quadratic
		extension, given as $K = \Q_p(\varpi)$ with $\varpi^2 = u p$, where $u \in \Z_p^\times$.
		Let $\sA$ and~$\sA^K$ denote the Néron models of~$A$ and of its quadratic
		twist~$A^K$ over~$\Z_p$.
		
		Then the image of~$\H^0(\sA^K, \Omega^1)$ in~$\H^0(A_K, \Omega^1)$ is
		$1/\varpi$ times the image of~$\H^0(\sA, \Omega^1)$.
	\end{lemma}
	
	We thank Kęstutis Česnavicius for help with the proof.
	
	\begin{proof}
		We consider
		the images of~$\H^0(\sA, \Omega^1)$ and
		of~$\H^0(\sA^K, \Omega^1)$ in~$V \defeq \H^0(A_K, \Omega^1)$, respectively;
		they are free $\Z_p$-submodules of~$V$ of rank $g = \dim_K V = \dim A$.
		The first image
		is the $\Z_p$-dual of $L \defeq \Lie(\sA) \inj \Lie(A_K)$ (i.e., it consists
		of the differentials~$\omega$ such that $\langle \lambda, \omega \rangle \in \Z_p$
		for all $\lambda \in L$ under the natural pairing between the Lie algebra
		(the tangent space at the origin) and $V$ (its dual)), and similarly
		for~$L^K \defeq \Lie(\sA^K) \inj \Lie(A_K)$ and the second image.
		To see this, note first that $\H^0(\sA, \Omega^1_{\sA})$ is identified with
		$\H^0(\Spec \Z_p, \epsilon^* \Omega^1_{\sA})$, where $\epsilon$ is the zero
		section of~$\sA$; see~\cite[\S\,4.2,~Prop.~1]{BLR}. Then by~\cite[\S\,2.2,~Prop.~7(b)]{BLR},
		$\epsilon^* \Omega^1_{\sA} = \cI/\cI^2$, where $\cI$ is the ideal sheaf of the
		zero section of~$\sA$. The functor of points definition of the Lie algebra
		then gives that $\Lie(\sA)$ is the $\Z_p$-dual of~$\cI/\cI^2$. These identifications
		are all compatible with base change to~$\Q_p$, so the duality is compatible with what
		is happening on the generic fiber.
		
		It therefore suffices to show that $L^K = \varpi \cdot L$. By our assumptions,
		$K|\Q_p$ is tamely ramified. By~\cite[Thm.~4.2]{Edixhoven1992}, the natural
		map induces an isomorphism of~$\sA^K$ with the subscheme of the restriction
		of scalars $R_{\Z_p[\varpi]/\Z_p} \sA_{\Z_p[\varpi]}$ fixed by the twisted
		action of~$\Gal(K|\Q_p)$. Taking the invariants under this action commutes
		with forming the Lie algebra, so we obtain that $L^K$ is obtained by taking
		the invariants under this twisted action on
		$\Lie(\sA_{\Z_p[\varpi]}) = L \otimes_{\Z_p} \Z_p[\varpi]$; this invariant
		space is exactly $\varpi \cdot L$.
	\end{proof}
	
	\begin{corollary} \label[corollary]{twisted Neron differentials}
		Let $A$ be an abelian variety over~$\Q$ and let $K$ be a quadratic number field
		of odd discriminant~$D_K$ such that all primes of bad reduction for~$A$ are unramified
		in~$K$. Let $\sA$ and~$\sA^K$ denote the Néron models of~$A$ and of its quadratic
		twist~$A^K$ over~$\Z$. Then the image of~$\H^0(\sA^K, \Omega^1)$ in~$\H^0(A_K, \Omega^1)$
		is $1/\sqrt{D_K}$ times the image of~$\H^0(\sA, \Omega^1)$.
	\end{corollary}
	
	\begin{proof}
		Fix a Néron basis $(\omega_1, \ldots, \omega_g)$ for~$A$, where
		$g = \dim A$. Identifying invariant $1$-forms on $A$ and on~$A^K$ with their
		images on~$A_K$, we see that $(\sqrt{D_K}^{-1} \omega_1, \ldots, \sqrt{D_K}^{-1} \omega_g)$
		is a $\Q$-basis of the space of invariant $1$-forms on~$A^K$. Since $K|\Q$
		is unramified at all places of bad reduction of~$A$, these $1$-forms will
		form a local Néron basis at all these places, and also at all places of
		good reduction for~$A$ at which $K|\Q$ is unramified. Finally, \cref{twisted differentials}
		(with $\varpi \leftarrow \sqrt{D_K}$) shows that they also form a local Néron basis
		at all places where $K|\Q$ is ramified. So we have obtained a (global) Néron
		basis for~$A^K$, and the claim follows.
	\end{proof}
	
	\begin{corollary} \label[corollary]{twisted period matrix}
		Let $A$ be an abelian variety over~$\Q$ and let $K$ be a quadratic number field
		of odd discriminant
		such that all primes of bad reduction for~$A$ are unramified in~$K$.
		Let $\Pi_A$ be a big period matrix for~$A$ with respect to a Néron basis
		of~$\H^0(A, \Omega^1)$. Then $\sqrt{D_K}^{-1} \Pi_A$ is a big period
		matrix for the quadratic twist~$A^K$ with respect to a Néron basis of~$\H^0(A^K, \Omega^1)$.
	\end{corollary}
	
	\begin{proof}
		We use the Néron bases described in the proof of~\cref{twisted Neron differentials}.
		Fixing an embedding $K \inj \C$ and a symplectic basis of~$\H_1(A(\C), \Z)$,
		we see for the resulting period matrices that $\Pi_{A^K} = \sqrt{D_K}^{-1} \Pi_A$.
	\end{proof}
	
	We can use this result, together with the following elementary statement
	about abelian groups with an involution, to relate the period of~$A/K$ to
	the real periods of $A/\Q$ and~$A^K/\Q$ when $K$ is an imaginary quadratic field.
	
	\begin{lemma} \label[lemma]{abelian groups with involution} \strut
		\begin{enumerate}[\upshape(1)]
			\item \label{invol:label1}
			Let $V$ be a finite dimensional $\F_2$-vector space and let $\iota \in \GL(V)$
			with $\iota^2 = \id_V$. Then $(V : (\id+\iota)(V)) = \#V^{\langle \iota \rangle}$.
			\item \label{invol:label2}
			Let $G$ be a finitely generated abelian group and let $\iota \in \Aut(G)$
			with $\iota^2 = \id_G$. Let $G_1 = \{g + \iota(g) : g \in G\}$
			and $G_2 = \{g - \iota(g) : g \in G\}$. Then
			\[ (G : G_1 + G_2) = \#\Bigl(\frac{G}{2G}\Bigr)^{\langle \iota \rangle} \,. \]
		\end{enumerate}
	\end{lemma}
	
	\begin{proof} \strut
		\begin{enumerate}[(1)]
			\item The map $\phi \defeq \id + \iota = \id - \iota$ has kernel $V^{\langle \iota \rangle}$.
			Since $V$ is finite, we have
			\[ (V : (\id+\iota)(V)) = \#\coker \phi = \#\ker \phi = \#V^{\langle \iota \rangle} \,. \]
			\item Since for each $g \in G$, $2g = (g + \iota(g)) + (g - \iota(g)) \in G_1 + G_2$,
			we have that $2G \subseteq G_1 + G_2$. Therefore, using part~\eqref{invol:label1},
			\[ (G : G_1 + G_2) = \Bigl(\frac{G}{2G} : \frac{G_1 + G_2}{2G}\Bigr)
			= \Bigl(\frac{G}{2G} : (\id + \iota)\Bigl(\frac{G}{2G}\Bigr)\Bigr)
			= \#\Bigl(\frac{G}{2G}\Bigr)^{\langle \iota \rangle} \,. \qedhere
			\]
		\end{enumerate}
	\end{proof}
	
	\begin{corollary} \label[corollary]{omega quotient}
		Let $A$ and~$K$ be as in~\cref{twisted period matrix}, with $D_K < 0$. Then
		\[ \frac{\Omega_{A/\Q} \Omega_{A^K/\Q} \sqrt{|D_K|}^{g}}{\Omega_{A/K}}
		= \frac{\#A(\R)[2]}{2^g}
		= \#\pi_0(A(\R)) \,. \]
	\end{corollary}
	
	\begin{proof}
		We use $\Pi_A$ and~$\Pi_{A^K}$ to denote big period matrices of $A$ and~$A^K$
		with respect to a Néron basis of the invariant $1$-forms.
		The period~$\Omega_{A/K}$ is $2^g$~times the covolume of the lattice $\Lambda \subseteq \C^g$
		(where, as usual, $g = \dim A$) generated by the columns of~$\Pi_A$
		(a Néron basis of~$\H^0(A/\Q, \Omega^1)$ gives a Néron basis
		of~$\H^0(A/K, \Omega^1)$, since at the bad places of~$A$, $K|\Q$ is unramified
		and Néron models are preserved by unramified base extension).
		The real periods~$\Omega_{A/\Q}$ and~$\Omega_{A^K/\Q}$ are the covolumes
		of the lattices in~$\R^g$ generated by the $\C|\R$-traces of the columns
		of $\Pi_A$ and~$\Pi_{A^K}$, respectively. The first lattice is
		$\Lambda_1 = \{\lambda + \bar{\lambda} : \lambda \in \Lambda\}$.
		By~\cref{twisted period matrix},
		$\Pi_{A^K} = D_K^{-1/2} \Pi_A$ (using suitable bases), which together with $D_K < 0$
		implies that $\Omega_{A^K/\Q}$ is $\sqrt{|D_K|}^{-g}$ times the covolume of
		the lattice in~$\R^g$ generated by the $\C|\R$-traces of the columns
		of~$\sqrt{-1} \cdot \Pi_A$. This is the same as the covolume of
		$\Lambda_2 \subseteq \sqrt{-1} \,\R^g$, where
		$\Lambda_2 = \{\lambda - \bar{\lambda} : \lambda \in \Lambda\}$.
		So $\Omega_{A/\Q} \Omega_{A^K/\Q} \sqrt{|D_K|}^{g}$ is the covolume
		of $\Lambda_1 + \Lambda_2 \subseteq \C^g$.
		Applying \hbox{\cref{abelian groups with involution}\,\eqref{invol:label2}}
		with $G = \Lambda$ and $\iota$ the restriction of complex conjugation, we finally obtain
		\[ \frac{\Omega_{A/\Q} \Omega_{A^K/\Q} \sqrt{|D_K|}^{g}}{\Omega_{A/K}}
		= \frac{(\Lambda : \Lambda_1 + \Lambda_2)}{2^g}
		= \frac{\#(\Lambda/2\Lambda)^+}{2^g}
		= \frac{\#A(\R)[2]}{2^g} \,,
		\]
		where $(\Lambda/2\Lambda)^+$ denotes the subgroup fixed under the induced action
		of complex conjugation. (Note that
		$A(\C)[2] \iso \tfrac{1}{2} \Lambda/\Lambda \iso \Lambda/2\Lambda$ as a $\Gal(\C|\R)$-module.)
		The last equality comes from the fact that $\#A(\R)^0[2] = 2^g$,
		since the connected component of the identity is a $g$-dimensional real torus.
	\end{proof}
	
	\begin{corollary} \label[corollary]{c_f c_pi under twist}
		Let $f \in \cN(N, g)$ and let $K$ be a quadratic number field such that $D_K$
		is coprime with~$2N$. Let further $\pi \colon A_f \to J$ be an isogeny
		defined over~$\Q$. We write $f^K$ and $\pi^K$ for the corresponding quadratic twists
		of~$f$ and~$\pi$, respectively. Then $c_{f^K} c_{\pi^K} = c_f c_\pi$.
	\end{corollary}
	
	\begin{proof}
		The map $\pi_f \colon J_0(N) \to A_f$ is geometrically defined, so $\pi_{f^K}$
		is the same as the quadratic twist $\pi_f^K$ of~$\pi_f$. Note that $c_f$ is the
		index of $\pi_f^* \H^0(\sA_f, \Omega^1_{\sA_f/\Z})$ in
		$\pi_f^*(A_f, \Omega^1) \cap \H^0(\sJ_0(N), \Omega^1_{\sJ_0(N)/\Z})$
		(where, as usual, $\sA_f$ and~$\sJ_0(N)$ denote the Néron models of~$A_f$
		and~$J_0(N)$ over~$\Z$). Applying~\cref{twisted Neron differentials}, we
		see that, considered inside~$\H^0(J_0(N)_K, \Omega^1)$, the images of both
		spaces are multiplied by~$1/\sqrt{D_K}$ by twisting, so the index stays the
		same. This shows that $c_{f^K} = c_f$. An analogous argument shows that
		$c_{\pi^K} = c_\pi$.
	\end{proof}

	\subsection{The case of $L$-rank $1$: Computing $\#\Sha(J/\Q)_\an$}
	
	We now show how we can compute $\#\Sha(J/\Q)_\an$ from $\#\Sha(J/K)_\an$
	as determined in Section~\ref{SS: Sha_an/K} and $\#\Sha(J^K/\Q)_\an$,
	which we can compute as in Section~\ref{ssec:Shaan rk 0} (and using
	Section~\ref{ssect: periods} to make the computation of the \enquote{$L$-ratio} for the
	quadratic twist feasible, which would otherwise be quite slow, as the level
	of the twisted newform tends to be fairly large)
	since the $L$-rank of~$J^K$ is zero by the Heegner hypothesis.
	
	From the induction formula $L(J/K,s) = L(J/\Q,s)L(J^K/\Q,s)$ we obtain
	\[
	L^{(g)}(J/K,1) = L^{(g)}(J/\Q,1) L(J^K/\Q,1) \,.
	\]
	
	Then one computes $\#\Sha(J/\Q)_\an$ from the relation
	\begin{align}
		\#\Sha(J/K)_\an
		&= \#\Sha(J/\Q)_\an \cdot \#\Sha(J^K/\Q)_\an \nonumber \\
		&\quad {} \cdot \frac{(\#J(K)_\tors)^2}%
		{(\#J(\Q)_\tors)^2 (\#J^K(\Q)_\tors)^2} \label{ShaShaSha} \\
		&\quad {} \cdot \frac{\Tam(J/\Q) \Tam(J^K/\Q)}{\Tam(J/K)}
		\cdot \frac{\Reg_{J/\Q}}{2^g \Reg_{J/K}}
		\cdot \frac{\Omega_J \Omega_{J^K} \sqrt{|D_K|}^{g}}{\Omega_{J/K}} \nonumber
	\end{align}
	that we obtain from~\eqref{eq:ShaAKan} and its analogues for $J/\Q$ and~$J^K/\Q$.
	
	As discussed at the beginning of Section~\ref{SS: Sha_an/K},
	both regulators are defined in terms of the normalized canonical height. This implies that
	\[ \frac{\Reg_{J/\Q}}{(\#J(\Q)_\tors)^2} \cdot \frac{(\#J(K)_\tors)^2}{\Reg_{J/K}}
	= (J(K) : J(\Q))^2 \,.
	\]
	Since we have computed~$J(K)$ already, we can easily determine this index.
	By~\cref{omega quotient}, the last factor is $\#J(\R)[2]/2^g$, assuming
	$D_K$ is odd.
	
	We can evaluate the factor involving Tamagawa numbers using the following result.
	This is not used in the proof of the formula in~\cref{formula for Sha_an rk 1}
	below, but will be useful for the example in Appendix~\ref{sec:7-torsion in Sha}.
	
	\begin{lemma} \label[lemma]{Tam quotient}
		We have
		\[ \Tam(J^K/\Q) = \Tam(J/\Q) \cdot \prod\limits_{p \mid D_K} c_p(J^K/\Q) \]
		and hence
		\[ \frac{\Tam(J/\Q) \Tam(J^K/\Q)}{\Tam(J/K)} = \prod\limits_{p \mid D_K} c_p(J^K/\Q) \,, \]
		where $c_p(J^K/\Q) = \#J(\F_p)[2]$ when $p \mid D_K$ is an odd prime.
	\end{lemma}
	
	\begin{proof}
		Since all bad primes~$p$ of~$J$ split in~$K$,
		the Tamagawa numbers of~$J^K$ at these primes are the same as those of~$J$ and
		also agree with the two Tamagawa numbers of~$J/K$ at the primes dividing~$p$.
		Since the only further primes of bad reduction for~$J^K$ are those
		dividing~$D_K$, we obtain the stated equalities.
		
		The last claim follows from the fact
		that $J^K$ has totally unipotent reduction at~$p$, which implies that
		there is no $2$-torsion in $\sJ^K(\F_p)^0$ (where $\sJ^K$ is the
		Néron model of~$J^K$), together with the fact that the component group
		is killed by~$2$ since $J^K$ obtains good reduction after a quadratic
		extension (see~\cite[Cor.~5.3.3.2]{HalleNicaise}). This gives an
		isomorphism between $J(\F_p)[2] \iso \sJ^K(\F_p)[2]$ and the $\F_p$-points
		of the component group.
	\end{proof}
	
	Recall that $\Tam(J/K) = \Tam(J/\Q)^2$ by~\eqref{Eq:Tam JK}. We then obtain the following.
	
	\begin{corollary} \label[corollary]{formula for Sha_an rk 1}
		Keep the notations and assumptions introduced so far. If $D_K$ is odd, then
		\begin{align*}
			\#\Sha(J/\Q)_\an
			&= \frac{\disc \Z[f]}{\disc \End_\Q(J)}
			\cdot \frac{4^g}{\#J(\R)[2] \cdot \Tam(J/\Q)} \\
			&\qquad{} \cdot \Bigl(\frac{I_{K,\pi}}{(J(K) : J(\Q)) \cdot u_K^{g}}\Bigr)^2
			\cdot \Bigl(\frac{L(J^K/\Q, 1)}{\Omega_{J^K}}\Bigr)^{-1} .
		\end{align*}
	\end{corollary}
	
	\begin{proof}
		Combine \eqref{ShaShaSha} with the remarks after it and
		with \cref{formula for Sha J/K an} and~\cref{formula for Sha_an rk 0},
		applied to~$J^K$.
	\end{proof}
	
	
	\section{Bounding the support of the Tate--Shafarevich group}
	\label{sec:finite support}
	
	Let $A/\Q$ be an absolutely simple $\GL_2$-type abelian variety with associated
	newform~$f$.
	In this section, we obtain an explicit bound on the support of the
	Tate--Shafarevich group coming from the Heegner point Euler system.
	This leads to an explicit description of a finite set of (regular) prime ideals~$\frp$
	of~$\Z[f]$ such that $\Sha(A/\Q)[\frp] = 0$ for all~$\frp$ not in this set.
	In the $L$-rank $0$ case we make the results of Kolyvagin--Logachëv~\cite{KolyvaginLogachev}
	explicit and in the $L$-rank $1$ case those of Nekovář~\cite{Nekovar2007}.
	We first prove a result on the vanishing of the first Galois cohomology group for
	irreducible~$\rho_\frp$ in Section~\ref{ssec:H1}. Specializing to the case where
	$A = J$ is a Jacobian for simplicity (so we do not have to deal with polarizations),
	we derive the explicit
	finite support for~$\Sha(J/\Q)$ in Section~\ref{ssec:finite support};
	see~\cref{thm:finite support,thm:finite support rank 1,Howard Euler system}.
	
	In the following subsection,
	$F$ will be a general number field and does not denote $\Frac \Z[f]$.
	
	\subsection{Vanishing of $\H^1(F(A[\frp])|F, A[\frp])$} \label{ssec:H1}
	
	We assume that $\frp$ is a regular prime ideal of~$\Z[f]$
	and set $p = p(\frp)$.
	The goal of this section is to show that the Galois cohomology group
	$\H^1(F(A[\frp])|F, A[\frp])$ vanishes when the mod~$\frp$ Galois representation
	is irreducible (and $p > 2$); see~\cref{irreducible implies trivial cohomology}.
	The vanishing
	of this group is an important input for~\cite[Proposition~5.10]{KolyvaginLogachev}.
	
	Let $F$ be a number field. 
	(This level of generality is needed later on and also useful for further applications,
	for example in our forthcoming work on the BSD conjecture over totally real fields with Pip Goodman.)
	Let $G \defeq \Gal(F(A[\frp])|F) \inj G_\frp^\max$ with $G_\frp^\max$ defined as in~\cref{def: Gfrpmax}.
	
	The main idea is that $\H^1(G,A[\frp]) = 0$ if $G$ contains a nontrivial homothety.
	Our arguments are purely group cohomological, without much arithmetic input.
	
	\begin{definition}
		A \emph{homothety} in the automorphism group of a vector space~$V$ over a
		field~$K$ is a map of the form $v \mapsto \lambda v$ with
		$\lambda \in K^\times$. It is \emph{nontrivial} if $\lambda \neq 1$.
	\end{definition}
	
	\begin{lemma}
		\label[lemma]{nontrivial homothety then H1 eq 0}
		Let $V$ be a finite-dimensional vector space over a finite field~$\F$
		and let $G$ be a subgroup of~$\GL(V)$.
		If $G$ contains a nontrivial homothety, then $\H^1(G, V) = 0$.
	\end{lemma}
	
	\begin{proof}
		(Compare~\cite[Lemma~3]{LawsonWuthrich2016}.)
		Let $g \in G$ be a nontrivial homothety; note that $\gen{g}$ is a normal subgroup
		of~$G$. Consider the associated inflation-restriction exact sequence:
		\[ \begin{tikzcd}
			0 \ar[r] & \H^1(G/\gen{g}, V^{\gen{g}}) \ar[r,"\inf"] & \H^1(G, V) \ar[r,"\res"]
			& \H^1(\gen{g}, V)
		\end{tikzcd}
		\]
		The left-hand group is trivial, since $V^{\gen{g}} = 0$
		(a nontrivial scalar matrix fixes no nontrivial element of a vector space),
		and the right-hand group is trivial because $\#\gen{g} \mid \#\F^\times$
		and $\#V = \#\F^{\dim{V}}$ are coprime. So the middle group must
		be trivial as well.
	\end{proof}
	
	\begin{lemma} \label[lemma]{trivial H1}
		Let $\F$ be a finite field of characteristic $p \ge 3$ and let \hbox{$G \subseteq \GL_2(\F)$}
		be such that $G$ does not fix a unique line in~$\F^2$. Then
		\[ \H^1(G, \F^2) = 0 \,. \]
	\end{lemma}
	
	\begin{proof}
		We proceed in a number of steps.
		\begin{enumerate}[(1)]
			\item \label{tH1_1} If $N$ is a normal subgroup in~$G$ and $N$ fixes a unique line,
			then so does~$G$. This is because $G$ acts on the lines fixed by~$N$.
			\item \label{tH1_2} If $N$ is a normal subgroup in~$G$ of index prime to~$p$,
			then (by inflation-restriction and since $\H^1(G/N, V) = 0$ for
			$V = (\F^2)^N$) $\H^1(N, \F^2) = 0$ implies $\H^1(G, \F^2) = 0$.
			\item \label{tH1_3} By \eqref{tH1_1} and~\eqref{tH1_2}, we can restrict to
			subgroups of~$\SL_2(\F)$, observing that $G \cap \SL_2(\F)$ is a normal
			subgroup of~$G$ of index dividing~$\#\F^\times$.
			\item \label{tH1_4} If $\#G$ is prime to~$p$, then $\H^1(G, \F^2) = 0$.
			We can therefore assume that $p$ divides~$\#G$ and therefore also $\#\P G$.
			\item \label{tH1_5} If $G$ contains~$-I$ (the unique nontrivial homothety
			in~$\SL_2(\F)$; note $p \geq 3$),
			then $\H^1(G, \F^2) = 0$ by~\cref{nontrivial homothety then H1 eq 0}.
			Since $-I$ is the unique element of order~$2$ in~$\SL_2(\F)$, this is
			the case whenever $\#G$ is even, so in particular when $\#\P G$ is even.
			\item \label{tH1_6} We consult~\cite[Thm.~2.1]{King2005}, which lists all
			subgroups of~$\PSL_2(\F)$. In cases (f), (i), (p) and~(u), $p = 2$.
			In cases (b)--(e), (g), (h), (j) and~(k), $p$ does not divide~$\#\P G$.
			In cases (n), (o), (q)--(t) and~(v), $\#\P G$ is even.
			In the remaining cases (a), (l) and~(m), $\P G$ is contained in a
			Borel subgroup and has order divisible by~$p$, so $G$ fixes a unique line.
			In each case, either the assumptions are violated, or we can conclude
			using \eqref{tH1_4} or~\eqref{tH1_5}.
			\qedhere
		\end{enumerate}
	\end{proof}
	
	\begin{proposition}[Irreducible implies trivial cohomology]
		\label[proposition]{irreducible implies trivial cohomology}
		Let $A$ be an absolutely simple abelian variety over~$\Q$ of $\GL_2$-type
		and let $\frp$ be a regular prime ideal of~$\End_\Q(A)$.
		Let $F$ be a number field that is a Galois extension of~$\Q$.
		We assume that $p(\frp) \ge 3$ and that $\rho_\frp|_{G_F}$ is irreducible.
		Then
		\[ \H^1(F(A[\frp])|F, A[\frp]) = 0 \,. \]
	\end{proposition}
	
	\begin{proof}
		Let $G$ be the image of~$\rho_\frp$ and let $G' \defeq \rho_\frp(G_F)$,
		which is a normal subgroup of~$G$. Since $\rho_\frp$ is irreducible,
		$G$ does not fix a unique line; by part~\eqref{tH1_1} of the proof
		of~\cref{trivial H1}, this implies that $G'$ also does not fix a
		unique line. Then \cref{trivial H1} says that
		\[ \H^1(F(A[\frp])|F, A[\frp]) = \H^1(G', \F_\frp^2) = 0 \,. \qedhere \]
	\end{proof}

	\subsection{Bounding the support of $\Sha(A/\Q)$} \label{ssec:finite support}
	
	Using our computations of $G_\frp$ from
	Section~\ref{sec:computation-of-the-residual-galois-representations}
	and the Heegner index from Section~\ref{sec:computation-of-the-heegner-points-and-index},
	we can improve~\cite{KolyvaginLogachev} and~\cite{Nekovar2007} to give an
	explicit finite bound for the support
	of $\Sha(A/\Q)$ considered as a $\Z$- or an $\cO$-module.
	
	We do not repeat the full proof, we only explain how to make the
	arguments explicit.
	
	\begin{assumption} \label[assumption]{assumption on K}
		Let $A$ be a modular abelian variety of level~$N$. We write $\cO \defeq \End_\Q(A)$
		and assume that $\cO$ is the maximal order of~$\Q(f)$ (via an isomorphism
		as in~\eqref{Eq:iso End0}). This is no essential restriction; compare~\cref{RM ab var}
		and note that the truth of the strong BSD~Conjecture is an isogeny invariant
		by~\cite[Theorem~I.7.3]{MilneADT}. (However, the support of~$\Sha$ can change
		under isogenies.)
		Let $K$ be a Heegner field of odd \emph{Heegner discriminant} $D_K \neq -3$
		(in particular, $D_K \notin \{-3, -4, -8\}$), i.e.,
		$K$ is an imaginary quadratic field such that all primes dividing the level~$N$
		split completely in~$K$. Then $y_{K,\pi} \in A(K)$ is a \emph{Heegner point}, and
		we assume that
		$\Lrk(A/K) = 1$, i.e., $y_{K,\pi}$ is non-torsion by the Gross--Zagier formula.
		Note that $y_{K,\pi}$ satisfies the Euler system relations
		from~\cite[§\,2]{KolyvaginLogachev} because the isogeny $\pi$ is
		equivariant with respect to the action of~$\Z[f] \subseteq \cO$;
		see~\eqref{Eq:pi equivariant}.
	\end{assumption}
	
	For several curves among our LMFDB examples, the endomorphism ring of the
	Jacobian is not maximal. However, in all these cases there is another curve in the
	database whose Jacobian is isogenous with endomorphism ring the maximal
	order; it then suffices to consider these other curves.
	
	In Table~\ref{tab:KL} we collect the most important objects and constants in~\cite{KolyvaginLogachev}.
	We specialize to the case that $A$ is the Jacobian~$J$ of a curve
	with its canonical principal
	polarization. In particular, since $J$ has RM, the Rosati involution
	associated to the polarization is the identity on $\End_\Q(J) \otimes_\Z \R$,
	which we need to use the results in~\cite[§\,2.1]{KolyvaginLogachev}.
	This implies that the polarization $\phi_\Lambda$ in Table~\ref{tab:KL} is principal.
	Recall that $\Tam(J/\Q) = \prod_\ell c_\ell(J/\Q)$ is the Tamagawa product of~$J$.
	The component group~$\pi_0(\sJ)$ of the Néron model~$\sJ/\Z$ of~$J/\Q$
	is an $\cO$-module. This allows us to consider its order~$\Tam(J/\Q)$
	as the corresponding characteristic ideal in~$\cO$ in the following.
	
	\begin{table}[t]
		{
			\begin{tabular}{rll}
				\toprule
				symbol & definition & properties\\
				\midrule
				$m_1$    & $\frp^{m_1}\cdot\Sel_{\frp^\infty}(J/\Q) = 0$ & $m_1 = m_3 + m_{10} + 2(m_9 + m_{11}) + m_{13}$ \\
				$m_2$    & $\ord_\frp(\Ann_{\cO} (\ker\phi_\Lambda))$ & $0$ if $\phi_\Lambda$ principal polarization \\
				$A$      & $\prod_v\Ann_\cO(\H^1(K_v^\nr|K_v, J))$ & divides $\Tam(J/K) = \Tam(J/\Q)^2$ \\
				$B$      & $\ord(h_K \cdot \lambda(j(\pi(0))))$ & divides $\#J(K)_\tors$ \\
				$x$      & $ABy_{K,\pi} \bmod \frp^n J(K)$ & \\
				$m_4$    & $[1]$ & $0$ if $\frp \nmid 2$ \\
				$m_6$    & $[3g] + m_2/2$ & $0$ if $\deg\phi_\Lambda = 1$ and $\frp \nmid 2$ \\
				$m_7$    & $gm_4 + m_6$ & $0$ if $\deg\phi_\Lambda = 1$ and $\frp \nmid 2$ \\
				$m_3$    & $[2] + 3m_7 + m_4 + m_2$ & \\
				& $= [12g+3] + 5 m_2/2$ & $0$ if $\deg\phi_\Lambda = 1$ and $\frp \nmid 2$ \\
				$m_9$    & Lemma~5.9 in~\cite{KolyvaginLogachev} & $0$ if $\frp \nmid 2$ and $\rho_\frp$ is irreducible \\
				$m_{10}$ & $\frp^{m_{10}} \H^1(K|\Q,J[\frp^n](K)) = 0$ & $0$ if $\frp \nmid 2$ \\
				$V$      & $K(J[\frp^{n + m_2(\frp)}])$ & \\
				$m_{11}$ & $\frp^{m_{11}} \H^1(V|K,J[\frp^n]) = 0$ & $0$ if $\rho_\frp$ irreducible \\
				$m_{13}$ & $r\cdot x \in \frp^n J(K) \implies r \in \frp^{n - m_{13}}\cO$ &
				$0$ if $\frp \nmid A B I_{K,\pi}$ \\
				\bottomrule
			\end{tabular}
		}
		\caption{The constants $m_k = m_k(\frp^{n})$ and important objects occurring
			in the proof of~\cite{KolyvaginLogachev}.
			The notation \enquote{$[m]$} denotes~$m$ when $p(\frp) = 2$ and~$0$ otherwise.}
		\label{tab:KL}
	\end{table}
	
	\begin{theorem}[Explicit finite support of $\Sha$ in the $L$-rank $0$ case]
		\label[theorem]{thm:finite support}
		Assume $\Lrk J = 0$. Suppose that $\frp$ is a maximal ideal of~$\cO$
		such that $\rho_\frp$ is irreducible and
		\[ \frp \nmid 2 \cdot \Tam(J/\Q) \cdot \gcd_K(\cI_{K,\pi}) \,, \]
		where $K$ runs through the Heegner fields for $J/\Q$. Then
		\[ \Sha(J/\Q)[\frp] = 0 \,. \]
	\end{theorem}
	
	\begin{proof}
		Note that the arguments in~\cite{KolyvaginLogachev} (there for a prime~$\ell$)
		also work for prime ideals~$\frp$; this is explained in~\cite[§\,7.1]{DograLeFourn}:
		annihilation of modules under~$\cO$ by~$p$ is translated to the annihilation
		under $\frp \mid p$ using the Chinese remainder theorem
		$\cO/p \isoto \bigoplus_{\frp \mid p}\cO/\frp^{e_\frp}$.
		We set $p \defeq p(\frp)$.
		
		Looking at Table~\ref{tab:KL}, all constants $m_i$ are $0$ for~$\frp$
		satisfying our hypotheses.
		
		\begin{enumerate}
			\item $m_3 = 0$ because $\frp \nmid 2$ and the polarization is principal.
			
			\item If $\frp \nmid \Tam(J/\Q)$, then $\frp \nmid A$: Let $p = p(\frp)$.
			Let $v \nmid p$ be a finite prime of~$K$ with residue field~$\F_v$.
			Let $\sJ$ be the Néron model of~$J/K$.
			By~\cite[Proposition~I.3.8]{MilneADT},\footnote{%
				Note the erratum at \url{https://jmilne.org/math/Books/index.html}.}
			\begin{align*}
				\H^1_\nr(K_v, J) &\iso \H^1(\F_v, \pi_0(\sJ)(\ol{\F}_v)) \\
				&\iso \varinjlim_n\H^1(\F_{v^n}|\F_v, \pi_0(\sJ)(\F_{v^n}))
			\end{align*}
			as $\cO$-modules.
			The last module is a subquotient of $\pi_0(\sJ)(\F_{v^n})$ since
			$\F_{v^n}|\F_v$ is cyclic~\cite[Proposition~1.7.1]{NSW2.3}.
			(Note that conjecturally, $\Tam(J/\Q)$ divides all Heegner indices~\cite[Conjecture~V.(2.2)]{GrossZagier1986};
			see also~\cref{formula for Sha J/K an}.)
			
			\item If $\frp \mid B$, then $J(\Q)[\frp] \neq 0$, so $\rho_\frp$ is reducible.
			
			\item $m_9 = 0$ if $\frp \nmid 2$ and $\rho_\frp$ is irreducible, because
			the irreducibility implies that the $p$-isogeny graph is reduced to a point;
			see~\cite[Lemma~5.9]{KolyvaginLogachev}.
			
			\item $\frp \nmid 2$ implies that $m_{10} = 0$, since $\#\Gal(K|\Q) = 2$
			is prime to~$\#A[\frp]$.
			
			\item $\rho_\frp$ irreducible implies $m_{11} = 0$
			by~\cref{irreducible implies trivial cohomology} with $F = K$.
			
			\item One can take $m_{13}$ to be
			\[ v_\frp(\cI_{K,\pi}) = v_\frp\big(\Char_{\cO}(J(K)/\cO y_K)\big) \,. \]
			This is because this choice of~$m_{13} = m_{13}(\frp^\infty)$
			satisfies~\cite[Proposition~5.12]{KolyvaginLogachev}.
		\end{enumerate}
		Hence $m_1 = 0$, so $\Sel_{\frp^\infty}(J/\Q) = 0$.
	\end{proof}
	
	To simplify notation below, we write $\sK_\frp(f)$ for the set of all Heegner fields~$K$
	such that $a_n(f) \not\equiv \epsilon_K(n) a_n(f) \pmod{\frp}$
	for some~$n$ coprime to $N$, where $\epsilon_K$
	is the nontrivial quadratic Dirichlet character associated with~$K|\Q$.
	In practice, we find a Heegner field $K$ such that $K \in \sK_\frp(f)$ for \emph{all} $\frp \nmid 2$
	by checking that for some small bound~$B$ the ideal
	\[ \bigl\langle a_n(f) - \epsilon_K(n) a_n(f) : (n,ND_K) = 1, \; n \leq B \bigr\rangle \]
	of~$\Z[f]$ has norm a power of~$2$ (note that the norm is always divisible by~$2$).
	Note that this ideal is non-zero
	if $f$ does not have CM by $\epsilon_K$ in the terminology of~\cref{def: CM form}.
	
	\begin{theorem}[Explicit finite support of $\Sha$ in the $L$-rank $1$ case]
		\label[theorem]{thm:finite support rank 1}
		Assume that $\Lrk J = 1$ and $J/\Q$ is simple and does not have~CM.
		Suppose that $\frp$ is a maximal ideal of~$\cO$
		such that $\rho_\frp$ is irreducible, $K \in \sK_\frp(f)$ and
		\[
		\frp \nmid 2 \cdot \Tam(J/\Q) \cdot \cI_{K,\pi} \,.
		\]
		Then
		\[ \Sha(J/\Q)[\frp] = \Sha(J/K)[\frp] = 0 \,. \]
	\end{theorem}
	
	\begin{proof}
		In the setting of~\cite{Nekovar2007}, $F = \Q$, $K$ is a Heegner field for $J/\Q$,
		the character $\alpha$ is trivial (see the main theorem at the beginning
		of~\cite{Nekovar2007}), therefore $K(\alpha) = K$ and hence $\beta = 1$ by~\cite[(3.1)]{Nekovar2007}.
		According to~\cite[7.3, 7.5, 7.5.3]{Nekovar2007}, we have to show that our hypotheses
		imply $C_i(\frp) = 0$ for $i = 1, \ldots, 6, 0$ (with $C_i(\frp)$ as defined in~\cite{Nekovar2007}):
		
		\begin{enumerate}[(1)]
			\item If $\frp \nmid \Tam(J/\Q)$, then $C_1(\frp) = 0$ because of the definition of~$C_1(\frp)$
			in~\cite[Proposition~5.12]{Nekovar2007} and by the argument in the proof
			of~\cref{thm:finite support}.
			
			\item If $\ker\big(\H^1(K,J[\frp]) \stackrel{\res}{\to} \H^1(K(J[\frp]),J[\frp])\big) = 0$,
			then $C_2(\frp) = 0$~\cite[Proposition~6.1.2]{Nekovar2007}.
			By~\cref{irreducible implies trivial cohomology} with $F = K$
			and the inflation-restriction sequence, this holds
			since $\rho_\frp|_{G_K}$ is irreducible by~\cref{rho restricted to K still irreducible}
			(here we use $D_K \neq -4, -8, -3$; compare~\cref{assumption on K}).
			
			\item Since $\alpha$ is trivial, $H = K(\alpha) = K$ in~\cite[§\,6]{Nekovar2007}.
			Hence there we only need to consider the Dirichlet character
			$\eta = \epsilon_K \colon \Gal(K|\Q) \to \{\pm1\}$
			of $\Gal(K|\Q)$ in~\cite[Proposition~6.2.2]{Nekovar2007}, thus~$C_3(\frp) = 0$ if
			$K \in \sK_\frp(f)$.
			
			\item We have $C_4(\frp) = 0$ because $\beta^2 = 1$.
			
			\item By definition~\cite[(7.4)]{Nekovar2007}, $C_5(\frp) = 0$ since $H = K(\alpha) = K$.
			
			\item One has $C_6(\frp) \defeq \ord_\frp \deg \phi$ with $\phi \colon J \to J^\dual$
			a polarization (for the Weil pairing)~\cite[(7.4)]{Nekovar2007}.
			Hence $C_6(\frp) = 0$ because $J$ is principally polarized.
			
			\item[(0)] Let $x \in J(K)$ be a Heegner point. One has
			\[ C_0(\frp) \defeq \max\{c \in \Z_{\geq 0} : x \in J(K)_\tors + \frp^c J(K)\} \,; \]
			see~\cite[(7.4)]{Nekovar2007}. Hence $C_0(\frp) = 0$ if $\frp \nmid \cI_{K,\pi}$.
		\end{enumerate}
		Note that these results imply $\Sha(J/K)[\frp] = 0$. Since $\frp \nmid 2$, $\Sha(J/\Q)[\frp] = 0$ follows.
	\end{proof}
	
	We use the refined information that is provided by considering $\Tam(J/\Q)$
	as an $\cO$-ideal in the following way.
	
	\begin{proposition} \label[proposition]{divisibility of Tamagawa numbers}
		We assume that $J$ is the Jacobian of a curve of genus~$2$.
		Fix an odd prime~$q$. Let $\sJ/\Z$ be the Néron model of $J/\Q$.
		If
		\begin{enumerate}[\upshape(i)]
			\item \label{item:Tam exactly one p}
			there is exactly one rational prime~$p$ with $v_q(c_p(J/\Q)) \geq 1$
			and we have $v_q(c_p(J/\Q)) = 1$ (then $a_p(f) \in \{\pm1,0\}$),
			
			\item \label{item:Tam exactly one reducible}
			$q \cO = \frq \frq'$ is split in $\cO$ with $\rho_{\frq'}$ irreducible,
			
			\item \label{item:Tam nontrival q-torsion over Q}
			$v_q(\exp(J(\Q)_\tors)) > v_q(p - a_p(f))$, where $\exp(J(\Q)_\tors)$
			denotes the exponent of the rational torsion subgroup of~$J$,
		\end{enumerate}
		then $v_{\frq'}(\Tam(J/\Q)) = 0$ and $v_\frq(\Tam(J/\Q)) = 1$.
	\end{proposition}
	
	\begin{proof}
		Since $q$ is odd, the $q$-primary part of~$J(\Q)_\tors$ injects
		into~$\sJ(\F_p)$. The group of~$\F_p$-points on the connected component
		of the identity of~$\sJ_{\F_p}$ has exponent $p - a_p(f)$ since it is a product
		of two copies of $\F_p^\times$ (when $a_p(f) = 1$) or of the norm~$1$
		subgroup of~$\F_{p^2}^\times$ (when $a_p(f) = -1$) or of~$\F_p$ (when $a_p(f) = 0$).
		So~\eqref{item:Tam nontrival q-torsion over Q} implies that the $q$-primary
		part of~$J(\Q)_\tors$ maps nontrivially into the $q$-primary part of the group of $\F_p$-points
		of the component group, $\pi_0(\sJ_{\F_p})(\F_p)$.
		By~\eqref{item:Tam exactly one p}, the latter is the same as the $q$-primary part
		of the ideal~$\Tam(J/\Q)$ and it has order~$q$, so by~\eqref{item:Tam exactly one reducible},
		this $q$-primary part of~$\Tam(J/\Q)$ is either $\frq$ or~$\frq'$.
		The map from the $q$-primary part of~$J(\Q)_\tors$ to~$\pi_0(\sJ_{\F_p})(\F_p)$
		respects the action of the endomorphism ring. Since $\rho_{\frq'}$ is irreducible,
		only~$\frq$ can occur in the characteristic ideal of~$J(\Q)_\tors$, so only~$\frq$
		can occur in~$\Tam(J/\Q)$.
	\end{proof}
	
	\begin{examples}
		There are three Jacobians of curves from the LMFDB for which we need to
		apply~\cref{divisibility of Tamagawa numbers} to show that the $\frp$-primary
		part of~$\Sha(J/\Q)$ is trivial for a degree~$1$ prime ideal~$\frp$ such that
		$\rho_\frp$ is irreducible. For one curve each at levels $N = 39$ and~$123$,
		the Tamagawa product is~$7$ and $7$ is split in the endomorphism ring. In both
		cases, there is rational $7$-torsion and $c_3 = 7$, so the proposition applies.
		In the last case, $N = 133$ and the Tamagawa product is~$3$, with $c_7 = 3$.
		We have $a_7 = 1$, so $v_3(7 - a_7(f)) = 1$, but luckily, $J(\Q) \isom \Z/9$,
		so condition~\eqref{item:Tam nontrival q-torsion over Q} above is still satisfied.
	\end{examples}
	
	Under stronger assumptions on~$\frp$ we can even get an upper bound for
	$\#\Sha(J/K)[\frp^\infty]$.
	
	\begin{theorem} \label[theorem]{Howard Euler system}
		Assume that $L'(f/K,1) \neq 0$. Let $\frp \mid p > 2$ be a regular prime ideal
		of~$\End_\Q(J)$ with $p \nmid h_K \cdot u_K \cdot N$.
		Suppose $\im\rho_{\frp^\infty} = G_{\frp^\infty}^\max$.
		Then
		\[
		\Sel_{\frp^\infty}(J/K) \isom (E_\frp/\cO_\frp) \oplus M \oplus M
		\]
		with $M$ finite of order bounded by $\#\cO_\frp/I_{K,\pi}$.
		Here, $E_\frp$ is the completion of $E = \End^0_\Q(J)$.
	\end{theorem}
	
	\begin{proof}
		This is~\cite[Theorem~A]{Howard2004}.
	\end{proof}
	
	Note that this implies $\Sha(J/K)[\frp^\infty] \isom M \oplus M$,
	so in particular,
	\[\#\Sha(J/K)[\frp^\infty] \le \#(\cO_\frp/I_{K,\pi})^2 \,. \]
	\cref{Howard Euler system} requires $p$ to be a prime of good reduction
	and $\rho_{\frp^\infty}$ to be surjective,
	which~\cref{thm:finite support,thm:finite support rank 1} do not.
	
	
	\section{Computing $\Sha(J/\Q)[\frp^\infty]$ using descent} \label{sec:descent}
	
	In this section, let $J$ be the Jacobian of a curve of genus~$2$ whose
	endomorphism ring $\cO = \End_\Q(J)$ is an order in a real quadratic field. Let $\frp$
	be a prime ideal of degree~$1$ of~$\cO$ of residue characteristic~$p$.
	Then $J$ is modular; let $N$ be the level.
	
	The results of the preceding section reduce the problem of showing that
	$\#\Sha(J/\Q) = \#\Sha_\an(J/\Q)$ to the verification that
	\[ v_p(\#\Sha(J/\Q)) = v_p(\#\Sha_\an(J/\Q)) \]
	for finitely many primes~$p$.
	The left hand side can be computed by studying $\Sha(J/\Q)[\frp^\infty]$
	for the prime ideals~$\frp$ of~$\cO$ dividing~$p$.
	In most cases, we need to show that $\Sha(J/\Q)[\frp^\infty] = 0$, for which
	it is sufficient to show that $\Sha(J/\Q)[\frp] = 0$. We can compute (the size of)
	$\Sha(J/\Q)[\frp]$ in principle by doing a \emph{descent}, i.e., by computing
	the $\frp$-Selmer group $\Sel_\frp(J/\Q)$ of~$J/\Q$. Recall that
	\[ \Sel_\frp(J/\Q) \defeq \ker\bigl(\H^1(\Q, J[\frp]) \to \prod_v \H^1(\Q_v, J(\bar{\Q}_v))\bigr) \]
	and that the Selmer group sits in the following exact sequence.
	\[ 0 \To \frac{J(\Q)}{\frp J(\Q)} \To \Sel_\frp(J/\Q) \To \Sha(J/\Q)[\frp] \To 0 \,. \]
	(There are analogous definitions with $p$ in place of~$\frp$.)
	This implies that
	\[ v_p\bigl(\#\Sha(J/\Q)[\frp]\bigr)
	= v_p(\#\Sel_\frp(J/\Q)) - (\deg \frp) \rk_{\cO}(J/\Q) - v_p(\#J(\Q)[\frp]) \,.
	\]
	So to verify that $\#\Sha(J/\Q)[\frp] = 0$, it is sufficient to show that
	\[ \dim_{\F_p} \Sel_\frp(J/\Q) \le (\deg \frp) \rk_{\cO}(J/\Q) + \dim_{\F_p} J(\Q)[\frp] \,. \]
	
	\subsection{Dealing with $p = 2$}\label{sec:dealing-with-p--2}
	
	We always need to determine~$\#\Sha(J/\Q)[2^\infty]$, since primes dividing~$2$
	are always excluded in~\cref{thm:finite support,thm:finite support rank 1}.
	Luckily, for Jacobians of hyperelliptic curves, the size of the $2$-Selmer group
	can be computed fairly easily. This is described in~\cite{Stoll2001} and is
	implemented in Magma. If this computation shows that $\Sha(J/\Q)[2] = 0$,
	then we know that $\#\Sha(J/\Q)$ is odd.
	
	The $2$-primary part of~$\Sha(J/\Q)$ is somewhat special, as its cardinality
	can be twice a square (the odd part, if finite, is always the square of some
	group). Using results of~\cite{PoonenStoll1999}, we can determine whether
	this is the case; in particular, the computation of $J(\Q)$ and the $2$-Selmer group is
	sufficient to detect that $\#\Sha(J/\Q)[2^\infty] = 2$. In all cases from our database where
	the $2$-part of $\#\Sha_\an(J/\Q)$ is $1$ or~$2$, this computation shows
	that $\#\Sha(J/\Q)[2^\infty]$ has the expected value.
	
	When $\#\Sha(J/\Q)[2] > 2$, but
	$\Sha(J/\Q)[4] = \Sha(J/\Q)[2]$, then this can be verified by computing the
	Cassels--Tate pairing on~$\Sel_2(J/\Q)$ (this is a symmetric bilinear form
	on the $\F_2$-vector space~$\Sel_2(J/\Q)$ whose kernel is the preimage
	of~$2 \Sha(J/\Q)[4]$). A method for doing this is described in the recent
	preprint~\cite{FisherYan} by Fisher and Yan, and an alternative approach
	will be detailed in forthcoming work by Shukla.
	
	There are ten cases in our database where $\#\Sha_\an = 4$ (in all cases,
	$\#\Sha_\an \in \{1,2,4\}$), corresponding to the levels and isogeny classes
	\begin{gather*}
		67c, \; 73a, \; 133e, \; 211a, \; 275a, \; 313a, \; 358a, \; 640a, \; 640b, \; 887a,
	\end{gather*}
	each of which occurs only once in the list. In each case, the $2$-descent
	computation shows that $\#\Sha(J/\Q)[2] = 4$ as expected.
	For the two curves at level~$640$
	(that are quadratic twists by~$-1$ of each other), there exists a Richelot
	isogenous Jacobian~$J'$, for which a $2$-descent shows that $\Sha(J'/\Q)[2] = 0$.
	This shows that in both cases, $\#\Sha(J/\Q)[2^\infty] = 4$, since elements
	of order~$4$ would have to survive the isogeny.
	
	To deal with the remaining eight cases, we need to compute the kernel
	of the Cassels--Tate pairing on~$\Sha(J/\Q)[2]$ and verify that this
	kernel is trivial.
	Fortunately, Fisher and Yan~\cite{FisherYan} have computed the pairing on the
	\hbox{$2$-Selmer} groups of all Jacobians of genus~$2$ curves in the LMFDB with
	even analytic order of~$\Sha$. In particular, they have verified that
	$\#\Sha(J/\Q)[2^\infty] = 4$ in all cases where $\#\Sha_\an(J/\Q) = 4$.
	This finishes the verification for the $2$-primary part of~$\Sha$
	in our LMFDB examples. There is one of the \enquote{Wang only} curves that also
	has $\#\Sha(J/\Q)[2] = \#\Sha(J/\Q)_\an = 4$ (the curve with label~125B);
	Tom Fisher has kindly checked for us using the code from~\cite{FisherYan}
	that $\#\Sha(J/\Q)[2^\infty] = 4$ for this curve as well.

	\subsection{Odd primes}\label{sec:odd-primes}
	
	We will now assume that $\frp$ is a prime ideal of~$\cO$ dividing an odd prime.
	We will also assume that $\deg \frp = 1$, as for primes of degree~$2$,
	the computation tends to get fairly involved. The situation is then
	analogous to that of a full $p$-descent (where $p = p(\frp)$) on an elliptic
	curve: the kernel of the isogeny is isomorphic to $\Z/p \times \Z/p$
	as a group, and it carries a Weil pairing. How to do a $p$-descent on an
	elliptic curve is analyzed in detail in~\cite{SchaeferStoll}; most of what
	is done below builds on this analysis. For the general theory of how to
	perform descent computations, see~\cite{BPS}. We assume that the prime ideal~$\frp$
	is principal, generated by a prime element~$\varpi$. This ensures that
	$J/J[\frp] \isom J$ via the multiplication-by-$\varpi$ map. This assumption
	is satisfied in all the examples that we have considered.
	
	We now assume in addition that $\rho_\frp$ is reducible; this is the most common
	situation when the general results do not allow us to conclude that
	$\Sha(J/\Q)[\frp] = 0$. We then have an exact sequence of Galois modules
	\[ 0 \To M_1 \To J[\frp] \To M_2 \To 0 \,, \]
	where $M_1$ and~$M_2$ are one-dimensional Galois modules corresponding to
	characters with values in~$\F_p^\times$.
	A frequently occurring case is that $J[\frp]$ contains a rational point
	of order~$p$; then $M_1 \isom \Z/p$, and the action on~$M_2$ is via the
	cyclotomic character~$\chi_p$. Let $A = J/M_1$ be the isogenous abelian
	surface; write $\phi \colon J \to A$ for the corresponding isogeny and
	$\psi \colon A \to J$ for the isogeny such that $\psi \circ \phi = \varpi$;
	its kernel is $\phi(J[\frp]) \iso M_2$. Then we have the associated Selmer
	groups $\Sel(\phi) \subseteq \H^1(\Q, M_1)$ and $\Sel(\psi) \subseteq \H^1(\Q, M_2)$
	and an exact sequence (see~\cite[Lemma~6.1]{SchaeferStoll})
	\[ 0 \To \frac{M_2(\Q)}{\phi(J[\frp](\Q))} \To \Sel(\phi) \To \Sel_\frp(J/\Q)
	\To \Sel(\psi) \,.
	\]
	This allows us to bound $\dim_{\F_p} \Sel_\frp(J/\Q)$ via
	\[ \dim_{\F_p} \Sel_\frp(J/\Q) \le \dim_{\F_p} \Sel(\phi) + \dim_{\F_p} \Sel(\psi)
	- \dim_{\F_p} \frac{M_2(\Q)}{\phi(J[\frp](\Q))} \,.
	\]
	
	Let $S$ be a finite set of primes. For a finite Galois module~$M$ such that
	$pM = 0$, we denote
	by $\H^1(\Q, M; S)$ the subgroup of cohomology classes unramified outside~$S$,
	i.e., mapping to zero in~$\H^1(I_q, M)$ for all primes $q \notin S$, where $I_q$
	is the inertia group at~$q$. (Since $p$ is odd, we can ignore the infinite place.)
	
	By Lemma~3.1 and the following text in~\cite{SchaeferStoll} (the arguments
	carry over from elliptic curves to abelian varieties), the Selmer group
	$\Sel(\theta)$ of an isogeny $\theta \colon A_1 \to A_2$ is contained
	in $\H^1(\Q, \ker \theta; S)$,
	where $S$ is the set of primes~$q$ such that the Tamagawa number $c_q(A_2)$
	is divisible by~$p$, together with~$p$. We set
	\[ S_J \defeq \{p\} \cup \{\text{$q$ prime} : p \mid c_q(J)\} \]
	and define $S_A$ in a similar way. Then
	\[ S_A \subseteq S'_J \defeq \{p\} \cup \{\text{$q$ prime} : q \mid N\} \,. \]
	By considering the reduction type of~$J$ at~$q$, we may be able to obtain
	a smaller upper bound for~$S_A$.
	
	Now let $M$ be a one-dimensional (over~$\F_p$) Galois representation given
	by the character $\chi \colon \GalQ \to \F_p^\times$. Let $M^\dual \defeq \Hom(M, \mu_p)$
	be the Cartier dual, with character $\chi_p \chi^{-1}$. Let $L$ be the
	fixed field of the kernel of~$\chi_p \chi^{-1}$. The degree of $L/\Q$ divides~$p-1$,
	so is prime to~$p$, so by inflation-restriction, we see that
	\[ \H^1(\Q, M) \iso \H^1(L, M)^{\Gal(L/\Q)} \iso \H^1(L, \mu_p)^{(1)} \iso (L^\times/L^{\times p})^{(1)} \,, \]
	where the superscript $(1)$ denotes the subspace on which the action of
	$\sigma \in \Gal(L/\Q)$ is given by multiplying by~$a_\sigma$ / raising to
	the $a_\sigma$th power, where $a_\sigma = \chi_p \chi^{-1}(\sigma) \in \F_p^\times$.
	We can restrict this isomorphism to the elements that are unramified
	outside~$S$. Here $\alpha L^{\times p} \in L^\times/L^{\times p}$ is considered to be
	unramified outside~$S$ when the extension $L(\sqrt[p]{\alpha})/L$ is unramified
	outside places above primes in~$S$; equivalently (when $p \in S$),
	$p$ divides all valuations $v_\frq(\alpha)$ for $\frq \mid q \notin S$.
	Denoting the subgroup of elements unramified outside~$S$ by~$L(S, p)$,
	this shows that $\Sel(\phi) \subseteq L_1(S_A, p)^{(1)}$ and $\Sel(\psi) \subseteq L_2(S_J, p)^{(1)}$,
	where $L_1$ and~$L_2$ are the fields associated to $M_1$ and~$M_2$, respectively.
	
	Since it occurs frequently, we give an explicit statement in the case
	that $J[\frp]$ contains a rational point of order~$p$ (then $M_1 = \Z/p$
	and $M_2 = \mu_p$). We will use the notation
	\[ [\cA] \defeq \begin{cases} 1, & \text{if $\cA$ is true} \\ 0, & \text{otherwise.} \end{cases} \]
	
	\begin{proposition} \label[proposition]{Sha_frp bound}
		Let $J$ be the Jacobian of a curve of genus~$2$ that has real multiplication;
		as before, let $N$ be its level. Let $\frp$
		be a prime ideal of degree~$1$ of~$\cO = \End_\Q(J)$ of residue characteristic~$p > 2$.
		Assume that $J[\frp](\Q) \isom \Z/p$ and that the class number of~$\Q(\mu_p)$
		is not divisible by~$p$. Then
		\begin{align*}
			\dim_{\F_p} \Sha(J/\Q)[\frp]
			&\le \#\{\text{$q$ prime} : \text{$q \mid N$ and $q \equiv 1 \bmod p$}\} \\
			&\qquad{} + \#\{\text{$q$ prime} : p \mid c_q(J)\} + [p \mid N] - \rk_{\cO} J(\Q) \,.
		\end{align*}
		If in addition there is a prime $\ell \equiv 1 \bmod p$ such that the natural map
		\[ r_\ell \colon \Q(\{q : \text{$p \mid c_q(J)$ or $q = p \mid N$}\}, p)
		\to \Q_\ell^\times/\Q_\ell^{\times p} \]
		is nontrivial and the map
		\[ r'_\ell \colon \Q(\mu_p)(\{p\} \cup \{q : \text{$q \mid N$ and $q \equiv 1 \bmod p$}\}, p)^{(1)}
		\to \Q_\ell^\times/\Q_\ell^{\times p} \]
		induced by any embedding $\Q(\mu_p) \to \Q_\ell$ is surjective,
		then the above inequality is strict.
	\end{proposition}
	
	\begin{proof}
		Since $M_1 = \langle P \rangle \isom \Z/p$, we have $L_1 = \Q(\mu_p)$.
		Similarly, since $M_2 \isom \mu_p$,
		we have $L_2 = \Q$. Let $F \in \Q(X)^\times$
		be a function whose divisor is $p D$, where the linear equivalence class of~$D$
		is~$P$; then the \enquote{descent map}
		$\delta \colon J(\Q) \to \H^1(\Q, M_2) \iso \Q^\times/\Q^{\times p}$
		is given by evaluating~$F$ on a representative divisor whose support is disjoint
		from that of~$D$. We have the analogous map $\delta_q \colon J(\Q_q) \to \Q_q^\times/\Q_q^{\times p}$
		for each prime~$q$. By the above, $\Sel(\psi \colon A \to J)$ is contained
		in the unramified outside~$S_J$ part of $\Q^\times/\Q^{\times p}$, which is
		the subgroup generated by the classes of the primes in~$S_J$.
		When $p \nmid N$, so that $J$ has good reduction at~$p$, then we can
		choose~$F$ in such a way that its reduction mod~$p$ is well-defined.
		When evaluating~$\delta_p$ on a point~$Q \in J(\Q_p)$, we can pick a representative
		divisor whose support is disjoint mod~$p$ from the support of~$D$;
		this shows that $F(Q)$ is in the image of~$\Z_p^\times$ and hence that
		$\Sel(\psi) \subseteq \langle q : p \mid c_q(J) \rangle$ in this case.
		So in any case, we have
		\begin{equation} \label{E:Sel psi}
			\dim_{\F_p} \Sel(\psi) \le \#\{\text{$q$ prime} : p \mid c_q(J)\} + [p \mid N] .
		\end{equation}
		On the other hand, $\Sel(\phi) \subseteq L_1(S'_J, p)^{(1)} = \Q(\mu_p)(S'_J, p)^{(1)}$.
		Since we assume that the class number of~$\Q(\mu_p)$ is not divisible by~$p$,
		the group $\Q(\mu_p)(S'_J, p)$ is generated by the images of a primitive $p$th
		root of unity~$\zeta_p$, a choice of fundamental units, $1 - \zeta_p$, and
		one element generating a suitable power of each prime ideal above a prime in~$S'_J$.
		Of these, only $\zeta_p$ and the totally split primes contribute to the relevant
		eigenspace (the fundamental units come, up to index prime to~$p$, from the maximal
		real subfield, and the ideal $\langle 1 - \zeta_p \rangle$ and the nonsplit
		prime ideals have nontrivial stabilizer), and the contribution of each totally
		split prime~$q$ is exactly~$1$. Since $q$ is totally split in~$\Q(\mu_p)$
		if and only if $q \equiv 1 \bmod p$, we obtain the bound
		\begin{equation} \label{E:Sel phi}
			\begin{array}{r@{\;}c@{\;}l}
				\dim_{\F_p} \Sel(\phi)
				&\le& \dim_{\F_p} \Q(\mu_p)(S'_J, p)^{(1)} \\[6pt]
				&=& 1 + \#\{\text{$q$ prime} : \text{$q \mid N$ and $q \equiv 1 \bmod p$}\} .
			\end{array}
		\end{equation}
		Finally, we have
		\begin{align*}
			\dim_{\F_p} \Sha(J/\Q)[\frp]
			&= \dim_{\F_p} \Sel(J/\Q)[\frp] - \dim_{\F_p} \frac{J(\Q)}{\frp J(\Q)} \\
			&\le \dim_{\F_p} \Sel(\phi) + \dim_{\F_p} \Sel(\psi) - 1 - \rk_{\cO} J(\Q) \\
			&\le \#\{\text{$q$ prime} : \text{$q \mid N$ and $q \equiv 1 \bmod p$}\} \\
			&\qquad {} + \#\{\text{$q$ prime} : p \mid c_q(J)\} + [p \mid N] - \rk_{\cO} J(\Q) ,
		\end{align*}
		where we have used $\dim J(\Q)/\frp J(\Q) = \rk_{\cO} J(\Q) + \dim J(\Q)[\frp]$
		and Equations \eqref{E:Sel psi} and~\eqref{E:Sel phi}.
		
		To show the refinement, consider the kernel-cokernel exact sequence associated
		to $J(\Q_\ell) \stackrel{\phi}{\To} A(\Q_\ell) \stackrel{\psi}{\To} J(\Q_\ell)$,
		\begin{align*}
			0 &\To \Z/p \To J(\Q_\ell)[\frp] \stackrel{\phi}{\To} \mu_p(\Q_\ell) \To \\
			&\To \frac{A(\Q_\ell)}{\phi(J(\Q_\ell))} \stackrel{\psi}{\To} \frac{J(\Q_\ell)}{\frp J(\Q_\ell)}
			\To \frac{J(\Q_\ell)}{\psi(A(\Q_\ell))} \To 0 \,.
		\end{align*}
		Since $\ell \equiv 1 \bmod p$, $\dim \mu_p(\Q_\ell) = 1$. Since $p \nmid \ell$,
		$\dim J(\Q_\ell)/\frp J(\Q_\ell) = \dim J(\Q_\ell)[\frp]$. These facts imply that
		\[ \dim \frac{A(\Q_\ell)}{\phi(J(\Q_\ell))} + \dim \frac{J(\Q_\ell)}{\psi(A(\Q_\ell))} = 2 \,. \]
		We note that the elements of~$\Sel(\phi)$ (respectively, $\Sel(\psi)$) map
		into the image of $A(\Q_\ell)/\phi(J(\Q_\ell))$ in $\H^1(\Q_\ell, \Z/p) \isom \Q_\ell^\times/\Q_\ell^{\times p}$
		(respectively, into the image of $J(\Q_\ell)/\psi(A(\Q_\ell))$
		in $\H^1(\Q_\ell, \mu_p) \iso \Q_\ell^\times/\Q_\ell^{\times p}$) under $r'_\ell$
		(respectively, $r_\ell$).
		If $\dim J(\Q_\ell)/\psi(A(\Q_\ell)) = 0$, then $\Sel(\psi)$ is contained in
		the kernel of~$r_\ell$; since $r_\ell$ is assumed to be nontrivial, this implies
		that the bound on $\dim \Sel(\psi)$ can be reduced by~$1$.
		Otherwise, $\dim A(\Q_\ell)/\phi(J(\Q_\ell)) \le 1 < 2 = \dim \Q_\ell^\times/\Q_\ell^{\times p}$.
		Since $r'_\ell$ is assumed to be surjective, this implies that the bound on
		$\dim \Sel(\phi)$ can be reduced by~$1$. So in all cases, the bound on
		$\dim \Sel(\phi) + \dim \Sel(\psi)$ is reduced by~$1$, which gives a corresponding
		improvement for the bound on~$\dim \Sha(J/\Q)[\frp]$.
	\end{proof}
	
	\begin{remark} \label[remark]{refinement for bad p}
		In some cases (for example in~\cref{ex desc 7}\,\eqref{Wang175}),
		we can improve the bound in~\cref{Sha_frp bound} by~$1$
		when $p \mid N$. Assume that we can find enough \enquote{descent functions} $F \in \Q(X)^\times$
		(i.e., whose divisor is $p$~times a divisor~$D$ such that $[D] \in J(\Q)$
		generates~$J[\frp]$) that reduce to well-defined functions~$\bar{F}$ on~$X_{\F_p}$
		and such that the support of~$\bar{F}$ consists of smooth points on~$X_{\F_p}$
		and that the following holds: if $D$ is a divisor of degree zero on~$X$ defined over~$\Q$ with
		reduction~$\bar{D}$ modulo~$p$, then we can find some~$F$ such that the divisor of~$\bar{F}$
		has support disjoint from~$\bar{D}$. Then the argument near the beginning
		of the proof of~\cref{Sha_frp bound} shows that $F(D) \in \Z_p^\times$, hence
		\[ \Sel(\psi) \subseteq \Q(\{\text{$q$ prime} : p \mid c_q(J)\}, p) \,, \]
		and so we can remove the term \enquote{$[p \mid N]$} in the final bound.
	\end{remark}
	
	We also note the following.
	
	\begin{lemma} \label[lemma]{dim Sha even}
		If we assume additionally that $\rk_{\cO} J(\Q) \le 1$ in the situation above,
		then $\dim_{\F_p} \Sha(J/\Q)[\frp]$ is even.
		In particular, $\dim_{\F_p} \Sha(J/\Q)[\frp] \le 1$ implies that
		$\Sha(J/\Q)[\frp] = 0$.
	\end{lemma}
	
	\begin{proof}
		By the general results on abelian surfaces with RM, $J$ is modular.
		The additional assumption then implies by~\cite{KolyvaginLogachev} that $\Sha(J/\Q)$ is finite.
		Therefore the Cassels--Tate pairing on~$\Sha(J/\Q)$ is perfect. It is also
		anti-symmetric, which implies (since $p$ is odd) that its restriction
		to~$\Sha(J/\Q)[\frp^\infty]$ is perfect and alternating. This in turn
		implies that $\Sha(J/\Q)[\frp^\infty] \isom M \times M$ for some finite
		$\cO_\frp$-module~$M$. In particular, $\Sha(J/\Q)[\frp] \isom M[\frp] \times M[\frp]$,
		and so $\dim_{\F_p} \Sha(J/\Q)[\frp] = 2 \dim_{\F_p} M[\frp]$.
	\end{proof}
	
	\begin{examples} \label[examples]{ex desc 1}
		For most pairs $(X, p)$ consisting of a curve~$X$ in our database
		of LMFDB curves and an odd prime~$p$ such that
		the table in Section~\ref{sec:examples-rho-frp} says that
		the semisimplification of~$\rho_\frp$ splits as $\mathbf{1} \oplus \chi_p$,
		where $\frp$ is a prime ideal of the endomorphism ring of degree~$1$,
		the non-strict bound in~\cref{Sha_frp bound} together with~\cref{dim Sha even} show
		that $\Sha(J/\Q)[\frp] = 0$. The exceptions are as follows. (We frequently
		indicate the isogeny class with a letter appended to the level~$N$.)
		\begin{enumerate}[(i)]
			\item Two of the four curves at level $N = 31$, where $p = 5$ (one has no rational
			point of order~$5$, the other gives a bound $\dim \Sha(J/\Q)[\frp] \le 2$).
			This is unproblematic, since these Jacobians are isogenous to the
			Jacobians of the two other curves at that level, for which the
			simple bound proves that $\Sha(J/\Q)[\frp] = 0$. By invariance of BSD
			in isogeny classes, it suffices to verify strong BSD for one of
			these Jacobians.
			\item \label{fail case ii}
			The pairs $(N, p) = (73a, 3)$ and $(85b, 3)$, where there
			is no rational point of order~$3$.
			\item \label{fail case iii}
			The curves at level~$133d$ with $p = 3$ and at level~$275a$ with $p = 5$,
			where the non-strict bound only gives $\dim \Sha(J/\Q)[\frp] \le 2$.
		\end{enumerate}
	\end{examples}
	
	\begin{examples} \label[examples]{ex desc 2}
		We consider the curves listed under~\eqref{fail case ii} above.
		In these cases, the $\frp$-torsion (where $\frp$ is an ideal of norm~$3$)
		sits in an exact sequence
		\[ 0 \To \mu_3 \To J[\frp] \To \Z/3\Z \To 0 \,. \]
		Using the fact that $M_2(\Q)/\phi(J[\frp]) \isom \Z/3\Z$ in these cases and that
		the Tamagawa numbers at all bad primes are not divisible by~$3$
		(and $3$ is not a bad prime, which implies that
		$\Sel(\psi) \subseteq \Q(\{\text{$q$ prime} : q \mid N\}, 3)$),
		this leads to a general bound of the form
		\begin{align*}
			\dim_{\F_3} \Sha(J/\Q)[\frp]
			&\le \dim_{\F_3} \Sel(\psi) - 1 + \dim_{\F_3} \Sel(\phi) - \rk_{\cO} J(\Q) \\
			&\le \dim_{\F_3} \Q(\{\text{$q$ prime} : q \mid N\}, 3) - 1 \\
			&\qquad {} + \dim_{\F_3} \Q(\mu_3)(\{3\}, 3)^{(1)} - \rk_{\cO} J(\Q) \\
			&\le \#\{\text{$q$ prime} : q \mid N\} - \rk_{\cO} J(\Q) ,
		\end{align*}
		which evaluates to $1$ and~$2$, respectively, for $N = 73$ and $85$.
		Using~\cref{dim Sha even}, this already shows that $\Sha(J/\Q)[\frp] = 0$
		for $N = 73$.
		
		To improve the bound for $N = 85$, we note that in this case
		\[ \Sel_\frp(J/\Q) \subseteq \H^1(\Q, J[\frp]; S)
		\subseteq \Q(\sqrt[3]{5})(\{3\}, 3) \times \Q(\mu_3)(\{3\}, 3) \]
		and that $\Sel(\psi) \subseteq \Q(\{5, 17\}, 3)$ maps into the first factor
		by the obvious map. (Compare Section~6 of~\cite{SchaeferStoll};
		one can check that the algebra that is called~$D$ there is a product
		of two copies of~$\Q(\sqrt[3]{5})$.) Since any element involving~$17$
		will have image ramified at~$17$, this shows that actually,
		$\Sel(\psi) \subseteq \Q(\{5\}, 3)$ (we even have equality here, since
		we know the map to $\Sel_\frp(J/\Q)$ has one-dimensional kernel),
		thus improving the bound by~$1$, which is sufficient to conclude.
	\end{examples}
	
	\begin{examples} \label[examples]{ex desc 3}
		We now consider the cases listed under~\eqref{fail case iii} above.
		For both pairs $(N, p)$, we can use the strict bound in~\cref{Sha_frp bound}.
		The non-strict bound for~$\dim \Sha(J/\Q)[\frp]$ is~$2$ in both cases.
		
		In the first case, $(7 \cdot 19, 3)$ for the curve in isogeny class~133d,
		we use $\ell = 7$. Here $c_7 = 3$ and $c_{19} = 1$, so the first condition
		is that $\Q(\{7\}, 3) \to \Q_7^\times/\Q_7^{\times 3}$ is nontrivial,
		which is clearly the case. The second condition is that
		$\Q(\mu_3)(\{3, 7, 19\}, 3)^{(1)} \to \Q_7^\times/\Q_7^{\times 3}$ is
		surjective, which follows from the fact that $\Q_7$ does not contain
		a primitive ninth root of unity. So both conditions are satisfied,
		and the bound can be improved to $\dim \Sha(J/\Q)[\frp] \le 1$,
		which is sufficient to conclude that $\Sha(J/\Q)[\frp] = 0$ by~\cref{dim Sha even}.
		
		We now consider $(5^2 \cdot 11, 5)$ for the curve in isogeny class~275a.
		Here we use $\ell = 11$. Both Tamagawa numbers are~$1$, so the first
		condition is that $\Q(\{5\}, 5) \to \Q_{11}^\times/\Q_{11}^{\times 5}$
		is nontrivial. This follows from the fact that $5$ is not a fifth
		power in~$\F_{11}$. The second condition is that
		$\Q(\mu_5)(\{5, 11\}, 5)^{(1)} \to \Q_{11}^\times/\Q_{11}^{\times 5}$
		is surjective. This follows in the same way as for the previous example.
		So we can again reduce the bound by~$1$ and obtain that $\Sha(J/\Q)[\frp] = 0$.
		
		We note that with some more work, one can show that in both these
		cases we have $\Sel(\psi) = 0$.
	\end{examples}
	
	The remaining cases (where $\rho_\frp^{\ss}$ splits into two nontrivial characters)
	come from the following pairs of level (+ isogeny class) and prime.
	\begin{gather*}
		(125a, 5), \; (147a, 7), \; (245a, 7), \; (250a, 5), \; (275b, 3) \\
		(289a, 3), \; (289a, 17), \; (375a, 5), \; (841a, 29) .
	\end{gather*}
	
	\begin{examples} \label[examples]{ex desc 4}
		We consider the curve at level~$N = 125 = 5^3$ and $\frp = \langle \sqrt{5} \rangle$.
		According to~\cref{ex reducible 1}, \eqref{ex_red_125}, we have
		$\rho_\frp \isom \chi_5^2 \oplus \chi_5^3$. This implies that
		\[ \Sel_\frp(J/\Q) \subseteq \Q(\mu_5)(\{5\}, 5)^{(-1)} \oplus \Q(\sqrt{5})(\{5\}, 5)^{(2)} \,, \]
		where the superscript~$(m)$ indicates that the action is via $\chi_5^m$.
		One finds easily that the first summand is trivial and the second has dimension~$1$.
		Since the $\cO$-rank of~$J(\Q)$ is~$1$, we obtain the bound
		\[ \dim_{\F_5} \Sha(J/\Q)[\sqrt 5] \le 0 + 1 - 1 = 0 \,, \]
		so $\Sha(J/\Q)[\sqrt 5] = 0$.
		
		Similarly, for the curves at levels $2 \cdot 5^3$ and~$3 \cdot 5^3$,
		we have (for $\frp = \langle \sqrt{5} \rangle$) that
		$\rho_\frp^{\ss} \isom \chi_5^2 \oplus \chi_5^3$. Since $2$ and~$3$ are primitive
		roots mod~$5$, we still have that $\Q(\mu_5)(\{q, 5\}, 5)^{(-1)}$ is trivial
		and $\Q(\sqrt{5})(\{q, 5\}, 5)^{(1)}$ is one-dimensional, where $q = 2$ or~$3$.
		This gives the bound
		\[ \dim_{\F_5} \Sha(J/\Q)[\sqrt 5] \le 0 + 1 - \rk_{\cO} J(\Q) \le 1 \]
		(the rank is zero for the curve at level $2 \cdot 5^3$ and one for the curve
		at level $3 \cdot 5^3$), which again suffices to conclude that
		$\Sha(J/\Q)[\sqrt 5] = 0$.
	\end{examples}
	
	The two pairs $(N, p) = (17^2, 17)$ and $(29^2, 29)$ can be dealt with in a similar
	way as $(N, p) = (5^3, 5)$.
	
	\begin{examples} \label[examples]{ex desc 5}
		We now consider $(N, p) = (3 \cdot 7^2, 7)$ and $(5 \cdot 7^2, 7)$.
		In both cases, $\cO = \Z[\sqrt{2}]$, and $\rho_\frp$ is reducible for exactly
		one of the two prime ideals above~$7$, with
		$\rho_\frp^{\ss} \isom \chi_7^3 \oplus \chi_7^4$. The two relevant groups are
		$\Q(\mu_7)^+(\{q,7\}, 7)^{(4)}$ and $\Q(\sqrt{-7})(\{q,7\}, 7)^{(3)}$,
		with $q = 3$ or $q = 5$. Since both are primitive roots mod~$7$, $q$ does
		not contribute to the relevant eigenspaces, and it is easy to see that the
		first group has dimension~$1$, whereas the second one is trivial.
		The $\cO$-rank is~$1$ in both cases, which directly shows that
		$\Sha(J/\Q)[\frp] = 0$.
	\end{examples}
	
	\begin{examples} \label[examples]{ex desc 6}
		For $(N, p) = (17^2, 3)$, we have
		$\rho_\frp^{\ss} \isom \epsilon_{17} \oplus \epsilon_{-3 \cdot 17}$,
		where $\epsilon_m$ denotes the quadratic character mod~$m$. This is one of
		two cases in our examples where the two characters are not powers of~$\chi_p$.
		The two relevant groups are
		$\Q(\sqrt{-3 \cdot 17})(\{3, 17\}, 3)^-$, which is trivial, and
		$\Q(\sqrt{17})(\{3, 17\}, 3)^-$, which has dimension~$1$. The $\cO$-rank is~$1$,
		which directly gives that $\Sha(J/\Q)[\frp] = 0$.
		
		The other similar case is $(N, p) = (5^2 \cdot 11, 3)$, where
		$\rho_\frp^{\ss} \isom \epsilon_{5} \oplus \epsilon_{-3 \cdot 5}$ and we obtain
		the same bound with a similar argument. (Note that $5$ and~$11$, like~$17$,
		are primitive roots mod~$3$.)
	\end{examples}
	
	\begin{examples} \label[examples]{ex desc 7}
		There are three \enquote{Wang only} curves for which a $\frp$-descent is necessary,
		namely
		\begin{enumerate}[(i)]
			\item \label{Wang117} Curve 117B with $\frp \mid 7$,
			\item \label{Wang125} Curve 125B with $\frp \mid 5$, and
			\item \label{Wang175} Curve 175 with $\frp \mid 5$.
		\end{enumerate}
		In all cases, we have to show that $\Sha(J/\Q)[\frp] = 0$.
		
		In case~\eqref{Wang117}, we have an exact sequence
		\[ 0 \To \epsilon_{-3} \To J[\frp] \To \epsilon_{-3} \cdot \chi_7 \To 0 \,. \]
		This case can be dealt with in a similar way as in~\cref{ex desc 5};
		we obtain $\dim \Sha(J/\Q)[\frp] \le 1$.
		
		In case~\eqref{Wang125}, $J[\frp] \isom \mathbf{1} \oplus \chi_5$ is split.
		It can be dealt with similarly to~\cref{ex desc 4}, leading again
		to $\dim \Sha(J/\Q)[\frp] \le 1$.
		
		Finally, in case~\eqref{Wang175}, we have a non-split exact sequence
		\[ 0 \To \mathbf{1} \To J[\frp] \To \chi_5 \To 0 \,, \]
		and $\frp = \langle \sqrt{5} \rangle$.
		The non-strict bound from~\cref{Sha_frp bound} gives us only $\dim \Sha(J/\Q)[\frp] \le 2$:
		we have $\Sel(\phi) \subseteq \Q(\mu_5)(\{5, 7\}, 5)^{(1)} = \Q(\mu_5)(\emptyset, 5)^{(1)}$,
		which has dimension~$1$, and $\Sel(\psi) \subseteq \Q(\{5, 7\}, 5)$ of
		dimension~$2$. We can improve the bound using~\cref{refinement for bad p}.
		The functions
		\[ \frac{(8 x^5 + 5 x^4 + 15 x^3 + 30 x^2 + 15 x + 4) \pm (6 x^2 + x - 2) y}{(x - a)^5} \]
		for $a \in \{0, \pm 1, \pm 2, \infty\}$ (where we set $1/(x - \infty)^5 \defeq 1$)
		form a suitable set of descent functions on the model
		\[ X \colon y^2 = x^6 - 2 x^5 - 3 x^4 - 6 x^3 - 14 x^2 - 8 x - 3 \]
		of the curve, so in fact $\Sel(\psi) \subseteq \Q(\{7\}, 5)$
		has dimension at most~$1$, which is enough to conclude that $\Sha(J/\Q)[\frp] = 0$.
		(In fact, $\dim \Sel(\psi) = 1$, since evaluating a suitable descent function
		on a point of order~$5$ gives~$7^2$, which is a fifth power in~$\Q_5$, but
		nontrivial in~$\Q_7^\times/\Q_7^{\times 5}$.)
	\end{examples}
	
	
	\section{Bounding $\Sha(J/\Q)[\frp^\infty]$ using Iwasawa Theory} \label{sec:Iwasawa theory}
	
	The Main Conjectures for modular forms in Iwasawa Theory imply the $p$-part
	of strong BSD under certain conditions when the $L$-rank is~$0$ or~$1$; see~\cref{Skinner_cor2}
	below and~\cite{CCSS}. For fixed~$p$ of good ordinary reduction, we can also compute an
	approximation to the $p$-adic $L$-function and use the known results on the
	$p$-adic BSD conjecture (see~\cref{p-adic BSD} below) to determine
	or at least bound the $p$-valuation of the order of~$\Sha$;
	this works even for higher rank, assuming Schneider's conjecture on the non-vanishing
	of the $p$-adic Schneider regulator (see~\cite{SteinWuthrich2013}).
	
	As we will only need these results in the case that $p$ is a prime of good
	ordinary reduction that is inert in~$\Z[f]$, we restrict to this situation in
	the following, although more general results are available,
	which apply more generally in the case when there is
	exactly one prime ideal of~$\Z[f]$ lying above $p$.
	
	This section generalizes~\cite{SteinWuthrich2013} from elliptic curves
	to modular abelian varieties of arbitrary dimension.
	The results in Sections~\ref{ssec:IGMC} and~\ref{ssec: Iwasawa theory}
	are not fully used in this paper, but
	can be used to extend the verification of strong BSD to examples
	not contained in our database.
	
	Let $f$ be a newform of level~$N$ with coefficient ring~$\Z[f]$.
	We assume that $a_p(f)$ is a $p$-adic unit,
	i.e., that $p$ is \emph{ordinary} for $f$. This implies that
	the Euler factor at~$p$ of~$f$ (equivalently, the characteristic polynomial of
	the Frobenius at~$p$ on the associated compatible system of $\ell$-adic Galois representations)
	has exactly one root that is a $p$-adic unit; we denote
	this root by~$\alpha$. We use the following notation. Let $f$ be a newform and
	$\frp \mid p$ a prime of $\Z[f]$.
	Let $\sL_\frp(f,T) \in \Q(f)_\frp\llbracket T\rrbracket$ be the \emph{$\frp$-adic $L$-function}
	of~$f$ constructed in~\cite[§\,2.2]{BMS}. If $A$ is the abelian variety associated to~$f$,
	$\sL_\frp(A,T) = \prod_{\sigma\colon \Z[f] \inj \R} \sL_\frp(f^\sigma, T)$; see~\cite[§\,2.3]{BMS}.
	
	\subsection{The Iwasawa--Greenberg Main Conjecture} \label{ssec:IGMC}
	
	We use the following known cases of the $\GL_2$ Iwasawa--Greenberg Main Conjecture.
	\begin{theorem}[Skinner--Urban, Skinner] \label[theorem]{Skinner_cor1}
		Let $f \in S_2(\Gamma_0(N))$ be a newform and $p > 2$ be a prime with
		\begin{enumerate}
			\item[\upshape{(ord)}] $v_p(N) \leq 1$ and $|a_p(f)|_p = 1$.
		\end{enumerate}
		Let $\frp \mid p$ be a finite place of $\Q(f)$ and $\Q_\infty$ be the cyclotomic $\Z_p$-extension of $\Q$ with Galois group $\Gamma \defeq \Gal(\Q_\infty|\Q)$. Let $\Lambda \defeq \Z[f]_\frp\llbracket\Gamma\rrbracket$ be the Iwasawa algebra.
		
		Assume that
		\begin{enumerate}
			\item[\upshape(irr)] $\rho_{f,\frp}$ is irreducible and
			\item[\upshape($\spadesuit$)] that there exists a prime $q \neq p$ with $v_q(N) = 1$ such that $\rho_{f,\frp}$ is ramified at $q$.
		\end{enumerate}
		Then one has an equality
		\[
		\Char_{\Q_\infty,\Z[f]_\frp}(f) = (\sL_\frp(f,T))
		\]
		of ideals in $\Lambda$. Here, $\Char_{\Q_\infty,\Z[f]_\frp}(f)$ is the characteristic ideal of the $\frp$-adic Selmer group of $f$ over $\Q_\infty$ and $\sL_\frp(f,T)$ is the $\frp$-adic $L$-function of $f$; both are defined in~\cite[§\,1.1]{SkinnerUrban2014}.
	\end{theorem}
	
	\begin{proof}
		See~\cite[Theorem~1]{SkinnerUrban2014} in the case $v_p(N) = 0$ and~\cite[Theorem~A]{Skinner2016} in the case $v_p(N) = 1$ (by reduction to~\cite{SkinnerUrban2014} using Hida theory). (Note the footnote~1 in~\cite[p.~172]{Skinner2016}, which says one can weaken the condition that there exists an $\Z[f]_\frp$-basis of $T_f$ with respect to which the image of $\rho_{f,\frp^\infty}$ contains $\SL_2(\Z_p)$ to condition ($\spadesuit$) for the Iwasawa Main Conjecture to hold integrally in~\cite[Theorem~1]{SkinnerUrban2014}.)
	\end{proof}
	
	\begin{example}
		We expect that
		\cref{Skinner_cor1} combined with the computation of $\sL_\frp(f,T)$ can be used to
		verify strong BSD for the Jacobian $J$ of level~$145$ of the curve
		\[
		C \colon y^2 = 20 x^5 - 19 x^4 + 118 x^3 - 169 x^2 + 50 x + 25
		\]
		with $\End_\Q(J) \cong \Z[\sqrt{2}]$. Our algorithm gives us
		that $\#\Sha(J/\Q)_\an = 1$ and $\Sha(J/\Q)[\frp] = 0$ except maybe for the two primes~$\frp$
		lying above $7$. For them, one can compute the $\frp$-adic $L$-function using Sage.
		(This curve is not contained in our dataset and is only an illustration
		how the results in this section can be used in the verification of the $p$-part of BSD
		when descent methods are impractical.
		It came up in forthcoming work of Kaya--Masdeu--Müller--van der Put that
		computes Schneider regulators of Mumford curves. A verification the $p$-adic BSD conjecture
		for this example up to high precision is likely to be included in their article.)
	\end{example}

	\subsection{Results on the $p$-part of BSD from Main Conjectures in Iwasawa Theory}
	\label{ssec: Iwasawa theory}
	
	We state a result that follows from Iwasawa Theory and shows
	that the $p$-part of~BSD holds in the $L$-rank~$0$ case under fairly mild assumptions
	when $p$ is inert in~$\Z[f]$. Note that an explicit descent computation is hard
	when $p$ is inert and~$\rho_{f,p}$ is irreducible, so this is useful to deal with
	such primes when our other methods do not show that $\Sha(A/\Q)[p^\infty]$ is trivial.
	
	\begin{theorem}[$p$-part of BSD in the $L$-rank $0$ case]
		\label[theorem]{Skinner_cor2}
		Let $A/\Q$ be a simple modular abelian variety associated to a newform~$f$ of
		level~$N$. Let $p > 2$ be a prime inert in~$\Z[f]$ with
		\begin{enumerate}
			\item[\upshape(ord)] $v_p(N) \leq 1$ and $|a_p(f)|_p = 1$.
		\end{enumerate}
		Assume that
		\begin{enumerate}
			\item[\upshape(irr)] $\rho_{f,p}$ is irreducible and
			\item[\upshape($\spadesuit$)] that there exists a prime $q \neq p$ with $v_q(N) = 1$ such that
			$\rho_{f,p}$ is ramified at~$q$.
		\end{enumerate}
		If $L(A/\Q,1) \neq 0$, then the $p$-part of BSD holds for~$A/\Q$, i.e.,
		\begin{align*}
			\biggl|\frac{L(A/\Q,1)}{\Omega_A}\biggr|_p
			&= \biggl|\frac{\Tam(A/\Q)
				\cdot \#\Sha(A/\Q)[p^\infty]}{\#A(\Q)_\tors \cdot \#A^\dual(\Q)_\tors}\biggr|_p \\
			&= \bigl|\Tam(A/\Q)\bigr|_p \cdot \bigl|\#\Sha(A/\Q)\bigr|_p .
		\end{align*}
	\end{theorem}
	
	\begin{proof}
		This is~\cite[Theorem~C]{Skinner2016} (the proof generalizes from elliptic curves
		to modular abelian varieties because~\cite[Theorem~B]{Skinner2016} is for general
		newforms and $p$ is inert in~$\Z[f]$). Note that
		$|\#A(\Q)_\tors \cdot \#A^\dual(\Q)_\tors|_p = 1$ since $\rho_{f,p}$ is irreducible.
	\end{proof}
	
	\begin{remark}\strut
		\begin{enumerate}[(i)]
			\item We have restricted ourselves to inert primes for simplicity.
			(The results hold more generally if there is exactly one prime above $p$.)
			The proofs for other primes $\frp$ would need to be adapted
			from the case of elliptic curves and one would need to define the
			algebraic factors in the strong BSD~formula as elements of~$\Z[f]_\frp$
			up to units. (It is likely that one can even define them up to squares
			of units over the Heegner field as most terms are in fact squares.)
			
			\item In the $L$-rank~$1$ case, there is~\cite[Theorem~10.3]{Zhang2014}
			for elliptic curves, which builds upon~\cite[Theorem~9.3]{Zhang2014},
			which is for general newforms. These theorems have stronger assumptions
			than~\cref{Skinner_cor2}.
			
			\item \cite{CGLS} proves the $p$-part in the $L$-rank~$1$ case for odd good
			ordinary primes~$p$ completely split in the Heegner field~$K$ with $\rho_p$ \emph{reducible} if
			$\rho_p|_{G_K}^{\text{ss}} \isom \phi \oplus \psi$ such that
			\hbox{$\phi|_{G_{K_v}},\psi|_{G_{K_v}} \neq \mathbf{1}, \chi_p$} for
			the places $v \mid p$ of~$K$.
			For example, this excludes the case $A(K)[p] \neq 0$.
			The work of the first author and Mulun Yin~\cite{KellerYin2024}
			removes the restriction on~$\phi$ and~$\psi$.
			Depending on~\cite{CGS} and~\cite{BSTW2024},
			it also treats the $L$-rank~$0$ case and
			the case of bad multiplicative reduction.
			The introduction of~\cite{KellerYin2024} contains an example
			of a modular abelian surface over~$\Q$ that illustrates the use of the
			results proved there to make
			a descent computation for a torsion prime unnecessary.
			
			\item The preprint~\cite{CCSS} proves the $p$-part of strong~BSD in the
			$L$-rank $0$ and~$1$ cases also for $p$ good \emph{non-ordinary},
			but it assumes that the level~$N$ is squarefree.
			
			\item \cref{thm:finite support,thm:finite support rank 1} give $\Sha(J/\Q)[\frp] = 0$
			under the assumptions there, but they do not prove the $p$-part of BSD.
			\cref{Skinner_cor2} does this,
			but not for non-inert primes or supersingular or bad additive primes.
			Note that by the Sato--Tate conjecture there are infinitely many~$p$
			with $a_p = 0$, and they can be treated with~\cite{CCSS} only when $N$
			is squarefree.
		\end{enumerate}
	\end{remark}

	\subsection{The $p$-adic BSD conjecture}
	
	We continue to assume that $p$ is inert in $\Z[f]$.
	The coefficient of its leading term of $\sL_p(A,T)$ at $T = 0$ is denoted by $\sL_p^*(A,0)$.
	Recall from the introduction to this section that $\alpha \in \Z[f]_\frp^\times$
	is the $p$-adic unit root of the Euler factor of~$f$ at~$p$.
	As in~\cite{BMS}, the \emph{$p$-adic multiplier} is
	\[ \epsilon_p(A/\Q) \defeq \Nm_{\Z[f]_\frp|\Z_p}(1 - \alpha^{-1})^2 \,. \]
	(Note that if $a_p = 1$ and $p \mid N$, then $\alpha = 1$
	and $\epsilon_p(A/\Q) = 0$, so conjecturally there will be an \emph{extra zero}.)
	The \emph{$p$-adic regulator} $\Reg_p(A/\Q)$ is the determinant of the $p$-adic
	height pairing on~$A(\Q)$ defined in~\cite[Definition~3.3]{BMS}.
	According to a Conjecture of Schneider, $\Reg_p(A/\Q)$ should be non-zero,
	but this is not known in general.
	Assume $p > 2$. Let $\gamma \in 1 + p\Z_p$ be a topological generator, the same used in the construction of $\sL_p(A,T)$. We take $\gamma = 1 + p$ in our computations. Let $\Reg_\gamma(A/\Q) \defeq \Reg_p(A/\Q)/\log_p(\gamma)^r$ where $r = \rk A(\Q)$.
	We then have the following $p$-adic version of the
	BSD~conjecture; see~\cite[Conjecture~1.4]{BMS}, generalizing~\cite{MTT1986}.
	
	\begin{conjecture}[$p$-adic BSD conjecture]
		Let $A/\Q$ be a principally polarized modular abelian variety and $p$ a prime of good ordinary reduction
		for~$A/\Q$. Then
		\[ \rk A(\Q) = \ord_{T=0}\sL_p(A,T) \]
		and
		\[
		\sL_p^*(A,0) = \epsilon_p(A/\Q) \cdot \frac{\#\Sha(A/\Q) \cdot \Tam(A/\Q) \cdot \Reg_\gamma(A/\Q)}%
		{(\#A(\Q)_\tors)^2} \,.
		\]
	\end{conjecture}
	
	We are interested in the size of the $p$-part of~$\Sha(A/\Q)$, so it
	is sufficient to compare the $p$-adic valuations of both sides.
	
	The following result due to Schneider shows that the conjecture holds in some cases
	at least up to a $p$-adic unit, but with~$\sL_p(A/\Q,T)$ replaced by the
	\emph{Iwasawa $L$-function}~$\sL_p^{(1)}(A/\Q,T)$ defined in~\cite[§\,2]{Schneider1985}.
	Note that the latter is defined without using modularity, but in a more algebraic way.
	
	\begin{theorem} \label[theorem]{p-adic BSD}
		Let $A/\Q$ be a simple principally polarized modular abelian variety with associated newform $f$. If
		\begin{enumerate}[\upshape(i)]
			\item $p$ is a prime of good ordinary reduction,
			\item such that the $p$-adic regulator $\Reg_\gamma(A/\Q)$ is non-zero, and
			\item $\Sha(A/\Q)[p^\infty]$ is finite,
		\end{enumerate}
		then the $p$-adic BSD conjecture holds (up to a $p$-adic unit) for $A/\Q$ and $\sL_p^{(1)}(A/\Q,T)$:
		
		The Iwasawa $L$-function $\sL_p^{(1)}(A/\Q,T)$ vanishes to order~$\rk A(\Q)$
		at \hbox{$T = 0$}, and its leading term has $p$-valuation equal to that of
		\[
		\epsilon_p(A/\Q) \cdot \frac{\#\Sha(A/\Q)[p^\infty] \cdot \Tam(A/\Q) \cdot \Reg_\gamma(A/\Q)}%
		{(\#A(\Q)_\tors)^2} \,.
		\]
	\end{theorem}
	
	\begin{proof}
		See~\cite[Theorem~$2'$]{Schneider1985}.
	\end{proof}
	
	We need to compare the two $L$-functions $\sL_p(A/\Q,T)$ and $\sL_p^{(1)}(A/\Q,T)$.
	
	\begin{theorem} \label[theorem]{Kato divisibility}
		Let $p > 2$ be a prime of good ordinary reduction for~$A/\Q$. Let $\frp$ be a prime ideal of $\Z[f]$ lying above $p$.
		Assume that
		the image of $\Gal(\Qbar|\Q(\mu_{p^\infty})) \to \Aut_{\Z[f]_\frp}(T_\frp A)$
		contains~$\SL_2(\Z_p)$
		(see Section~\ref{sec:the-image-of-the-pfr-adic-galois-representation}).
		Then
		\[
		\sL_p^{(1)}(A/\Q,T) \mid \sL_p(A/\Q,T) \in \Z[f]_p\llbracket T \rrbracket \,.
		\]
	\end{theorem}
	
	\begin{proof}
		See~\cite[Theorem~17.4\,(3)]{Kato2004}.
	\end{proof}
	
	\begin{theorem} \label[theorem]{ord sLpf}
		Let $p > 2$ such that $a_p(f) \in \Z[f]_p^\times$. Then
		\[
		\ord_{T=0} \sL_p(f,T) \geq \cork_{\Z[f]_p} \Sel_{p^\infty}(A_f/\Q) \ge \rk_{\Z[f]} A_f(\Q) \,.
		\]
		Here, the corank $\cork_{\Z[f]_p} M$ of a discrete torsion $\Z[f]_p$-module is the $\Z[f]_p$-rank of its $\Z[f]_p$-Pontrjagin dual.
	\end{theorem}
	
	\begin{proof}
		See~\cite[Theorem~18.4]{Kato2004}.
	\end{proof}
	
	\begin{corollary} \label[corollary]{bound of Sha p using padic L}
		Let $A/\Q$ be a simple principally polarized modular abelian variety with associated newform~$f$. If
		\begin{enumerate}[\upshape(i)]
			\item $p > 2$ is a prime of good ordinary reduction
			\item such that the $p$-adic regulator $\Reg_\gamma(A/\Q)$ is non-zero,
			\item $\Sha(A/\Q)[p^\infty]$ is finite,
			\item the image of $\Gal(\Qbar|\Q(\mu_{p^\infty})) \to \Aut_{\Z[f]_p}(T_p A)$
			contains~$\SL_2(\Z_p)$, and
			\item $\ord_{T=0} \sL_p(f^\sigma,T) \le \rk_{\cO} A(\Q)$ for all~$\sigma$, \label{item:ordLp le rk}
		\end{enumerate}
		then equality holds in~\eqref{item:ordLp le rk} and
		\[
		v_p(\#\Sha(A/\Q)[p^\infty])
		\le v_p\biggl(\frac{\prod_\sigma \sL_p^*(f^\sigma, 0) \cdot (\#A(\Q)_\tors)^2}%
		{\epsilon_p(A/\Q) \cdot \Reg_\gamma(A/\Q) \cdot \Tam(A/\Q)}\biggr) \,.
		\]
	\end{corollary}
	
	\begin{proof}
		Combine \cref{p-adic BSD,Kato divisibility,ord sLpf}.
	\end{proof}
	
	Note that~\cite[Theorems~A and~B]{Skinner2016} (with the case of good ordinary
	reduction coming from~\cite{SkinnerUrban2014}) establishes equality up to units
	in~$\Z[f]_p\llbracket T\rrbracket$ in the ordinary case under some conditions like
	$\rho_p$ being surjective.
	
	So if we can compute the $p$-adic valuations of~$\sL_p^*(A,0)$ and of~$\Reg_\gamma(A/\Q)$,
	this result allows us to bound the order of the $p$-part of~$\Sha(A/\Q)$ from above.
	Note that we know by the results of Kolyvagin--Logachëv and their extensions that
	$\Sha(A/\Q)[p^\infty]$ is finite in the cases of interest,
	so the main assumption is that $p$ is a prime of good ordinary reduction.
	When the $L$-rank is zero, then the $p$-adic regulator is~$1$, so it remains to
	compute~$\sL_p(A,T)$. We explain in the next subsection how to do this. When
	the $L$-rank is~$1$, we also need to compute the $p$-adic height pairing.
	In the case that $p$ is a prime of good ordinary reduction for the Jacobian of
	a genus~$2$ curve, this is accomplished in~\cite[§\,3.4]{BMS};
	see also~\cite{GajovicMueller}.

	\subsection{Computing approximations to $p$-adic $L$-functions}
	
	We use Greenberg's improvement~\cite{Greenberg2007} of the Pollack--Stevens
	algorithm~\cite{PollackStevens2011} to compute the $p$-adic $L$-function of
	a newform~$f$ with $a_p(f)$ a $p$-adic unit.
	Our Magma implementation is based on that of Darmon--Pollack~\cite{DarmonPollack2006}
	with their Magma code available at~\cite{DarmonPollack2006Code}.
	We modified the code so it also works with newforms with arbitrary coefficient rings,
	accepts newforms as input and outputs the $p$-adic $L$-function as a $p$-adic power
	series to any specified precision for the uniformizer~$\pi$ of~$\cO_\frp$ and~$T$.
	For performance reasons, we specialized to weight $k = 2$.
	
	\subsubsection{Computing the $p$-stabilization of $\phi_f$ if $p \nmid N$.}
	\label{sec:p-stabilization}
	The construction of the $p$-adic $L$-function via the overconvergent modular symbol
	algorithm of~\cite{Greenberg2007, PollackStevens2011}
	needs $v_p(N) = 1$. In the case that $p \nmid N$ is a good ordinary prime, one has to
	\emph{$p$-stabilize} $f$ to a form of level~$Np$ with the same $T_\ell$-eigenvalues
	as $f$ and $U_p$-eigenvalue the unit root of $T^2 - a_p(f)T + p$.
	Note that the $p$-adic $L$-function is defined with respect to that ($p$-stabilized if $p \nmid N$) lift;
	see~\cite[§§\,4.1, 4.2, 4.4]{BelabasPerrinRiou2021}.
	
	Recall that $a_p(f) \in \Z[f]_\frp^\times$. If $p \nmid N$, let $\alpha \in \Z[f]_\frp^\times$ be
	the unit root of $T^2 - a_p(f)T + p$. Let $\beta$ be the other root.
	We then $p$-stabilize in the sense that we
	replace $f$ by its \emph{$p$-stabilization} $f_\alpha(z) = f(z) - \beta f(pz)$
	of level~$Np$ with the same Hecke eigenvalues away from~$p$.
	Otherwise, i.e., if $v_p(N) = 1$, we can directly use the
	algorithm of~\cite{Greenberg2007, PollackStevens2011} to
	compute a lift.
	(The reason for the $p$-stabilization is that the distribution property
	for the modular symbol follows if it is an eigenvector under the $U_p$-operator,
	but not under the $T_p$-operator.)
	
	\subsubsection{Modular symbols.} Let $\Delta_0$ be the left $\Z[\GL_2(\Q)]$-module $\Div^0(\bP^1(\Q))$, where $\GL_2(\Q)$
	acts via fractional linear transformations on $\bP^1(\Q)$.
	We write an element $[s] - [r]$ of $\Delta_0$ with $r,s \in \bP^1(\Q)$ as $\{r \to s\}$.
	
	Consider the monoid
	\[
	\Sigma_0(p) \defeq \biggl\{\begin{pmatrix} a & b \\ c & d \end{pmatrix} \in \Mat_2(\Z)
	: (a,p) = 1, \; p \mid c, \; ad-bc \neq 0\biggr\} \,.
	\]
	Note that because $p \mid N$, $\Gamma_0(N) \subseteq \Sigma_0(p)$. Let $V$ be a right $\Sigma_0(p)$-module with the action of~$S \in \Sigma_0(p)$ on $v \in V$ denoted
	by~$v|S$.
	Then the Hecke operators
	\[
	T_\ell \defeq \begin{pmatrix} \ell & 0 \\ 0 & 1 \end{pmatrix}
	+ \sum_{a = 0}^{\ell-1} \begin{pmatrix} 1 & a \\ 0 & \ell \end{pmatrix} \in \Z[\Sigma_0(p)] \,,
	\quad \text{$\ell \neq p$ prime,}
	\]
	and
	\[
	U_p \defeq \sum_{a = 0}^{p-1} \begin{pmatrix} 1 & a \\ 0 & p \end{pmatrix} \in \Z[\Sigma_0(p)]
	\]
	act on~$V$ on the right, making it into a Hecke module.
	(The matrices constituting $T_\ell, U_\ell$ are in $\Sigma_0(p)$, but not in $\Gamma_0(N)$.)
	In Magma, the elements of~$\Delta_0$ are called the \emph{modular symbols}
	of~$\Gamma_0(N)$, but we refer to the elements of the module
	\[ \Symb_{\Gamma_0(N)}(V) \defeq \Hom_{\Z[\Gamma_0(N)]}(\Delta_0, V) \]
	as the \emph{$V$-valued modular symbols}. The abelian group $\Hom_{\Z}(\Delta_0, V)$
	is a \emph{right} $\Z[\Gamma_0(N)]$-module and we write the action similarly
	as~$\varphi|S$. Spelled out explicitly, the $\Gamma_0(N)$-equivariance means that
	for $\phi \in \Symb_{\Gamma_0(N)}(V)$ one has $\phi|S = \phi$, i.e.,
	$\phi(S\{r \to s\})|S = \phi\{r \to s\}$ for $\{r \to s\} \in \Delta_0$ and $S \in \Gamma_0(N)$.
	Then $\Symb_{\Gamma_0(N)}(V)$ is a right Hecke module via the left action
	on~$\Delta_0$ and the right action on~$V$.
	
	\subsubsection{Computing the canonical periods $(\Omega_{f^\sigma}^\pm)_\sigma$ attached to $f$.}
	We have to find canonical periods of~$\{f^\sigma\}$ as
	defined in~\cite[§\,2]{BMS}, unique up to~$\Z[f]^\times$. We compute an approximation of
	the periods~$(\Omega_{f^\sigma}^\pm)_\sigma$ with $\sigma\colon \Q(f) \inj \R$ as follows:
	We compute the period integrals
	\[
	p_{f^\sigma}^\pm(r) \defeq \bpi i \Big(\int_r^{i\infty} f^\sigma(z) \,dz \pm \int_{-r}^{i\infty} f^\sigma(z) \,dz \Big)
	\]
	(i.e., over the path $\{r \to \infty\} \in \H_1(X_0(N)(\C),\{\text{cusps}\};\Z)$)
	for a finite set of $r \in \Q$ such that the $\{r \to \infty\}$ generate the
	$\Delta_0$ as a $\Gamma_0(N)$-module
	to a high enough precision.
	By~\cite[§\,2.1, Theorem~2.2 and~2.4]{BMS}, there are \emph{canonical periods} $(\Omega_{f^\sigma}^\pm)_\sigma$
	attached to $f$ such that $p_{f^\sigma}^\pm/\Omega_{f^\sigma}^\pm \in \Q(f)$
	and \enquote{being compatible with twists} (for a precise formulation see~\cite[Theorem~2.2]{BMS}).
	
	By~\cite[text after Remark~2.6]{BMS}, there is a $b \in \Q(f)^\times$ such that
	$(\Omega_{f^\sigma}^\pm)_\sigma = (\sigma(b) p_{f^\sigma}^\pm(r))_\sigma$.
	We approximate a representative vector of an equivalence class $(\Omega_{f^\sigma}^\pm)_\sigma$.
	To compute the period integrals $p_{f^\sigma}^\pm(r)$, we combine~\eqref{Eq:LJ/Om}
	(and the equation before it) and~\cref{relation period over omegafsigma and Neron differential}
	and get
	\begin{equation} \label{Eq:period rel Lp}
		\Omega_{J/\Q} = \frac{c_fc_\pi}{\sqrt{\disc\Z[f]}} \cdot \frac{\#\coker\pi_\R}{\#\ker\pi_\R} \cdot \Big|\prod_\sigma \Omega_{f^\sigma}^+\Big| \,.
	\end{equation}
	We compute $\Omega_{J/\Q}$ as in~\cref{computation of Omega_J}, $c_fc_\pi$ as in~\cref{computation of c_f c_pi}, and $\#\coker\pi_\R$, $\#\ker\pi_\R$ as described in Section~\ref{sec:period matrices}. In fact, since the latter two constants are powers of $2$, we do not need to compute them if $p \neq 2$. Since we work with a $\Z$-basis $(g_i)$ of $S_2(f, \Z)$, there is no factor $\sqrt{\disc\Z[f]}$ when we replace $\prod_\sigma \Omega_{f^\sigma}^+$ by the corresponding product for the $(g_i)$.
	
	Assuming that Equation~\eqref{Eq:period rel Lp} holds,
	the canonical periods are unique up to $\sigma(b)$ for some $b \in \Z[f]_\frp^\times$.
	So the $\frp$-adic valuation of the leading term of the $p$-adic $L$-function
	is uniquely determined.
	
	\subsubsection{The modular symbol $\phi_f$.}
	
	The \emph{modular symbol associated with $f$} is $\phi_f \in \Hom_{\Z[\Gamma_0(N)]}(\Delta_0, \Q(f))$,
	where $\Q(f)$ is a trivial $\Gamma_0(N)$-module.
	It is defined as $\phi_f\{r \to \infty\} = p_f^+(r)/\Omega_f^+ + p_f^-(r)/\Omega_f^- \in \Q(f)$
	and by extension $\phi_f\{r \to s\} = \phi_f\{r \to \infty\} - \phi_f\{s \to \infty\}$
	to all paths between cusps $r,s \in \bP^1(\Q)$.
	It is enough to know $\phi_f$ on a finite set of generators
	of $\Delta_0$ as a $\Z[\Gamma_0(N)]$-module,
	which can be obtained via \emph{Manin symbols}~\cite[§§\,2.2, 2.3]{Cremona1997}.
	
	\subsubsection{Overconvergent modular symbols and $p$-adic $L$-functions.}
	Let $\frp \mid p$ be a prime ideal of~$\Z[f]$ inducing an embedding of~$\Q(f)$ into
	the completion~$\Q(f)_\frp$. To approximate the distribution~$L_\frp(f)$
	(see~\cite[§\,6.2]{PollackStevens2011}) and its associated power series~$\sL_\frp(f,T)$,
	we implemented Greenberg's improvement~\cite{Greenberg2007} of Pollack--Stevens's computation of an
	overconvergent modular eigenlift of the modular symbol
	$\phi_f \in \Symb_{\Gamma_0(N)}(\Q(f))$ (see~\cite[§\,6.3]{PollackStevens2011}).
	
	For the following, see~\cite[§\,3.1]{PollackStevens2011} and~\cite[§\,1]{Greenberg2007}.
	For $r \in |\C_p^\times|_p$, define $\sA[r](\Q(f)_\frp)$ to be the
	$\Q(f)_\frp$-Banach space of $\Q(f)_\frp$-affinoid functions on
	\[ B[\Z_p,r] \defeq \{z \in \C_p : \exists a \in \Z_p \text{ with } |a - z|_p \leq r\} \]
	endowed with the supremum norm. Let
	\[ \sA^\dagger(\Z_p) \defeq \varinjlim_{s > 1} \sA[s](\Q(f)_\frp) \,, \]
	endowed with the colimit topology, denote the algebra of $\Q(f)_\frp$-over\-con\-ver\-gent
	functions on~$B[\Z_p,1]$. The overconvergent distributions~$\sD^\dagger(\Q(f)_\frp)$
	are defined as its continuous $\Q(f)_\frp$-linear dual endowed with the strong topology.
	
	\subsubsection{The strategy to compute the distribution $L_\frp(f)$.}
	We abbreviate $\Gamma_0(N)$ by~$\Gamma$ in the following.
	The method of constructing the $\frp$-adic $L$-function~$L_\frp(f)$ of~$f$ as
	a distribution is summarized in the following diagram. Here $\Symb_{\Gamma}(\Q(f))[f]$
	denotes the subspace of~$\Symb_{\Gamma}(\Q(f))$ on which the Hecke operators~$T_n$
	with $p \nmid n$ act as multiplication by~$a_n(f)$ and the superscript \enquote{$U_p = \alpha$}
	means the subspace where the Hecke operator~$U_p$ acts as multiplication by the unit
	root eigenvalue~$\alpha$ of the characteristic polynomial at~$p$.
	
	\[
	\begin{tikzcd}
		\Symb_{\Gamma}(\sD^\dagger(\Q(f)_{\frp}))[f]^{U_p = \alpha}
		\ar[rrr, "\Phi_f \mapsto \Phi_f\{0 \to \infty\} = L_\frp(f)"]
		\ar[d, dashed, "\rho^*","\iso" left] &&&
		\sD^\dagger(\Q(f)_{\frp})
		\ar[d,"L_\frp(f) \mapsto \int_{\Z_p^\times}\mathbf{1} L_\frp(f)"] \\
		\Symb_{\Gamma}(\Q(f))[f]^{U_p = \alpha}
		\ar[rrr,"\phi_f \mapsto \phi_f\{0 \to \infty\}"] &&&
		\Q(f)_{\frp}
	\end{tikzcd}
	\]
	
	To compute the distribution~$L_\frp(f)$, we start with~$\phi_f$ in the lower
	left corner and lift it via the algorithm described below to a distribution valued
	modular symbol~$\Phi_f$; this is the left vertical isomorphism. Then $L_\frp(f)$
	is the evaluation~$\Phi_f\{0 \to \infty\}$, i.e., the image under the top
	horizontal morphism.
	
	\subsubsection{Computing the modular symbol $\phi_f^\pm$ attached to $f$.}
	Knowing the canonical periods, we can compute~$\phi_f^\pm\{r \to s\}$ with $r, s \in \bP^1(\Q)$ as described in~\cite{Wuthrich2018}
	by computing the corresponding period integrals for all $f^\sigma$ and dividing by~$\Omega_{f^\sigma}^\pm$
	and recognizing the elements in~$\Q(f)$ using that
	\[
	\prod_\sigma (X - \phi_{f^\sigma}^\pm\{r \to s\}) \in \Q[X]
	\]
	has rational coefficients
	with the denominator of the coefficient of $X^{[\Q(f):\Q]-n}$
	bounded by $(4^g \cdot c_fc_\pi \cdot \#J(\Q)_\tors)^n$;
	see Section~\ref{ssec: LJ/Om} and~\cite[Prop.~1]{Wuthrich2018}.
	
	One can compute the modular symbol associated to~$f$ up to a factor in $\Q(f)$ with Magma.
	(Note that one has to take a suitable $\Q(f)$-linear combination because Magma takes
	a basis of $S_2(f,\Z)$.)
	Hence, alternatively, one can compute this factor by comparing with our above computation
	for some $\phi_f\{r \to s\} \neq 0$ and scale.
	
	\subsubsection{Determining the required $p$-adic precision for the desired $T$-adic precision of $\sL_p(f,T)$.} \label{sssec:p-adic precision for Lp}
	Note that a distribution $\mu \in \sD^\dagger(\Q(f)_\frp)$ is uniquely
	determined by its moments~$\mu(x^j)$, $j \geq 0$. We represent a
	distribution in an approximation module $\bA^M\sD^\dagger(\Z[f]_\frp)$
	by storing its $j$-th moment up to precision $M+1 - j$, i.e., as an element
	of~$\Z[f]_\frp/\frp^{M+1-j}$. To be able to compute lifts, we store
	the $0$-th moment up to the final precision.
	
	To ensure the moments we consider are integral, we do the following:
	According to~\cite[Corollary~7.6, §\,8.3\,3]{PollackStevens2011},
	to obtain a precision of~$n$ $p$-adic digits, one needs to compute with
	a precision~$M$ satisfying
	\[ M - \lceil\log (M+2)/\log p\rceil - 1 \ge n \,. \]
	We take $m$ minimal with $p^m > M+1$ and scale~$\phi_f$ by~$p^{m+1}$.
	Then we perform $M$~steps of the algorithm to obtain the approximation
	to $p^{m+1} \Phi_f$ in~$\bA^M\sD(\Z[f]_\frp)$. Finally, we divide by~$p^{m+1}$.
	
	\subsubsection{Computing the lift $\Phi$ of $\phi$.}
	This is the key step in the computation of the $p$-adic $L$-function
	and a simplification of~\cite{PollackStevens2011} due to Greenberg~\cite{Greenberg2007},
	with the additional non-critical slope assumption, which is satisfied if there is
	a unit root of the characteristic polynomial of $p$-Frobenius.
	
	Let $\phi \in \Symb_{\Gamma}(\Q(f))[f]^{U_p = \alpha}$ be a Hecke eigensymbol.
	The action of $\gamma = \left(\begin{smallmatrix} a & b \\ c & d \end{smallmatrix}\right)$
	with $p \mid c$ and $p \nmid a$ on a distribution~$\mu$ is given by
	$(\mu | \gamma)(f) = \mu(\gamma \cdot f)$, where $\gamma$ acts on the function~$f$
	as $(\gamma \cdot f)(z) = (a + cz)^2\cdot f(\frac{b+dz}{a+cz})$.
	
	We start with~$(\Phi)_0 \defeq \phi$. To lift~$(\Phi)_M$ from
	$\Symb_{\Gamma_0(N)}(\bA^M\sD^\dagger(\Z[f]_\frp))$ to precision~$M+1$,
	we first lift all values on the finitely many generators of~$\Delta_0$
	(computed using Manin symbols) arbitrarily to precision~$M+1$.
	The resulting function will usually not be additive, $\Gamma_0(N)$-equivariant
	or an $U_p$-eigensymbol anymore, but it will be after we apply the operator~$U_p$.
	
	We stop after we have reached the precision from~\ref{sssec:p-adic precision for Lp}.
	
	\subsubsection{Computing an approximation to the $p$-adic $L$-function from the distribution $L_\frp(f)$.} To go from an approximation of~$L_\frp(f)$ to the power series~$\sL_\frp(f,T)$,
	we use the formulas in~\cite[§\,9]{PollackStevens2011}. The computation depends on the choice of a topological generator $\gamma \in 1 + p\Z_p$.
	We take $\gamma = 1 + p$ in our computations.
	
	\subsection{Examples}
	
	The example below represents the only case where our other methods
	are not sufficient to compute $\#\Sha(A/\Q)[p^\infty]$, so that we need
	to use $p$-adic $L$-functions. (Note that a $\frp$-descent would amount
	to a full $3$-descent in this case, which is not really feasible with
	current methods.)
	
	\begin{example} \label[example]{Sha3 for N 188}
		For the Jacobian~$J$ of the
		\href{https://www.lmfdb.org/Genus2Curve/Q/35344/a/565504/1}{curve $X$ of level $188 = 2^2 \cdot 47$}
		in our data set, the Tamagawa product is $9$ (we have $c_2 = 9$
		and $c_{47} = 1$); the prime~$3$ is inert in the endomorphism ring $\Z[(1+\sqrt{5})/2]$.
		To show that $\Sha(J/\Q)[3^\infty] = 0$ in this case, we therefore
		compute the $3$-adic $L$-function. Since $3 \nmid N$,
		we have to $3$-stabilize, see~\ref{sec:p-stabilization}. Note that
		$a_3(f) = -\frac{3+\sqrt{5}}{2}$ is a $3$-adic unit,
		hence the reduction at~$3$ is good ordinary. Then $\epsilon_3(J/\Q)$
		is a unit as well (since $\alpha \equiv a_3(f) \not\equiv 1 \bmod 3$).
		
		We verify using~\cref{p-adic image,mod-9 image} that
		$\SL_2(\Z_3)$ is contained in the image of~$\rho_{f,3^\infty}$.
		
		The $L$-rank is $1$ and $J(\Q) \isom \Z^2$.
		A computation of $3$-adic heights shows that $v_3(\Reg_3(J/\Q)) = 0$. Therefore, one has
		\[
		v_3(\Reg_\gamma(J/\Q)) = 0 - r\cdot v_3(\log_3(1+3)) = -2 \,.
		\]
		We thank Steffen Müller for computing it for us using the code
		described in~\cite{GajovicMueller}.
		Using the algorithm sketched above, we find that
		\begin{align*}
			\sL_3(J,T) &= \sL_3(f,T) \cdot \sL_3(f^\tau,T) \\
			&= (O(3^3) + uT + O(T^2)) \cdot (O(3^3) + u'T + O(T^2))
		\end{align*}
		with $3$-adic units~$u, u'$, where  $\tau$ is the nontrivial automorphism of~$\Z[f]$.
		(The computation took 30 minutes and 214\,MiB RAM on a AMD Ryzen 7 PRO 6850U.)
		Since $\rk J(\Q) = 2$ and $\Sha(J/\Q)[3^\infty]$ is finite, $\cork_{\Z[f]_3} \Sel_{3^\infty}(J/\Q) = 1$, so~\cref{bound of Sha p using padic L} shows that the vanishing order of $\sL_3(J,T)$ at $T = 0$ must be exactly~$2$ and that
		\begin{align*}
			v_3\bigl(\#\Sha(J/\Q)[3^\infty]\bigr) &\le v_3\Biggl(\frac{\prod_\sigma 1 \cdot 1 \cdot 1}%
			{1 \cdot 3^{-2} \cdot 9}\Biggr) = 0 \,,
		\end{align*}
		hence $\Sha(J/\Q)[3^\infty] = 0$.
	\end{example}
	
	
	\section{Examples} \label{sec:examples}
	
	\subsection{Jacobians of genus $2$ Atkin--Lehner quotients}
	
	Tables displaying the results for the Hasegawa curves can be found in~\cite{KellerStoll2022}.

	\subsection{All genus $2$ curves with absolutely simple modular Jacobian from the LMFDB}
	
	We compute the analytic order of $\Sha$ using the results from Sections~\ref{sec:computation-of-the-heegner-points-and-index} and~\ref{sec:Sha_an}. It turns out that all of them are $1$, $2$, or $4$. We also discover some twists $J^K$ which have analytic order of $\Sha$ divisible by $3^2$, $5^2$, or~$7^2$.
	
	It turns out that combining the information about the images of the residual Galois representations from Section~\ref{sec:computation-of-the-residual-galois-representations} with the Heegner indices from Section~\ref{sec:computation-of-the-heegner-points-and-index} (for a few examples, we have to use two Heegner fields) and the Euler system from Section~\ref{sec:finite support} prove that $\Sha(J/\Q)[\frp] = 0$ except in the following cases:
	\begin{enumerate}[$\bullet$]
		\item $p(\frp) = 2$: these are dealt with in Section~\ref{sec:dealing-with-p--2}.
		\item Odd primes $\frp$ with $\rho_\frp$ reducible: these are dealt with in Section~\ref{sec:odd-primes}.
		\item One example with $N = 188$ for which $3 \mid \Tam(J/\Q)$ and $\rho_3$ is irreducible: see~\cref{Sha3 for N 188}.
	\end{enumerate}
	
	This completes the verification of strong BSD for all the $97$ absolutely simple modular Jacobians in the LMFDB.
	
	The Heegner discriminants used in the computation were typically $\leq 51$ in absolute value.
	The largest Heegner discriminant used was $-131$ for the example of level $165$. 
	We needed two Heegner discriminants in the examples of level~$523$ and~$621$. 
	The Heegner point used in one of the examples of level~$275$ has unusually large height
	($\approx 83.863$).
	Computing the Heegner point and the Mordell--Weil group of the 
	Jacobian over the Heegner field are the most time-consuming computations in our verification.
	
	\subsection{Jacobians of the four remaining Wang curves}
	
	The Jacobians $J$ of the four Wang curves of levels 65A, 117B, 125B, and~175
	(in the notation of~\cite{FLSSSW}) also have analytic order of Sha in~$\{1, 2, 4\}$.
	The remaining descent computations that are necessary to finish the proof
	that $\#\Sha(J/\Q) = \#\Sha(J/\Q)_\an$ in these cases are sketched in~\cref{ex desc 7}.
	
	
	\appendix
	\section{An example with $7$-torsion in $\Sha$} \label{sec:7-torsion in Sha}
	
	In this appendix, which heavily relies on contributions by
	Sam Frengley, we verify strong BSD for a genus~$2$ Jacobian~$J$ such that
	$\#\Sha(J/\Q) = 7^2$. This contrasts with the examples in our database,
	where we always have $\#\Sha(J/\Q) \mid 4$.
	
	\subsection{Visibility}
	
	We first briefly recall some generalities about visibility of elements of Tate--Shafarevich groups of abelian varieties, following e.g., \cite{CM_VEITSTG,AS_VOSTGOAV,AgasheStein2005,F_VEOO7ITTSGOAEC}.
	
	Let $A_1/K$ and $A_2/K$ be abelian varieties defined over a number field $K$. Suppose that there exist finite $\Gal(\ol{K}|K)$-submodules $\Delta_1 \subset A_1(\ol{K})$ and $\Delta_2 \subset A_1(\ol{K})$
	equipped with an isomorphism $\phi \colon \Delta_1 \isoto \Delta_2$ of $\Gal(\ol{K}|K)$-modules. Write $A_1' = A_1/\Delta_1$ and $A_2' = A_2/\Delta_2$ and let $\psi_1 \colon A_1 \to A_1'$ and $\psi_2 \colon A_2 \to A_2'$ be the quotient morphisms. The isomorphism $\phi$ induces an isomorphism on Galois cohomology $\H^1(K, \Delta_1) \isoto \H^1(K, \Delta_2)$. We might hope that if enough \enquote{local coincidences} occur at the bad primes, then the $\psi_1$ and $\psi_2$-Selmer groups $\Sel(\psi_1)$ and $\Sel(\psi_2)$ may be isomorphic. In this case if the order of the group $A_1(K) / \psi_1 A_1'(K)$ is smaller than the order of $A_2(K) / \psi_2 A_2'(K)$, the latter will contribute to a discrepancy between the order of $A_1(K) / \psi_1 A_1'(K)$ and $\Sel(\psi_1)$, thereby \enquote{explaining} some nontrivial elements of $\Sha(A_1/K)[\psi_1]$.
	
	To make this idea precise let $\Delta \subset (A_1 \times A_2)(\ol{K})$ be the graph of the isomorphism $\phi$ and consider the abelian variety $B = (A_1 \times A_2)/\Delta$. We have a pair of morphisms
	\begin{equation*}
		\iota \colon A_1 \to A_1 \times A_2 \to B.
	\end{equation*}
	
	\begin{definition}
		We say that an element or subgroup of $\H^1(K, A_1)$ is \emph{visible in $B$} if it is contained in the kernel of the homomorphism $\iota_* \colon \H^1(K, A_1) \to \H^1(K, B)$. We write $\operatorname{Vis}_B \H^1(K, A_1)$ (respectively $\operatorname{Vis}_B \Sha(A_1/K)$) for the subgroup of $\H^1(K, A_1)$ (respectively $\Sha(A_1/K)$) consisting of elements visible in $B$.
	\end{definition}
	
	We will use the following theorem (which is proved in~\cite[Theorem~2.2]{F_VEOO7ITTSGOAEC}; see also~\cite{AS_VOSTGOAV} and~\cite[Appendix]{AgasheStein2005}), which we state in the case when $K = \Q$.
	\begin{theorem}
		\label{thm:when-visible-sha}
		If $A_1(\Q) / \phi_1 A_1'(\Q) = 0$, then the subgroup $\operatorname{Vis}_B \H^1(\Q,A)$ is isomorphic to $A_2(\Q) / \phi_2 A_2'(\Q)$. Moreover if $\#\Delta$ is odd, and
		\begin{enumerate}[\upshape(i)]
			\item all Tamagawa numbers of $A_1/\Q$ and $A_2/\Q$ are coprime to $\#\Delta$, and
			\item the abelian variety $B$ has good reduction at all primes dividing $\#\Delta$,
		\end{enumerate}
		then $\operatorname{Vis}_B \Sha(A_1/\Q) \isom A_2(\Q)/\phi_2 A_2'(\Q)$.
	\end{theorem}
	
	\subsection{The example}
	
	Let $C/\Q$ be the genus~$2$ curve given by the Weierstrass equation
	\[ C \colon y^2 = -10(x^6 - 10 x^5 + 32 x^4 - 40 x^3 + 38 x^2 - 20 x + 4) \,. \]
	Its Jacobian~$J$ is of $\GL_2$-type; the level is $N = 3200$.
	This genus $2$ curve is obtained, up to quadratic twist, by specializing the family
	of genus~$2$ curves with $\sqrt{2}$-multiplication given by Bending~\cite[Theorem~4.1]{B_COG2W2M}
	at $(A,P,Q) = (-10,1,-5)$.
	This example was found by computing (for many such specializations) twists $X_{J[\frp]}^\pm(7)$
	of the modular curve~$X(7)$ whose $K$-points parameterize elliptic curves~$E$ equipped with an
	isomorphism of $\Gal(\ol{K}|K)$-modules $J[\mathfrak{p}] \cong E[7]$ (see e.g.,
	\cite[Section~4.4]{PSS_TOX7APSTX2Y3Z7} for the construction of these twists).
	The details of these calculations and further examples
	appear in the PhD~thesis of Sam Frengley~\cite{FrengleyPhD} and in~\cite{Frengley2024}.
	
	We begin by showing that the support of~$\Sha(J/\Q)$ is contained in~$\{7\}$.
	
	\begin{proposition}
		\label[proposition]{properties of Sha 49 curve}
		The Jacobian $J$ has $\cO \defeq \End_\Q(J) \isom \Z[\sqrt{2}]$ and
		$r_\an = 0 = r$, $J(\Q) = 0$, $\Sel_2(J/\Q) = 0$,
		$c_2 = c_5 = 1$, and $I_{\Q(\sqrt{-31}),\pi} = 7$.
		
		In particular, $\#\Sha(J/\Q)_\an = 7^2$ and $\#\Sha(J/\Q)$
		is a power of~$7$.
	\end{proposition}
	
	\begin{proof}
		The endomorphism ring contains~$\Z[\sqrt{2}]$ by construction
		of the curve. Since $\Z[\sqrt{2}]$ is the maximal order of~$\Q(\sqrt{2})$,
		it follows that $\cO \isom \Z[\sqrt{2}]$.
		A computation of the $2$-Selmer group shows that $\rk J(\Q) = 0$ and
		$\Sha(J/\Q)[2^\infty] = 0$. We check that $L(J/\Q, 1) \neq 0$ by computing
		$L(J/\Q, 1)/\Omega_J$ as described in Section~\ref{ssec: LJ/Om}. The torsion
		subgroup of~$J(\Q)$ turns out to be trivial. The Tamagawa number at~$5$ can
		be determined using van Bommel's Magma code~\cite[\S\,4.4]{vanBommel2019}.
		However, Magma is unable to compute a regular model at~$2$. So we computed a
		regular model by hand (see \texttt{Sha7-curve.m}
		for the computation of $c_2$) and found that the reduction type is $[\text{III}_2^*]$
		in~\cite{NamikawaUeno}. Note that this is also needed to determine the power
		of~$2$ in the \enquote{compensation factor}~$C$ in~\cref{computation of Omega_J}, which
		we need for the computation of $\#\Sha(J/\Q)_\an = 7^2$.
		Using the approach in Section~\ref{sec:computation-of-the-heegner-points-and-index},
		we find that $I_{K,\pi} = 7$ for $K = \Q(\sqrt{-31})$.  As all residual Galois
		representations are irreducible, the claim now follows from~\cref{thm:finite support}.
	\end{proof}
	
	\begin{lemma} \label[lemma]{upper bound on Sha}
		We have that $\#\Sha(J/\Q)[7^\infty] \mid 7^2$.
	\end{lemma}
	
	\begin{proof}
		Note that $7$ is split in the endomorphism ring~$\cO \isom \Z[\sqrt{2}]$,
		so the Heegner index as an ideal of~$\cO$ equals~$\frp$ with $\frp$ one of the two prime
		ideals above~$7$. Using~\cref{p-adic image}, we find that the $\frp$-adic Galois representations
		for these two primes have image~$\GL_2(\Z_7)$. Furthermore, $7 \nmid h_{\Q(\sqrt{-31})} N$.
		Hence~\cref{Howard Euler system} shows $\Sha(J/K) \inj (\Z/7)^2$.
		Since $[K : \Q] = 2$ is coprime to~$7$, we get $\Sha(J/\Q) \inj (\Z/7)^2$.
	\end{proof}
	
	Applying~\cref{thm:when-visible-sha}, we can show the following.
	
	\begin{proposition}[Sam Frengley]
		\label[proposition]{prop:7-div-sha}
		Let $E$ be the elliptic curve with LMFDB label \LMFDBLabel{3200.a1} and Weierstrass equation
		\[ E \colon y^2 = x^3 - 100x + 400 \,. \]
		There exists a prime $\frp \mid 7$ in $\cO$ such that $J[\frp] \isom E[7]$ as
		$\Gal(\Qbar|\Q)$-modules. Moreover, if $\Delta \subset J \times E$ is the graph
		of this isomorphism, then $\Sha(J/\Q)$ contains a subgroup isomorphic to $(\Z/7)^2$
		which is visible in the abelian threefold $(J \times E)/\Delta$.
	\end{proposition}
	
	\begin{proof}
		We first show that $J[\frp]$ is isomorphic to~$E[7]$ as a $\Gal(\Qbar|\Q)$-module,
		following \cite[Theorem~6.3]{F_VEOO7ITTSGOAEC}. Let $\cK \defeq J/\{\pm 1\}$ be the
		Kummer surface of~$J$ given by the model in~\cite[Chapter~3]{CF_PTAMAOCOG2} and let
		$x_J \colon J \to \cK$ be the quotient morphism. Similarly let
		$x_E \colon E \to E/\{\pm 1\} \cong \bP^1$ be the $x$-coordinate morphism.
		Since the mod-$7$ Galois representation attached to~$E/\Q$ is surjective, by
		\cite[Proposition~6.1]{F_VEOO7ITTSGOAEC} to show that there exists a quadratic twist~$E^d/\Q$
		of~$E/\Q$ such that $J[\frp]$ is isomorphic to~$E^d[7]$ as a $\Gal(\Qbar/\Q)$-module
		it suffices to show that there exist points $P \in J[\frp]$ and $Q \in E[7]$ such that
		$\Q(x_J(P))$ and~$\Q(x_E(Q))$ are isomorphic. Using the approach detailed in
		\cite[Theorem~6.3]{F_VEOO7ITTSGOAEC} we give an explicit degree~$24$ number field~$L/\Q$
		and equations for points $x_J(P) \in \cK(L)$ and $x_E(Q) \in \bP^1(L)$ which generate
		$L/\Q$. For the computations see the file~\texttt{congruence.m}.
		
		Finally, if $d \in \Z$ is chosen to be squarefree, then $d$ is divisible only by bad
		primes of $E$ and~$C$. As discussed in~\cite[(5.2)]{F_VEOO7ITTSGOAEC},
		by \cite[Section~2.1]{FLSSSW} or~\cite[Lemma~3]{MS_COG2WGRAF2WARWP} an isomorphism of
		Galois modules $J[\frp] \cong E^d[7]$ induces a congruence
		\begin{equation} \label{eq:trace-cond}
			a_p(E^d)^2 - t_p a_p(E^d) + n_p \equiv 0 \bmod{7} .
		\end{equation}
		Here $t_p = p + 1 - N_1 $ and $n_p = (N_1^2 + N_2 )/2 - (p+1) N_1 - p$,
		where $N_1 = \# C(\F_p)$ and $N_2 = \# C(\F_{p^2})$. Note that $E$ and~$C$ have bad
		reduction at $2$ and~$5$ and good reduction at all other primes, and for each integer $d \neq 1$
		dividing $10$ the congruence in~\eqref{eq:trace-cond} fails to hold at one of
		$p = 11$, $17$, or $23$. It follows that $J[\frp]$ and~$E[7]$ are isomorphic
		as $\Gal(\Qbar|\Q)$-modules.
		
		To show that $\Sha(J/\Q)$ contains a subgroup isomorphic to~$(\Z/7)^2$, first note that
		the Tamagawa numbers of~$E/\Q$ are coprime to~$7$ and by~\cref{properties of Sha 49 curve}
		so are the Tamagawa numbers of~$J/\Q$. The torsion subgroups of $E/\Q$ and~$J/\Q$ are trivial,
		the rank of~$E/\Q$ is~$2$, and the rank of~$J/\Q$ is~$0$, again by~\cref{properties of Sha 49 curve}.
		It follows from Theorem~\ref{thm:when-visible-sha} that $\Sha(J/\Q)[7]$ contains a subgroup
		isomorphic to $(\Z/7)^2$ which is visible in $(J \times E)/\Delta$.
	\end{proof}
	
	Combining these results, we obtain the following.
	
	\begin{theorem}
		For $J$ as above, we have $\#\Sha(J/\Q)_\an = 7^2 = \#\Sha(J/\Q)$.
	\end{theorem}
	
	The computations with precision~462 took 57~hours and $3.3$\,GiB~RAM
	on an MIT server running Magma V2.28-3 provided to us by Andrew Sutherland.
	The log of \texttt{N3200.m} can be found in \texttt{3200.log}. The bottleneck
	was the computation of the Heegner point, which required 783700 Fourier coefficients
	of the newform. (The remaining computations take less than 2~minutes.)
	Note that one must use $c_2 = 1$
	from~\cref{properties of Sha 49 curve} to obtain the correct value $\#\Sha(J/\Q)_\an = 7^2$;
	Magma cannot compute a regular semistable model of the curve at $p = 2$.
	
	\subsection{Further examples}
	
	One expects the existence of elements of order~$p$ in~$\Sha$ in quadratic twists
	of~$A/\langle P \rangle$ where $0 \neq P \in A[p](K)$ by~\cite{ShnidmanWeiss}.
	
	\medskip
	
	Our computations of the analytic orders of Sha for the Jacobians~$J$ of the LMFDB~examples
	with $L$-rank~$1$ yield the following examples of twists of~$J$ by the first Heegner field
	such that there is nontrivial $p$-torsion in Sha for some $p \in \{3,5,7\}$.
	The number in parentheses indicates the index of the curve in the list of LMFDB examples.
	\begin{enumerate}[(i)]
		\item[(30)] The twist of the second curve with $N = 133$
		by $D_K = -31$ has $I_K = 3^2$ and $\#\Sha(J^K/\Q)_\an = 2^2\cdot3^2$.
		
		\item[(55)] The twist of the first curve with $N = 275$
		by $D_K = -19$ has $I_K = 3$ and $\#\Sha(J^K/\Q)_\an = 2^2\cdot3^2$.
		
		\item[(57)] The twist of the curve with $N = 289$
		by $D_K = -15$ has $I_K = 3$ and $\#\Sha(J^K/\Q)_\an = 2^2\cdot3^2$.
		
		\item[(74)] The twist of the curve with $N = 523$
		by $D_K = -35$ has $I_K = 3^2$ and $\#\Sha(J^K/\Q)_\an = 2^2\cdot3^4$.
		
		\item[(77)] The twist of the first curve with $N = 621$
		by $D_K = -11$ has $I_K = 7$ and $\#\Sha(J^K/\Q)_\an = 7^2$.
		
		\item[(82)] The twist of the curve with $N = 647$
		by $D_K = -11$ has $I_K = 5$ and $\#\Sha(J^K/\Q)_\an = 2^2\cdot5^2$.
	\end{enumerate}
	To obtain these values, we computed $\#\Sha(J/K)_\an$ for a Heegner field $K = \Q(\sqrt{D})$ and use
	that $\#\Sha(J/\Q)_\an$ is a power of~$2$ for the LMFDB~examples,
	hence $\#\Sha(J^K/\Q)_\an$ and~$\#\Sha(J/K)_\an$ differ by a power of~$2$.
	By~\cref{rho restricted to K still irreducible}, the set of $\frp \nmid 2$ with $\rho_{J/\Q,\frp}$
	irreducible remains the same when $\rho$ is restricted to~$G_K$.
	By~\cref{Tam quotient} and~\eqref{Eq:Tam JK}, the odd prime factors of~$\Tam(J^K/\Q)$
	are those of~$\Tam(J/\Q)$.
	
	To also show that $\#\Sha(J^K/\Q)$ agrees with the analytic order of Sha,
	one has to compute $\Sha(J^K/\Q)[2^\infty]$, which our computation
	predicts to be~$(\Z/2)^2$ in all of the above cases except the curve with $N = 621$.
	We verified this using the code from~\cite{FisherYan}.
	Thus, using~\cref{thm:finite support rank 1,Howard Euler system}, the remaining
	tasks for verifying strong BSD are as follows.
	\begin{enumerate}[(i)]
		\item[(30)] Show that $\dim_{\F_\frp} \Sha(J^K/\Q)[\frp] = 0$ and $= 2$ for the two prime ideals $\frp \mid 3$, respectively.
		
		\item[(55)] Show that $\dim_{\F_\frp} \Sha(J^K/\Q)[\frp] = 2$ for one of the two prime ideals $\frp \mid 3$. Since the Heegner index as an $\cO$-ideal has norm $3$, this shows $\Sha(J^K/\Q)[\frp] = 0$ for the other $\frp \mid 3$.
		
		\item[(57)] Show that $\dim_{\F_\frp} \Sha(J^K/\Q)[\frp] = 2$ for one of the two prime ideals $\frp \mid 3$. Since the Heegner index as an $\cO$-ideal has norm $3$, this shows $\Sha(J^K/\Q)[\frp] = 0$ for the other $\frp \mid 3$. For one of them, one can perform an isogeny descent.
		
		\item[(74)] Show that $\Sha(J^K/\Q)$ has a subgroup isomorphic to $(\Z/3)^4$; note that the prime $3$ is inert in $\Z[f]$.
		
		\item[(77)] Show that $\dim_{\F_\frp} \Sha(J^K/\Q)[\frp] = 2$ for one of the two prime ideals $\frp \mid 7$. Since the Heegner index as an $\cO$-ideal has norm $7$, this shows $\Sha(J^K/\Q)[\frp] = 0$ for the other $\frp \mid 7$.
		
		\item[(82)] Show that $\dim_{\F_\frp} \Sha(J^K/\Q)[\frp] = 2$ for one of the two prime ideals $\frp \mid 5$.
	\end{enumerate}


		\section*{Acknowledgments}
		We thank Francesc Castella, Kęstutis Česnavicius, Maarten Derickx, Michael Dettweiler, Tom Fisher,
		Sam Frengley, Pip Goodman, Walter Gubler, Peter Hum\-phries, David Loeffler, Steffen Müller, Yusuf Mustopa,
		Jakob Stix, Chris Williams, and Chris Wuthrich for hints and helpful discussions.
		We thank the referee for their careful reading and suggestions to improve the article.
		
		\section*{Funding statement}
		The authors were supported by the Deutsche Forschungsgemeinschaft (DFG),
		STO~299/18-1 (AOBJ: 667349) (MS and TK) and~STO~299/18-2 (AOBJ: 686837) (MS)
		while working on this article.
		TK was partially supported by the 2021 MSCA Postdoctoral Fellowship
		01064790 -- ExplicitRatPoints. MS was supported by the National Science Foundation
		under Grant No.~DMS-1928930 while he was in residence at the Simons Laufer
		Mathematical Sciences Institute in Berkeley, California, during the Spring 2023 semester.
		This stay was also supported through DFG grants STO~299/13-1 (AOBJ: 635353) and
		STO~299/17-1 (AOBJ: 662415).
		
		\newcommand{\etalchar}[1]{$^{#1}$}
		\providecommand{\bysame}{\leavevmode\hbox to3em{\hrulefill}\thinspace}
		\providecommand{\MR}{\relax\ifhmode\unskip\space\fi MR }
		\providecommand{\MRhref}[2]{%
			\href{http://www.ams.org/mathscinet-getitem?mr=#1}{#2}
		}
		\providecommand{\href}[2]{#2}

\end{document}